\documentclass[12pt,a4paper]{article}
\usepackage{amsmath,amsthm,amsfonts,graphicx,psfrag,amssymb,dsfont,mathrsfs}
\usepackage{minitoc,color}

\usepackage[colorlinks]{hyperref}
\usepackage[alphabetic,initials]{amsrefs}
\usepackage{bbm}
\usepackage{dsfont}
\usepackage[ngerman,english]{babel}
\usepackage{stackengine}
\usepackage{tikz}
\usepackage{enumitem}
\usepackage[ngerman,english]{babel}

\usetikzlibrary{matrix,arrows}

\newcommand\dhxrightarrow[2][]{%
  \mathrel{\ooalign{$\xrightarrow[#1\mkern4mu]{#2\mkern4mu}$\cr%
  \hidewidth$\rightarrow\mkern4mu$}}
}

\newcommand{\R}{\mathbb{R}}
\newcommand{\C}{\mathbb{C}}
\newcommand{\Z}{\mathbb{Z}}
\newcommand{\N}{\mathbb{N}}
\newcommand{\J}{\mathcal{J}}

\newcommand{\wind}{\operatorname{wind}}
\newcommand{\A}{\mathbf{A}}
\newcommand{\crit}{\operatorname{crit}}
\newcommand{\sign}{\operatorname{sign}}
\newcommand{\F}{\mathcal{F}}
\newcommand{\ind}{\operatorname{ind}}

\def\br#1{\left\{#1\right\}}
\def\bp#1{\left(#1\right)}

\def\fl#1{\lfloor#1\rfloor}
\def\ceil#1{\lceil#1\rceil}
\def\abs#1{\left|#1\right|}

\theoremstyle{plain}
\newtheorem{thm}{Theorem}[section]

\newtheorem{conjeorem}[thm]{Conjeorem}
\newtheorem{prop}[thm]{Proposition}
\newtheorem{cor}[thm]{Corollary}
\newtheorem{lemma}[thm]{Lemma}

\theoremstyle{definition}
\newtheorem{example}[thm]{Example}
\newtheorem*{notation}{Notation}
\newtheorem{defn}[thm]{Definition}

\newtheorem{disclaimer}[thm]{Disclaimer}
\newtheorem{remark}[thm]{Remark}
\newtheorem{conj}[thm]{Conjecture}
\newtheorem{assumptions}[thm]{Assumptions}

\usepackage{scalerel}

\newcommand\reallywidehat[1]{\arraycolsep=0pt\relax%
\begin{array}{c}
\stretchto{
  \scaleto{
    \scalerel*[\widthof{\ensuremath{#1}}]{\kern-.5pt\bigwedge\kern-.5pt}
    {\rule[-\textheight/2]{1ex}{\textheight}} 
  }{\textheight} %
}{0.8ex}\\           
#1\\                 
\rule{-1ex}{0ex}
\end{array}
}

\date{}

\begin{document}
\begin{titlepage}
\begin{center}
\LARGE\textbf{Algebraic torsion in higher-dimensional contact manifolds} \\
\vspace{0.8cm}
\Large\textbf{Dissertation} \\ 
\vspace{0.2cm}

zur Erlangung des akademischen Grades \\
\vspace{0.2cm}

\textbf{Doctor rerum naturalium} \\
(Dr. rer. nat.) \\
\vspace{0.2cm}

\textbf{Im Fach: Mathematik} \\
eingereicht an der \\
\vspace{0.2cm}

\textbf{Mathematisch-Naturwissenschaftlichen Fakult\"{a}t} \\
\textbf{der Humboldt-Universit\"{a}t zu Berlin} \\
\vspace{0.2cm}
von \\ 

\textbf{M.Sc Agustin Moreno} \\
geboren am 17. August 1990 in Montevideo, Uruguay \\

\vspace{0.2cm}
\textbf{Pr\"{a}sident der Humboldt-Universit\"{a}t zu Berlin} \\
\textbf{Prof. Dr.-Ing. Dr. Sabine Kunst} \\

\vspace{0.2cm}
\textbf{Dekan der Mathematisch-Naturwissenschaftlichen Fakult\"{a}t}\\
\textbf{Prof. Dr. Elmar Kulke}\\ 
\end{center}

\begin{flushleft} \large \textbf{Gutachter:} \\ 
\hspace{0.5cm}\textbf{1.} \\
\vspace{0.1cm}
\hspace{0.5cm}\textbf{2.} \\
\vspace{0.1cm}
\hspace{0.5cm}\textbf{3.} \\
\vspace{0.2cm}
\textbf{Tag der m\"{u}ndlichen Pr\"{u}fung:}
\end{flushleft}
\end{titlepage}

\thispagestyle{empty}
\newpage

\vspace*{\fill}
\begingroup 
\centering
\Large \textbf{Selbst\"{a}ndigkeitserkl\"{a}rung}
\vspace{0.2cm}

\normalsize Hiermit erkl\"{a}re ich, dass ich die vorliegende Dissertation selbst\"{a}ndig und nur unter Verwendung der angegebenen Literatur und Hilfsmittel angefertig habe.
\endgroup
\vspace*{\fill}

\newpage

\vspace*{\fill}
\begin{abstract} We construct examples in any odd dimension of contact manifolds with finite and non-zero algebraic torsion (in the sense of \cite{LW}), which are therefore tight and do not admit strong symplectic fillings.
We prove that Giroux torsion implies algebraic $1$-torsion in any odd dimension, which proves a conjecture in \cite{MNW}. We construct infinitely many non-diffeomorphic examples of $5$-dimensional contact manifolds which are tight, admit no strong fillings, and do not have Giroux torsion. We obtain obstruction results for symplectic cobordisms, for which we give a proof not relying on SFT machinery. We give a tentative definition of a higher-dimensional spinal open book decomposition, based on the $3$-dimensional one of \cite{LVHMW}. An appendix written in co-authorship with Richard Siefring gives a basic outline of the intersection theory for punctured holomorphic curves and hypersurfaces, which generalizes his $3$-dimensional results of \cite{Sie11} to higher dimensions. From the intersection theory we obtain an application to codimension-$2$ holomorphic foliations, which we use to restrict the behaviour of holomorphic curves in our examples.  
\end{abstract}
\vspace*{\fill}

\addtocontents{toc}{\protect\sloppy}

\newpage

\vspace*{\fill}
\selectlanguage{ngerman}
\begin{abstract}
Wir konstruieren Beispiele von Kontaktmannigfaltigkeiten in jeder
ungeraden Dimension, welche endliche nicht-triviale algebraische
Torsion (im Sinne von \cite{LW}) aufweisen, somit straff sind und keine
starke symplektische F\"ullung haben. Wir beweisen, dass Giroux Torsion
algebraische $1$-Torsion in jeder ungeraden Dimension impliziert, womit
eine Vermutung aus \cite{MNW} bewiesen wird. Wir konstruieren unendlich
viele nicht diffeomorphe Beispiele von $5$-dimensionalen
Kontaktmannigfaltigkeiten, welche straff sind, keine starke
symplektische F\"ullung zulassen und keine Giroux Torsion haben. Wir
erhalten Obstruktionen f\"ur symplektische Kobordismen, ohne f\"ur
deren Beweis die SFT Maschinerie zu verwenden. Wir geben eine
provisorische Definition eines spinalen offenen Buchs in h\"oherer
Dimension an, basierend auf der vom $3$-dimensionalen Fall
aus \cite{LVHMW}. In einem Anhang geben wir in gemeinsamer
Autorenschaft mit Richard Siefring eine wesentliche Zusammenfassung
der Schnitttheorie f\"ur punktierte holomorphe Kurven und
Hyperfl\"achen an, welche die $3$-dimensionalen Resultate aus
\cite{Sie11} auf h\"ohere Dimensionen verallgemeinert. Mittels der
Schnitttheorie erhalten wir eine Anwendung f\"ur holomorphe
Bl\"atterungen von Kodimension zwei, die wir benutzen um das Verhalten
von holomorphem Kurven in unseren Beispielen einzuschr\"anken.
\end{abstract}
\selectlanguage{english}
\vspace*{\fill}

\addtocontents{toc}{\protect\sloppy}

\newpage

\tableofcontents

\newpage

\section{Introduction}

\subparagraph*{Background and motivation.}

The full power of holomorphic curves, so far, has been exploited most successfully in dimensions 3 and 4, where low-dimensional techniques and Gauge theory are also available. In higher dimensions, on the one hand, the \emph{flexible} side of contact and symplectic topology, the one concerning phenomena governed by h-principles and algebraic topology, has seen some spectacular recent advances. Most notably, with the introduction of \emph{loose Legendrians} \cite{Murphy}, and the generalization of the classically 3-dimensional notion of \emph{overtwistedness} \cite{Eli89,BEM}. The classification problem for overtwisted contact manifolds is thus reduced to a classification problem of \emph{almost} contact manifolds, which belongs to the world of topology/obstruction theory. Moreover, being overtwisted is an obstruction to being \emph{symplectically fillable}, that is, of being the convex boundary of a symplectic manifold \cite{Gro85,Eli90,Nie,BEM}. On the other hand, the \emph{rigid} side, the one concerning holomorphic curves and Floer-type invariants, seems, in comparison and in higher dimensions, much less explored. The results are more topological in spirit (e.g.\ \cite{McD,BGZ,GNW}), and the theory is less systematic, consisting more of a collection of interesting examples and phenomena. Rather unfortunate is the absence of Gauge theory in higher dimensions, although we still have holomorphic curves as an available tool at our disposal.

Whereas the classification of overtwisted contact manifolds is a matter of topology, the understanding of \emph{tight} contact manifolds (i.e.\ not overtwisted) is much more delicate. In dimension 3, the theory of \emph{convex surfaces} has been specially fruitful \cite{Hoa, Hob,CGH,Et}.  
In higher dimensions, while there have been attempts at a theory of convex hypersurfaces (notably by Honda and his student Yang Huang---see \cite{HH}), not much is known to be true. In fact, many of the classification results obtained in dimension 3 are known to be false. For instance, the uniqueness of tight contact structures on spheres, disproved by the computations in \cite{U99}, or the contact prime decomposition theorem of \cite{Col97}, disproved in \cite{GNW}. Another example of $3$-dimensional phenomena, which does not as of yet have a higher-dimensional analogue, is the coarse classification of \cite{CGH}. Loosely speaking, it tells us that in order to find examples of infinitely many non-isotopic contact structures in the same homotopy class, then adding Giroux torsion is the only way (see below for a definition). One can wonder if the same holds in higher dimensions.

One could try to establish hierarchical structures in the category of all contact manifolds, and to relate it to the classification problem. An attempt to do this is an invariant called \emph{algebraic torsion}. In order to address such (very hard) questions, first one needs a better repertoire of methods for constructing potentially interesting contact structures in higher dimensions, together with developing techniques to handle them. This will be the main direction we will be pursuing in this document. 

\subparagraph*{Goals.} We will construct higher-dimensional examples of contact manifolds $M$, with ``interesting'' properties, and we will develop techniques in order to compute rigid invariants. They will present a geometric structure which we call a \emph{spinal open book decomposition} or SOBD (based on the 3-dimensional version of \cite{LVHMW}), of a certain type which can be ``detected'' algebraically by algebraic torsion, a holomorphic-curve contact invariant. Though we will not formalize it in such form, these type of SOBDs, which one could call \emph{partially planar}, mimics the notion of \emph{planar $m$-torsion} domains as defined in \cite{Wen2}. Indeed, they contain a portion, the \emph{planar piece}, which is a fibration resembling an open book, whose fibers are surfaces of genus zero and $m+1$ boundary components; and another similar portion, the \emph{non-planar piece}, which is not diffeomorphic to the former. For suitable data, these surfaces become holomorphic, and are leaves of a foliation of $\mathbb{R}\times M$. The isolated ones may be counted in a suitable way, and the result is an invariant which ``recovers'' the number $m$. This is the idea inspiring algebraic $m$-torsion. One would recover the fact that a---suitably defined--- planar $m$-torsion (a geometric structure) implies algebraic $m$-torsion (an algebraic fact), in higher dimensions.

We exhibit two main constructions (model A and model B) of contact forms in the same underlying manifold, in some sense dual to each other, and which hold in any odd dimension. Both are ``supported'' by dual SOBDs, having the same geometric decomposition, but where the role of each piece is reversed. One can do such reversion, since the ``monodromy'' is trivial. Therefore one can view these contact forms as ``Giroux'' forms for each decomposition, and in fact they are isotopic to each other. In general, we expect to associate, to any SOBD (with suitable extra data), a unique isotopy class of supported contact structures, which generalizes a part of Giroux's correspondence to this setting. 

We will aim at computing the algebraic torsion of these examples, which we show is finite, and, in certain cases, non-zero. In those cases, they are tight and admit no strong symplectic fillings. In order to estimate this invariant, we need a detailed understanding of holomorphic curves in their symplectization, and the SOBD structure is very helpful towards this end. 

Whereas for model A one obtains a foliation of its symplectization by holomorphic curves, for model B, one gets a foliation by holomorphic hypersurfaces (i.e.\ real codimension-2). By a careful and non-trivial asymptotic analysis, one can show that holomorphic curves with certain prescribed asymptotics are restricted to lie in the leaves of the hypersurfaces. This follows from an application of the intersection theory between curves and hypersurfaces which are asymptotically cylindrical, in the sense of Richard Siefring. We include an outline of the basic theory, in an appendix written in co-authorship with Siefring, in which we exhibit this application to codimension-2 holomorphic foliations. This is a prequel of his upcoming work \cite{Sie2}, generalizing his results of \cite{Sie11}.

We will also relate algebraic torsion with a geometric condition, \emph{Giroux torsion}. While this is a classical notion in dimension 3, the higher-dimensional version was introduced in \cite{MNW}. We will show that the geometric presence of certain torsion domains inside a contact manifold can be detected algebraically by SFT. More concretely, Giroux torsion implies algebraic $1$-torsion, in any odd dimension. Furthermore, we shall restrict our study to certain $5$-dimensional cases of our construction, for which we can show do not admit Giroux torsion, and for which we can obstruct the existence of certain exact cobordisms, in a way which does not use the abstract perturbation scheme for SFT promised by polyfold theory. We shall derive a corollary from this, preventing the existence of exact symplectic cobordisms between one of these examples, and the contact manifold resulting from applying a \emph{Lutz-Mori twist} to it. This is an interesting result, since the underlying manifolds are diffeomorphic, and the contact structures are homotopic as almost contact structures.

We will also give a tentative definition of an SOBD in higher-dimensions, since the language will prove rather useful, and discuss some particular cases. This notion was originally defined in \cite{LVHMW} in dimension 3, and it was motivated by the study of the structure one gets in the boundary of a 4-dimensional Lefschetz fibration having any surface (not just a disk) as a base. The case of the disk corresponds to the classical notion of an open book decomposition, which are thus a particular case of SOBD. By a result in \cite{DiGe}, which generalizes old work by Lutz in dimension 3, such structures appear naturally on $S^1$-principal bundles over higher-dimensional symplectic bases, defined by an Euler class which vanishes along a set of dividing hypersurfaces. We call these \emph{prequantization SOBDs}, and we have a particular case where the $S^1$-bundles are trivial, which corresponds to the case of \emph{Giroux SOBDs}. For these SOBDs we give a notion of a ``Giroux form'', giving a  canonical isotopy class of contact structures.

\subparagraph*{On the invariant.} The invariant we will use, algebraic torsion, has connections to questions about symplectic fillability. As defined in \cite{LW}, it is a contact invariant taking values in $\mathbb{Z}^{\geq 0}\cup \{\infty\}$, which was introduced, using the machinery of \emph{Symplectic Field Theory}, as a quantitative way of measuring how non-fillable a given contact structure is. In this sense, the philosophy is: the less order of algebraic torsion our manifold has, the ``less fillable'' it is, giving rise to a ``hierarchy of fillability obstructions'', cf.\ \cite{Wen2}. At least morally, $0$-torsion should correspond to overtwistedness, whereas $1$-torsion is implied by Giroux torsion (the converse is not true). 

Having $0$-torsion is actually equivalent to being \emph{algebraically overtwisted}, which means that the contact homology, or equivalently its SFT, vanishes (Proposition 2.9 in \cite{LW}). This is well-known to be implied by overtwistedness, but the converse is still open. In fact, there have been claims on results that could be interpreted as evidence that it does not hold in higher dimensions \cite{Ek}, but these have been recently withdrawn. The fact that Giroux torsion implies algebraic $1$-torsion was already known in dimension 3, and its higher-dimensional analogue was conjectured in \cite{MNW}.

The key fact about this invariant is that it behaves well under exact symplectic cobordisms, which implies that the concave end inherits any order of algebraic torsion that the convex end has. Thus, algebraic torsion may be also thought of as an obstruction to the existence of exact symplectic cobordisms. In particular, it serves as an obstruction to symplectic fillability. 

Moreover, there are connections to dynamics: any contact manifold with finite torsion satisfies the Weinstein conjecture (i.e.\ there exist closed Reeb orbits for every contact form).  

One should mention that there are other notions of algebraic torsion in the literature which do not use SFT, but which are only 3-dimensional (see \cite{KMvhMW} for the version using Heegard Floer homology, or the appendix in \cite{LW} by Hutchings, using ECH).

\subparagraph*{Statement of results.}

For the SFT setup, we follow \cite{LW}, where we refer the reader for more details. We will take the SFT of a contact manifold $(M,\xi)$ (with coefficients) to be the homology $H^{SFT}_*(M, \xi; \mathcal{R})$ of a $\mathbb{Z}_2$-graded unital $BV_\infty$-algebra $(\mathcal{A}[[\hbar]], \mathbf{D}_{SFT})$ over the group ring $R_{\mathcal{R}}:=\mathbb{R}[H_2(M;\mathbb{R})/\mathcal{R}]$, for some linear subspace $\mathcal{R}\subseteq H_2(M;\mathbb{R})$. Here, $\mathcal{A}=\mathcal{A}(\lambda)$ has generators $q_\gamma$ for each good closed Reeb orbit $\gamma$ with respect to some nondegenerate contact form $\lambda$ for $\xi$, $\hbar$ is an even variable, and the operator $$\mathbf{D}_{SFT}: \mathcal{A}[[\hbar]] \rightarrow\mathcal{A}[[\hbar]]$$ is defined by counting rigid solutions to a suitable abstract perturbation of a $J$-holomorphic curve equation in the symplectization of $(M, \xi)$. It satisfies

\begin{itemize}
 \item $\mathbf{D}_{SFT}$ is odd and squares to zero,
 \item $\mathbf{D}_{SFT}(1) = 0$, and
 \item $\mathbf{D}_{SFT} = \sum_{k \geq 1} D_k\hbar^{k-1},$
\end{itemize}

where $D_k : \mathcal{A} \rightarrow \mathcal{A}$ is a differential operator of order $\leq k$, given by
$$
D_k = \sum_{\substack{\Gamma^+,\Gamma^-, g,d\\ |\Gamma^+|+g=k}}\frac{n_g(\Gamma^+,\Gamma^-,d)}{C(\Gamma^-,\Gamma^+)}q_{\gamma_1^-}\dots q_{\gamma_{s^-}^-}z^d\frac{\partial}{\partial q_{\gamma_1^+}}\dots \frac{\partial}{\partial q_{\gamma_{s^+}^+}}
$$
The sum ranges over all non-negative integers $g \geq 0$, homology classes $d \in H_2(M; \mathbb{R})/\mathcal{R}$ and ordered (possibly empty) collections of good closed Reeb orbits $\Gamma^\pm = (\gamma_1^{\pm},\dots,\gamma_{s^\pm}^\pm)$ such that $s^+ + g = k$. After a choice of spanning surfaces as in \cite{EGH} (p. 566, see also p. 651), the projection to $M$ of each finite
energy holomorphic curve $u$ can be capped off to a 2-cycle in $M$, and so it gives rise to a homology class $[u]\in H_2(M)$, which we project to define $\overline{[u]} \in H_2(M; \mathbb{R})/\mathcal{R}$. The number $n_g(\Gamma^+,\Gamma^-,d) \in \mathbb{Q}$ denotes the count of (suitably perturbed) holomorphic curves of genus $g$ with positive asymptotics $\Gamma^+$ and negative asymptotics $\Gamma^-$ in the homology class $d$, including asymptotic markers as explained in \cite{EGH}, or \cite{Wen3}, and including rational weights arising from automorphisms. $C(\Gamma^-, \Gamma^+) \in \mathbb{N}$ is a combinatorial factor defined as $C(\Gamma^-, \Gamma^+) = s^-!s^+!\kappa_{\gamma_1^-}\dots\kappa_{\gamma_{s^-}^-}$, where $\kappa_\gamma$ denotes the covering multiplicity of the Reeb orbit $\gamma$.

The most important special cases for our choice of linear subspace $\mathcal{R}$ are $\mathcal{R} = H_2(M; \mathbb{R})$ and $\mathcal{R} = \{0\}$, called the \emph{untwisted} and \emph{fully twisted} cases respectively, and $\mathcal{R} = \ker \Omega$ with $\Omega$ a closed 2-form on $M$. We shall abbreviate the latter case as $H^{SFT}_*(M, \xi;\Omega) := H^{SFT}_*(M, \xi; \ker \Omega)$, and the untwisted case simply by $H^{SFT}_*(M, \xi):=H^{SFT}_*(M, \xi;H_2(M;\mathbb{R}))$.

\begin{defn} Let $(M, \xi)$ be a closed manifold of dimension $2n+1$ with a positive, co-oriented contact structure. For any integer $k \geq 0$, we say that $(M, \xi)$ has $\Omega$-twisted algebraic torsion of order $k$ (or $\Omega$-twisted $k$-torsion) if $[\hbar^k] = 0$ in $H^{SFT}_*(M, \xi;\Omega)$. If this is true for all $\Omega$, or equivalently, if $[\hbar^k] = 0$ in $H^{SFT}_*(M, \xi;\{0\})$, then we say that $(M, \xi)$ has fully twisted algebraic $k$-torsion. 
\end{defn}

We will refer to \emph{untwisted} $k$-torsion to the case $\Omega=0$, in which case $R_\mathcal{R}=\mathbb{R}$ and we do not keep track of homology classes. Whenever we refer to torsion without mention to coefficients we will mean the untwisted version. We will say that, if a contact manifold has algebraic $0$-torsion for every choice of coefficient ring, then it is \emph{algebraically overtwisted}, which is equivalent to the vanishing of the SFT, or its contact homology. By definition, $k$-torsion implies $(k+1)$-torsion, so we may define its \emph{algebraic torsion} to be
$$
AT(M,\xi;\mathcal{R}):=\min\{k\geq 0:[\hbar^k]=0\}\in \mathbb{Z}^{\geq 0}\cup\{\infty\},
$$
where we set $\min \emptyset=\infty$. We denote it by $AT(M,\xi)$, in the untwisted case.

This construction is well-behaved under symplectic cobordisms: 
Any exact symplectic cobordism $(X, \omega=d\alpha)$ with positive end $(M_+, \xi_+)$ and negative end $(M_-, \xi_-)$ gives rise to a natural $\mathbb{R}[[\hbar]]$-module morphism on the untwisted SFT,
$$\Phi_X: H^{SFT}_*(M_+, \xi_+) \rightarrow H^{SFT}_*(M_-, \xi_-),$$ a \emph{cobordism map}. This implies that if $(M_+, \xi_+)$ has $k$-torsion, then so does $(M_-, \xi_-)$. There is also a version with coefficients for the case of non-exact cobordisms (see \cite{LW}), but we will use only the version for fillings, which we now state. In the following, $\overline{R_{\mathcal{R}}}$ will denote the \emph{Novikov completion} of the ring $R_{\mathcal{R}}$ (see \cite{EGH}). 

\begin{prop}\cite{LW}\label{SFTfill} Suppose $(X, \omega)$ is a strong filling of $(M, \xi)$ (see next section for a definition), and $\mathcal{R}(X) \subseteq \ker \omega \subseteq H_2(X; \mathbb{R})$ and $\mathcal{R} \subseteq H_2(M; \mathbb{R})$ are linear subspaces for which the natural map from $H_2(M; \mathbb{R})$ to $H_2(X; \mathbb{R})$ takes $\mathcal{R}$ into $\mathcal{R}(X)$. Then there is a natural $\mathbb{R}[[\hbar]]$-module morphism
$$
\Phi_X : H^{SFT}_*(M, \xi; \mathcal{R}) \rightarrow \overline{R_{\mathcal{R}(X)}}[[\hbar]],
$$
which acts on $R_{\mathcal{R}}[[\hbar]] \subseteq H^{SFT}_*(M, \xi; \mathcal{R})$ as the natural map to $R_{\mathcal{R}(X)}[[\hbar]]$ induced
by the inclusion $M \hookrightarrow X$. In particular, the untwisted SFT of $(M, \xi)$ admits a nonzero $\mathbb{R}[[\hbar]]$-module morphism
$$
\Phi_X : H^{SFT}_*(M, \xi)\rightarrow \overline{\mathcal{R}_{\ker \omega}}[[\hbar]]
$$
\end{prop}

\begin{remark} The existence of the above maps is actually valid under a weaker notion of fillability, called \emph{stable fillability}, which in dimension 3 is equivalent to weak fillability (see \cite{LW}). Any weak filling can be deformed to a stable filling, up to, if necessary, slightly perturbing the cohomology class of the symplectic form so that it is rational along the boundary (\cite{MNW}, corollary 1.12). It follows that if a contact manifold $(M,\xi)$ has finite algebraic torsion with coefficients in $R_\mathcal{R}$, then it does not admit any weak fillings $(W,\omega)$ for which $[\omega\vert_{M}]$ annihilates every element in $\mathcal{R}$, provided that this cohomology class is rational. This rationality condition is because the annihilation condition is not preserved under perturbations (unless $\mathcal{R}=0$, in which case it has finite fully twisted torsion, and does not admit weak fillings whatsoever). It is known to be unnecessary in dimension $3$, and probably also is in higher dimensions.
\end{remark}

Examples of 3-dimensional contact manifolds with any given order of torsion $k-1$, but not $k-2$, were constructed in \cite{LW}. The underlying manifold is the product manifold $M_g:=S^1\times \Sigma$, for $\Sigma$ a surface of genus $g$ which is divided into two pieces $\Sigma_+$ and $\Sigma_-$ along some dividing set of simple closed curves $\Gamma$ of cardinality $k$, where the latter has genus $0$, and the former has genus $g-k+1$. The contact structure $\xi_k$ is $S^1$-invariant and may be obtained, for instance, by a construction originally due to Lutz (see \cite{Lutz}). Its isotopy class is characterized by the fact that every section $\{pt\}\times\Sigma$ is a convex surface with dividing set $\Gamma$. The behaviour of algebraic torsion under cobordisms then implies that there is no exact symplectic cobordisms having $(M_g,\xi_k)$ and $(M_{g^\prime},\xi_{k^\prime})$ as convex and concave ends, respectively, if $k<k^\prime$.

The existence of the analogue higher dimensional contact manifolds was conjectured in \cite{LW}. We will consider a modified version of their examples. The modification we do here consists in taking the $S^1$-factor and replacing it by a \emph{closed} $(2n-1)$-manifold $Y$, having the special property that $Y\times I$ admits the structure of a Liouville domain (here, $I$ denotes the interval $[-1,1]$). This means that it comes with an exact symplectic form $d\alpha$, and has disconnected contact-type boundary $\partial(Y \times I,d\alpha)=(Y_-,\xi_-=\ker\alpha_-)\bigsqcup (Y_+,\xi_+=\ker\alpha_+)$, where $Y_\pm$ coincide with $Y$ as manifolds, but $(Y_\pm,\xi_\pm)$ are \emph{not contactomorphic} to each other. In fact, $Y_\pm$ have different orientations, and so they might not even be \emph{homeomorphic} to each other (not every manifold admits an orientation-reversing homeomorphism). A Liouville domain of the form $(Y\times I,d\alpha)$ is what we will call a \emph{cylindrical Liouville semi-filling} (or simply a cylindrical semi-filling). Their existence in every odd dimension was established in \cite{MNW}. We immediately see that this generalizes the previous 3-dimensional example, since $S^1$ admits the \emph{Liouville pair} $\alpha_\pm=\pm d\theta$ (which means that the 1-form $e^{-s}\alpha_- +e^s\alpha_+$ is Liouville in $S^1 \times \mathbb{R}$). We prove that the manifold $Y\times \Sigma$ indeed achieves $(k-1)$-torsion (Theorem \ref{thm5d}), for a suitable contact structure which we now describe. 

First, for once and for good, we will fix the following notation:

\begin{notation}\label{not} Throughout this document, the symbol $I$ will be reserved for the interval $[-1,1]$.
\end{notation}

\begin{figure}[t]\centering \includegraphics[width=0.70\linewidth]{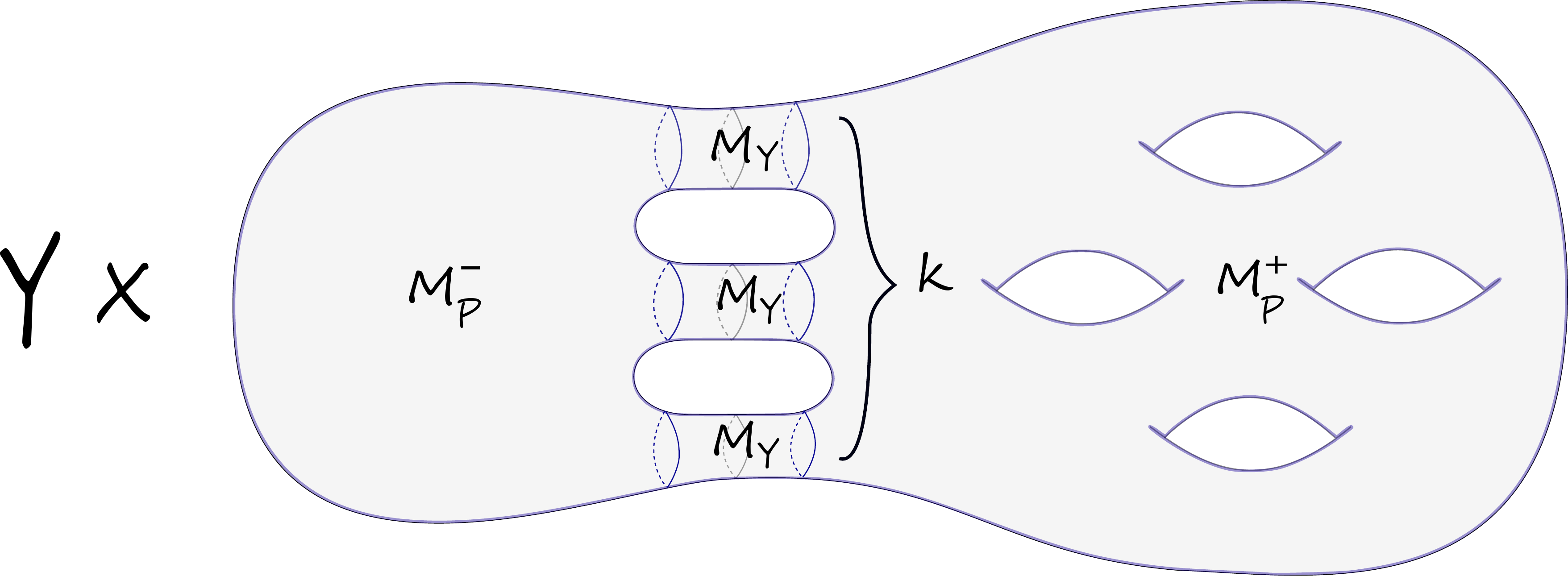}
\caption{\label{themanifold} The SOBD structure in $M$.}
\end{figure}

We can adapt the construction of the contact structures in \cite{LW} to our models. The starting idea is to decompose the manifold $	M=M_g=Y\times \Sigma$ into three pieces
$$
M_g=M_Y \bigcup M_P^\pm,
$$
where $M_Y = \bigsqcup^k Y \times I \times S^1$, and $M_P^\pm = Y \times \Sigma_\pm$ (see Figure \ref{themanifold}). 
We have natural fibrations $$\pi_Y: M_Y \rightarrow Y \times I$$$$ \pi_P^\pm: M_P^{\pm} \rightarrow Y_\pm,$$ with fibers $S^1$ and $\Sigma_\pm$, respectively, and they are compatible in the sense that
$$
\partial((\pi_P^\pm)^{-1}(pt))=\bigsqcup^k \pi_Y^{-1}(pt)
$$

While $\pi_Y$ has a Liouville domain as base, and a contact manifold as fiber, the situation is reversed for $\pi_P^\pm$, which has contact base, and Liouville fibers. This is a prototypical example of a SOBD. 

Using this decomposition, we can construct a contact structure $\xi_k$ which is a small perturbation of the stable Hamiltonian structure $\xi_\pm \oplus T\Sigma_\pm$ along $M_P^\pm$, and is a \emph{contactization} for the Liouville domain $(Y\times I,\epsilon d\alpha)$ along $M_Y$, for some small $\epsilon>0$. This means that it coincides with $\ker(\epsilon\alpha + d\theta)$, where $\theta$ is the $S^1$-coordinate. We will do this in detail in Section \ref{modelAsec}.

For these contact manifolds, we can estimate their algebraic torsion. First, recall that a contact structure is hypertight if it admits a contact form without contractible Reeb orbits (which we call a hypertight contact form). In particular, there are no holomorphic disks in their symplectization, which implies that there is no 0-torsion. By a well-known theorem by Hofer and its generalization to higher dimensions by Albers--Hofer (in combination with \cite{BEM}), hypertight contact manifolds are tight.  

\begin{thm}\label{thm5d}
For any $k\geq 1$, and $g\geq k$, the $(2n+1)$-dimensional contact manifolds $(M_g=Y\times \Sigma,\xi_k)$ satisfy $AT(M_g,\xi_k)\leq k-1$. Moreover, if $(Y,\alpha_\pm)$ are hypertight, and $k\geq 2$, the corresponding contact manifold $(M_g,\xi_k)$ is also hypertight. In particular, $AT(M_g,\xi_k)>0$, and it is tight.  
\end{thm}

In fact, the examples of Theorem \ref{thm5d} admit $\Omega$-twisted $k-1$-torsion, for $\Omega$ defining a cohomology class in $\mathcal{O}:=\mbox{Ann}(\bigoplus_k H_1(Y;\mathbb{R})\otimes H_1(S^1; \mathbb{R}))$, the annihilator of $\bigoplus_k H_1(Y;\mathbb{R})\otimes H_1(S^1; \mathbb{R}) \subseteq H_2(M_g;\mathbb{R})$. Here, we take the homology of the subregion $\bigsqcup^k Y\times \{0\}\times S^1$, lying along the region $M_Y$ where $\Sigma_\pm$ glue together. Using Proposition \ref{SFTfill} and the ensuing remark, we obtain:

\begin{cor}\label{thm5d2}
The examples of Theorem \ref{thm5d} do not admit weak fillings $(W,\omega)$ for which $[\omega\vert_{M_g}]$ is rational and lies in $\mathcal{O}$. In particular, they are not strongly fillable. 
\end{cor} 

\begin{remark} \label{mit}
\begin{itemize}[wide, labelwidth=!, labelindent=0pt] $\;$
\item By a result of Mitsumatsu in \cite{Mit}, any 3-manifold $Y$ which admits a smooth Anosov flow preserving a smooth volume form satisfies that $Y \times I$ can be enriched with a cylindrical Liouville semi-filling structure. Therefore any of these $3$-manifolds can be used in the construction of $5$-dimensional contact models with $AT\leq k-1$, for any $k\geq 1$.

\item The examples of Liouville cylindrical semi-fillings of \cite{MNW} satisfy the hypertightness condition. Then we have a doubly-infinite family of contact manifolds with $0< AT(M_g,\xi_k)\leq k-1$, in any dimension. These are then an instance of higher-dimensional tight but not strongly fillable contact manifolds, since they have non-zero and finite algebraic torsion. For $k=2$, this precisely computes the algebraic torsion.
\end{itemize}
\end{remark}

The authors of \cite{MNW} define a generalized higher-dimensional version of the notion of Giroux torsion. This notion is defined as follows: consider $(Y,\alpha_+,\alpha_-)$ a Liouville pair on a closed manifold $Y^{2n-1}$ (see Definition \ref{LP} below), and consider the \emph{Giroux $2\pi$-torsion domain} modeled on $(Y,\alpha_+,\alpha_-)$ given by the contact manifold $(GT,\xi_{GT}):=(Y\times [0,2\pi]\times S^1,\ker\lambda_{GT})$, where 
\begin{equation}
\begin{split}
\lambda_{GT}&=\frac{1+\cos(r)}{2}\alpha_++\frac{1-\cos(r)}{2}\alpha_-+\sin (r)d\theta\\
\end{split}
\end{equation} and the coordinates are $(r,\theta)\in [0,2\pi] \times S^1$. Say that a contact manifold $(M^{2n+1},\xi)$ has \emph{Giroux torsion} whenever it admits a contact embedding of $(GT,\xi_{GT})$. In this situation, denote by $\mathcal{O}(GT)\subseteq H^2(M;\mathbb{R})$ the annihilator of $\mathcal{R}_{GT}:=H_1(Y;\mathbb{R})\otimes H_1(S^1;\mathbb{R})$, viewed as a subspace of $H_2(M;\mathbb{R})$. The following was conjectured in \cite{MNW}:

\begin{thm}\label{AlgvsGiroux}
If a contact manifold $(M^{2n+1},\xi)$ has Giroux torsion, then it has $\Omega$-twisted algebraic 1-torsion, for every $[\Omega] \in \mathcal{O}(GT)$, where $GT$ is a Giroux $2\pi$-torsion domain embedded in $M$. 
\end{thm}

The proof uses the same techniques as Theorem \ref{thm5d}, and the main idea is to interpret Giroux torsion domains in terms of a specially simple kind of SOBD, which we call \emph{Giroux SOBD}. 

\vspace{0.5cm}

A natural corollary is the following:

\begin{cor}\label{GTfill} If a contact manifold $(M^{2n+1},\xi)$ has Giroux torsion, then it does not admit weak fillings $(W,\omega)$ with $[\omega\vert_M]\in \mathcal{O}(GT)$ and rational, where $GT$ is a Giroux $2\pi$-torsion domain embedded in $M$. In particular, it is not strongly fillable.
\end{cor}

This is essentially corollary 8.2 in \cite{MNW}, which was obtained with different methods. Observe that if $\mathcal{R}_{GT}=0$ then $(M,\xi)$ does not admit weak fillings at all. This is in fact the condition used in \cite{MNW} to obstruct weak fillability.

Conversely, one could ask whether there exist examples of non-fillable contact manifolds without Giroux torsion. Examples of 3-dimensional weakly but not strongly fillable contact manifolds without Giroux torsion were given in \cite{NW11}. In higher dimensions, the following theorem can be proved without appealing to the abstract perturbation scheme for SFT (see disclaimer \ref{disc} below). We use the fact that the unit cotangent bundle of a hyperbolic surface fits into a cylindrical semi-filling \cite{McD}.

\begin{thm}\label{nocob1} Let $(M_0^{5},\xi_0)$ be a $5$-dimensional contact manifold with Giroux torsion, and let $Y$ be the unit cotangent bundle of a hyperbolic surface. If $(M=\Sigma \times Y,\xi_k)$ is the corresponding $5$-dimensional contact manifold of Theorem \ref{thm5d} with $k\geq 3$, then there is no exact symplectic cobordism having $(M_0,\xi_0)$ as the convex end, and $(M,\xi_k)$ as the concave end.
\end{thm}


In particular, we obtain

\begin{cor}\label{corGT}
If $Y$ is the unit cotangent bundle of a hyperbolic surface, and $(M=\Sigma \times Y,\xi_k)$ is the corresponding $5$-dimensional contact manifold of Theorem \ref{thm5d} with $k\geq 3$, then $(M,\xi_k)$ does not have Giroux torsion.
\end{cor}

In the other direction of Theorem \ref{AlgvsGiroux}, examples of contact 3-manifolds which have $1$-torsion but not Giroux torsion where constructed in \cite{LW}. In higher dimensions, we have fairly strong geometric reasons for the following:

\begin{conj}\label{1tc} The examples of Corollary \ref{corGT} have untwisted algebraic 1-torsion (for \emph{any} $k\geq 1$).
\end{conj}

For $k\geq 3$, this would give an infinite family of contact manifolds with no Giroux torsion and algebraic 1-torsion in dimension $5$.

The interesting thing about this conjecture is that it builds on the relationship between SFT and string topology as discussed in \cite{CL09}. There are certain non-zero counts of punctured tori in the symplectization of $(Y,\alpha)$, where $\alpha$ is the standard Liouville form, which survive in $M$. They are given in terms of the coefficients of the Goldman--Turaev cobracket operation on strings in the underlying hyperbolic surface. One can find elements in the SFT algebra of $(Y,\alpha)$ whose differential in $\mathbb{R} \times Y$ has non-zero contributions from these tori and other pair of pants configurations, but by index considerations, only the former survive in $M$. This is how 1-torsion should arise. In fact, out of 35 possibilities, consisting of both cylinders and 1-punctured tori, this is the \emph{only} possible configuration from which 1-torsion can arise (see Section \ref{proooof}). We will therefore refer to it as \emph{sporadic}. This conjecture relies on being able to show the existence of obstruction bundles for certain building configurations, which is technically challenging. However, we provide a heuristic argument as to why we expect it to be true.

The above conjecture may be also given the following (also conjectural) interpretation: Theorem \ref{nocob1} would be then a result which is beyond the scope of algebraic torsion, since the presence of such a cobordism would only yield the already known fact that both ends have algebraic 1-torsion. In fact, in the proof of Theorem \ref{nocob1}, only holomorphic \emph{cylinders} play a role, whereas 1-punctured tori do not. In contrast, both are taken into consideration for $1$-torsion. This suggests the existence of an invariant more subtle than algebraic torsion, which we suspect would be obtained from the \emph{rational} SFT, rather than the whole SFT, and yet needs to be discovered. In this sense, in Theorem \ref{nocob1} we would morally be exploiting the properties of this hypothetical invariant. Also, the computations in \cite{LW} could also be interpreted in this way.

Putting Theorem \ref{thm5d} (and Remark \ref{mit}), together with Corollaries \ref{thm5d2} and \ref{corGT}, we obtain the following:

\begin{cor}\label{corG}
There exist infinitely many non-diffeomorphic $5$-dimensional contact manifolds $(M,\xi)$ which are tight, not strongly fillable, and which do not have Giroux torsion.
\end{cor}

To our knowledge, there are no other known examples of higher-dimensional contact manifolds as in Corollary \ref{corG}. Moreover, to our understanding, our proposed methods are so far the only way to show the non-existence of Giroux torsion in higher dimensions. Also, according to Conjecture \ref{1tc}, we expect the above examples to have algebraic 1-torsion. 

\vspace{0.5cm}

One can twist the contact structure of Theorem \ref{thm5d} close to the dividing set, by performing the \emph{$l$-fold Lutz--Mori twist} along a hypersurface $H$ lying in $\partial(\bigsqcup^k Y \times I \times S^1)$. This notion was defined in \cite{MNW}, and builds on ideas by Mori in dimension 5 \cite{Mori09}. It consists in gluing copies, along $H$, of a \emph{$2\pi l$-Giroux torsion domain} $(GT_l,\xi_{GT}):=(Y \times [0,2\pi l]\times S^1,\ker \lambda_{GT})$, the contact manifold obtained by gluing $l$ copies of $GT=GT_1$ together. In general, to be able to twist, $H\subseteq M$ is required to satisfy that it is transverse to $\xi$ and admits an orientation preserving identification with $S^1\times Y$, for some contact manifold $(Y,\xi_0)$, such that $\xi \cap TH = TS^1 \oplus \xi_0$. This is what the authors of \cite{MNW} call a \emph{$\xi$-round hypersurface modeled on $(Y,\xi_0)$}. The resulting contact structures are, in general, all homotopic as almost contact structures, but in our case they are distinguishable by a suitable version of cylindrical contact homology. This is proven in the same way as for the aforementioned hypertight examples in \cite{MNW}, which are also obtained by twisting (cf.\ Appendix \ref{LMtwists}). By construction, all of these have Giroux torsion, so by Theorem \ref{AlgvsGiroux} they have $\Omega$-twisted $1$-torsion, for $[\Omega] \in \mathcal{O}=\mathcal{O}(GT)$.

As a corollary of Theorem \ref{nocob1}, we get:

\begin{cor}\label{nocob}
Let $Y$ be the unit cotangent bundle of a hyperbolic surface, and let $(M=\Sigma \times Y,\xi)$ be the corresponding $5$-dimensional contact manifold of Theorem \ref{thm5d}, with $k\geq 3$. If $(M,\xi_l)$ denotes the contact manifold obtained by an $l$-fold Lutz--Mori twist of $(M,\xi)$, then there is no exact symplectic cobordism having $(M,\xi_l)$ as the convex end, and $(M,\xi)$ as the concave end (even though the underlying manifolds are diffeomorphic, and the contact structures are homotopic as almost contact structures).  
\end{cor}

In the proof of Theorem \ref{nocob1}, we make use of obstruction bundles, in the sense of Hutchings--Taubes \cite{HT1,HT2}. We prove they exist in our setup, under the condition that they exist inside the leaves of a codimension-2 holomorphic foliation. In the process of showing this condition for some particular examples, we get a result for super-rigidity of holomorphic cylinders in $4$-dimensional symplectic cobordisms, which might be of independent interest. This is a natural adaptation of the results of section 7 in \cite{Wen6} to the punctured setting. Recall that a somewhere injective holomorphic curve is super-rigid if it is immersed, has index zero, and the normal component of the linearized Cauchy--Riemann operator of every multiple cover is injective.  

\begin{thm}\label{superig} For generic $J$, every somewhere injective holomorphic cylinder in a $4$-dimensional symplectic cobordism, with index zero, and vanishing Conley--Zehnder indices (in some trivialization), is super-rigid.
\end{thm}

The original goal of this project was to generalize to higher dimensions the results by Latschev--Wendl, i.e.\ construct examples in any odd dimension with any given algebraic torsion. With this objective in mind, the candidates of Theorem \ref{thm5d} were suggested to us by C. Wendl, and are inspired by examples in \cite{MNW}. However, the following result remains open: 

\begin{conj}\label{conj}
For any $n\geq 1$, any $k \geq 1$, there exist a suitable $(2n-1)$-dimensional $Y$ for which the corresponding example of Theorem \ref{thm5d} satisfies that $AT(M_g,\xi_k)=k-1$ (for some genus $g$). 
\end{conj}

Let us remark that the $5$-dimensional examples discussed above were our original candidates for $AT=2$ (for $k=3$), but the existence of the sporadic configuration, discovered a posteriori, very likely prevents this. In this direction, we will develop general techniques for obtaining lower bounds on the order of torsion for any contact manifold, and illustrate their use in the proof of Theorem \ref{nocob1}. We will address this conjecture in Section \ref{conj14} below, where we state a lemma for the slightly more general setup of Section \ref{app}, which we hope will be of use when we pursue this conjecture in further work. We remark that, in all likelihood, the difficulty of settling it for the case $k=3$, at least in our approach, is comparable to that of Theorem \ref{nocob1}.



\begin{disclaimer}\label{disc} Since the statements of our results make use of machinery from Symplectic Field Theory, they come with the standard disclaimer that they assume that its analytic foundations are in place. They depend on the abstract perturbation scheme promised by the polyfold theory of Hofer--Wysocki--Zehnder. We shall assume that it is possible to achieve transversality by introducing an arbitrarily small abstract perturbation to the Cauchy-Riemann equation, and that the analogue of the SFT compactness theorem still holds as the perturbation is turned off. In practice, this means that, in order to study curves for the perturbed data, we need to also study holomorphic building configurations for the unperturbed one. However, we have taken special care in that the approach taken not only provides results that will be fully
rigorous after the polyfold machinery is complete, but also gives several direct results that are \emph{already} rigorous.
\end{disclaimer}

\subparagraph*{On tight and non-fillable contact manifolds.} Examples of tight but not fillable contact manifolds are rather well-known objects in dimension $3$. The first examples of ones which are tight but not weakly fillable can be found in \cite{EtHo}, which involves a combination of Giroux convex surface theory, and previous results by Lisca using Seiberg--Witten theory and Donaldson's theorem on the intersection form for 4-manifolds. Using that hypertightness implies tightness, and old work by Lutz, it is not hard to construct examples of such which do not go through SW technology. Examples in higher dimensions can be found in \cite{MNW}, and the underlying manifold is very similar to our manifold $M=Y\times \Sigma$, with the difference that both $\Sigma_\pm$ have genus $g$, and the contact structures are obtained by Lutz--Mori twists of an exactly fillable one. In fact, these examples are hypertight, and so it follows from theorem \ref{AlgvsGiroux} that they satisfy $AT=1$. 

\subparagraph*{On cylindrical Liouville semi-fillings.} Recall that for our results we have used the existence of Liouville domains of the form $Y\times I$, which we call ``cylindrical Liouville semi-fillings'', in every even dimension. The history of examples of such Liouville domains is interesting in itself. Their existence answers negatively a question of Calabi, which asks whether the fact from complex geometry, that the pseudoconvex boundary of a compact complex manifold is necessarily connected, still holds true in the symplectic category. They were first constructed in dimension 4 by McDuff in \cite{McD}, where $Y$ was taken to be the unit cotangent bundle of a hyperbolic surface, by using a suitable deformation of the standard symplectic form by a magnetic field 2-form. We have used this example in our construction of $5$-dimensional models. Further examples in dimension 4, which generalize McDuff's, and also dimension 6 were found by Geiges in \cite{Geiges1} and \cite{Geiges2}. There, $Y$ is a compact left-quotient of certain Lie group $G$ by a lattice, and the exact symplectic form is constructed using pairs of left-invariant contact forms on $G$ which satisfy a suitable algebraic condition, which makes them into what has been named a ``Geiges pair''. The same set of ideas is used by the 4-dimensional examples by Y.Mitsumatsu in \cite{Mit}. The examples are the same as in \cite{Geiges1} but from a slightly different point of view. The emphasis is put on unimodular 3-dimensional Lie algebras and a notion of a bilinear and symmetric pairing called the \emph{linking pairing}, such that the \emph{self linking} of a 1-form is exactly a measure of the failure of the contact condition to hold (cf.\ Definition \ref{LKparinig} below). It is indefinite only for $sl(2,\mathbb{R})$ and $solv$, the Lie algebras corresponding to the $SL(2,\mathbb{R})$ and $Solv$ geometries respectively, according to Thurston's geometrization. The construction of such Liouville domains in any even dimension was carried out in \cite{MNW}, using, again, Lie groups and left-invariant 1-forms, and in a wider context. The price to pay in is that now there is some algebraic number theory involved (fundamentally, Dirichlet's unit theorem), and the main body of work lies in showing the existence of co-compact lattices in certain Lie groups of affine transformations. 

\subparagraph*{Some of the difficulties.} There are several technical difficulties along the way, which come mainly from the lack of higher-dimensional tools/techniques. One of them is proving transversality for the curves in the foliation, so that indeed one has a space of isolated curves to count. A standard technique in dimension 3 is to use automatic transversality \cite{Wen1}, which consists in checking a fairly straightforward numerical inequality involving topological data associated to the curve. In higher dimensions, things become cumbersome. The way we approached this difficulty was to prove first regularity for suitable Morse-Bott data (which has more symmetries), from which regularity for sufficiently close non-degenerate (Morse) data follows, by the implicit function theorem. In the former case, we explicitly compute the normal component of the linearization of the Cauchy--Riemann equation, after making suitable choices of coordinates and Morse-Bott function, which hold without loss of generality. By splitting the vector bundles in question into sums of invariant line bundles for which one can use automatic transversality, we manage to compute precisely the elements in the kernel of the normal operator. They correspond to nearby holomorphic curves in a finite energy foliation. We obtain that the dimension of the kernel coincides with the Fredholm index of the operator, from which regularity follows.

Another difficulty is the lack of a higher-dimensional intersection theory between punctured holomorphic \emph{curves}, in the sense of Siefring. In dimensions 3 and 4, this is a good tool for proving that holomorphic curves with certain prescribed asymptotics are unique, and therefore one knows exactly what to count. In higher dimensions, we resorted to varied methods, the combination of which, to our knowledge, is new. We solved this difficulty by using the SOBD structure suitably to obtain a branched cover, together with energy estimates to reduce the problem to dimension 4 (which draws inspiration from \cite{BvK10}), as well as holomorphic cascades and gluing in the sense of \cite{Bo} (Theorem \ref{uniqueness}). 

Moreover, we need to understand the number of ways certain holomorphic building configurations may glue to honest curves, which we intend to count, after making $J$ generic. For this, we make use of obstruction bundles, a fairly non-trivial gadget to deal with in practice. The idea is to count the number of honest curves obtained by gluing buildings, by algebraically counting the zero set of a section of an obstruction bundle. While we will not explicitly compute those numbers (presumably a fairly hard thing to accomplish), and just prove existence results in suitable cases, the symmetries in the setup imply that there are cancellations which can be exploited to obtain our results. Under the assumption that obstruction bundles exist leaf-wise, we show their existence in our examples by careful analytical considerations, exploiting the fact that curves lies in hypersurfaces of the foliation.

\section*{Guide to the document}

In Section \ref{Liouvilledoms}, we look at examples of cylindrical Liouville semi-fillings $Y \times I$, both in dimension 3  and higher. The main result (Corollary \ref{explicitJ}) involves the construction of explicit compatible almost complex structures that are cylindrical away from a small neighbourhood of the \emph{non-contact central slice} $Y \times \{0\}$, to which the Liouville vector field is tangent. These almost complex structures have the nice (but rather incidental) property that certain hypersurface inside this slice becomes holomorphic. This will serve as a building block for the construction of our model contact forms of further sections, the contact manifolds of Theorem \ref{thm5d}, and perhaps are suitable for constructions of holomorphic curves in the Liouville completion of $Y\times I$, something which might be needed in further work.

\vspace{0.5cm}

The main construction is dealt with in two main sections (\ref{modelAsec} and \ref{RulingOut}), which correspond to the building of models A and B. Model A is used to prove Theorem \ref{thm5d}, done in Section \ref{Torsion}, and model B was originally intended to be used to give lower bounds on the algebraic torsion. However, we use it to prove theorem \ref{nocob1}, in Section \ref{5dmodel}. We discuss obstruction bundles in Section \ref{obbundles}. 

\vspace{0.5cm}

The proof of Theorem \ref{AlgvsGiroux} is dealt with in Section \ref{GirouxTorsion1-tors}, which is basically a reformulation of the previous sections, with the key input being an adaptation of the uniqueness Theorem \ref{uniqueness}.
  
\vspace{0.5cm}

We include three Appendixes: Appendix \ref{LMtwists} describes the Lutz--Mori twists, and applies it to our model A, to obtain the contact structures of Corollary \ref{nocob}. Appendix \ref{SOBDap} contains a definition of a SOBD in higher dimensions, which is adapted to the examples we have dealt with in this document. These SOBDs will be a useful tool for proving Theorem \ref{AlgvsGiroux}. The last appendix, Appendix \ref{AppSiefring}, is written in co-authorship with Richard Siefring, and gives a basic outline of his intersection theory for punctured holomorphic curves and hypersurfaces. We prove the aforementioned application to codimension-2 holomorphic foliations, needed in the main section of the model B construction (Section \ref{curvesvshyp}), and used also in Section \ref{prequantcurves} to prove a uniqueness result for curves over certain prequantization spaces.

\newpage

\vspace*{\fill}
\begin{center}
\begingroup
\centering 
\large \textbf{Acknowledgements}
\endgroup 
\end{center} 

First of all, my thanks go to my supervisor, Chris Wendl, for introducing me to this project and for his support and patience throughout its duration. To Richard Siefring, for very helpful conversations and for co-authoring an appendix in this document. To Patrick Massot, Sam Lisi, Jonny Evans, Yanki Lekili, and Janko Latschev, for helpful conversations/correspondence on different topics. To Gabriel and Miguel Paternain, for supporting and believing in me, and providing a role-model, during my times in Uruguay and Cambridge. To friend and colleague Momchil Konstantinov, for endless and very enlightening conversations about maths, books, the lifting of heavy things, and general nonsense, and for putting up with me and my (occasional?) fits of insanity. To the rest of the London crew: Tobias Sodoge, my German bro, thanks for sitting next to me and not end up killing me; Brunella Torricelli, the Swiss lady of the many languages, thanks for implicitly teaching me German and Italian; Emily Maw, my favourite English-woman, and my source of gossip, thanks for being you; Navid Nabijou, thanks for all the maths shared and for having the hair that you have; Jacqui Espina, thanks for helping out when my Conley-Zehnder indices were running amiss, and for great home-made chocolate; Oldrich Spacil, from whom I learned the \cite{BEM} paper. To my Berlin office mate Alex Fauck, for interesting maths conversations, but more importantly for opening our door when I forget my keys (which happens more often than I care to admit), and lending me money when I forget my wallet (which luckily happens less often). To the rest of the Berlin Symplectic group: Klaus Mohnke, Felix Schm\"{a}schke (for also lending me money when Alex is not around), Viktor Fromm. To Daniel \'{A}lvarez-Gavela, for sharing maths and music, but mostly for driving me around California in Ernesto, \emph{Il Furgone Presto}. To all of those involved in the Kylerec workshop, specially the organizers, which helped make it the success that it was. To \'{A}lvaro del Pino, for bearing with me back when I was trying to understand the painful details behind the classification of overtwisted contact 3-manifolds. To my friends of our Cambdrige-born collective chamber music ensemble, String Theory: Franca, Lay, Will, Philip \& Teresa, Jumpei, Jozefien, Kin; let us continue to reunite and play music wherever our fingers land on the map. To all of those whom I am forgetting and we crossed paths along the way, thank you.     

And last but not least, to my family and friends back home, specially my mom and brother, whose support and love has always been their gift to me, even though they have no clue of what these pages can possibly be used for (save, perhaps, fuel for a good BBQ).

This research has been partly carried out in (and funded by) University College London (UCL) in the UK, and by the Berlin Mathematical School (BMS) in Germany.

\vspace*{\fill}

\newpage

\paragraph*{Basic notions}\label{basicn}

A contact form in a $(2n+1)$-dimensional manifold $M$ is a $1$-form $\alpha$ such that $\alpha \wedge d\alpha^{n}$ is a volume form, and the associated contact structure is $\xi=\ker \alpha$ (we will assume all our contact structures are co-oriented). The Reeb vector field associated to $\alpha$ is the unique vector field $R_\alpha$ on $M$ satisfying $$\alpha(R_\alpha)=1,\;i_{R_\alpha}d\alpha=0$$

A $T$-periodic Reeb orbit is $(\gamma,T)$ where $\gamma: \mathbb{R}\rightarrow M$ is such that $\dot{\gamma}(t)=TR_\alpha(\gamma(t))$, $\gamma(1)=\gamma(0)$. We will often just talk about a Reeb orbit $\gamma$ without mention to $T$, called its period, or action. If $\tau>0$ is the minimal number for which $\gamma(\tau)=\gamma(0)$, and $k \in \mathbb{Z}^+$ is such that $T=k\tau$, we say that the covering multiplicity of $(\gamma,T)$ is $k$. If $k=1$, then $\gamma$ is said to be simply covered (otherwise it is multiply covered). A periodic orbit $\gamma$ is said to be non-degenerate if the restriction of the time $T$ linearised Reeb flow $d\varphi^T$ to $\xi_{\gamma(0)}$ does not have $1$ as an eigenvalue. More generally, a Morse--Bott submanifold of $T$-periodic Reeb orbits is a closed submanifold $N \subseteq M$ invariant under $\varphi^T$ such that $\ker(d\varphi^T - \mathds{1})=TN$, and $\gamma$ is Morse--Bott whenever it lies in a Morse--Bott submanifold, and its minimal period agrees with the nearby orbits in the submanifold. The vector field $R_\alpha$ is non-degenerate/Morse--Bott if all of its closed orbits are non-degenerate/Morse--Bott.

A stable Hamiltonian structure (SHS) on $M$ is a pair $\mathcal{H}=(\Lambda,\Omega)$ consisting of a closed $2$-form $\Omega$ and a $1$-form $\Lambda$ such that 

$$\ker \Omega \subseteq \ker d\Lambda,\;\mbox{and}\;\Omega\vert_{\xi} \mbox{ is non-degenerate, where } \xi=\ker \Lambda$$ 

In particular, $(\alpha,d\alpha)$ is a SHS whenever $\alpha$ is a contact form. The Reeb vector field associated to $\mathcal{H}$ is the unique vector field on $M$ defined by $$\Lambda(R)=1,\; i_{R}\Omega=0$$ There are analogous notions of non-degeneracy/Morse--Bottness for SHS.

A symplectic form in a $2n$-dimensional manifold $W$ is a $2$-form $\omega$ which is closed and non-degenerate. A Liouville manifold (or an exact symplectic manifold) is a symplectic manifold with an exact symplectic form $\omega=d\lambda$, and the associated Liouville vector field $V$ is defined by the equation $i_{V}d\lambda=\lambda$. Any Liouville manifold is necessarily open. A boundary component $M$ of a Liouville manifold (endowed with the boundary orientation) is convex if the Liouville vector field is positively transverse to $M$, and is concave, if it is so negatively. An exact cobordism from a (co-oriented) contact manifold $(M_+,\xi_+)$ to $(M_-,\xi_-)$ is a compact Liouville manifold $(W,\omega=d\alpha)$ with boundary $\partial W=M_+\bigsqcup M_-$, where $M_+$ is convex, $M_-$ is concave, and $\ker \alpha\vert_{M_\pm}=\xi_\pm$. Therefore, the boundary orientation induced by $\omega$ agrees with the contact orientation on $M_+$, and differs on $M_-$. A Liouville filling (or a Liouville domain) of a --possibly disconnected-- contact manifold $(M,\xi)$ is a compact Liouville cobordism from $(M,\xi)$ to the empty set. A strong symplectic cobordism and a strong filling are defined in the same way, with the difference that $\omega$ is exact only in a neighbourhood of the boundary of $W$ (so that the Liouville vector field is defined in this neighbourhood, but not necessarily in its complement).

The symplectization of a contact manifold $(M,\xi=\ker \alpha)$ is the symplectic manifold $(\mathbb{R}\times M, \omega=d(e^a\alpha))$, where $a$ is the $\mathbb{R}$-coordinate. In particular, it is a non-compact Liouville manifold. Similarly, the symplectization of a stable Hamiltonian manifold $(M,\Lambda,\Omega)$ is the symplectic manifold $(\mathbb{R}\times M, \omega^\varphi)$, where $\omega^\varphi=d(\varphi(a)\Lambda)+\Omega$, and $\varphi$ is an element of the set $$\mathcal{P}=\{\varphi \in C^\infty(\mathbb{R},(-\epsilon,\epsilon)):\varphi^\prime>0\}$$ Here, $\epsilon>0$ is chosen small enough so that $\omega^\varphi$ is indeed symplectic. An $\mathcal{H}$-compatible (or simply cylindrical) almost complex structure on a symplectization $(W=\mathbb{R}\times M,\omega^\varphi)$ is $J \in \mbox{End}(TW)$ such that $$J \mbox{ is } \mathbb{R}\mbox{-invariant}, J^2=-\mathds{1},\;J(\partial_a)=R,\;J(\xi)=\xi,\;J\vert_{\xi} \mbox{ is }\Omega \mbox{-compatible }$$ The last condition means that $\Omega(\cdot,J\cdot)$ defines a $J$-invariant Riemannian metric on $\xi$. If $J$ is $\mathcal{H}$-compatible, then it is easy to check that it is $\omega^\varphi$-compatible, which means that $\omega^\varphi(\cdot,J\cdot)$ is a $J$-invariant Riemannian metric on $\mathbb{R}\times M$.

To any closed $T$-periodic Reeb orbit $(\gamma,T)$ one can associate an asymptotic operator $\mathbf{A}_\gamma$. To write it down, choose a symmetric connection $\nabla$ on $M$, and a $\mathcal{H}$-compatible almost complex structure $J$, and define
$$
\mathbf{A}_\gamma=\mathbf{A}_{\gamma,J}:W^{1,2}(\gamma^*\xi)\rightarrow L^2(\gamma^*\xi)
$$
$$
\mathbf{A}_\gamma\eta =-J(\nabla_t \eta -T\nabla_\eta R)
$$

Alternatively, one has the expression 
$$
\A_{\gamma}\eta(t)=-J\left.\frac{d}{ds}\right|_{s=0}d\varphi^{-Ts}\eta(t+s),
$$
for $\eta \in W^{1,2}(\gamma^*\xi)$, where $\varphi^s$ again is the time-$s$ Reeb flow.

Morally, this is the Hessian of a certain action functional on the loop space of $M$ whose critical points correspond to closed Reeb orbits. It is symmetric with respect to a suitable $L^2$-product. A periodic orbit $\gamma$ is non-degenerate if and only if $0$ does not lie in the spectrum of $\mathbf{A}_\gamma$, and more generally, if $\gamma$ is Morse--Bott and lies in a Morse--Bott submanifold $N$, then $\dim \ker \mathbf{A}_\gamma=\dim N - 1$. Under a choice of unitary trivialization $\tau$ of $\gamma^*\xi$, this operator looks like 
$$
\mathbf{A}_\gamma:W^{1,2}(S^1,\mathbb{R}^{2n})\rightarrow L^2(S^1,\mathbb{R}^{2n})
$$
$$
\mathbf{A}_\gamma =-i\partial_t-S(t),
$$
where $S$ is a smooth loop of symmetric matrices (the coordinate representation of $-TJ\nabla R$), which comes associated to a trivialization of $\mathbf{A}_\gamma$. When $\gamma$ is non-degenerate, its Conley--Zehnder index with respect to $\tau$ is defined to be the Conley--Zehnder index of the path of symplectic  matrices $\Psi(t)$ satisfying $\dot{\Psi}(t)=iS\Psi(t)$, $\Psi(0)=\mathds{1}$. We denote this by $\mu_{CZ}^\tau(\gamma)=\mu_{CZ}(\Psi)$.

We will consider, for cylindrical $J$, punctured $J$-holomorphic curves $u:(\dot{\Sigma},j)\rightarrow (\mathbb{R}\times M,J)$ in the symplectization of a stable Hamiltonian manifold $M$, where $\dot{\Sigma}=\Sigma \backslash \Gamma$, $(\Sigma,j)$ is a compact connected Riemann surface, and $u$ satisfies the nonlinear Cauchy--Riemann equation $du \circ j = J \circ u$. We will also assume that $u$ is asymptotically cylindrical, which means the following. Partition the punctures into \emph{positive} and \emph{negative} subsets $\Gamma = \Gamma^+\cup \Gamma^-$, and at each $z \in \Gamma^\pm$, choose a biholomorphic identification of a punctured neighborhood of $z$ with the half-cylinder $Z_\pm$, where $Z_+ = [0, \infty) \times S^1$ and $Z_- = (-\infty, 0] \times S^1$. Then writing $u$ near the puncture in cylindrical coordinates $(s, t)$, for $|s|$
sufficiently large, it satisfies an asymptotic formula of the form 
$$u \circ\phi(s, t) = exp_{(Ts,\gamma(T t))} h(s, t)$$
Here $T > 0$ is a constant, $\gamma : \mathbb{R} \rightarrow M$ is a $T$-periodic Reeb orbit, the exponential map is defined with respect to any $\mathbb{R}$-invariant metric on $\mathbb{R} \times M$, $h(s, t) \in \xi_{\gamma(T t)}$ goes to $0$ uniformly in $t$ as $s \rightarrow \pm \infty$ and $\phi : Z_\pm \rightarrow Z_\pm$ is a smooth embedding such that $\phi(s, t) - (s + s_0, t + t_0) \rightarrow 0$ as $s \rightarrow \pm \infty$ for some constants $s_0 \in \mathbb{R}$, $t_0 \in S^1$. We will refer to punctured asymptotically cylindrical $J$-holomorphic curves simply as $J$-holomorphic curves.

Observe that, for any closed Reeb orbit $\gamma$ and cylindrical $J$, the trivial cylinder over $\gamma$, defined as $\mathbb{R}\times \gamma$, is $J$-holomorphic. 

The Fredholm index of a punctured holomorphic curve $u$ which is asymptotic to non-degenerate Reeb orbits in a $(2n+2)$-dimensional symplectization $W^{2n+2}=\mathbb{R} \times M$ is given by the formula 

\begin{equation} \label{index} 
ind(u)=(n-2)\chi(\dot{\Sigma})+2c_1^\tau(u^*TW)+\mu^\tau_{CZ}(u)
\end{equation} 

Here, $\dot{\Sigma}$ is the domain of $u$, $\tau$ denotes a choice of trivializations for each of the bundles $\gamma_z^*\xi$, where $z \in \Gamma$, at which $u$ approximates the Reeb orbit $\gamma_z$. The term $c_1^{\tau}(u^*TW)$ is the relative first Chern number of the bundle $u^*TW$. In the case $W$ is $2$-dimensional, this is defined as the algebraic count of zeroes of a generic section of $u^*TW$ which is asymptotically constant with respect to $\tau$. For higher-rank bundles, one determines $c_1^\tau$ by imposing that $c_1^\tau$ is invariant under bundle isomorphisms, and satisfies the Whitney sum formula (see e.g.\ \cite{Wen8}). The term $\mu^\tau_{CZ}(u)$ is the total Conley--Zehnder index of $u$, given by 

$$
\mu^\tau_{CZ}(u)=\sum_{z \in \Gamma^+}\mu^\tau_{CZ}(\gamma_z) - \sum_{z \in \Gamma^-}\mu^\tau_{CZ}(\gamma_z)
$$ 

Given a $\mathcal{H}$-compatible $J$, and a $J$-holomorphic curve $u$ in $\mathbb{R}\times M$, the expression $u^*\Omega$ is a non-negative integrand, and one can define its $\Omega$-energy $$\mathbf{E}(u)=\int u^*\Omega$$ It is non-negative, and vanishes if and only if $u$ is a (multiple cover of) a trivial cylinder.

\section{Cylindrical Liouville semi-fillings revisited}\label{Liouvilledoms}

\begin{figure}[!ht]\centering
\includegraphics[width=0.92\linewidth]{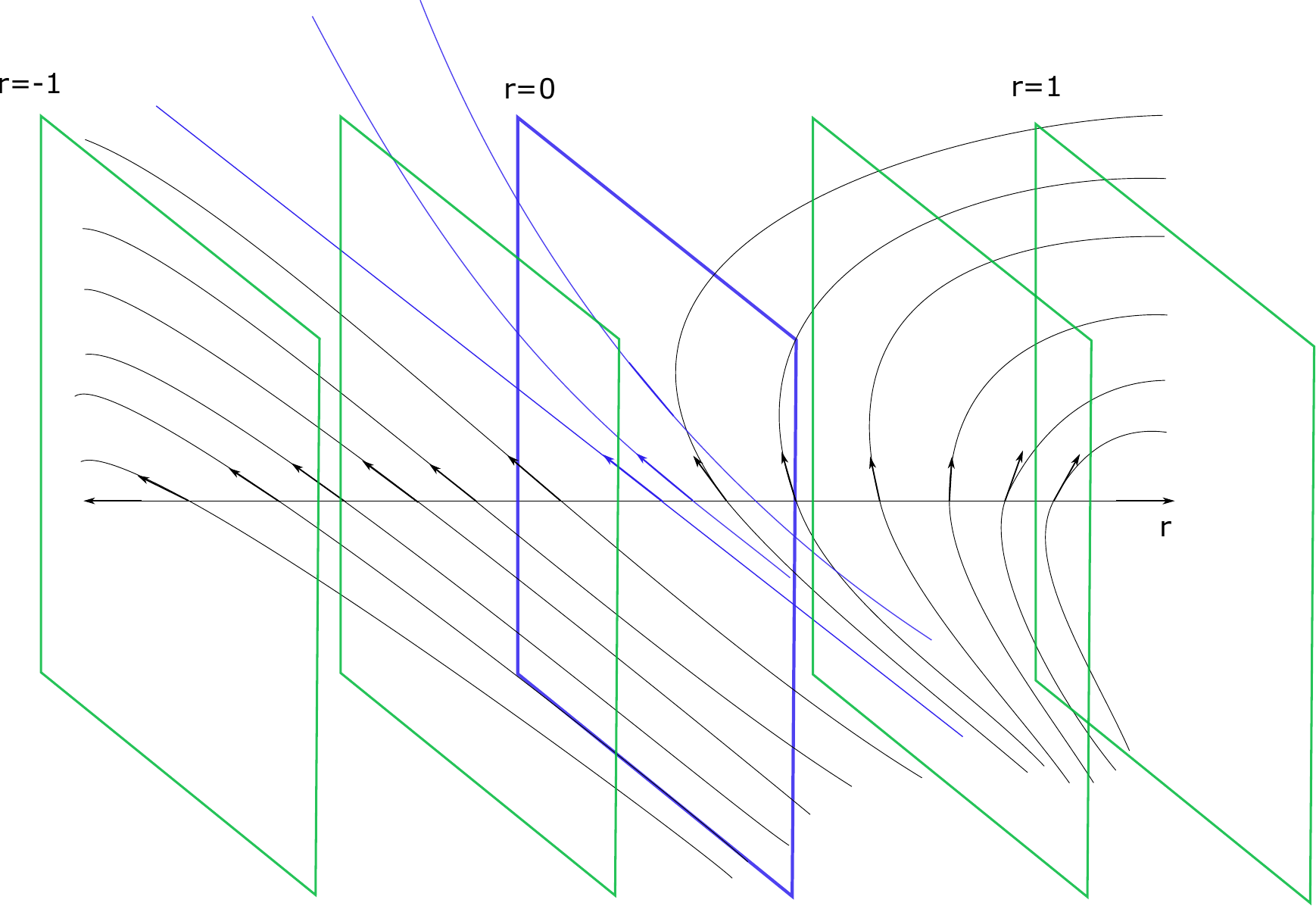}
\caption{\label{LouvilleDomain} The picture illustrates the qualitative behaviour of the flow of the Liouville vector field $V$ on any cylindrical Liouville semi-filling, for which the central slice ($r=0$) is invariant. If the coefficient in the $\partial_r$ component of $V$ depends only on $r$ (as in the examples here), one has a diffeomorphism between $Y \times [-1,0)$ and $Y \times (0,1]$ with the negative symplectizations $\mathbb{R}^{\leq 0} \times Y_-$ and $\mathbb{R}^{\leq 0} \times Y_+$ respectively. Indeed, if $\varphi_t^{V}$ denotes the time $t$ flow of $V$, one can relate $t$ (away from $r=0$) to the interval parameter $r$ by diffeomorphisms $f_-:(-\infty,0]\rightarrow [-1,0)$ and $f_+: (-\infty,0] \rightarrow (0,1]$ defined by $\varphi_t^{V}(Y \times \{-1\})=\{f_-(t)\} \times Y_-$ and $\varphi_t^{V}(Y \times \{1\})=\{f_+(t)\}\times Y_+$, satisfying $\alpha_{f_\pm(t)}=e^{t}\alpha_\pm$. One may therefore informally think of such a Liouville domain as being obtained by gluing two negative symplectizations along a ``non-contact hypersurface''.}
\end{figure}

Let us begin by revisiting some of the simplest of the low dimensional cylindrical Liouville semi-fillings in a very concrete way, following Mitsumatsu's exposition. This helps to illustrate some aspects and difficulties of the setup we will be dealing with, and, perhaps more importantly, gain some intuition (for which we direct the reader's attention to Figure \ref{LouvilleDomain}).

Based on the $3$-dimensional notion of \cite{Mit}, we introduce some piece of terminology:

\begin{defn}\label{LKparinig} Given a compact manifold $Y^{2n-1}$, and a volume form $dvol \in \Omega^{2n-1}(Y)$, for $\alpha, \beta \in \Omega^1(Y)$, we will denote by $link(\alpha,\beta) \in C^{\infty}(M)$ the function defined by the equation
$$
\alpha \wedge d\beta^{n-1}=link(\alpha,\beta)dvol,
$$
and we will refer to it as the \emph{linking function} between $\alpha$ with $\beta$. The case of interest is when $\alpha=\beta$, in which we denote the \emph{self-linking function} by $link(\alpha):=link(\alpha,\alpha)$.
\end{defn}  

\begin{example} Let us look first at the 3-dimensional Lie algebra $solv$, which has basis $X,W,T$ with the relations $[T,X]=X$, $[T,W]=-W$, $[X,W]=0$. This corresponds to manifolds of the form $Y \times I$ for $Y$ a compact quotient of the unique simply connected Lie group corresponding to $solv$ (which may be constructed as a semidirect product $\mathbb{R}\rtimes \mathbb{R}^2$ for a suitable action of $\mathbb{R}$ in $\mathbb{R}^2$ \cite{MNW,Mit,Geiges1}). Let $X^*,T^*,W^*$ denote the dual basis in $solv^*$, which may be viewed as left-invariant differential forms on $Y$ dual to the left-invariant vector fields $X,T,W$ under the usual evaluation pairing. Recalling the formula $d\omega(V,W)=-\omega([V,W])$ for $\omega$ a left-invariant 1-form and $V,W$ left-invariant vector fields, we obtain the equations 
$$dX^*=-T^*\wedge X^*,\; dW^*=T^*\wedge W^*,\;dT^*=0$$

If we let $\alpha_r=X^*+rW^*$ for $r \in I$, we then have that
$$d\alpha_r=X^*\wedge T^* +rT^*\wedge W^*$$ and therefore
$$\alpha_r\wedge d\alpha_r=2rX^*\wedge T^* \wedge W^*,$$ so that, the self-linking function of $\alpha_r$ is given by $link(\alpha_r)=2r$. This means that $\alpha_r$ is a contact form for the contact structure 
\begin{equation}\label{xirsolv}
\xi_r=\ker\alpha_r=\langle rX-W,T\rangle,
\end{equation}
for every $r \neq 0$, where we take this basis to be positive so that $d\alpha_r$ is a positive/negative symplectic form for non-zero $r$. This depends on the sign of $r$, since $d\alpha_r(rX-W,T)=2r$. Its Reeb vector field is $$R_r=\frac{1}{2}\left(X+\frac{W}{r}\right),$$ which explodes to order 1 at $r=0$. In this case $\xi_0=\langle T,W\rangle$ integrates to a foliation, by Frobenius' theorem, which has $X$ as a trivialisation of the transverse direction. This can in fact be interpreted as an Anosov splitting, where the Anosov flow is generated by $T$, the stable/unstable foliations are $\langle T,W \rangle$ and $\langle T, W\rangle$, and $\xi_\pm=\xi_{\pm 1}$ form a $\pi/4$ angle with the foliation $\xi_0$---see \cite{ETh98,Mit}. 

If we now view this as a 1-form $\alpha=\{\alpha_r\}_{r\in I}$ on $Y \times I$ (which satisfies $\alpha(\partial_r)=0$), we obtain that the symplectic form is given by $$d\alpha=dr\wedge \frac{\partial\alpha_r}{\partial r} + d\alpha_r=dr\wedge W^*-T^*\wedge X^* + rT^*\wedge W^*$$ 

The associated Liouville vector field is $V=2r\partial_r-T$, which becomes tangent to the central slice $Y \times \{0\}$, and hence this is not a contact-type hypersurface. It is manifestly positively/negatively transverse to the slices $Y_\pm:=Y\times \{\pm 1\}$. 

It is interesting to see what happens when $H:Y \times I \rightarrow \mathbb{R}$ is chosen to be a Hamiltonian depending only on $r$. In this case, the slices $Y_r:=Y \times \{r\}$ are level sets for $H$ and hence, for $r\neq 0$, the Reeb vector field $R_r$ and $X_H$ are colinear, spanning the characteristic line field of this contact type hypersurface. Indeed, we have the formula $R_r=X_H/\alpha(X_H)$, which even makes sense outside of the set of critical points of $H$ which are different from zero, since $R_r$ does. Even though $R_r$ is not well-defined at $r=0$, we have that $X_H$ always is, and thus, if $H$ is taken to have no critical points, $R_r$ gives a well-defined \emph{line field} at $r=0$. We shall refer to it as the \emph{characteristic line field of $Y_0$}. In this case, it is spanned by $W$. We can check that $$X_H=H^\prime(r)(rX+W),$$ and $$\alpha(X_H)=2rH^\prime(r),$$ where we stick to the sign convention $i_{X_H}d\alpha=-dH$ for Hamiltonian vector fields. 

An interesting feature is that $\alpha_{-r}(R_r)=0$ for every $r \in I$, i.e.\ the Reeb vector field of $\xi_r$ actually lies in the mirror $\xi_{-r}$, so that when $r \rightarrow 0$ the vector fields $R_r$ and $R_{-r}$ approach $\xi_0$ whilst pointing to opposite directions. One also has that the 1-form $i_{R_r}d\alpha_{-r}=T^*$ is independent of $r$, so that in particular the equation $i_{R_r}d\alpha_{-r}=i_{R_{-r}}d\alpha_r$ holds. Also, observe that $\xi_r$ is a trivial bundle, and so $c_1(\xi_r)=0$ for any choice of compatible $J$.

\end{example}

\begin{example}
\label{sl2ex}

Another related example, which will be important for Section \ref{5dmodel}, is when the Lie algebra is $sl(2,\mathbb{R})$, having a basis $h,l,k$ with brackets $[h,l]=-k, [h,k]=-l,[l,k]=h$. Compact quotients of $SL(2,\mathbb{R})$ include unit cotangent bundles of hyperbolic surfaces, which can be written as $PSL(2,\mathbb{R})/\Gamma$, where $\Gamma$ is identified with the fundamental group of the surface. 

We obtain $$dh^*=k^*\wedge l^*,\; dl^*=h^*\wedge k^*,\; dk^*=h^*\wedge l^*$$ Letting $\alpha_r=\frac{1}{2}((1+r)k^*+(1-r)h^*)$, so that $$\xi_r=\langle l, (1-r)k-(1+r)h\rangle,$$ where the basis is taken to be positive, one sees that $link(\alpha_r)=r$, and $$d\alpha=\frac{1}{2} (dr\wedge k^*-dr \wedge h^*+(1+r)h^*\wedge l^*+(1-r)k^*\wedge l^*)$$
We then find $$V=2r\partial_r-l,\; R_r=\frac{1}{2r}((r-1)h+(r+1)k),$$ and, for $H$ a $Y$-independent Hamiltonian, that $$X_H=H^{\prime}(r)((r-1)h+(r+1)k)$$$$\alpha(X_H)=2H^\prime(r)r$$  One also has the equation $\alpha_{-r}(R_r)=0$, and the fact that $i_{R_r}d\alpha_{-r}=l^*$ is $r$-independent, as can be seen explicitly, as well as $c_1(\xi_r)=0$ for any compatible $J$. The characteristic line field of $Y_r$ at $r=0$ is spanned by $k-h$.

\end{example}

The following is a way of getting cylindrical Liouville semi-fillings in higher dimensions:

\begin{defn}\label{LP}\cite{MNW} A \emph{Liouville pair} on an oriented $(2n-1)$-dimensional manifold $Y$ is $(\alpha_+,\alpha_-)$ a pair of contact forms such that $\pm \alpha_\pm \wedge d\alpha_\pm^{n-1}>0$, and the $1$-form 
$$
\beta:=\frac{1}{2}(e^{-s}\alpha_-+e^s\alpha_+)
$$
on $\mathbb{R}\times Y$ is a Liouville form (i.e.\ satisfies $d\beta^n>0$).
\end{defn}

\begin{example}
We now have a look now at the higher-dimensional examples, and do some explicit computations. In \cite{MNW}, the authors provide 1-forms $\alpha_\pm$ in $\mathbb{R}^{2n+1}$, given by the formulas
$$
\alpha_\pm = \pm e^{\sum_{i=1}^n t_i} d\theta_0 + e^{-t_1}d\theta_1+\dots + e^{-t_n}d\theta_n,
$$
with coordinates $(\theta_0,\dots, \theta_n,t_1,\dots, t_n) \in \mathbb{R}^{2n+1}$, such that the pair $(\alpha_+,\alpha_-)$ is a Liouville pair on $\mathbb{R}^{2n+1}$. By taking the orientation-reversing diffeomorphism 
$$
\varphi: (0,\pi) \rightarrow \mathbb{R}
$$
$$
\varphi(r)=ln\left(\frac{1+\cos(r)}{\sin(r)}\right)
$$
we obtain the formula 
\begin{equation}
\begin{split}
\alpha:=&\varphi^*\beta=\frac{1}{2}\left(\frac{1+\cos(r)}{\sin(r)} \alpha_+ + \frac{1-\cos(r)}{\sin(r)}\alpha_-\right)\\
&=\frac{1}{\sin(r)}\left(e^{\sum_{i=1}^n t_i}\cos(r)d\theta_0+e^{-t_1}d\theta_1+\dots+e^{-t_n}d\theta_n\right),
\end{split}
\end{equation}
which is a Liouville form in $Y \times (0,\pi)$, whenever $Y$ is a compact quotient of $\mathbb{R}^{2n+1}$ to which $\alpha_\pm$ descends to. Finding such quotient is the main work behind this story in the aforementioned paper. We see from the expression for $\alpha$ that it is also given by a 1-parameter family $\alpha=\{\alpha_r\}$. It actually does not extend to $Y \times [0,\pi]$ unless we multiply by $\sin(r)$. This is what is called, in Giroux terminology, an \emph{ideal Liouville domain} (see Definition \ref{idLD} in Section \ref{GirouxTorsion1-tors}). The singular slice is then $Y_{\pi/2}$, to which the Liouville vector field becomes tangent. Since the Liouville form is so explicit, we may thus compute all the associated data, which will be of use later when defining an almost complex structure.

One can check that  

\begin{equation}
\begin{split}
\alpha_r \wedge d\alpha_r^n&=2(-1)^{\frac{n(n+1)}{2}}\frac{\cos(r)n!}{\sin(r)^{n+1}}d\theta_0 \wedge dt_1 \wedge \dots \wedge dt_n\wedge d\theta_1 \wedge \dots \wedge d\theta_n 
\end{split}
\end{equation}

So indeed the only non-contact type slice happens at $r=\pi/2$, so that $d\alpha_r\vert_{\xi_r}$ is a symplectic form for every other value of $r$, where 
\begin{equation}
\begin{split}
\xi_r&=\ker\alpha_r\\
&= \langle \partial_{t_1}, \dots, \partial_{t_n}, e^{-t_i}\partial_{\theta_1}-e^{-t_1}\partial_{\theta_i} \;( i=2,\dots,n), e^{-t_1}\partial_{\theta_0}-e^{\sum_{j=1}^n t_j}\cos(r)\partial_{\theta_1} \rangle
\end{split}
\end{equation}
Observe that again this is a trivial bundle, so $c_1(\xi_r)=0$ for any compatible $J$.

\vspace{0.5cm}

A routine computation shows that the Reeb and Liouville vector fields are

\begin{equation}
\label{LouvReeb}
\begin{split}
R_r&=\frac{\sin(r)}{n+1}\left(\frac{e^{-\sum_{i=1}^n t_i}}{\cos(r)}\partial_{\theta_0}+e^{t_1}\partial_{\theta_1}+\dots+ e^{t_n}\partial_{\theta_n}\right)\\
V&= - \frac{\sin(r)}{\sin(r)^2+(n+1)\cos(r)^2}\left(\sin(r)\sum_{i=1}^n \partial_{t_i}+(n+1)\cos(r)\partial_r\right)
\end{split}
\end{equation} 
We then clearly see that $R_r$ explodes at $r=\pi/2$, whereas $V$ is everywhere defined,  and is tangent to $Y_{\pi/2}$. $R_r/\sin(r)$ and $V/\sin(r)$ extend to $Y \times [0,\pi]$ as non-vanishing vector fields, and the latter coincides with the outwards-pointing vector field $\pm \partial_r$ along the boundary. One may also check that $$i_{R_r}d\alpha_{\pi-r}=\frac{2}{n+1}\sum_{i=1}^ndt_i,$$ which is $r$-independent. The characteristic line field of $Y_r$ is spanned by $cos(r)R_r$, which extends to $r=\pi/2$, being spanned by $\partial_{\theta_0}$.

\vspace{0.5cm}

Observe that all of the above computations hold for the compact quotients $Y$ of $\mathbb{R}^{2n+1}$ which inherit the Liouville pair $\alpha_\pm$, which will be $\mathbb{T}^{n+1}$-bundles over $\mathbb{T}^n$ where the coordinates in the fibers are $(\theta_0,\dots,\theta_n)$, and the coordinates in the base are $(t_1,\dots,t_n)$. 

\end{example}

\subparagraph*{Almost complex structures adapted to cylindrical semi-fillings}

We now turn to constructing explicit and manageable $d\alpha$-compatible almost complex structure on cylindrical Liouville semi-fillings, which will serve for future constructions and computations. 

\begin{prop}
\label{explicitJ}
Given a cylindrical Liouville semi-filling $(Y\times I,d\alpha)$, for any small $\delta>0$, we can construct an oriented hyperplane distribution $\widetilde{\xi}_r$ on $Y_r=Y\times \{r\}$, depending smoothly on $r$, and satisfying:

\begin{itemize}
\item $\widetilde{\xi}_r$ coincides with $\xi_r$ as distributions for $r$ in the complement of $(-\delta,\delta)$, with orientation which agrees for $r>\delta$, and differs for $r<-\delta$. Here $\xi_r$ has the  orientation which coincides or differs by the one induced by contact form $\alpha_r$, depending on the sign of $r$.
\item $d\alpha_r\vert_{\widetilde{\xi}_r}$ is a positive symplectic form for every $r \in I$. 

\end{itemize}

Moreover, in all the examples discussed above, we further have:
\begin{itemize}
\item $\widetilde{\xi}_0$ is integrable.
\end{itemize}

\end{prop}

\begin{proof}

The distribution $\widetilde{\xi}_r$ is defined, in general, as the $d\alpha$-symplectic complement in $Y_r$ of its characteristic line field. The first two claims are then immediate.

For the \emph{solv} case, one can show $\widetilde{\xi}_0=\langle X,T\rangle$, which is integrable by Frobenius' theorem, since $[X,T]=-X$ lies in it.







For the $sl(2,\mathbb{R})$ case, one has $\widetilde{\xi}_0=\langle h + k,l \rangle$, which is integrable since $[l,h+k]=h+k$.




For the higher-dimensional examples, $\widetilde{\xi}_0= \langle \partial_{\theta_1},\dots, \partial_{\theta_n}, \partial_{t_1},\dots,\partial_{t_n} \rangle$,  which is again integrable. 






\end{proof} 
 
\begin{cor}
\label{compJ0}
Given a cylindrical Liouville semi-filling $(Y \times I,d\alpha)$, for any small $\delta>0$, it is possible to make a choice of $d\alpha$-compatible almost complex structure $J_0$, so that 
\begin{itemize}
\item $J_0$ is compatible with the symplectization structure away form $(-\delta,\delta)$, that is, it preserves $\xi_r$, where it is $d\alpha_r$-compatible, it is invariant under the flow of $V$, and maps $V$ to $R_r$.
\item There exist globally defined smooth vector fields $\widetilde{V}$ and $\widetilde{R}_r$, which satisfy:
\begin{enumerate}
\item $\widetilde{V}$ and $\widetilde{R}_r$ span $\widetilde{\xi_r}^{d\alpha}$, the $d\alpha$-symplectic complement of $\widetilde{\xi}_r$. 
\item They coincides with $\mp V$ and $\mp R_r$ in the complement of $(-\delta,\delta)$ (the upper sign corresponds to the region $\{r\geq \delta\}$, and the lower to $\{r \leq -\delta\}$).
\item $d\alpha(\widetilde{V},\widetilde{R}_r)=1$.

\item $dr(\widetilde{V})<0$.

\item $\widetilde{R}_r$ and $R_r$ are colinear. More explicitly,  $\widetilde{R}_r=\alpha(\widetilde{R}_r)R_r$, where the function $\alpha(\widetilde{R}_r)$ vanishes at $r=0$ and equals $\mp 1$ in the complement of $(-\delta,\delta)$. In particular, $dr(\widetilde{R}_r)=0$.

\item $J_0(\widetilde{V})=\widetilde{R}_r$.
\end{enumerate}
\item $c_1(Y\times I,J_0)=0$ (which implies that $c_1(Y\times I,J)=0$ for \emph{any} $d\alpha$-compatible $J$).
\end{itemize}

Moreover, in all the examples here discussed, we can arrange:

\begin{itemize}
\item $\alpha(\widetilde{V})=0$.
\item There is a $J_0$-holomorphic hypersurface contained in $Y \times \{0\}$.  
\end{itemize}
\end{cor}
\begin{proof} We write $\widetilde{\xi}_r^{d\alpha}=\langle \widetilde{V}, \widetilde{R}_r\rangle$, which can be chosen as in the statement. We choose $J_0$ to preserve $\widetilde{\xi}_r$ where it is $d\alpha_r$-compatible (which is the same as being $d\alpha$-compatible since $dr(\widetilde{\xi}_r)=0$), so that it also preserves its symplectic complement $\widetilde{\xi_r}^{d\alpha}$, mapping $\widetilde{V}$ to $\widetilde{R}_r$. We see that 
$$
c_1(Y\times I,J_0)=c_1(\widetilde{\xi}_r)+c_1(\langle \widetilde{V},\widetilde{R}_r\rangle)=0,
$$
since $\langle \widetilde{V},\widetilde{R}_r\rangle$ is a trivial line bundle, and $c_1(\widetilde{\xi}_r)=0$. The condition $dr(\widetilde{V})<0$ follows from the computations in Remark \ref{remark1} below. The last condition comes from the integrability of $\widetilde{\xi}_0$ in the examples we have seen. 

\vspace{0.5cm}

To see that we can take $\widetilde{V}$ such that $\alpha(\widetilde{V})=0$ in the examples discussed, we simply compute.

In the \emph{solv}-case, we may take 

$$ \widetilde{R}_r= \frac{1}{2\beta_0(r)r+\beta_1(r)}(rX+Y), \;\widetilde{V}=(2\beta_0(r)r+\beta_1(r))\partial_r-\beta_0(r)T,
$$ which satisfy the desired conditions. 

In the $sl(2,\mathbb{R})$-case,
we take

$$\widetilde{R}_r= \frac{1}{2\beta_0(r)r+\beta_1(r)}((r-1)h+(r+1)k), \;\widetilde{V}=(2\beta_0(r)r+\beta_1(r))\partial_r-\beta_0(r)l,
$$

For the higher-dimensional case, one may take 
$$
\widetilde{V}= \Phi\bigg\{\left(\frac{\cos(r)}{\sin(r)}\beta_1(r)+e^{\sum_{i= 1}^nt_i} \beta_0(r)\sin(r)\right)\sum_{i=1}^n \partial_{t_i}
$$
$$
+ \left((n+1)\beta_0(r)e^{\sum_{i= 1}^nt_i}\cos(r)-\beta_1(r)\right)\partial_r\bigg\}
$$
where the non-vanishing function $\Phi$ has the expression
$$
\Phi=-\frac{\sin(r)}{e^{\sum_{i=1}^n t_i} (\sin(r)^2+(n+1)\cos(r)^2)},
$$

and

$$
\widetilde{R}_r=\frac{1}{\Psi}\left(\partial_{\theta_0}+\cos(r)e^{\sum_{i=1}^nt_i}\sum_{i=1}^n e^{t_i}\partial_{\theta_i}\right)
$$
where 
$$
\Psi=\frac{1}{\sin(r)}\left((n+1)\beta_0(r)e^{\sum_{i= 1}^nt_i}\cos(r)-\beta_1(r)\right)
$$

\end{proof}

Let us make some further technical observations:

\begin{remark} \label{remark1} 
\begin{enumerate}[wide, labelwidth=!, labelindent=0pt] $\;$
\item If $H:Y \times I \rightarrow \mathbb{R}$ is a Hamiltonian which only depends on $r$, we can show that 
\begin{equation}\label{gradham}
\nabla H=dr(\widetilde{V}) H^\prime \widetilde{V},\; X_H= dr(\widetilde{V}) H^\prime \widetilde{R}_r, 
\end{equation}
where the gradient is the one corresponding to the metric $g=d\alpha(\cdot,J_0 \cdot)$. In particular, $\nabla H$ always has a non-vanishing $\partial_r$-component away from its critical points, where the sign coincides with that of $H^\prime$, and, in our examples, we have 
\begin{equation} \label{grad}
\alpha(\nabla H)=0
\end{equation}

\item Using that the splitting $T(Y \times I)=\widetilde{\xi}_r \oplus \langle \widetilde{V},\widetilde{R}_r\rangle$ is symplectic, we may write 

$$d\alpha=d\alpha_r\vert_{\widetilde{\xi}_r} + d\widetilde{V}\wedge d\widetilde{R}_r=d\alpha_r + dr \wedge \frac{\partial \alpha_r}{\partial r}$$

We then obtain that 
$$d\widetilde{V}=-i_{\widetilde{R}_r}d\alpha=\frac{\partial\alpha_r}{\partial r}(\widetilde{R}_r)dr,$$
where we have used that $\widetilde{R}_r=\alpha(\widetilde{R}_r)R_r$.

Since $\partial_r=\frac{\widetilde{V}}{dr(\widetilde{V})}+W$, for some $W \in \widetilde{\xi}_r$, we have that $d\widetilde{V}(\partial_r)=\frac{1}{dr(\widetilde{V})}$, so plugging $\partial_r$ into the above equation we get 

\begin{equation} 
\label{hey}
\frac{\partial \alpha_r}{\partial r}(R_r)=\frac{1}{\alpha(\widetilde{R}_r)} \frac{\partial \alpha_r}{\partial r}(\widetilde{R}_r)=\frac{1}{dr(\widetilde{V})\alpha(\widetilde{R}_r)},
\end{equation} which has a pole at $r=0$. The function $\frac{\partial \alpha_r}{\partial r}(R_r)$ is positive for $r>0$, and negative for $r<0$, and so it follows that $dr(\widetilde{V})<0$.

\item The vector field $R_r$ has a well-defined pole at $r=0$, so that $R_r$ may be written in a neighbourhood of $r=0$ as $R_r=\frac{\widehat{R_r}}{r^p}$ where $\widehat{R_0}$ is nowhere-vanishing.

For any Hamiltonian function $H:Y \times I \rightarrow \mathbb{R}$ depending only on $r$, we have $R_r=\frac{X_H}{\alpha(X_H)}$ for $r\neq 0$, by recalling that $Y_r$ is a contact-type hypersurface. In particular, if we take $H$ so that it does not have critical points, we have that $\alpha(X_H)$ has to have a zero at $r=0$ with some order $p$. This order is independent of $H$, since $R_r$ is. Moreover, we have $i_{X_H}d\alpha_r=0$, which is clearly true for $r \neq 0$, and so true everywhere by continuity. Therefore, we obtain
$$
H^\prime dr = dH=-i_{X_H}d\alpha=-i_{X_H}\left(d\alpha_r+dr\wedge \frac{\partial \alpha_r}{\partial r}\right)=\frac{\partial \alpha_r}{\partial r}(X_H)dr,
$$
so that $$H^\prime = \frac{\partial \alpha_r}{\partial r}(X_H),$$ and therefore 
\begin{equation}\label{somexpr}
\frac{\partial\alpha_r}{\partial r}(R_r)= \frac{H^\prime}{\alpha(X_H)},
\end{equation}
and this expression is independent of $H$. In particular, this shows that $\frac{\partial\alpha_r}{\partial r}(R_r)$ has a pole at $r=0$ of the same order as $R_r$, which by (\ref{hey}) coincides with that of $1/\alpha(\widetilde{R}_r)$.

\end{enumerate}
\end{remark}

\section{Model A. Achieving torsion}\label{modelAsec}

\subsection{Construction of model A}
\label{model}

In this section, we shall start by introducing model A on the underlying manifold $M$ (to be specified below in detail), which will make use of the cylindrical Liouville semi-fillings discussed in the previous section. We will use the ``double completion'' construction, originally appearing in \cite{LVHMW}, which has the effect of endowing $M$ with a contact structure and an explicit deformation to a SHS, by viewing it as a contact-type hypersurface in a non-compact Liouville manifold. It will also be apparent that our model carries a higher-dimensional version of a spinal open book decomposition, which in dimension $3$ is the subject of \cite{LVHMW} (see appendix B for a definition). This observation suggests that perhaps the setup of many things in this paper can be carried over to an abstract setting of a manifold carrying a more general spinal open book decomposition, but we have decided to stay modest. 

As is usual in explicit constructions, the contact form thus obtained will be degenerate, in this case along the \emph{spine}, and so a standard Morse function technique as in \cite{Bo} will be necessary if we want to compute anything in SFT, for which nondegeneracy is crucial. As usual, our function will be defined in the space of Reeb orbits, which can be identified with $Y \times I$. One has two approaches here: either choose a Morse function on this space (which makes the contact structure automatically non-degenerate along this region and such that low-action Reeb orbits correspond to its critical points), or choose a Morse--Bott function only depending on the interval parameter $r$ first (which throws away only one degree of degeneracy), and then perturb again by a Morse function on $Y$. The result in the end, in terms of the foliation by holomorphic curves that we get, will be the same as long as we choose a Morse function of the form $H(y,r)=f(r)+\gamma(r)g(y)$ for $g$ and $f$ Morse in $Y$ and $I$ respectively, and $\gamma$ a bump function supported in a small neighbourhood of $r=0$. The difference is that the Morse--Bott approach simplifies the proof of uniqueness and regularity significantly, and the transition from Morse--Bott to Morse is understood via the results in \cite{Bo}.

Another technicality we have to deal with is the smoothening of some corners, in order that the Liouville vector field in the double completion is actually transverse to $M$. This is a rather standard technique, which has the unfortunate effect of making the expressions for the contact form and associated data rather unattractive. 

\vspace{0.5cm}

Without further ado, let $Y$ be a closed $(2n-1)$-manifold such that $(Y \times I,d\alpha)$ is a Liouville domain, for some exact symplectic form $d\alpha \in \Omega^2(Y \times I)$ (recall that throughout this document, $I$ will denote the interval $[-1,1]$). We will assume that the Liouville form $\alpha=\{\alpha_r\}_{r \in I}$ is given by a 1-parameter family of $1$-forms in $Y$. In particular, we get that $\alpha(\partial_r)=0$. We can write the symplectic form as $$d\alpha=d\alpha_r + dr\wedge \frac{\partial \alpha_r}{\partial r}$$ The Liouville vector field $V$, defined to be $d\alpha$-dual to $\alpha$, points outwards at each boundary component, and hence, using its flow, we can choose our coordinate $r \in I$ so that $V$ agrees with $\pm \partial_r$ near the boundary $\partial(Y \times I)=Y \times \{\pm 1\}=:Y_\pm$. Therefore, we can assume that $\alpha= e^{\pm r-1}\alpha_\pm$ on $Y \times [-1,-1+\delta)$ and $Y \times (1-\delta,1]$, respectively, for some small $\delta>0$. Then $Y_\pm$ carries a contact structure $\xi_\pm=\ker\alpha_\pm$, where $\alpha_\pm = i_Vd\alpha|_{T(Y_\pm)}=\alpha|_{T(Y_\pm)}$. Observe that the behaviour of $V$ near the ends necessarily implies that there are values $r \in I$ such that $V|_{Y\times \{r\}}$ lies in $TY$, and hence $Y_r:=Y \times \{r\}$ is not a contact type hypersurface. This means that $(Y_\pm, \xi_\pm)$ will in general not be contactomorphic (perhaps $Y_-$ and $Y_+$ are not even \emph{homeomorphic}, as observed before), which justifies the notation. The slices $Y_r$ which are of contact type inherit a contact structure $\xi_r=\ker \alpha_r$ and the resulting Reeb vector field $R_r$ satisfies $R_r=e^{1\mp r}R_\pm$ in the respective components of $\{|r| >1-\delta\}$, where $R_\pm$ is the Reeb vector field of $\alpha_\pm=\alpha_{\pm 1}$. We shall assume throughout that the only non-contact type slice is $Y_0$, so that $\alpha_r$ is a contact form for every $r \neq 0$. Also, we shall make the convention that whenever we deal with equations involving $\pm$'s and $\mp$'s, one has to interpret them as to having a different sign according to the region (the ``upper'' sign denotes the ``plus'' region, and the ``lower'', the ``minus'' region).

\vspace{0.5cm}

Let now $M= Y \times \Sigma$ be a product $(2n+1)$-manifold, where $\Sigma$ is the orientable genus $g$ surface obtained by gluing a connected genus $0$ surface with $k$ boundary components $\Sigma_-$, to a connected genus $g-k+1>0$ surface with $k$ boundary components $\Sigma_+$ along the boundary, by an orientation preserving map. The surface $\Sigma$ then inherits the orientation of $\Sigma_-$, which is opposite to the one in $\Sigma_+$. On each boundary component of $\partial \Sigma_\pm$, choose collar neighbourhoods $\mathcal{N}(\partial \Sigma_\pm)=(-\delta,0]\times S^1$ (for the same $\delta$ as before), and coordinates $(t_\pm,\theta_\pm) \in \mathcal{N}(\partial \Sigma_\pm)$, so that $\partial \Sigma_\pm=\{t_\pm = 0\}$. 

In order to use the setup above, we will consider $\Sigma_-$ and $\Sigma_+$ to be attached at each of the $k$ boundary components by a cylinder $I\times S^1$, so that $M$ at this region is the disjoint union of $k$ copies of $Y \times I\times S^1$, with the $Y \times I$ identified with the Liouville domain above. We can then write the points of $M$ here as $(y,r,\theta)$, where the $\theta \in S^1$ coordinate can be chosen to coincide with $\theta_\pm$ where the gluing takes place. We shall therefore drop the subscript $\pm$ when talking about the $\theta$ coordinate. Denote also $$\mathcal{N}(-Y):=Y \times [-1,-1+\delta)\times S^1,\; \mathcal{N}(Y):=Y \times (1-\delta,1]\times S^1,$$ in the above identification. 

Therefore we may write our manifold $M$ as $$M=M_Y \cup M_P^{\pm},$$ where $M_Y=\bigsqcup^k Y \times I\times S^1$ is a region gluing $M_P^{\pm}=Y_\pm \times \Sigma_{\pm}$ together (recall Figure \ref{themanifold}). We shall refer to them as the \emph{spine} or \emph{cylindrical region}, and the \emph{positive/negative paper}, respectively. Observe that we have fibrations $$\pi_Y: M_Y \rightarrow Y \times I$$$$ \pi_P^\pm: M_P^{\pm} \rightarrow Y_\pm,$$ with fibers $S^1$ and $\Sigma_\pm$, respectively, and hence can be given the structure of a SOBD (see Appendix B for a definition). We shall use this fact to construct a suitable contact structure on $M$, and a stable Hamiltonian structure arising as a deformation of this contact structure. 

\vspace{0.5cm}

We construct now the following $(2n+2)$-manifold: $$E=[-1,-1+\delta) \times M_P^- \;\bigsqcup\; (1-\delta,1] \times M_
P^+ \;\bigsqcup\; (-\delta,0]\times  M_Y/\sim,$$ where we identify $(r,y,t_-,\theta) \in [-1,-1+\delta) \times Y \times \mathcal{N}(\partial \Sigma_-)$ with $(t,y,r,\theta)\in (-\delta,0] \times \mathcal{N}(-Y)$ if and only if $t=t_-$, and $(r,y,t_+,\theta) \in (1-\delta,1] \times Y \times \mathcal{N}(\partial \Sigma_+)$ with $(t,y,r,\theta)\in (-\delta,0] \times \mathcal{N}(Y)$ if and only if $t=t_+$.

We shall also construct an open manifold, which can be considered to be an open completion of $E$, as follows. Denote by $\Sigma_\pm^\infty$ the open manifolds obtained from $\Sigma_\pm$ by attaching cylindrical ends of the form $(-\delta,+\infty) \times S^1$ at each boundary component, where the subset $(-\delta,0] \times S^1$ coincides with the collar neighbourhoods chosen above. The coordinates $t_\pm$ and $\theta$ extend to these ends in the obvious way, and we shall refer to the cylindrical ends as $\mathcal{N}(\partial \Sigma_\pm^\infty)$. We also consider the cylinder $\mathbb{R} \times S^1$ obtained by enlarging the cylindrical region $I\times S^1$ we had above. Denote then $$M_P^{\pm,\infty}=Y \times \Sigma_\pm^\infty$$$$M_Y^\infty= Y \times \mathbb{R} \times S^1$$$$ \mathcal{N}^\infty(-Y)=Y\times (-\infty, -1+\delta)\times S^1,\;\mathcal{N}^\infty(Y)=Y \times (1-\delta,+\infty) \times S^1$$ and define the \emph{double completion} of $E$ to be $$E^{\infty,\infty}=(-\infty,-1+\delta) \times M_P^{-,\infty} \;\bigsqcup\; (1-\delta,+\infty) \times M_
P^{+,\infty}\; \bigsqcup \; (-\delta,+\infty)\times M_Y^\infty /\sim,$$ where we identify $(r,y,t_-,\theta)\in (-\infty,-1+\delta)\times Y \times \mathcal{N}(\partial\Sigma_-^\infty)$ with $(t,y,r,\theta) \in (-\delta,+\infty) \times \mathcal{N}^\infty(-Y)$ if and only if $t=t_-$, and $(r,y,t_+,\theta)\in (1-\delta,+\infty)\times Y \times \mathcal{N}(\partial\Sigma_+^\infty)$ with $(t,y,r,\theta) \in (-\delta,+\infty) \times \mathcal{N}^\infty(Y)$ if and only if $t=t_+$  (see Figure \ref{doublecompletion2}). Observe that, by definition, the $t$ coordinate coincides with the $t_\pm$ coordinates, where these are defined, so we shall again drop the $\pm$ subscripts from the variables $t_\pm$. Note also that the $r$ coordinate is globally defined, whereas $t$ is not. Denote then by $E^{\infty,\infty}(t)$ the region of $E^{\infty,\infty}$ where the coordinate $t$ is defined.

Choose now $\lambda_\pm$ to be Liouville forms on the Liouville domains $\Sigma^\infty_\pm$, such that $\lambda_\pm=e^td\theta$ on $\mathcal{N}(\partial \Sigma_\pm^\infty)$. Then observe that this last expression makes sense in the region of $E^{\infty,\infty}$ where both $\theta$ and $t$ are defined, and where they are not, the form $\lambda_\pm$ makes sense. So this yields a globally defined 1-form $\lambda \in \Omega^1(E^{\infty,\infty})$ which coincides with $\lambda_\pm$ where these are defined. Also, the same argument works for $\alpha$, since it looks like $e^{r-1}\alpha_+$ in $\{r> 1-\delta\}$, and like $e^{-r-1}\alpha_-$ in $\{r<-1+\delta\}$, and in $\{|r|\leq 1-\delta\}$ we can just use the expression $\alpha$, so that we get a global $\alpha \in \Omega^1(E^{\infty,\infty})$. More concretely, we have that the manifold $E^{\infty,\infty}$ ``fibers'' (as a trivial ``fibration'', but with fibers which change topology) over the symplectic base $(Y \times \mathbb{R},d\alpha)$. Its symplectic fibers which can be identified with $(\Sigma_\pm^\infty, d\lambda_\pm)$ on the suitable components of the region $\{|r|\geq 1-\delta\}$, and with $([-\delta,+\infty)\times S^1 ,d(e^td\theta))$ on the region $\{|r|< 1-\delta\}$. This global $\lambda$ may be considered as a fiberwise Liouville form for this ``fibration'', whereas the form $\alpha \in \Omega^1(E^{\infty,\infty})$ is just the pullback of the form $\alpha \in \Omega^1(M_Y^\infty)$ under the projection map $E^{\infty,\infty}\rightarrow M_Y^\infty$. 

For $K\gg 0$ a big constant, $\epsilon>0$ a small one, and $L\geq 1$, choose a smooth function $$\sigma=\sigma^{L}_{\epsilon,K}:\mathbb{R}\rightarrow \mathbb{R}^+$$ satisfying
\begin{itemize}
\item $\sigma \equiv K$ on $\mathbb{R}\backslash[-L,L]$.
\item $\sigma \equiv \epsilon$ on $[-L+\delta,L-\delta]$.
\item $\sigma^\prime(r)<0$, for $r\in (-L,-L+\delta)$, $\sigma^\prime(r) > 0$ for $r\in(L-\delta,L)$. 
\end{itemize}

We have that the 1-form $\sigma \alpha$ is Liouville on $Y \times \mathbb{R}$. Indeed, if $dvol$ is a positive volume form in $Y$ with respect to the $\alpha_-$-orientation, we may write 
$$
d\alpha^n=dvol\wedge dr,\; \alpha_r\wedge d\alpha_r^{n-1}=link(\alpha_r)dvol,
$$
where the sign of the \emph{self-linking} function $r \mapsto link(\alpha_r)$ (Def.\ \ref{LKparinig}) is opposite to that of $r \in \mathbb{R}$. Observe that the first equation holds in $Y \times \mathbb{R}$, and the second one in $Y_r$.
We then compute that
$$
d(\sigma\alpha)^n=\sigma^{n-1}(\sigma - n \sigma^\prime link(\alpha_r))dvol\wedge dr
$$
Since the only region where $\sigma^\prime$ is non-zero is $\{L>|r|>L-\delta\}$, and on the two distinct components of this region its sign is is opposite with that of $link(\alpha_r)$, we have that the above expression is positive.

Moreover, one has that the associated Liouville vector field is given by 
\begin{equation} \label{Vsig}
V_\sigma:=\frac{\sigma}{\sigma + \sigma^\prime dr(V)}V=\left\{\begin{array}{ll} V,& \mbox{ on } (\mathbb{R}\backslash (-L,L))\cup [-L+\delta,L-\delta]\\
\frac{\sigma}{\sigma + \sigma^\prime}\partial_r,& \mbox{ on } (L-\delta,L)\\-\frac{\sigma}{\sigma - \sigma^\prime}\partial_r,& \mbox{ on } (-L,-L+\delta)
\end{array}\right.
\end{equation} 

Observe that $V_\sigma$ is everywhere \emph{positively} colinear with $V$. 

After extending the form $\sigma \alpha$ to $E^{\infty,\infty}$ in the natural way, we can use Thurston's trick \cite{Th} for the case of a trivial fibration to conclude that the form $$\lambda_\sigma:=\lambda_\sigma^{L}:=\sigma\alpha + \lambda$$ is a Liouville form on $E^{\infty,\infty}$, that is, the 2-form $$\omega_\sigma:=\omega_\sigma^{L}:=d\lambda_\sigma=\sigma^\prime dr \wedge \alpha + \sigma d\alpha + d\lambda$$ is symplectic. 

We then can define a Liouville vector field $X_\sigma$, given by the formula $$i_{X_\sigma}\omega_\sigma=\lambda_\sigma$$ If $X_\pm$ denotes the Liouville vector field on $\Sigma_\pm^\infty$ which is $d\lambda_\pm$-dual to $\lambda_\pm$, coinciding with $\partial_t$ in $\mathcal{N}(\partial \Sigma_\pm^\infty)$, we can define a smooth vector field on $E^{\infty,\infty}$ by 
\begin{equation}
\label{X}
X=\left\{\begin{array}{ll}
  X_+,& \mbox{ on } \{r>1-\delta\} \\
  X_-, &\mbox{ on } \{r<-1+\delta\} \\
  \partial_t,& \mbox{ on } \{|r|<1\}\\
  \end{array}\right.
\end{equation}
We then have that $$X_\sigma= X+V_\sigma,$$ as is straightforward to check.  

Denote now by $$E^{L,Q}=E^{\infty,\infty}/(\{|r|>L\} \cup \{t > Q\}),$$ for $Q \geq 0$, and $L \geq 1$. We have its ``horizontal'' and ``vertical'' boundaries $$\widetilde{M}_Y^{L,Q}:=\partial_hE^{L,Q} := \{t=Q\} \cap \{r \in [-L,L]\}$$$$\widetilde{M}_P^{L,Q}:=\partial_vE^{L,Q} := \widetilde{M}_P^{-,L,Q} \bigsqcup \widetilde{M}_P^{+,L,Q},$$ where $$\widetilde{M}_P^{-,L,Q} :=\{r=-L\}\cap \{t \leq Q\}$$$$\widetilde{M}_P^{+,L,Q}:= \{r=L\}\cap \{t \leq Q\}$$More rigourously, these regions should actually be understood to contain the region of $\{r=\pm L\}$ where $t$ is not defined, which are copies of $\Sigma_\pm/\mathcal{N}(\partial\Sigma_\pm)\times Y$, respectively. We write them as above for the sake of simplicity. The manifold $$\widetilde{M}^{L,Q}:=\partial E^{L,Q}:= \partial_h E^{L,Q} \cup \partial_v E^{L,Q},$$ is then a manifold with corners $$\partial_h E^{L,Q} \cap \partial_v E^{L,Q}=\{|r|=L\}\cap \{t=Q\}$$  
  
The above expression for $X_\sigma$ implies that $$X_\sigma=\pm\frac{\sigma}{\sigma \pm \sigma^\prime}\partial_r + \partial_t$$ in the corresponding components of the region $\{|r|>L-\delta\}\cap \{t\geq -\delta\}$. This means that $X_\sigma$ will be transverse to the smoothening of $\partial E^{L,Q}$ that we shall now construct, so that this smoothened manifold will be a contact type hypersurface in $(E^{\infty,\infty},\omega_\sigma)$.

\vspace{0.5cm}

\begin{figure}[t] \centering \includegraphics[width=0.60\linewidth]{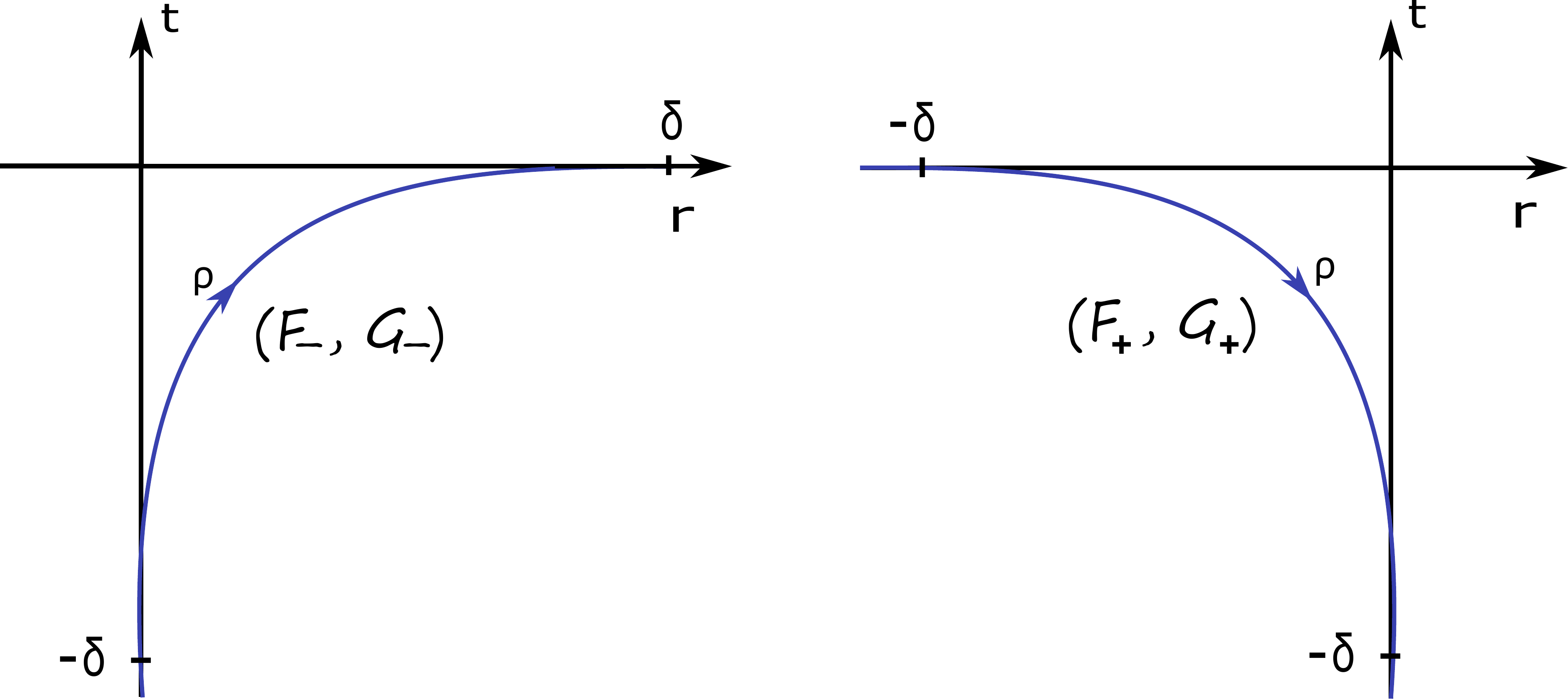}
\caption{\label{SmoothendCorner} This picture shows the paths $\rho \mapsto (F_\pm(\rho),G_\pm(\rho))$, for our choice of smoothening functions for the corners.}
\end{figure}

\begin{figure}[t]\centering \includegraphics[width=0.68\linewidth]{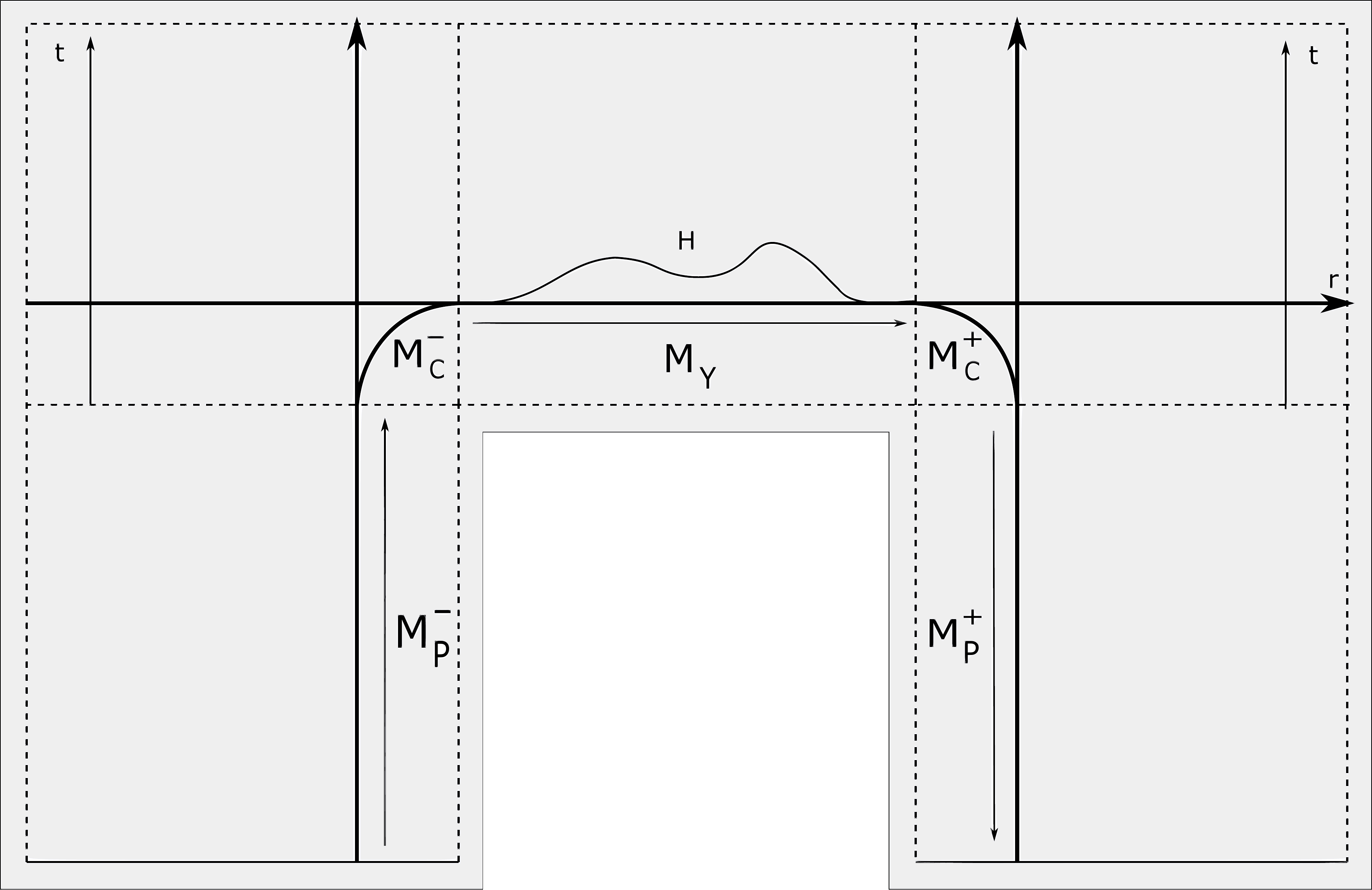}
\caption{\label{doublecompletion2} The double completion, and a Morse function $H$ along the spine.}
\end{figure}

Choose smooth functions $F_\pm,G_\pm:(-\delta,\delta)\rightarrow (-\delta,0]$ such that $$\left\{\begin{array}{ll}
 (F_+(\rho),G_+(\rho))=(\rho,0),\;(F_-(\rho),G_-(\rho))=(0,\rho),& \mbox{ for } \rho \leq -\delta/3\\
 G_+^\prime(\rho)<0, F_-^\prime(\rho)>0,& \mbox{ for } \rho > -\delta/3\\
 G_-^\prime(\rho)>0, F_+^\prime(\rho)>0,& \mbox{ for } \rho < \delta/3\\
 (F_+(\rho),G_+(\rho))=(0,-\rho),\;(F_-(\rho),G_-(\rho))=(\rho,0),& \mbox{ for } \rho \geq \delta/3\\
\end{array}\right.$$ See Figure \ref{SmoothendCorner}.

We now smoothen the corner $\partial_h E^{L,Q} \cap \partial_v E^{L,Q}$ by substituting the region $$\partial E^{L,Q} \cap \left(\{t \in (Q-\delta,Q]\}\cup \{|r|\in (L-\delta,L]\}\right),$$ which contains the corners, with the smooth manifold $$M_{C}^{\pm,L,Q}:=\{(r=\pm L+F_\pm(\rho), t = Q + G_\pm(\rho),y,\theta): (\rho,y,\theta) \in (-\delta,\delta) \times Y \times S^1\}$$

The smoothened boundary can then be written as $$M^{L,Q}:= M_P^{\pm,L,Q} \bigcup M_C^{\pm,L,Q} \bigcup M_Y^{L,Q},$$ where $$M_P^{\pm,L,Q} = \widetilde{M}_P^{\pm,L,Q} \cap \{t<Q-\delta/3\}$$$$M_Y^{L,Q}= \widetilde{M}_Y^{L,Q} \cap \{|r|<L-\delta/3\},$$ with the same caveat as above in the definitions of $\widetilde{M}_P^{\pm,L,Q}$. The Liouville vector field $X_\sigma$ is transverse to this manifold, so that we get a contact structure on $M^{L,Q}$ given by $$\xi^{L,Q}= \ker(\lambda^{L}_\sigma|_{TM^{L,Q}})$$ Observe that $M^{L,Q}$ is diffeomorphic to $M$. Geometrically, it can be thought of as a suitable deformation of the manifold obtained from ``stretching'' the $\Sigma$ component of $M$ at both boundary components of the cylindrical region, by an $L-1$ amount, and the collar neighbourhoods, by a $Q$ amount. Therefore, this actually yields a contact structure on our original manifold $M$. Observe that, by construction, we have non-empty intersections $$M_P^{\pm,L,Q}\cap M_C^{\pm,L,Q} = \{r=\pm L\}\cap \{t \in (Q-\delta,Q-\delta/3)\}$$$$ M_Y^{L,Q}\cap M_C^{\pm,L,Q}=\{t=Q\}\cap \{|r| \in (L-\delta,L-\delta/3)\}$$

We shall construct a stable Hamiltonian structure on $M^{L,Q}$ which arises as a deformation of the above contact structure, such that both coincide on $M_Y^{L,Q}$, as follows.

Choose a smooth function $\beta: \mathbb{R}\rightarrow [0,1]$ such that $\beta(t)=0$ for $t \leq -\delta+\delta/9$, $\beta(t)=1$ for $t\geq -2\delta/3-\delta/9$, and $\beta^\prime \geq 0$. Set 
\begin{equation}\label{Z}
Z=\left\{\begin{array}{ll} V_\sigma+\beta(t)X,& \mbox{ in the region } E^{\infty,\infty}(t) \mbox{ where } t \mbox{ is defined}\\ 
\frac{\sigma}{\sigma+\sigma^\prime}\partial_r,& \mbox{ in } E^{\infty,\infty}(t)^c \cap \{r > 1-\delta\}\\
-\frac{\sigma}{\sigma-\sigma^\prime}\partial_r,& \mbox{ in } E^{\infty,\infty}(t)^c \cap \{r < -1+\delta\},
\end{array}\right. 
\end{equation}
which yields a smooth vector field on $E^{\infty,\infty}$, a deformation of $X_\sigma$. We have that $Z$ is still transverse to $M^{L,Q}$, and we claim that the pair $$\mathcal{H}:=\mathcal{H}^{L,Q}:=(\Lambda^{L,Q}=i_Z\omega_\sigma|_{TM^{L,Q}},\Omega^{L,Q}=\omega_\sigma|_{TM^{L,Q}})$$ yields a stable Hamiltonian structure on $M^{L,Q}$. Observe that in $M_Y^{L,Q}$, the above pair is contact and coincides with $(\lambda_\sigma,\omega_\sigma)|_{TM^{L,Q}}$, given by the formulas $$\Lambda^{L,Q}=\epsilon\alpha + e^Qd\theta$$$$\Omega^{L,Q}=\epsilon d\alpha$$ 

In other words, for $Q=0$, $(M^{L,0}_Y, \Lambda^{L,0})$ can be seen as the \emph{contactization} of the Liouville domain $(Y\times I,\epsilon d\alpha)$. 







In order to prove that $M^{L,Q}$ is stabilized by $Z$, one needs to show that $$\ker(\Omega^{L,Q})\subseteq \ker(d\Lambda^{L,Q}),\; \Lambda^{L,Q}\wedge (\Omega^{L,Q})^n>0$$ This is a straightforward check.

\vspace{0.5cm}

Along $M_Y^{L,Q}$ the Reeb vector field is given by $R^{L,Q}= \frac{\partial_\theta}{e^Q},$ which is totally degenerate, so we need to perturb our spine $M_Y^{L,Q}$, using a standard Morse function method. We have two approaches: the ``Morse'' approach, and the ``Morse--Bott'' one, depending on the choice of perturbing function. 

Observe first that the space of Reeb orbits is clearly identified with $Y \times [-L,L]$. The Morse approach consists in choosing $H_L: Y \times [-L,L] \rightarrow \mathbb{R}^{\geq 0}$ which is Morse outside of $Y \times \{|r|>L-\delta\}$, and such that in these regions depends only in $r$ and satisfies $\partial_rH_L\leq 0$ near $r=L$, $\partial_rH_L\geq 0$ near $r=-L$, vanishing in small neighbourhoods of $r=\pm L$, but strictly positive in $\{L-\delta<|r|<L-2/3 \delta\}$. The effect of the perturbation that we shall produce with this approach is that the Reeb orbits (with action up to a large action threshold) will correspond to the critical points of $H_L$, and will therefore be finite and non-degenerate. The second approach consists simply of choosing $H_L$ to depend only on $r$ \emph{globally}, with respect to which it is a Morse function, and not only near the boundary. This will have the effect of introducing a finite number of Morse--Bott families of Reeb orbits, one for each of the critical points of the function, each of which may be identified with a copy of $Y$. In both cases, the previous function may be extended to a function on $M^{L,Q}$, by setting $\widehat{H}_L(y,r,\theta)=H_L(y,r)$ in $M_Y^{L,Q}$, $\widehat{H}_L(y,\rho,\theta)=H_L(y,\pm L + F_\pm(\rho))$ in $M_C^{\pm,L,Q}$, and zero in $M_P^{\pm,L,Q}$. In what follows we shall usually take $H_L$ to be arbitrarily $C^1$-small as needed. We shall make no distinction between $\widehat{H}_L$ and $H_L$.

\vspace{0.5cm}

If $t\mapsto \Phi^t_Z$ denotes the flow of $Z$, choose $\epsilon > 0$ sufficiently small so that the manifold $$M^{\epsilon; L,Q}:=\{\Phi_Z^{\epsilon H_L(x)}(x) \in E^{\infty,\infty}: x \in M^{L,Q}\}$$ is still transverse to (and therefore stabilized by) $Z$. Geometrically, this manifold is just a small perturbation of $M^{L,Q}$, where how much we perturb is given by the function $H_L$ (see Figure \ref{doublecompletion2}). Again we have a stable Hamiltonian structure $$\mathcal{H}_\epsilon:=\mathcal{H}_\epsilon^{L,Q}:=(\Lambda_\epsilon^{L,Q},\Omega_\epsilon^{L,Q}):=(i_Z\omega_\sigma|_{TM^{\epsilon;L,Q}},\omega_\sigma|_{TM^{\epsilon;L,Q}}),$$ and a decomposition $$M^{\epsilon;L,Q}=M_Y^{\epsilon;L,Q} \cup M_C^{\epsilon;\pm,L,Q} \cup M_P^{\epsilon;\pm,L,Q},$$ where each component is the perturbation of the corresponding component of $M^{L,Q}$. Since $H_L$ vanishes on $M_P^{\pm,L,Q}$, we actually have that $$M_P^{\epsilon;\pm,L,Q}=M_P^{\pm,L,Q},$$ along which we have

$$\Lambda_\epsilon^{L,Q}=\Lambda^{L,Q}=K \alpha_\pm $$$$\Omega_\epsilon^{L,Q}=\Omega^{L,Q}=Kd\alpha_\pm + d\lambda$$ Since $Z$ is Liouville along the spine $M_Y^{L,Q}$, and therefore its flow exponentiates the contact form $\Lambda^{L,Q}$, we have, in $M_Y^{\epsilon;,L,Q}$, the expressions $$\Lambda^{L,Q}_\epsilon = e^{\epsilon H_L}\Lambda^{L,Q}=e^{\epsilon H_L}(\epsilon\alpha + e^Qd\theta)$$$$\Omega_\epsilon^{L,Q}=d\Lambda_\epsilon^{L,Q}=\epsilon e^{\epsilon H_L}dH_L \wedge (\epsilon\alpha + e^Qd\theta)+\epsilon e^{\epsilon H_L}d\alpha$$ Now, the region $M_C^{\epsilon;\pm,L,Q}$ is obtained by flowing along the vector field $Z=V_\sigma+\partial_t=\pm \frac{\sigma}{\sigma\pm\sigma^\prime}\partial_r + \partial_t$, which has time-$\epsilon H$ flow given in coordinates by $$r \mapsto \Phi_{V_\sigma}(\epsilon H_L,r)$$$$t \mapsto t + \epsilon H_L,$$ where $\Phi_{V_\sigma}(s,\cdot)$ is the time $s$ flow of $V_\sigma$. Recalling that $H_L$ does not depend on $y$ here, the $r$ and $t$ coordinates on the perturbation of the smoothened corners are $$r=F_\pm^{\epsilon;L}(\rho)=\Phi_{V_\sigma}(\epsilon H_L(\cdot,\pm L+F_\pm(\rho)),\pm L+F_\pm(\rho))$$$$ t=G_\pm^{\epsilon;L,Q}(\rho)=Q+G_\pm(\rho)+ \epsilon H_L(\cdot,\pm L+F_\pm(\rho)),$$ so that 

\begin{equation}
\label{lambdaC}
\begin{split}
\Lambda^{L,Q}_\epsilon &=  \sigma(F_\pm^{\epsilon;L}(\rho)) e^{\pm F_\pm^{\epsilon;L}-L} \alpha_\pm + \beta(G_\pm^{\epsilon;L,Q}(\rho))e^{G_\pm^{\epsilon;L,Q}(\rho)}d\theta\\
\Omega_\epsilon^{L,Q}=&\;(\sigma^\prime\pm \sigma)(F_\pm^{\epsilon;L}(\rho))e^{\pm F^{\epsilon;L}_\pm(\rho)-L} (F_\pm^{\epsilon;L})^\prime(\rho)d\rho \wedge \alpha_\pm\\& + \sigma(F_\pm^{\epsilon;L}(\rho)) e^{\pm F_\pm^{\epsilon;L}(\rho)-L}d\alpha_\pm\\& + e^{G_\pm^{\epsilon;L,Q}(\rho)}(G_\pm^{\epsilon;L,Q})^\prime(\rho) d\rho \wedge d\theta,
\end{split}
\end{equation}

in this region.

\vspace{0.5cm}

We now compute the Reeb vector field $R^{L,Q}_\epsilon$ associated to this stable Hamiltonian structure. 

In $M_P^{\epsilon;\pm,L,Q}$, $$R_\epsilon^{L,Q}=R^{L,Q}=\frac{R_\pm}{K}$$ 

In $M_Y^{\epsilon;L,Q}$, we have the expression 
\begin{equation}\label{ReebinY}
R_\epsilon^{L,Q}=e^{-\epsilon H_L-Q}\left((1+\epsilon \alpha(X_{H_L}))\partial_\theta-e^QX_{H_L}\right),
\end{equation} 

where $X_{H_L}$ is the Hamiltonian vector field on $Y \times I$ associated to $H_L$. The sign convention for Hamiltonians we are taking is $i_{X_{H_L}}d\alpha = -dH_L$. We see then by the formulas above that $R^{L,Q}_\epsilon$ converges to $R^{L,Q}$, so the previous is indeed a perturbation of the latter. Observe that if we evaluate at a critical point $(y,r)$ of $H$, the vector field $X_{H_L}$ vanishes, and so we have that the Reeb is given by $e^{-\epsilon H_L(y,r) - Q}\partial_\theta$. Thus, its orbits will wrap around the $S^1$ factor of the spine with period $2\pi e^{\epsilon H_L(y,r)+Q}$. If we are taking the Morse approach, the set crit$(H_L)$ of its critical points is finite, and hence we have only a finite number of Reeb orbits of the form $\gamma_p:=\{p\} \times S^1 \subseteq \mbox{crit}(H_L)\times S^1$ in this region. 

Observe that choosing $H_L$ to be $C^1$-small has the effect of making the vector field $X_{H_L}$ also small, so that the closed orbits which do not arise from critical points of $H_L$ have large period. So, taking any large (but fixed) $T\gg 0$, we can choose $H_L$ small enough so that all the periodic Reeb orbits in this region up to period $T$ are of the form $\gamma_p^l$, for $p \in \mbox{crit}(H_L)$, and $l\leq N$, for some covering threshold $N$ depending on $T$. A similar situation holds for the Morse--Bott case, when the low action Reeb orbits are still of the same form, only that $p$ varies in $Y$, and there are as many copies of these families as critical points of $H_L$. 

\vspace{0.5cm}

Finally, in $M_C^{\epsilon;\pm,L,Q}$, the Reeb is 

\begin{equation}
\label{ReebinC}
R_\epsilon^{L,Q}=\frac{1}{\Phi^{\epsilon; L,Q}_\pm(\rho)}\left(e^{\mp F_\pm^{\epsilon;L}(\rho)+L}(G_\pm^{\epsilon;L,Q})^\prime(\rho) R_\pm- e ^{-G^{\epsilon;L,Q}_\pm(\rho)}(\sigma^\prime\pm \sigma)(F_\pm^{\epsilon;L}(\rho))(F_\pm^{\epsilon;L})^\prime(\rho)\partial_\theta\right),
\end{equation}

where
 
\begin{equation} \label{Phieps} 
\Phi^{\epsilon; L,Q}_\pm(\rho)= \sigma(F_\pm^{\epsilon;L}(\rho)) (G_\pm^{\epsilon;L,Q})^\prime(\rho) - \beta(G_\pm^{\epsilon;L,Q}(\rho))(\sigma^\prime \pm \sigma)(F_\pm^{\epsilon;L}(\rho))(F_\pm^{\epsilon;L})^\prime(\rho)
\end{equation}

One can check that $\Phi^{\epsilon; L,Q}_\pm$ has sign which is \emph{opposite} to its subscript. 

\begin{remark} \label{remark11}

\begin{enumerate}[wide, labelwidth=!, labelindent=0pt] $\;$
\item We observe that for any fixed large $T\gg 0$ action threshold, by making $K$ large enough and $H_L$ small enough, and choosing $\sigma$ suitably, \emph{all} the closed Reeb orbits of action less than $T$ (not only the ones contained in $M_Y^\epsilon$, but in the whole model) will be covers of Reeb orbits of the form $\{p\} \times S^1$, for $p \in \mbox{crit}(H_L)$. Therefore they are \emph{non-degenerate}/Morse--Bott in the Morse/Morse--Bott approach.

\item One can check that $\lambda_\sigma(R_\epsilon)=1$ (recall that $\lambda_\sigma\vert_{M^{\epsilon;L,Q}}$ is the primitive of $\Omega_\epsilon=\omega_\sigma\vert_{M^{\epsilon;L,Q}}$). Therefore, if we have an almost complex structure which is compatible with our SHS, as we will in upcoming sections, and an asymptotically cylindrical $J$-holomorphic curve $u$ with positive/negative punctures $\Gamma^\pm$, then the $\Omega_\epsilon$-energy of $u$ can be computed as
\begin{equation}
\label{action}
\int u^*\Omega_\epsilon=\int u^*d\lambda_\sigma=\sum_{z \in \Gamma^+}T_z-\sum_{z \in \Gamma^-}T_z,
\end{equation}
where $T_z$ is the action of the Reeb orbit corresponding to the puncture $z$. In particular, if the positive punctures correspond to critical points of $H_L$, then so will the negative ones, by the previous remark.

\item By inspecting the expressions of the Reeb vector field we see that there are no contractible closed Reeb orbits for the SHS, if we assume this same condition for $R_\pm$. Moreover, the direction of the Reeb vector field does not change after perturbing back to sufficiently close contact data (cf.\ Section \ref{SHStoND} below), and so this also holds for the latter data. It follows that the isotopy class defined by model A is hypertight, and this shows the hypertightness condition of Theorem \ref{thm5d}.

\end{enumerate}
\end{remark}

\subsection{Compatible almost complex structure}
\label{ACS}

\paragraph*{Construction}

We will from now on set $L=1$, and $Q=0$, dropping the superscripts $L$ and $Q$ from all of the notation. We will proceed to define a suitable, though non-generic, almost complex structure $J=J_\epsilon$ on the symplectization $W^\epsilon=\mathbb{R} \times M^\epsilon$, where $$M^\epsilon=M^{\epsilon;1,0}=M^\epsilon_Y \bigcup M_P^{\epsilon;\pm} \bigcup M^{\epsilon;\pm}_C$$ It will be compatible with the stable Hamiltonian structure $\mathcal{H}_\epsilon$, and the fibers $\Sigma_\pm$ of our fibration $\pi_P^\pm$, the ``pages'', will lift as holomorphic curves. We will blur the distinction between $M=M^{0;1,0}$ and its diffeomorphic perturbed copy $M^\epsilon$ (as well as for $W$ and $W^\epsilon$), so that we are actually working on a fixed $M$ with a SHS which depends on $\epsilon$. 

\vspace{0.5cm}

Denote by $\xi_\epsilon:= \ker \Lambda_\epsilon$. We will define $J$ on $\xi_\epsilon$, in an $\mathbb{R}$-invariant way, and then simply set $J(\partial_a):=R_\epsilon$.

Choose a $d\alpha$-compatible almost complex structure $J_0$ on $Y\times I$ as in Corollary \ref{compJ0}. Observe that the vector field $\partial_\theta$ is transverse to $\xi_\epsilon$ along $M_Y$. Therefore, we may then define $J_Y=\pi_Y^* J_0$ on $\xi_\epsilon\vert_{M_Y}$.

\vspace{0.5cm}

Recalling that in $\Sigma_\pm/\mathcal{N}(\partial\Sigma_\pm) \times Y\subseteq M_P^{\pm}$ we have that $\Lambda_\epsilon=K\alpha_\pm$, whereas in $\{t \in (-\delta,-\delta/3)\}\times Y\subseteq M_P^{\pm}$, that $\Lambda_\epsilon=K\alpha_\pm +\beta(t)e^td\theta$, we get that $$\xi_\epsilon=\xi_\pm \oplus T\Sigma_\pm$$$$\xi_\epsilon = \xi_\pm \oplus \langle \partial_t , \beta(t)e^tR_\pm -K\partial_\theta \rangle,$$ respectively. Therefore the restriction of the projection $\pi_P^\pm: M_P^\pm\rightarrow \Sigma_\pm$ to each piece induces an isomorphism $d\pi_P^\pm: \xi_\epsilon/\xi_\pm \stackrel{\simeq}{\longrightarrow} T\Sigma_\pm$. Choose $j_\pm$ to be a $d\lambda$-compatible almost complex structures on $\Sigma_\pm$, so that $j_\pm(\partial_t)= K \partial_\theta$ in $\mathcal{N}(\partial \Sigma_\pm)$. Choose also a $d\alpha_\pm$-compatible almost complex structure $J^\pm$ on $\xi_\pm$. Define $$J=J^\pm \oplus (\pi_P^\pm)^*j_\pm$$ on $\xi_{\epsilon}\vert_{M_P^\pm}$.

\vspace{0.5cm}

Now, in $M_C^{\pm}$, where $\Lambda_\epsilon=\sigma(F_\pm^\epsilon) e^{\pm F^\epsilon_\pm-1}\alpha_\pm + \beta(G_\pm^\epsilon)e^{G_\pm^\epsilon}d\theta$, we have that 
$$\xi_\epsilon=\xi_\pm \oplus \langle \partial_\rho , T_\pm \rangle,$$ where 
$$T_\pm = \beta(G_\pm^\epsilon)e^{G^\epsilon\pm}R_\pm - \sigma(F_\pm^\epsilon) e^{\pm F^\epsilon_\pm-1}\partial_\theta$$ 

Then, $$\Omega_\epsilon(\partial_\rho,T_\pm)=e^{\pm F^\epsilon+G^\epsilon_\pm-1}(\beta(G_\pm^\epsilon)(\sigma^\prime\pm \sigma)(F_\pm^\epsilon)(F_\pm^\epsilon)^\prime - \sigma(F_\pm^\epsilon)(G_\pm^\epsilon)^\prime)=-e^{\pm F_\pm^\epsilon+G^\epsilon_\pm-1}\Phi_\pm^\epsilon,$$ where we have used expression (\ref{Phieps}). Its sign matches the subscripts (by recalling that $\Phi_\pm^\epsilon$ has sign opposite to the subscripts). We can then set 

\begin{equation}\label{wpm}
v_1=\partial_\rho,\;v_2^\pm=\frac{T_\pm^\epsilon}{\Omega_\epsilon(\partial_\rho,T_\pm^\epsilon)}=a_\pm(\rho)R_\pm + b_\pm(\rho) \partial_\theta,
\end{equation}
where $$a_\pm=-\frac{\beta(G_\pm^\epsilon)}{e^{\pm F^\epsilon_\pm-1}\Phi_\pm^\epsilon}$$ 
are functions whose signs match the subscript, whereas the signs of 
$$b_\pm=\frac{\sigma(F_\pm^\epsilon)}{ e^{G_\pm^\epsilon}\Phi_\pm^\epsilon}$$ 
are opposite to the subscripts. We have $\Omega_\epsilon(v_1,v_2^\pm)=1$ and 
$$\xi_\epsilon=\xi_\pm \oplus \langle v_1 , v_2^\pm \rangle,$$ a manifestation of the fact that $\Omega_\epsilon$ is non-degenerate when restricted to $\xi_\epsilon=\ker \Lambda_\epsilon$. Moreover, this splitting is $\Omega_\epsilon$-orthogonal, as is easy to check. Now, in the overlaps $M_C^{\pm}\cap M_Y$, we have that the vector $v_2^\pm$ is sent to $a_\pm(\rho)R_\pm$ under the linearised projection $d\pi_Y$. Moreover, since $r=\rho$, we have $v_1=\partial_r$, and therefore we get 
\begin{equation}
\begin{array}{lll}
J(v_1)&=(d\pi_Y|_{\xi_\epsilon})^{-1}J_0(\partial_r)&=\pm (d\pi_Y|_{\xi_\epsilon})^{-1}J_0(V)\\
&=\pm (d\pi_Y|_{\xi_\epsilon})^{-1}( R_{F^\epsilon_\pm(\rho)})&=\pm e^{\mp F^\epsilon_\pm(\rho)+1} (d\pi_Y|_{\xi_\epsilon})^{-1}(R_\pm)\\
&=\pm e^{\mp F^\epsilon_\pm(\rho)+1} \left(R_\pm +\frac{b_\pm(\rho)}{a_\pm(\rho)}\partial_\theta\right)&=\pm \frac{e^{\mp F^\epsilon_\pm(\rho)+1}}{a_\pm(\rho)} v_2^\pm\\
&= g_\pm^Y(\rho) v_2^\pm &\\
\end{array}
\end{equation} where we set $g_\pm^Y:=\pm \frac{e^{\mp F^\epsilon\pm(\rho)+1}}{a_\pm(\rho)}=\mp\frac{\Phi_\pm^\epsilon}{\beta(G_\pm^\epsilon)}$, which is always positive.

Moreover, in $M_P^{\pm}\cap M_C^{\pm}$, we have $t=G_\pm(\rho)=G_\pm^\epsilon(\rho)=\mp \rho$, $\Phi^\epsilon_\pm = G_\pm^\prime = \mp 1$, $\beta(G_\pm^\epsilon)=0$, and so $v_1=\mp \partial_t$ and $v_2^\pm=\mp K e^{\pm \rho}\partial_\theta$, hence 

$$J(v_1)=j_\pm(\mp \partial_t)=\mp K\partial_\theta= e^{\mp \rho} v_2^\pm$$ 

Since $g_\pm^Y$ and $e^{\mp\rho}$ are both positive, we can now take any smooth positive functions $h_\pm:(-\delta,\delta) \rightarrow \mathbb{R}^+$ which coincide with $e^{\mp \rho}$ near $\rho=\pm \delta$ and with $g_\pm^Y$ near $\rho=\mp\delta$. We glue the two definitions by setting 
\begin{equation}
\label{w}
J(v_1)=h_\pm(\rho)v_2^\pm:=w^\pm,
\end{equation}
and we make $J$ agree with $J^\pm$ on $\xi_\pm$.

This gives a well-defined cylindrical $J$ in $\mathbb{R}\times M$.

\paragraph*{Compatibility}

We check that $J$ is $\mathcal{H}_\epsilon$-compatible.

For $v \in TW$, write $$v=v^{\xi}+v^RR_\epsilon + v^a \partial_a,$$ where $v^\xi \in \xi_\epsilon$, so that $$Jv=Jv^\xi-v^R\partial_a+v^aR_\epsilon$$ 
Therefore $\Omega_\epsilon(v,Jv)=\Omega_\epsilon(v^\xi,Jv^\xi)$, and $\Omega_\epsilon(Jv,Jw)=\Omega_\epsilon(Jv^\xi,Jw^\xi)$, for $w \in TW$.

Outside of the cylindrical region, where $\Omega_\epsilon= K d\alpha_\pm+d\lambda_\pm$, $\xi_\epsilon=\xi_\pm \oplus T\Sigma_\pm$ and $J\vert_{\xi_\epsilon}=J^\pm \oplus j_\pm$ for some $d\alpha_\pm$-compatible $J^\pm$, we can write $v^\xi=v^{\xi_\pm}+v^{\Sigma\pm}$, and we have that $$\Omega_\epsilon(v,Jv)=K d\alpha_\pm(v^{\xi_\pm},J^\pm v^{\xi_\pm})+ d\lambda_\pm(v^{\Sigma\pm},j_\pm v^{\Sigma\pm})\geq 0,$$ with equality if and only if $v^\xi=0$. If $w \in TW$ is another vector, we get 
\begin{equation}
\begin{split}
\Omega_\epsilon(Jv,Jw)=&K d\alpha_\pm(J^\pm v^{\xi_\pm},J^\pm w^{\xi_\pm})+ d\lambda_\pm(j_\pm v^{\Sigma\pm},j_\pm w^{\Sigma\pm})\\
=&K d\alpha_\pm(v^{\xi_\pm}, w^{\xi_\pm})+ d\lambda_\pm( v^{\Sigma\pm},w^{\Sigma\pm})\\
=&\Omega_\epsilon(v,w)
\end{split}
\end{equation}

In the cylindrical region, where $\Omega_\epsilon=d\Lambda_\epsilon=d\left(e^{\epsilon H}(\epsilon\alpha + d\theta)\right)$, if we let $\pi:=\pi_Y$ be the projection to $Y\times I$, by definition of $J$ we have that $\pi Jv^\xi=J_0\pi v^\xi$. Thus we obtain 

$$\Omega_\epsilon(v,Jv)=\epsilon e^{\epsilon H}d\alpha(\pi v^\xi,J_0 \pi v^\xi)\geq 0$$

We have equality if and only if only $\pi v^\xi=0$, in which case $v^\xi=0$ since $\pi\vert_{\xi_\epsilon}$ is an isomorphism. Similarly, 

\begin{equation}
\begin{split}
\Omega_\epsilon(Jv,Jw)=&\epsilon e^{\epsilon H}d\alpha(J_0\pi v^\xi,J_0 \pi w^\xi)\\
=&\epsilon e^{\epsilon H}d\alpha(\pi v^\xi,\pi w^\xi)\\
=&\Omega_\epsilon(v,w)
\end{split}
\end{equation}

Over the corners $M_P^\pm$, where $\xi_\epsilon=\xi_\pm\oplus\langle v_1,w^\pm \rangle$ is $\Omega_\epsilon$-orthogonal and $\Omega_\epsilon(v_1,w^\pm)=h_\pm$, we can write 

$$v^\xi=v^{\xi_\pm}+v^1v_1+v^2 w^\pm$$
$$Jv^\xi=J^\pm v^{\xi_\pm}+v^1 w^\pm - v^2 v_1,$$
and then
$$
\Omega_\epsilon(v,Jv)= \sigma(F_\pm^\epsilon)e^{\pm F_\pm^\epsilon-1}d\alpha_\pm(v^{\xi_\pm},J^\pm v^{\xi_\pm})+h_\pm((v^1)^2+(v^2)^2)\geq 0,
$$ again with equality if and only if $v^\xi=0$.
Also,
\begin{equation}
\begin{split}
\Omega_\epsilon(Jv,Jw)=& \sigma(F_\pm^\epsilon)e^{\pm F_\pm^\epsilon-1}d\alpha_\pm(J^\pm v^{\xi_\pm},J^\pm v^{\xi_\pm})+h_\pm(v^1w^2-v^2w^1)\\
=& \Omega_\epsilon(v,w)
\end{split}
\end{equation}

\begin{remark}\label{dlambda}$\;$
We observe that over $\mathbb{R}\times M_P^\pm$, where $\Lambda_\epsilon=K\alpha_\pm$, we have $d\Lambda_\epsilon(v,Jv)\geq 0$ , with equality if and only if $v^{\xi_\pm}=0$, so that the projection to $TY$ of $v$ lies in the span of $R_\pm$.
\end{remark}

\subsection{Finite energy foliation}
\label{Curves}

We will now consider the symplectization of our stable Hamiltonian manifold $(M, \mathcal{H}_\epsilon)$, given by $(W=\mathbb{R}\times M,\omega^\varphi_\epsilon=d(\varphi(a)\Lambda_\epsilon)+\Omega_\epsilon)$, where $\varphi \in \mathcal{P}$. We will construct a finite-energy foliation of $W$ by $J$-holomorphic curves, consisting of three distinct types, which we describe in the next theorem. 

\begin{thm}
There exists a finite energy foliation of the symplectization $(W,\omega_\epsilon^\varphi)$ by simple $J$-holomorphic curves of the following types:

\begin{itemize}
\item \textbf{trivial cylinders $C_p$}, corresponding to Reeb orbits of the form $\gamma_p=\{p\} \times S^1$ for $p \in$ crit$(H)$, and which may be parametrized by 

$$C_p: \mathbb{R}\times S^1 \rightarrow \mathbb{R} \times M_Y =\mathbb{R} \times (Y \times I) \times S^1 $$$$C_p(s,t)=(s,p,e^{-\epsilon H(p)}t)$$ 

\item \textbf{flow-line cylinders $u^a_\gamma$}, parametrized by $$u^a_\gamma: \mathbb{R}\times S^1 \rightarrow \mathbb{R} \times M_Y =\mathbb{R} \times (Y \times I) \times S^1$$
\begin{equation}\label{trivcyl}
u^a_\gamma(s,t)=(a(s),\gamma(s),\theta(s)+t)
\end{equation}
for a proper function $a: \mathbb{R}\rightarrow \mathbb{R}$, a function $\theta:\mathbb{R}\rightarrow S^1$, and a map $\gamma: \mathbb{R}\rightarrow Y\times I$ satisfying 

\begin{equation}
\label{floweqs}
\begin{array}{l} 
\dot{\gamma}= \nabla H(\gamma(s))\\
\dot{a}(s)= e^{\epsilon H(\gamma(s))},\; a(0)=a\\
\dot{\theta}(s)=-\epsilon\alpha(\nabla H (\gamma(s))),\; \theta(0)=0
\end{array}
\end{equation}

Here, the gradient is computed with respect to the metric $g_{d\alpha,J_0}:= d\alpha(\cdot, J_0\cdot)$.

They have for positive/negative asymptotics the Reeb orbits corresponding to $$p_\pm:=\lim_{s\rightarrow \pm \infty} \gamma(s) \in \mbox{crit}(H)$$ 

\item \textbf{positive/negative page-like holomorphic curves $u^\pm_{y,a}$}, which consist of a trivial lift at symplectization level $a$ of a page $P_y^\pm:= \left(\Sigma_\pm\backslash\mathcal{N}(\partial\Sigma_\pm)\right) \times \{y\}$, for $y \in Y$, glued to $k$ cylindrical ends which lift the smoothened corners and then enter the symplectization of $M_Y$, asymptotically becoming a flow-line cylinder. They have $k$ positive asymptotics at Reeb orbits of the form $\gamma_p$, exactly one for each component of $M_Y$. The positive curves have genus $g-k+1>0$ and $k$ punctures, whereas the negative curves have genus $0$, and also $k$ punctures.  
\end{itemize}
\end{thm}

\begin{figure}[!ht]
\centering
\includegraphics[width=0.9\linewidth]{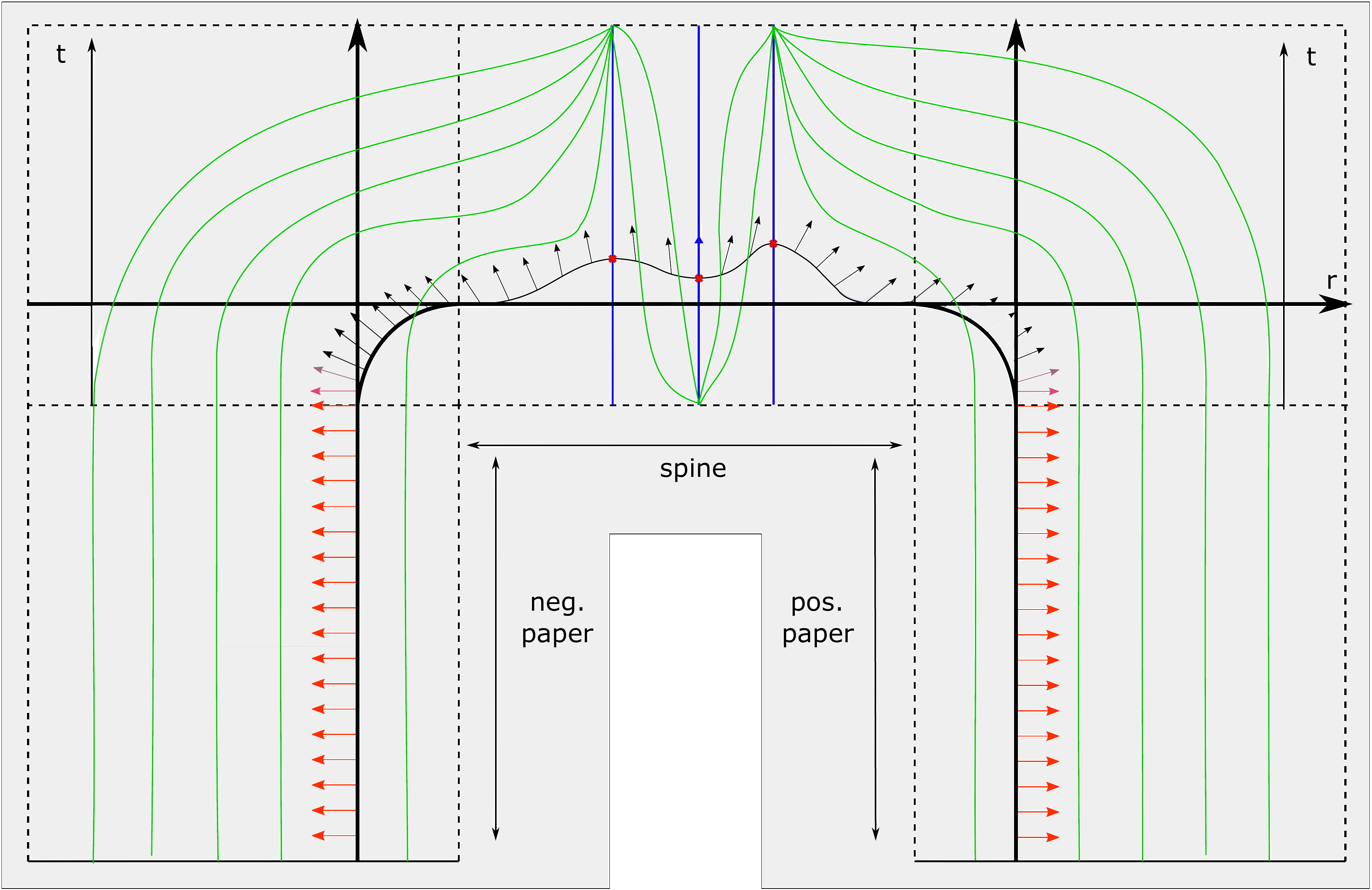}
\caption{\label{DoubleCompletion} The picture shows the double completion $E^{\infty,\infty}$, where one can see both $M$ and its perturbed version $M^\epsilon$ sitting as contact type hypersurfaces, such that their positive symplectizations embed in $E^{\infty,\infty}$ (one has to follow the flow of $X_\sigma$ inside this manifold to get the $a$-coordinate). We draw the stabilizing vector field $Z$. It coincides with the Liouville vector field $X_\sigma$ only over the spine and a portion of the smoothened corners, where it is drawn in black. We draw it blue at $r=0$, where $V$ becomes tangent to $Y$, so that it becomes transverse to the picture. The transition in color from red to black happens where the homotopy takes place. The Morse/Morse--Bott function shown has three critical points in the interval direction, 2 local maxima and a minimum, giving rise to three $Y$-families of Reeb orbits in the Morse--Bott case. The foliation by holomorphic curves is shown in green (for the non-trivial curves) and blue (for the trivial cylinders). The ones on the left correspond to the genus zero $u^-_{y,a}$'s, built by gluing three components, each corresponding to the dotted quadrants. The ones in the right to the positive genus $u^+_{y,a}$'s, and the ones in the middle are the holomorphic cylinders $u^a_\gamma$'s. The latter move in the direction of the \emph{positive} Morse flow lines, so their behaviour is as shown.}
\end{figure}

\begin{figure}[!ht] \centering
\includegraphics[width=0.70\linewidth]{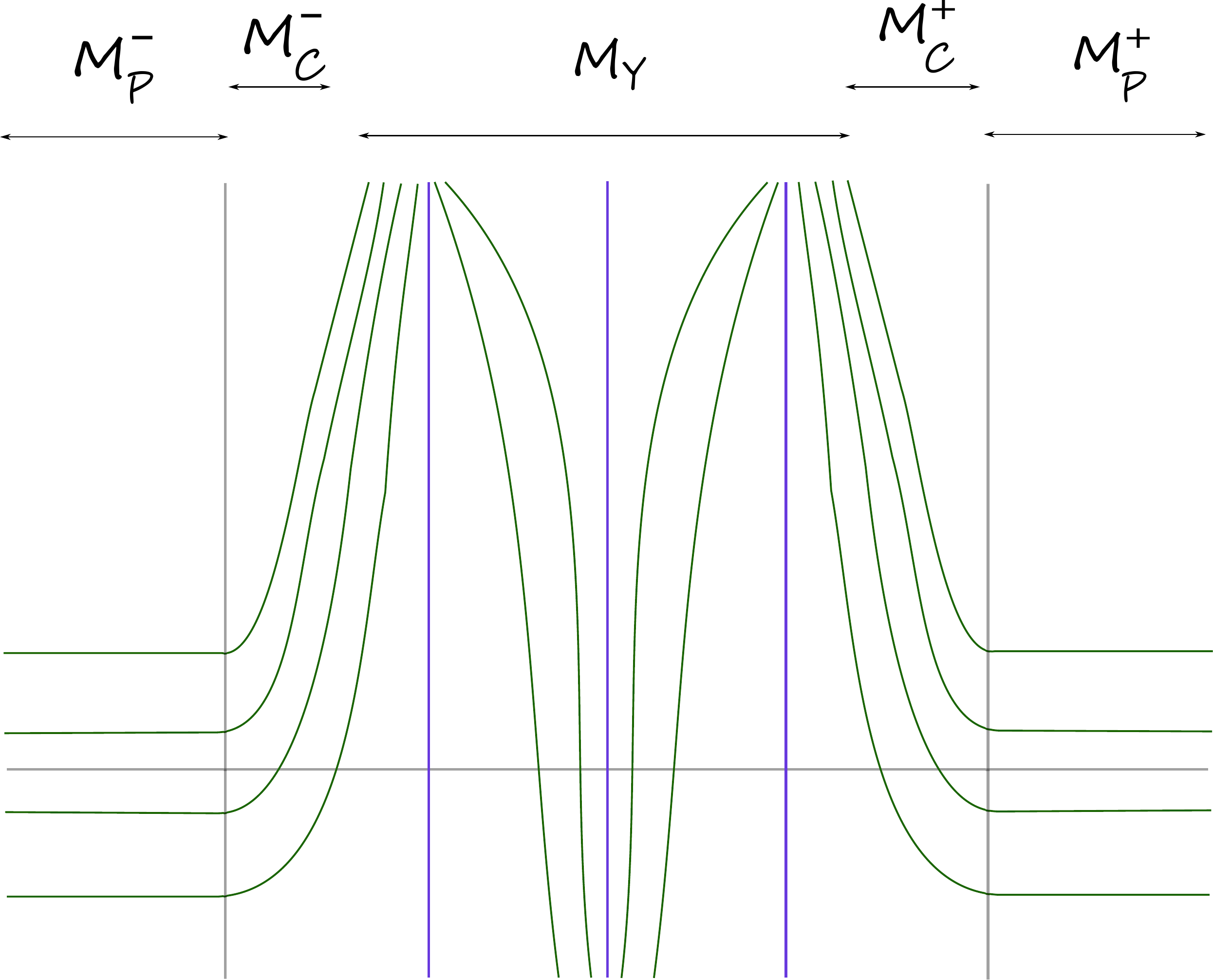}
\caption{\label{fol} Here is a sketch of the finite energy foliation from the point of view of the actual symplectization. The curves on the left and right (corresponding to $u_{y,a}^-$ and $u_{y,a}^+$ respectively) start at a constant symplectization level and glue to cylindrical ends which asymptote the trivial cylinders depicted in blue. The foliation in the middle correspond to the flow-line cylinders, with asymptotic behaviour determined by the function $H$.}
\end{figure}

Observe that in the Morse-Bott case the function $\theta$ vanishes identically, since $\alpha(\nabla H)=0$ (Eq.\ \ref{grad}).

Figure \ref{DoubleCompletion} summarizes the situation. We shall not distinguish the curves $u_{y,a}^{\pm}$ and $u^a_\gamma$ from the simple holomorphic curves that they parametrize, and we will drop the $a$ for the notation whenever we wish to refer to the equivalence class of the curves under $\mathbb{R}$-translation. 

We begin by checking that the flow-line cylinders are indeed holomorphic. This is the same computation as e.g.\ is done in \cite{Sie}, from where we borrow some notation.

\vspace{0.5cm}

Take on $Y\times I$ the metric $d\alpha(\cdot,J_0 \cdot)$, and denote by $\nabla H$ the resulting gradient vector field associated to $H$, and $X_H=J_0 \nabla H$ its Hamiltonian vector field. For every $X \in T(Y\times I)$, we denote by $$\widehat{X}= X-\epsilon \alpha(X)\partial_\theta$$ its lift to $\xi_\epsilon$, so that we have $J\widehat{X}=\widehat{J_0 X}$, by definition of $J$ along $M_Y$. Therefore, expression (\ref{ReebinY}) can be rewritten as 
$$
R_\epsilon=e^{-\epsilon H}\left(\partial_\theta-\widehat{X_H}\right)
$$ 

We then check that the Cauchy--Riemann equation holds for $u_\gamma^a$:

\begin{equation}
\begin{split}
\partial_su_\gamma^a(s,t)=&\dot{a}(s)\partial_a+\dot{\gamma}(s)+\dot{\theta}(s)\partial_\theta\\
=&e^{\epsilon H(\gamma(s))}\partial_a+\nabla H (\gamma(s)) - \epsilon \alpha(\nabla H(\gamma(s)))\partial_\theta\\
=&e^{\epsilon H(\gamma(s))}\partial_a+\widehat{\nabla H}(\gamma(s))
\end{split}
\end{equation}

\begin{equation}
\begin{split}
\partial_tu_\gamma^a(s,t)=& \partial_\theta\\
=&\partial_\theta-\widehat{X_H}(\gamma(s)) +\widehat{X_H}(\gamma(s))\\
=&e^{\epsilon H(\gamma(s))}R_\epsilon+\widehat{X_H}(\gamma(s)),
\end{split}
\end{equation}

so that clearly $J\partial_su_\gamma^a=\partial_tu_\gamma^a$.

The fact that the flow lines plus the critical points give a foliation of $Y \times I$ implies that the $u_\gamma^a$'s together with the trivial cylinders over the Reeb orbits corresponding to critical points foliate the whole of the region $\mathbb{R}\times M_Y \subseteq W$. 

\vspace{0.5cm}

We now construct the page-like curves.

The pages $P^\pm_y:=(\Sigma_\pm\backslash\mathcal{N}(\partial \Sigma_\pm)) \times \{y\}$, $y \in Y$, clearly lift to a foliation of the region $\mathbb{R} \times M_P^{\pm} \subseteq W$, which takes the form $\{\{a\} \times P^\pm_y :a \in \mathbb{R}, y \in Y\}$. It consists of $J$-invariant curves, since $J$ was taken to extend a page-wise almost complex structure $j_\pm$. We wish now to glue flow-line cylindrical ends to this lifts, along the corners $M_C^{\pm}$ and a portion of $M_Y$.

\vspace{0.5cm}

We have that $Jv_1=h_\pm v_2^\pm$ and $R_\epsilon$ are both linear combinations of the vector fields $R_\pm$ and $\partial_\theta$ along $M_C^\pm$, with the coefficients only depending on $\rho$. Since these are not colinear, we have smooth functions $B,C:(-\delta,\delta)\rightarrow \mathbb{R}$ such that $$\partial_\theta = B R_\epsilon + CJv_1$$ One can in fact compute that the following expressions hold: 
\begin{equation}\label{BC}
\begin{array}{l}
B=e^{G_\pm^\epsilon}\beta(G_\pm^\epsilon)\\
C=\frac{e^{G_\pm^\epsilon}(G_\pm^\epsilon)^\prime}{h_\pm}
\end{array}
\end{equation}

In particular, $B$ is non-negative and $C$ is positive/negative in the minus/plus regions (as follows by looking at the summands in $(G_\pm^\epsilon)^\prime= G_\pm^\prime+\epsilon \partial_rH F_\pm^\prime$). 

We have that 
\begin{equation}\label{holcorner}
J\partial_\theta = -B\partial_a  -Cv_1=-B\partial_a-C\partial_\rho
\end{equation}
We conclude that 
$$\langle \partial_\theta , B\partial_a+C\partial_\rho\rangle = \langle \partial_\theta ,J\partial_\theta \rangle
$$

It follows that the distribution above has integral submanifolds which are unparametrized holomorphic curves. We can actually find holomorphic parametrizations given by 

$$w^\pm_{y,a}: (-\delta,\delta) \times S^1 \rightarrow  \mathbb{R} \times M_C^{\pm} =\mathbb{R}\times  Y \times (-\delta,\delta)\times S^1$$$$ (s,t)\mapsto (b(s),y,\rho(s),t),$$ 

for some fixed $y \in Y$, and functions $b,\rho:(-\delta,\delta)\rightarrow \mathbb{R}$ satisfying 

$$\dot{\rho}(s)=C(\rho(s)),\; \rho(-\delta)=-\delta$$
$$\dot{b}(s)=B(\rho(s)),\;b(-\delta)=a$$

Equation (\ref{holcorner}) shows that $w^\pm_{y,a}$ is indeed holomorphic. The fact that $B$ is always positive implies that $b$ is increasing with respect to the $s$-parameter, and the fact that $C$ is positive/negative implies that $\rho$ is increasing/decreasing. Hence, by the orientation of  $\rho$, we have that $w_{y,a}^\pm$ will glue to a flow-line cylinder which has only \emph{positive} ends (one for each boundary component), for both the ``plus'' and ``minus'' curves. 

\vspace{0.5cm}

Now, the curves $w^\pm_{y,a}$ glue with curves $u^a_\gamma$ for a flow line $\gamma$ which is constant in the $Y$-direction close to $r=\pm 1$. This is because we have chosen $H$ to depend only on $r$ in the region $\{|r|>1-\delta\}$. This means that the curves $u^a_\gamma$ which start there look like $u^a_\gamma(s,t)=(a(s),y(s),r(s),t)$ for some $y(s) \in Y$ such that $y(s)\equiv y$ near $r=\pm 1$ and $\lim_{s\rightarrow +\infty}y(s)=y^+$, and some $r:\mathbb{R}\rightarrow I$ with $\lim_{s\rightarrow +\infty}r(s)=r^+$, so that $(y^+,r^+)\in crit(H)$. This shows that indeed we have smooth $J$-holomorphic curves $$u^\pm_{y,a}:=P^\pm_{y,a} \bigcup w^\pm_{y,a} \bigcup u^a_\gamma$$ which asympote $k$ Reeb orbits $\gamma_{y,a;i}^\pm=\{p_{y,a;i}^\pm\} \times S^1$, where $p_{y,a;i}^\pm \in crit(H)$, for $i=1,\dots,k$, and which have genus $g(u^-_{y,a})=0$ and $g(u^+_{y,a})=g-k+1$.  

\subsection{Index computations}
\label{indexcomp}

In this section, we compute the Fredholm index of the curves in the foliation.

\begin{thm} $\;$

\begin{enumerate}[wide, labelwidth=!, labelindent=0pt] 
\item After a sufficiently small Morse perturbation making Reeb orbits along $M_Y$ non-degenerate, we can find a natural trivialization $\tau$ of the contact structure along $\gamma_p$ (inducing a trivialization $\tau^l$ along all of its covers $\gamma_p^l$), and $N\in \mathbb{N}$, which depends on $H$ and grows as $H$ gets smaller, such that the Conley--Zehnder index of $\gamma_p^l$ is given by $$\mu^{\tau^l}_{CZ}(\gamma^l_p)=\mbox{ind}_p(H)-n,$$ for $l\leq N$.
  
\item In the Morse approach, the Fredholm indexes of the curves in our finite energy foliation are given by 
\begin{equation}
\begin{split}
\mbox{ind}(u_{y,a}^-)&=2n(1-k)+ \sum_{i=1}^k \mbox{ind}_{p^-_{y,a;i}}(H)\\ 
\mbox{ind}(u_{y,a}^+)&=2n(1-g-k)+\sum_{i=1}^k ind_{p_{y,a;i}^+}(H)\\
\mbox{ind}(u_\gamma^a)&=\mbox{ind}_{p_+}(H)-\mbox{ind}_{p_-}(H)
\end{split}
\end{equation}
\end{enumerate}
\end{thm} 

\subparagraph*{Computation of Conley--Zehnder indices}

We compute the Conley--Zehnder index for covers of the Reeb orbits $\gamma_p$, for a covering bound depending on $H$. This is a rather standard computation in Floer theory (see also \cite{vK}, for the computation on general prequantization spaces). We include our version of it, for completeness.

The asymptotic operator associated to the Reeb orbit $\gamma_p^l$ is
$$
\mathbf{A}_{\gamma_p^l}:W^{1,2}((\gamma_p^l)^*\xi_\epsilon)\rightarrow L^2((\gamma_p^l)^*\xi_\epsilon)
$$
$$\mathbf{A}_{\gamma_p^l}\eta=-J(\nabla_t \eta - T_{\gamma_p^l}\nabla_\eta R_\epsilon),
$$
where $\nabla$ is any symmetric connection in $M$, and $T_{\gamma_p^l}=\Lambda_\epsilon(\dot{\gamma}^l_p)=2\pi l e^{\epsilon H(p)}$ is the action of $\gamma_p^l$, which we assume parametrized by $$\dot{\gamma_p^l}=T_{\gamma_p^l} R_\epsilon(\gamma_p^l)=2\pi l\partial_\theta$$ 

Extend the metric $g_{d\alpha,J_0}$ so that $g_{d\alpha,J_0}(\partial_\theta,\partial_\theta)=1$ and $g_{d\alpha,J_0}(\partial_\theta,T(Y\times I))=0$, and take $\nabla$ to be the associated Levi--Civita connection, which is symmetric and satisfies $\nabla \partial_\theta=0$.

Under the holomorphic identification $(\gamma_p^l)^*\xi_\epsilon\simeq [0,1]\times T_p(Y\times I)$, a trivial bundle, we can find a unitary trivialization by choosing a symplectic basis with respect to $d\Lambda_\epsilon$, which makes $J=J_0=i$, the standard almost complex structure. From this, we can assume that $J$ is parallel with respect to $\nabla$ (i.e.\ $J\nabla = \nabla J$). Using formula (\ref{ReebinY}), we can compute
\begin{equation}
\begin{split}
\nabla R_\epsilon =& d\left(e^{-\epsilon H}\right) e^{\epsilon H}\otimes R_\epsilon +e^{-\epsilon H}\left(\epsilon d\left(\alpha(X_H)\right)\otimes \partial_\theta - \nabla X_H\right)\\
=& -\epsilon dH \otimes R_\epsilon + e^{-\epsilon H}\left(\epsilon d\left(\alpha(X_H)\right)\otimes \partial_\theta - J\nabla^2 H\right)
\end{split}
\end{equation}

Along $\gamma_p$ we have $dH=0$, and the $\partial_\theta$-term is lost when projecting to $T_p(Y\times I)$, so that the restriction to $(\gamma_p^l)^*\xi_\epsilon$ of the above expression reduces to
$$
\nabla R_\epsilon=-e^{-\epsilon H(p)}J\nabla^2 H
$$

Therefore, the symmetric matrix associated to $\mathbf{A}_{\gamma_p^l}$ is 
$$
S:=-T_{\gamma_p^l}J\nabla R_\epsilon =-2\pi l \nabla^2 H,
$$
so that the equation $\mathbf{A}_{\gamma_p^l}\eta=0$ is equivalent to $\dot{\eta}=iS\eta$. The matrix of solutions $\Psi(t)$, satisfying $\Psi(0)=\mathds{1}$, $\dot{\Psi}(t)=iS\Psi(t)$, is the path of symplectic matrices given by $\Psi(t)=e^{itS}$. Observe that if $H$ is chosen so that $||H||_{C^2}<\frac{1}{l}$, then $||S||<2\pi$, and $\Psi(t)$ will not have $1$ as an eigenvalue for all $t \in [0,1]$. Therefore, by making $H$ sufficiently $C^2$-small, we can ensure that this inequality holds, and the smaller we make $H$, the larger we are allowed to take $l$. For a fixed $C^2$-small $H$, we have an upper bound $l\leq N$ as in the statement. In this case, using Proposition 35.6 (signature axiom) in \cite{Gutt}, we have
\begin{equation}
\begin{split}
\mu^{\tau^l}_{CZ}(\gamma_p^l)=&\mu_{CZ}(\Psi)\\
=&\frac{\mbox{Sign}(S)}{2}\\
=&n-\mbox{ind(S)}\\
=&n-(2n-ind_p(H))\\
=&ind_p(H)-n
\end{split},
\end{equation}
where Sign$(S)$ is the signature of $S$ (the number of positive eigenvalues of $S$ minus the number of negative ones). 

One can also compute these indices by using the relationship between the Conley--Zehnder index and the spectral flow of asymptotic operators, as is done for instance in \cite{Wen3}; or one may appeal to lemma 2.4 in \cite{Bo}, relating the Conley--Zehnder index of orbits arising from perturbation to the Robbin--Salamon index of the Morse--Bott submanifolds.

\subparagraph*{Computation of the Fredholm index} If we look at one of the $u^a_\gamma$, which approaches Reeb orbits $\gamma_\pm=\{p_\pm\}\times S^1$, we have trivializations $\tau=(\tau_-,\tau_+)$ of $\gamma_\pm^*\xi_\epsilon$ which arise from a global trivialization of $(u_\gamma^a)^*\xi_\epsilon$ given by the $S^1$-invariance of $\xi_\epsilon$. This implies that $c_1^\tau((u_\gamma^a)^*\xi_\epsilon)=0$, and since the Euler characteristic of a cylinder is zero, we get $$ind(u^a_\gamma)=\mu^\tau_{CZ}(\gamma_+)-\mu^\tau_{CZ}(\gamma_-)=ind_{p_+}(H)-ind_{p_-}(H)$$

Observe that we have analogous trivializations of $(\beta_\pm)^*\xi_\epsilon$ (still called $\tau$), where $\beta_\pm$ is any $S^1$-fiber contained in $M_P^\pm\cap E^{\infty,\infty}(t)$, where $\xi_\epsilon$ is also $S^1$-invariant, obtained from the $\gamma_p$'s by ``sliding'' along the collar direction. This trick then shows that $c_1^\tau((u_{y,a}^\pm)^*\xi_\epsilon)=c_1^\tau((P^\pm_{y,a})^*\xi_\epsilon)$. This is convenient since $\xi_\epsilon=\xi_\pm \oplus T\Sigma_\pm$ over $M_P^{\pm}$, and hence, given that $P^\pm_{y,a}$ is constant in $Y$, we have that 

$$(P^\pm_{y,a})^*\xi_\epsilon=E^n\oplus T\Sigma_\pm,$$ 

where $E^n$ is the trivial complex bundle of rank $n$ with fiber $(\xi_\epsilon)_y$. Therefore

$$c_1^\tau((P^\pm_{y,a})^*\xi_\epsilon)=c_1^\tau(T\Sigma_\pm)=\chi(\Sigma_\pm),$$ 

so that 
\begin{equation}
\begin{split}
ind(u_{y,a}^-)&=(n-2)(2-k)+2c_1^\tau((u_{y,a}^-)^*\xi_\epsilon)+\sum_{i=1}^k \mu^\tau_{CZ}(\gamma^-_{y,a;i})\\
&=(n-2)(2-k)+2(2-k)+\sum_{i=1}^k (ind_{p_{y,a;i}^-}(H)-n)\\
&=2n(1-k)+\sum_{i=1}^k ind_{p_{y,a;i}^-}(H)
\end{split}
\end{equation}

Similarly, one computes $ind(u_{y,a}^+)$.

\begin{remark} 
\begin{enumerate}[wide, labelwidth=!, labelindent=0pt] $\;$

\item Observe that, if $ind_{p_{y,a;i}^-}(H)=2n$ for every $i$, i.e.\ the curve approaches maxima at each boundary component (which is the generic behaviour), the index becomes $ind(u_{y,a}^-)=2n$.
Also we see that since $ind_{p^+_{y,a;i}}(H)\leq 2n$ for every $i$, then $ind(u_{y,a}^+)\leq 2n(1-g)\leq 0$, since $g\geq 1$. This tells us that these curves cannot possibly achieve transversality, and hence, after a perturbation making $J$ generic, they will disappear.

\item When $n=1$, since the only possible values for the Morse indexes are $0$ (elliptic minima), $1$ (hyperbolic critical point) and $2$ (elliptic maxima), the index formula becomes $$ind(u^-_{y,a})=2-2k+\#\{\mbox{hyperbolic}\}+2\#\{\mbox{maxima}\}=$$$$=2-2k+\#\{\mbox{hyperbolic}\}+2(k-\#\{\mbox{hyperbolic}\}-\#\{\mbox{minima}\})=$$$$=2-\#\{\mbox{hyperbolic}\}-2\#\{\mbox{minima}\}$$ Moreover, since the manifold $Y \times I$ has non-empty boundary, it is always possible to choose our generic (and small) Morse function $H$ to have no index $0$ critical points. Hence in the $3$-dimensional case (when $Y$ can be taken to be $S^1$), the index formula is just $ind(u_{y,a}^-)=2-\#\{\mbox{hyperbolic}\}$.

\end{enumerate}
\end{remark}

\subsection{Fredholm regularity}
\label{Regularity}

In this section, we shall prove that the curves we have constructed are Fredholm regular. It can be skipped on a first read, since it is quite technical.

In the Morse case, regularity of the flow-line cylinders can be basically reduced to the Morse--Smale condition for $H$, as we will show. This fact is known to experts, though does not seem to appear in full details in the literature, so we shall include them. The ideas are mainly taken from chapter 10 in \cite{Wen3}. We perturb the data along the spine $M_Y$ to a SHS for which holomorphic cylinders correspond to solutions to the Floer equation in $Y\times I$. We may then appeal to a standard result from Floer theory, which tells us that for a sufficiently small Morse Hamiltonian $H$, any element in the kernel of the corresponding linearised Floer operator is $S^1$-independent, the linearised version of the fact that every solution to the Floer equation is a Morse flow line. The Morse--Smale condition gives the surjectivity of these operators (i.e.\ the regularity of the perturbed holomorphic cylinders), and a standard perturbation argument using the implicit function theorem finishes the proof.    

For regularity of the other curves, we will assume the Morse--Bott situation, and prove regularity of the genus zero curves in the foliation (the postive genus ones have negative index and hence clearly cannot be regular). We will use the fact from \cite{Wen1} that Fredholm regularity is equivalent to the surjectivity of the normal component of the linearized Cauchy-Riemann operator, which is again a Fredholm operator. After some assumptions on our choice of coordinates and the Morse-Bott function $H$ (which we can always assume hold), and a suitable choice of normal bundle, we will explicitly write down an expression for this operator. We will obtain a set of PDEs whose solutions are precisely the elements in its kernel, for which we can check that curves in the foliation which are nearby a fixed leaf correspond to solutions. By splitting the operator, and using automatic transversality \cite{Wen1}, we show that these are all possible solutions. This will imply that the index of the normal operator coincides with the dimension of its kernel, from which surjectivity follows. From the implicit function theorem, we also obtain regularity for Morse data chosen \emph{sufficiently close} to Morse-Bott data, which is enough for our purposes.

In order to do computations with linearized operators, we will choose a suitable connection on $W=\mathbb{R}\times M$. It will be given by the Levi--Civita connection of a suitable metric $g$ and hence symmetric. 

\paragraph*{Constructing a symmetric connection}

We now proceed to define a metric on $W$, which shall induce our preferred connection for computations in this section. Given an almost complex structure $J$ which is compatible with a symplectic form $\omega$, we will denote by $g_{J,\omega}=\omega(\cdot,J\cdot)$ the associated Riemannian metric. Given a metric $g$, the symbol $\nabla^g$ will denote the Levi--Civita connection associated to $g$. It is worth recalling the fact that a product metric on a product manifold yields a product Levi--Civita connection (this can be proven, for instance, by looking at the classical formula for its Christoffel symbols), as we shall be using it throughout. 

\vspace{0.5cm}

Define, in the regions $\mathbb{R} \times Y\times \Sigma_\pm\backslash\mathcal{N}(\partial\Sigma_\pm)$, the metric 

$$g=\left(\begin{array}{cccc} 1  & 0 & 0 & 0 \\
0& g_{J_P^\pm,d\alpha_\pm}  &0 & 0 \\
0 & 0 & 1 & 0 \\
0 & 0 & 0 & g_{j_\pm,d\lambda_\pm} 
\end{array}\right),$$ where we are using the splitting $$T(\mathbb{R} \times Y\times \Sigma_\pm\backslash\mathcal{N}(\partial\Sigma_\pm))=\langle \partial_a \rangle \oplus \xi_\pm \oplus \langle R_\epsilon\rangle \oplus T\Sigma_\pm=\langle \partial_a \rangle \oplus TY \oplus T\Sigma_\pm$$

This yields a connection of the form $\nabla=d \times \nabla^Y\times \nabla_\pm^\Sigma$, where $\nabla_\pm^\Sigma=\nabla^{g_{j_\pm}}$, $\nabla^Y$ is a connection in $Y$, and $d$ is the standard flat connection. Since in the cylindrical ends the metric $g_{j_\pm}$ is standard, we have that $\nabla_\pm^\Sigma=d$ there (that is, $\nabla \partial_t=\nabla \partial_\theta = 0$). Then, we may extend this connection to the symplectization of each smoothened corner $\mathbb{R}\times M_C^\pm=\mathbb{R}\times Y \times (-\delta,\delta)\times S^1$ as $d\times \nabla^Y \times d \times d$. This connection also arises from the metric obtained by replacing $g_{j_\pm,d\lambda_\pm}$ by the identity in the basis $\{\partial_\rho,\partial_\theta\}$ in the above matrix, which glues smoothly to $g$, and is therefore symmetric.

Extend now $g$ over the region $\mathbb{R} \times Y\times I \times S^1$ as 
$$g=\left(\begin{array}{ccc} 1  & 0 & 0\\
0& g_{J_0,d\alpha} & 0\\
0 & 0 & 1  
\end{array}\right) = \left(\begin{array}{ccccc} 1  & 0 & 0 & 0 & 0 \\
0& g_{\widetilde{\xi}_r}  &0 & 0 & 0 \\
0 & 0 & 1 & 0 & 0 \\
0 & 0 & 0 & 1 & 0 \\
0 & 0 & 0 & 0 & 1 
\end{array}\right),$$ where we use the splitting $$T(\mathbb{R} \times Y\times I \times S^1)=\mathbb{R} \times T(Y \times I) \times \langle \partial_\theta\rangle= \langle \partial_a \rangle \oplus \widetilde{\xi}_r \oplus \langle \widetilde{R}_r, \widetilde{V}\rangle \oplus \langle \partial_\theta \rangle,$$ and where $g_{\widetilde{\xi}_r}=g_{J_0,d\alpha}\vert_{\widetilde{\xi}_r}$ (so that $g\vert_{Y\times I}=g_{d\alpha,J_0}$). Here, we use the $d\alpha$-symplectic splitting of $T(Y\times I)$ provided by Corollary \ref{compJ0}.

This gives a smooth metric $g$ over $W$, and $\nabla=\nabla^g$ is its Levi--Civita connection, and is thus symmetric. This shall be the connection we will use to write down all our linearised Cauchy--Riemann operators.  

\subsubsection{Regularity for unbranched covers of flow-line cylinders in the Morse case}
\label{regcyls} 

Let us recall that $\Lambda_\epsilon=e^{\epsilon H}(\epsilon \alpha + d\theta)$ along $M_Y$. We then have  
$$
d\Lambda_\epsilon=\epsilon \widehat{\Omega}_\epsilon,
$$
where we define 
$$
\widehat{\Omega}_\epsilon:=e^{\epsilon H}\left(H^\prime dr \wedge (\epsilon \alpha + d\theta)+d\alpha\right)=\frac{1}{\epsilon}\Omega_\epsilon
$$
 
We claim that $\widehat{\mathcal{H}}_\epsilon:=(\Lambda_\epsilon, \widehat{\Omega}_\epsilon)$ is a SHS for sufficiently small $\epsilon\geq 0$, which is contact for $\epsilon>0$. Indeed, clearly one has $\ker\widehat{\Omega}_\epsilon\subseteq \ker d\Lambda_\epsilon$, and 

\begin{equation}
\begin{split}
\Lambda_\epsilon\wedge\widehat{\Omega}_\epsilon^n&=e^{n \epsilon H}\Lambda_\epsilon\wedge (d\alpha^n+n H^\prime dr \wedge d\theta \wedge d\alpha^{n-1}+n \epsilon H^\prime dr\wedge \alpha \wedge d\alpha^{n-1})\\
&=e^{(n+1)\epsilon H}d\theta \wedge d\alpha^n>0
\end{split}
\end{equation}

This means that we may consider the contact form $\Lambda_\epsilon$ as a small perturbation of the SHS $\widehat{\mathcal{H}}_0$. 

The Reeb vector field associated to $\widehat{\mathcal{H}}_\epsilon$ is the same as before, given by the formula
$$
R_\epsilon= e^{-\epsilon H}((1+\epsilon\alpha(X_H))\partial_\theta-X_H)
$$
Therefore, we see that 
$$
\lim_{\epsilon\rightarrow 0}R_\epsilon=\partial_\theta - X_H=R_0,
$$
the Reeb vector field of $\widehat{\mathcal{H}}_0$.

Observing that by taking $\epsilon$ to go to zero we also have a perturbation from our original $J=J_\epsilon$ defined on $\xi_\epsilon=\ker \Lambda_\epsilon$, to $J_0$ defined on $\xi_0=\ker \Lambda_0=T(Y\times I)$.  

Let us now consider a $J_0$-holomorphic cylinder
$$
u=(f,v,g):\mathbb{R}\times S^1\rightarrow \mathbb{R}\times M_Y=\mathbb{R}\times (Y \times I)\times S^1
$$
Its partial derivatives are given by

\begin{equation}\label{parders}
\begin{split}
u_s=&f_s\partial_a+v_s+g_s\partial_\theta\\
u_t=&f_t\partial_a+v_t+g_t\partial_\theta
\end{split}
\end{equation}

Since $J_0(\partial_a)=\partial_\theta-X_H$, the holomorphic curve equation applied to $u$ reads 
\begin{equation}
\label{holiceq}
\begin{split}
0&=u_s+J_0(u_t)\\
&=(f_s-g_t)\partial_a+(f_t+g_s)\partial_\theta+v_s+J_0(v_t+g_tX_H)-f_tX_H
\end{split}
\end{equation}
This implies that $f_s=g_t$, $f_t=-g_s$, so that we may view the map $(f,g)$ as a holomorphic map from the cylinder (the punctured plane) to itself, preserving the punctures. Every such map is of the form $z\mapsto a z^b$ for $a,b \in \mathbb{C}^*$, so that, if we assume that $u$ is simply covered, we can then reparametrize $u$ so that $(f,g)$ is the inclusion map. Assuming this, from (\ref{holiceq}) we obtain
\begin{equation}\label{Floereq}
v_s+J_0(v_t+X_H)=0
\end{equation}
so that $v$ satisfies the Floer equation with respect to the Hamiltonian $H$.

\vspace{0.5cm}

Observe that if $v_t=0$, then $v(s,t)=\gamma(s)$ is a flow-line for $H$, and $u(s,t)=(s,\gamma(s),t)$. Assume this is the case. Let us compare the linearisation of the Floer equation for $v$ with the linearisation of the Cauchy--Riemann equation for $u$, and also with the linearisation of the Morse equation for $\gamma$. If we take $\nabla$ to be the symmetric connection we already constructed, the linearised Floer operator can be written down as
$$
\mathbf{D}_v^{\textbf{Floer}}:W^{1,2}(v^*T(Y\times I))\rightarrow L^2(v^*T(Y\times I))
$$
$$
\mathbf{D}_v^{\textbf{Floer}}\eta=\nabla_s \eta+J_0(v)\nabla_t \eta+\left(\nabla_\eta J_0\right)v_t-\nabla_\eta \nabla H
$$

On the other hand, the linearised Cauchy--Riemann operator for $u$ is the operator
$$
\mathbf{D}_u:W^{1,2,\delta_0}(u^*TW)\oplus V_{\Gamma} \rightarrow L^{2,\delta_0}(\overline{\mbox{Hom}}_\mathbb{C}(\mathbb{R}\times S^1,u^*TW)),
$$
$$
(\mathbf{D}_u\eta)\partial_s=\nabla_s \eta+J_0(u)\nabla_t \eta+\left(\nabla_\eta J_0\right)u_t
$$
Here, $\delta_0>0$ is a sufficiently small weight making the asymptotic operator non-degenerate, and hence making $\mathbf{D}_u$ Fredholm, and $V_\Gamma$ is a real $4$-dimensional vector space of smooth sections asymptotic at the punctures to constant linear combinations of $R_0$ and $\partial_a$ (there is no \emph{Teichm\"uller slice}---see \cite{Wen1}---since Teichm\"uller space of a cylinder is zero-dimensional). The index of $u$ is the Fredholm index of this operator, and $u$ is said to be Fredholm regular when this operator is surjective. Observe that $u^*TW=v^*(T(Y\times I))\oplus \langle \partial_a,R_0\rangle,$ a $J_0$-invariant splitting, and observe that we can identify the normal bundle of $u$ with $N_u:=v^*T(Y\times I)$. The normal component of the Cauchy--Riemann operator $\mathbf{D}_u$ is then 
$$
\mathbf{D}_u^N:W^{1,2,\delta_0}(N_u)\rightarrow L^{2,\delta_0}(\overline{\mbox{Hom}}_\mathbb{C}(\mathbb{R}\times S^1,N_u)),
$$
defined as $\mathbf{D}_u^N=\pi_N \mathbf{D}_u,$ where $\pi_N:u^*TW \rightarrow N_u$ is the natural projection. A result from \cite{Wen1} tells us that the dimension of the cokernel of the operator $\mathbf{D}_u^N$  coincides with that of $\mathbf{D}_u$, so Fredholm regularity of $u$ is equivalent to surjectivity of $\mathbf{D}_u^N$. 

Finally, the linearised Morse operator for $\gamma$ is given by
$$
\mathbf{D}_\gamma^{\textbf{Morse}}:W^{1,2}(\gamma^*T(Y\times I))\rightarrow L^2(\gamma^*T(Y\times I))
$$
$$
\mathbf{D}_\gamma^{\textbf{Morse}}\eta= \nabla_s \eta - \nabla_\eta \nabla H
$$

All of this operators are Fredholm, and 
\begin{equation}\label{eqind}
ind\left(\mathbf{D}_u\right)=ind\left(\mathbf{D}_v^{\textbf{Floer}}\right)=ind\left(\mathbf{D}_v^{\textbf{Morse}}\right)=ind_{p_+}(H)-ind_{p_-}(H)
\end{equation}

Equations (\ref{parders}) and (\ref{Floereq}) now read
\begin{equation}\label{parders2}
\begin{array}{rl}
u_s=&\partial_a+v_s\\
u_t=&\partial_\theta+v_t\\
v_s+J_0v_t=&\nabla H
\end{array}
\end{equation}
 
If $\eta \in W^{1,2}(v^*T(Y\times I))$, using (\ref{parders2}), and that $\nabla \partial_a=\nabla \partial_\theta=0$, we can compute 

\begin{equation}\label{coincide}
\begin{split}
\mathbf{D}_v^{\textbf{Floer}}\eta=&\nabla_s\eta+J_0(v)\nabla_t\eta+\nabla_\eta(J_0(v)v_t)-J_0(v)\nabla_\eta v_t-\nabla_\eta \nabla H\\
=&\nabla_s\eta+J_0(v)\nabla_t\eta+\nabla_\eta(\nabla H-v_s)-J_0(v)\nabla_\eta v_t-\nabla_\eta
\nabla H\\
=&\nabla_s\eta+J_0(v)\nabla_t\eta-\nabla_\eta v_s - J_0(v)\nabla_\eta v_t\\
=&\nabla_s\eta+J_0(u)\nabla_t\eta-\nabla_\eta u_s - J_0(u)\nabla_\eta u_t\\
=&\nabla_s\eta+J_0(u)\nabla_t\eta+\nabla_\eta (J_0(u)u_t) - J_0(u)\nabla_\eta u_t\\
=&(\mathbf{D}^N_u\eta)\partial_s
\end{split}
\end{equation}

From Floer theory we know that, if $H$ is sufficiently $C^2$-small, then every closed Hamiltonian orbit is a critical point of $H$, the solutions to the corresponding Floer equation are $S^1$-independent Morse flow lines, and every element in the kernel of $\mathbf{D}_v^{\textbf{Floer}}$ is $t$-independent, so it actually defines an element in the kernel of $\mathbf{D}_v^{\textbf{Morse}}$ (e.g.\ see \cite{SZ}). This means that $\dim \ker \mathbf{D}_v^{\textbf{Floer}}=\dim \ker \mathbf{D}_v^{\textbf{Morse}}$, so (\ref{eqind}) implies that $\mathbf{D}_v^{\textbf{Floer}}$ is surjective if and only if $\mathbf{D}_v^{\textbf{Morse}}$ is, and the latter is equivalent to the Morse--Smale condition. Moreover, the above computation is also true if we replace $u$ with an unbranched multiple cover of the form $u^m(s,t)=(s,\gamma(s),mt)$. This, combined with (\ref{coincide}) (and the fact that two $(0,1)$-forms coincide iff they coincide on $\partial_s$), we conclude that:

\begin{itemize}
\item Every $J_0$-holomorphic cylinder coincides with one of the $u_\gamma$ up to multiple covers and to reparametrization (and the parametrization $u_\gamma(s,t)=(s,\gamma(s),t)$ is $J_0$-holomorphic, which is compatible with equations (\ref{floweqs}) by setting $\epsilon=0$).

\item The Morse--Smale condition implies that every unbranched cover of $u_\gamma$ is Fredholm regular.
\end{itemize}

To obtain regularity for sufficiently close contact data, we use the implicit function theorem. First, observe that we have constructed a smooth path of compatible almost complex structures $\{J_\epsilon\}$. Then, the associated $\widehat{\mathcal{H}}_\epsilon$-energy $\mathbf{E}(u)=\sup_{\varphi \in \mathcal{P}}\int u^* \omega_\epsilon^\varphi$ of a $J_\epsilon$-holomorphic curve $u$  depends smoothly on the parameter $\epsilon$. And second, that for $p$ a critical point of $H$, the Reeb vector field $R_\epsilon$ coincides with $e^{-\epsilon H(p)}\partial_\theta$, so that $\gamma_p=\{p\}\times S^1$ is a $\widehat{\mathcal{H}}_\epsilon$-Reeb orbit of action $2\pi e^{\epsilon H(p)}$, for \emph{every} $\epsilon\geq 0$. Therefore, if $u_{k}$ is an unbranched cover of a $J_{\epsilon_k}$-holomorphic flow-line cylinder in our foliation where $\epsilon_k\rightarrow 0$, all having the same asymptotics, Stokes' theorem gives bounds on the energy, and Gromov compactness tells us that $u_k$ admits a subsequence which converges to a $J_0$-holomorphic building $\mathbf{u}_\infty$. Its levels consist of unbranched covers of flow line cylinders, by the uniqueness for $\epsilon=0$. Since these are Fredholm regular for $\epsilon=0$, by the implicit function theorem we may perturb to $\epsilon>0$ small. Since Fredholm regularity is an open condition, $u_k$ will be also regular for sufficiently large $k$. We conclude that the unbranched flow line cylinders $u^m_\gamma$ are regular for sufficiently small $\epsilon>0$. 

\subsubsection{Regularity for genus zero page-like curves in the Morse--Bott case}\label{reggenuszero}

In this section, we fix a genus zero curve $u:=u_{y,a}^-$ in the foliation, and denote by $u_{Y}:=u^{-1}(\mathbb{R}\times M_Y)$, $u_C:=u^{-1}(\mathbb{R}\times M_C^-)$, and $u_P:=u^{-1}(\mathbb{R}\times M_P^-)$. We will show that $\dim\ker \mathbf{D}_u^N=\mbox{ind}\;\mathbf{D}_u^N$, where $\mathbf{D}_u^N$ is the normal component of the linearized Cauchy-Riemann operator, from which regularity follows, and from which we see that the only solutions to the normal linearized Cauchy-Riemann equations correspond to nearby curves to $u$ in the foliation, which are already known to lie in $\ker \mathbf{D}_u^N$. Morally, regularity is a ``linearized'' version of uniqueness.

We will deal with the Morse--Bott case, where $H$ depends only on $r$. In this case, the operator we need to look at is given by 
$$
\mathbf{D}_u: W^{1,2,\delta_0}(u^*TW)\oplus  V_\Gamma \oplus X_\Gamma \rightarrow L^{2,\delta_0}(\overline{Hom}_{\mathbb{C}}((\dot{\Sigma}_-,j_-),(u^*TW,J))
$$
$$
\mathbf{D}_u\eta=\nabla \eta + J(u)\circ \nabla \eta \circ j_- + (\nabla_\eta J)\circ du \circ j_-
$$
Here, $\delta_0$ is a small weight making the operator Fredholm, $V_\Gamma$ is a $2k$-dimensional vector space of smooth sections asymptotic to constant linear combinations of $R_\epsilon$ and $\partial_a$, and $X_\Gamma$ is a $k(2n-1)$-dimensional vector space of smooth sections which are supported along $u_C\cup u_Y$ (a disjoint union of its cylindrical ends, which we denote $u_i$, $i=1,\dots, k$), and are constant equal to a vector in $T_yY$ along $u_Y$. We also have the operator
$$
\mathbf{D}_{(j_-,u)}\overline{\partial}_J: T_{j_-}\mathcal{T} \oplus W^{1,2,\delta_0}(u^*TW)\oplus V_\Gamma \oplus X_\Gamma \rightarrow L^{2,\delta_0}(\overline{Hom}_{\mathbb{C}}((\dot{\Sigma}_-,j_-),(u^*TW,J))
$$
$$
(Y,\eta)\mapsto J \circ du \circ Y + \mathbf{D}_u \eta,
$$ 
where $\mathcal{T}$ is a Teichm\"uller slice through $j_-$ (see e.g.\ \cite{Wen1} for a definition of this). The curve $u$ is said to be regular whenever $\mathbf{D}_{(j_-,u)}\overline{\partial}_J$ is surjective. By a result in \cite{Wen1}, this is equivalent to the surjectivity of its normal component. 
In this case, this operator is
$$
\mathbf{D}_u^N:W^{1,2,\delta_0}(N_u)\oplus X_\Gamma \rightarrow L^{2,\delta_0}(\overline{Hom}_{\mathbb{C}}(T\dot{\Sigma}_-,N_u))
$$
$$
\eta \mapsto \pi_N \mathbf{D}_u=\pi_N(\overline{\partial}_J\eta + (\nabla_\eta J)\circ du \circ j_-), 
$$
where $\pi_N$ is the orthogonal projection to the normal bundle $N_u$. Recall that the latter is any choice of $J$-invariant complement to the tangent space to $u$, which coincides with the contact structure at infinity. Observe that $\mbox{ind}\;\mathbf{D}_u^N=2n,$ which morally follows by the fact that $(y,a)$ varies in a $2n$-dimensional space, and it's backed by the Riemann-Roch formula in the Morse-Bott case (using the weight to perturb the Hessian of $H$):
\begin{equation}
\begin{split}
\mbox{ind}\;\mathbf{D}_u^N&=\mbox{ind}(u)\\
&=(n-2)(2-k)+2(2-k)+\sum_{i=1}^k \mbox{ind}(2\pi\nabla_{-2\delta}^2H - \delta_0 id)-n\\
&=n(2-k)+\sum_{i=1}^k n\\
&=2n
\end{split}
\end{equation}

The operator $\mathbf{D}_u$ and $\mathbf{D}_u^	N$ are of Cauchy-Riemann type, which means in particular that they satisfy the Leibnitz rule 
\begin{equation}\label{Leib}
\mathbf{D}_u(f \eta)=f \mathbf{D}_u \eta+ \overline{\partial}f \otimes \eta
\end{equation}

\subparagraph*{Splitting over the paper}
      
We will think of the punctured surface $\dot{\Sigma}_-$ as being obtained abstractly from the surface with boundary $\Sigma_-\backslash\mathcal{N}(\partial\Sigma_-)$ by attaching cylindrical ends. Observe first that over the region $Y \times \Sigma_-\backslash\mathcal{N}(\partial\Sigma_-)$, we have a splitting $$u^*T(\mathbb{R}\times M)=T\Sigma_-\oplus T_{(a,y)}(\mathbb{R} \times Y)$$ 
Since $J$ preserves this splitting, this gives an identification $N_u=T_{(a,y)}(\mathbb{R}\times Y)=(\xi_-)_y\oplus\langle \partial_a, R_\epsilon(y)=R_-(y)/K \rangle$ of the normal bundle of $u$ along this region. We also have that $E_p\:=\mathbb{R}\times Y \times \{p\}$ is a holomorphic hypersurface, for every $p \in \Sigma_-\backslash\mathcal{N}(\partial \Sigma_-)$, all of which may be canonically identified with each other. This implies that the linearised Cauchy--Riemann operator $\mathbf{D}_u$ preserves this splitting, as can be seen as follows. 

Indeed, fix any such $p$ and let $E:=E_p$. Let $u^\tau=(u^\tau_\Sigma,u^\tau_E): \dot{\Sigma}_- \rightarrow \Sigma \times E$ be a 1-parameter family of maps of class $W^{1,2}$, for $\tau \in (-\epsilon,\epsilon)$, defined on the domain of $u$, such that $u^0=u$ and $\partial_\tau\vert_{\tau=0}u^\tau=(\partial_\tau\vert_{\tau=0} u^\tau_\Sigma,\partial_\tau\vert_{\tau=0} u^\tau_E)=(v,w) \in W^{1,2}(T\Sigma_-\oplus T_{(a,y)}E)$ over $\dot{\Sigma}_-\backslash\mathcal{N}(\partial \Sigma_-)$. Since over this region we have $u=(\mbox{id},(y,a))$, We obtain 
\begin{equation*}
\begin{split}
\mathbf{D}_u(v,w)&= \nabla_{\tau}(\overline{\partial}_J u^\tau)\vert_{\tau=0} \\
& = \nabla_{\tau}(\overline{\partial}_{j_-} u^\tau_\Sigma + \overline{\partial}_J u^\tau_E)\vert_{\tau=0}\\
& = (\mathbf{D}_{id}\overline{\partial}_{j_-}) v + (\mathbf{D}_{(y,a)}\overline{\partial}_J)w
\end{split}
\end{equation*} 

This yields a splitting 

$$\mathbf{D}_u=\left(
\begin{array}{cc}
\mathbf{D}_{id}\overline{\partial}_{j_-} & 0\\
0 & \mathbf{D}_{(y,a)}\overline{\partial}_J
\end{array}\right)$$

Moreover, the trivial bundle with fiber $T_{(a,y)}(\mathbb{R}\times Y)$ also splits into complex subbundles as $T_{(a,y)}(\mathbb{R}\times Y)=(\xi_-)_y \oplus \langle \partial_a,R_\epsilon(y)=R_-(y)/K\rangle$. If we take a Hermitian basis $v_1,Jv_1,\dots,v_{n},Jv_n$ of $(\xi_-)_y$, and we extend it by $\partial_a$ and $R_\epsilon(y)$, all of these vector fields (and constant linear combinations of these) generate a 1-parameter family of holomorphic maps of the form $\{u_{a^\prime(\tau),y^\prime(\tau)}^-\}_{\tau \in (-\epsilon,\epsilon)}$ for some $a^\prime(\tau), y^\prime(\tau)$, so that $\mathbf{D}_u$ vanishes on this basis. If we write any section of $T_{(a,y)}(\mathbb{R}\times Y)$ as a linear combination, we see, by the Leibnitz rule (\ref{Leib}), that $\mathbf{D}_u$ acts on this bundle as the standard Cauchy--Riemann operator. So, we may also further split the operator, by using the basis, as  

\begin{equation} 
\label{splitting}
\mathbf{D}_u=\left(
\begin{array}{cccc}
\mathbf{D}_{id}\overline{\partial}_{j_-} & 0 & \hdots & 0\\
0 & \overline{\partial} & \hdots & 0 \\
\vdots & \vdots  & \ddots & \vdots\\
0 & 0  & \hdots & \overline{\partial}
\end{array}\right),
\end{equation}
where $\overline{\partial}$ is the standard Cauchy--Riemann operator on the trivial line bundle over $\dot{\Sigma}_-$.  

In terms of the normal operator, (\ref{splitting}) may be rewritten as 
$$\mathbf{D}_u=\left(
\begin{array}{cc}
\mathbf{D}_{id}\overline{\partial}_{j_-} & 0\\
0 & \mathbf{D}_u^N
\end{array}\right),$$

where the normal Cauchy--Riemann operator $\mathbf{D}_u^N$ splits as $\mathbf{D}_u^N=\bigoplus \overline{\partial}$.

\subparagraph*{Some technical assumptions} In order to be able to write down a manageable expression for $\mathbf{D}_u^N$ over the rest of the regions, we will assume, without loss of generality, that:
\begin{assumptions}\label{Assumptions}
\begin{enumerate}$\;$
\item[A.] $H$ has a unique critical point \emph{away} from the neighbourhood $(-\delta,\delta)$ where $J_0$ is non-cylindrical. Away from this neighbourhood, $\widetilde{\xi}_r=\xi_-$ (according to the sign of $r$). Choose, say, $H^\prime(-2\delta)=0$.
\item[B.] Choose our coordinate $r$ so that the Liouville vector field $V$ coincides with $\pm\partial_r$ on the complement of $(-\delta,\delta)$.
\end{enumerate}
\end{assumptions}

In other words, we ``bias'' the foliation towards the ``minus'' side. Remember that we already know that the curves in the ``plus'' side are not regular. 

These assumptions are only used in this section to show regularity, and do not affect other sections. Therefore we will, for simplicity, lift them in the rest of the sections.

\subparagraph*{Choosing a suitable normal bundle.} We now specify how we will extend our normal bundle $N_u$ to $u$ along its cylindrical ends. 

Recall that the contact structure $\xi_\epsilon$ is given by $\xi_\epsilon=\ker(\epsilon \alpha + d\theta)$ along $M_Y$, which is isomorphic to $T(Y\times I)$ by projection along the $\theta$-direction. Moreover, $J$ was chosen on $\xi_\epsilon\vert_{M_Y}$ so that the identification $(\xi_\epsilon,J)\rightarrow (T(Y\times I),J_0)$ is holomorphic, where $J_0$ may be \emph{any} $d\alpha$-compatible almost complex structure in $Y\times I$ which is cylindrical at the ends (in particular, the ones that we constructed explicitly in Chapter 2). So we will identify these two bundles, so that $J\vert_{\xi_\epsilon}=J_0$ along $M_Y$. Observe that $J$ is always $\theta$-independent. And here we use assumptions A and B: we can choose $J_0$ so that it is \emph{cylindrical} in the complement of $(-\delta,\delta)$. That is, we take $J_0=J^-\oplus i$ with respect to the splitting $T(Y\times I)=\xi_-\oplus\langle V=\partial_r, R_-\rangle$, for some $d\alpha$-compatible $J^-$ in $\xi_-$ (Here, $V$ is the Liouville vector field associated to $d\alpha$). And the important thing to observe is that $J_0$ is \emph{$r$-independent} ($r$ is the Liouville coordinate). This means that $J$ is also, under the above complex indentification.

We then choose $N_u$ as follows:

\begin{itemize}
\item $N_u=\langle \partial_a,R_-(y) /K\rangle \oplus (\xi_-)_y = T_{(a,y)}(\mathbb{R}\times Y)$ along 
$u_P$ (recall $u\vert_{u_P}$ is constant in the $\mathbb{R}$-direction). 

\item $N_u=T_{(y,-2\delta)}(Y\times I)=\langle \partial_r, R_-(y)/K \rangle \oplus (\xi_-)_y$ along $u_{Y}$. 

\end{itemize}

Since $u_{Y}$ is contained in the cylindrical ends of $u$, the second expression for $N_u$ should be considered its \emph{asymptotic} expression. Also, we can find a suitable interpolation $(-\delta,\delta)\ni\rho \mapsto v_\rho:=v$ along the corner $M_C^-=Y\times (-\delta,\delta) \times S^1$, between $\partial_a=v_{-\delta}$ and $\partial_r=v_{\delta}$, such that $v_\rho$ is always normal to $u$ for every $\rho \in (-\delta,\delta)$, and $N_u=\langle v_\rho, R_-(y)/K\rangle \oplus (\xi_-)_y$ along $u_C$. Indeed, just take $v=v_\rho = -J(R_-(y)/K)$ ($J$ depends on $\rho$ along $M_C^-$). So we just write a global expression
$$
N_u=\langle v,R_-(y)/K\rangle \oplus (\xi_-)_y=\langle v \rangle \oplus T_yY,
$$
where the first splitting is $J$-complex. Observe that $N_u$ is a trivial $J$-complex bundle. 

\subparagraph*{Writing down the normal Cauchy-Riemann operator globally.} We will write down an expression for $\mathbf{D}_u^N$, exploiting our choice of normal bundle. We have already shown that, since $J$ splits along $M_P^-$, we get $\mathbf{D}_u^N=\bigoplus \overline{\partial}$ along $u_P$. We will compute an asymptotic expression for $\mathbf{D}_u^N$. First, we observe that we may parametrize each cylindrical end $u_i$ of $u$ ($i=1,\dots,k$) along $\mathbb{R}\times M_Y$ as
$$
u_i:[0,+\infty)\times S^1 \rightarrow \mathbb{R}\times M_Y=\mathbb{R}\times Y \times I \times S^1 
$$
$$
u_i(s,t)=(a(s),y,r(s),t),
$$
where $\dot{r}(s)=H^\prime(r(s))$, $r(0)=-1$, $\dot{a}(s)=e^{\epsilon H(r(s))}$, $a(0)=a$. In the Morse-Bott case, this parametrization is holomorphic, which follows from equations (\ref{floweqs}), together with Equation (\ref{grad}) (since $H=H(r)$).

From our Conley-Zehnder index computations of Section \ref{indexcomp}, we know that, in the Morse-Bott case, the asymptotic operator at a Reeb orbit $\gamma_p=\{p\}\times S^1$ for $p=(-2\delta,y)$, $y \in Y$, is:
$$
\mathbf{A}_{\gamma_p}=-i\partial_t+2\pi \nabla^2_{-2\delta} H
$$
This operator is degenerate, but the weight $\delta_0$ in the Sobolev spaces fixes that. The operator $\mathbf{D}_u^	N$ is \emph{conjugate} to a Fredholm operator $\mathbf{D}_{\delta_0}$ defined on a space of sections of Sobolev class $W^{1,2}$, something we have used implicitly already in the computation of $\mbox{ind}\; \mathbf{D}_u^N$, which is nothing else than $\mbox{ind}\;\mathbf{D}_{\delta_0}$. The operator $\mathbf{D}_{\delta_0}$ is asymptotic to the slightly perturbed asymptotic operator 
$$
\mathbf{A}_{\gamma_p}^{\delta_0}=\mathbf{A}_{\gamma_p}-\delta_0id,
$$
where $\delta_0$ needs to be smaller than the \emph{spectral gap} of $\mathbf{A}_{\gamma_p}$. This means that $\mathbf{D}_{\delta_0}$ that can be written asymptotically (i.e in the cylindrical coordinates $(s,t)$ along $u_{Y}$) as
$$
\mathbf{D}_{\delta_0}\eta(s,t)= \overline{\partial}\eta(s,t) + S^{\delta_0}_H(s,t)\eta(s,t), 
$$
where $S^{\delta_0}_H(s,t)$ is a symmetric matrix such that $S^{\delta_0}_H(s,\cdot)\rightarrow -2\pi \nabla_{-2\delta}^2H+\delta_0 id$ as $s\rightarrow +\infty$, uniformly in the second variable. We have a similar asymptotic expression for $\mathbf{D}_u^N$, which corresponds to the matrix $S_H:=S^0_H$ (for $\delta_0=0$), which is nothing else that the coordinate representation of $\nabla J$.

Moreover, the symmetric connection $\nabla$ we constructed along $\mathbb{R}\times Y \times [-1,-\delta)\times S^1$, under assumptions A and B, becomes the Levi-Civita connection of the metric
$$
g=\left(\begin{array}{ccc} 1  & 0 & 0\\
0& g_{J_0,d\alpha} & 0\\
0 & 0 & 1  
\end{array}\right) = \left(\begin{array}{ccccc} 1  & 0 & 0 & 0 & 0 \\
0& g_{\xi_-}  &0 & 0 & 0 \\
0 & 0 & 1 & 0 & 0 \\
0 & 0 & 0 & 1 & 0 \\
0 & 0 & 0 & 0 & 1 
\end{array}\right),$$ where we use the splitting $$T(\mathbb{R} \times Y\times [-1,-\delta) \times S^1)= \langle \partial_a \rangle \oplus \xi_- \oplus \langle R_-, \partial_r \rangle \oplus \langle \partial_\theta \rangle,$$ and where $g_{\xi_-}=g_{J_0,d\alpha}\vert_{\xi_-}$. 

The takeaway is that $\nabla$, as well as $J$, are both \emph{independent of the coordinates $a,r,$ and $\theta$} along $\mathbb{R}\times M_Y$. It follows that $S_H=\nabla J$ is independent of both $s$ and $t$, and so
$$
S_H=-2\pi H^{\prime\prime}(-2\delta)e_{11},
$$
where $e_{11}$ denotes the matrix with a $1$ at the $(1,1)$-entry, and zero everywhere else. Since $-2\delta$ is a local maximum of $H$, we can write $-2\pi H^{\prime\prime}(-2\delta)=: \nu>0$, which may be chosen arbitrarily small, but positive. Along the corner $M_C^-$, while $\nabla$ is independent on $\theta, \rho$ and $a$, and $J$ is independent on $\theta$ and $a$, the latter, as we have observed already, is not $\rho$-invariant. So, as we traverse the smoothened corner $u_C$, we pick up a smooth path of symmetric matrices $(-\delta,\delta)\ni\rho \mapsto S_H(\rho)$. By choosing $v_\rho=-J(R_-(y)/K)$, it takes the form $S_H(\rho)=s_H(\rho)e_{11}$, where $s_H(-\delta)=0$ and $s_H(\delta)=\nu$, so that 
$$
\mathbf{D}_u^N\eta(s,t)=\overline{\partial}\eta(s,t)+S_H(\rho(s))\eta(s,t)
$$ along $M_C^-$. 

The upshot of the discussion is that 
$$
\mathbf{D}_u^N=\overline{\partial} + s_He_{11},
$$
where $s_H$ is supported along the cylindrical ends, and is asymptotically constant, equal to a small positive number $\nu$.

\subparagraph*{Computing the kernel.} We now investigate $\ker \mathbf{D}_u^N$, and check that all nearby holomorphic curves in the foliation indeed give solutions to $\mathbf{D}_u^N \eta=0$. The latter curves give a $2n$-dimensional subspace of $\ker \mathbf{D}_u^N$, and we will show that it coincides with $\ker \mathbf{D}_u^N$, from which regularity follows. 

We write any section as
$$
\eta=(\eta_0,\eta_\xi, \eta_\Gamma) \in W^{1,2,\delta_0}(\langle v, R_-(y)/K\rangle) \oplus W^{1,2,\delta_0}((\xi_-)_y) \oplus X_\Gamma, 
$$ 
with $\eta_0=(\eta_0^1,\eta_0^2) \in W^{1,2,\delta_0}(\langle v \rangle) \oplus  W^{1,2,\delta_0}(\langle R_-(y)/K\rangle)$. We denote 
$$
X_\Gamma^R:=\pi_R X_\Gamma
$$
$$
X_\Gamma^\xi:=\pi_\xi X_\Gamma,
$$
where 
$$\pi_R:T_yY \rightarrow \langle R_-(y)/K \rangle
$$
$$
\pi_\xi: T_yY \rightarrow (\xi_-)_y
$$ are the orthogonal projections with respect to the metric $g$. Then we write
$$
\eta_\Gamma=(\eta_\Gamma^R,\eta_\Gamma^\xi) \in X_\Gamma^R \oplus X_\Gamma^\xi,
$$ where $\eta_\Gamma^R=\pi_R \eta_\Gamma,$ $\eta_\Gamma^\xi=\pi_\xi \eta_\Gamma$. 

Then, $\eta \in \ker \mathbf{D}_u^N$ if and only if $\overline{\partial}\eta =-s_He_{11}\eta$, which is a small perturbation in the $e_{11}$-direction of the classical holomorphic equation, and coincides with it along $u_P$. In holomorphic coordinates $(s,t)$, which we always take cylindrical in the cylindrical ends, it is equivalent to the set of equations 
\begin{equation}\label{kerD}
\left\{\begin{array}{l} \partial_s\eta_0^1-\partial_t\left(\eta_0^2+\eta_\Gamma^R\right)=-s_H\eta_0^1\\
\partial_t\eta_0^1+\partial_s\left(\eta_0^2+ \eta_\Gamma^R\right)=0\\
\overline{\partial}\left(\eta_\xi + \eta_\Gamma^\xi\right)=0
\end{array}\right.
\end{equation}

Observe that the $\eta_\Gamma$ terms all disappear away from $u_C$. 

Let us check that indeed all the nearby holomorphic curves in the foliation satisfy the above equations. Observe first that, in general, the section $\eta_\xi + \eta_\Gamma^\xi$ is always bounded, and holomorphic. It follows that necessarily $\eta_\xi+\eta_\Gamma^\xi\equiv const$, which is the $(\xi_-)_y$ component of $\eta$. This accounts for constant push-offs of $u$ in the $(\xi_-)_y$-direction, which already correspond to curves in our foliation. Moreover, any small constant push-off of $u$ in the $Y$-directions correspond to a holomorphic curve $u_{y^\prime,a}$, which at a linear level corresponds to a constant section $\eta\equiv Y_0 \in T_yY\subset N_u$. Write $Y_0=(Y_0^R,Y_0^\xi)=(\pi_R Y_0, \pi_\xi Y_0)$. Let $\beta$ be a smooth bump function which vanishes along $u_P$, and is constant equal to $1$ along $u_{Y}$. We can write
$$
Y_0=\beta Y_0 + (1-\beta)Y_0 \in X_\Gamma \oplus W^{1,2,\delta_0}(T_yY), 
$$   
so it defines an element in the domain of $\mathbf{D}_u^N$, with $\eta_0^1=0$, $\eta_0^2=(1-\beta)Y_0^R$, $\eta_\xi=(1-\beta)Y_0^\xi$, $\eta_\Gamma^R= \beta Y_0^R$, $\eta_\Gamma^\xi=\beta Y_0^\xi$. Since $\eta_0^2+\eta_\Gamma^R=Y^R_0$ and $\eta_\xi + \eta_\Gamma^\xi=Y_0^\xi$ are constant, $\eta$ satisfies equations (\ref{kerD}). 

On the other hand, we can describe asymptotically all solutions $\eta$ with $\eta_\Gamma=\eta_0^2=\eta_\xi=0$, in which case (\ref{kerD}) simplifies to $\partial_s\eta_0^1=-s_H\eta_0^1$, and $\partial_t\eta_0^1=0$. Asymptotically they are of the form 
\begin{equation}\label{decays}
\eta=(e^{-\nu s}\eta_0^1(0),0,0),
\end{equation} which decay exponentially since $\nu>0$. A translation of $u$ by an amount $A_0$ in the $\mathbb{R}$-direction  corresponds to such linear solutions $\eta$, satisfying $\eta\equiv\eta_0^1\equiv A_0$ along $u_P$, and (\ref{decays}) with $\eta_0^1(0)=A_0$ along $M_Y$. 

Now, we show that indeed the unique solutions to equations (\ref{kerD}) are the ones we described. 

We have shown that the operator $\mathbf{D}_u^N$ splits into a direct sum
$$
\mathbf{D}_u^N:=\mathbf{D}_u^\xi \oplus \mathbf{D}_u^R,
$$
where
$$
\mathbf{D}_u^\xi: W^{1,2,\delta_0}((\xi_-)_y)\oplus X_\Gamma^\xi \rightarrow L^{2,\delta_0}(\overline{Hom}_{\mathbb{C}}(\dot{\Sigma}_-,(\xi_-)_y))
$$
and
$$
\mathbf{D}_u^R: W^{1,2,\delta_0}(\langle v, R_-(y)/K\rangle)\oplus X_\Gamma^R \rightarrow L^{2,\delta_0}(\overline{Hom}_{\mathbb{C}}(\dot{\Sigma}_-,\langle v, R_-(y)/K\rangle))
$$

We will show that both operators $\mathbf{D}_u^\xi$ and $\mathbf{D}_u^R$ are surjective, and this finishes the proof.

We have already seen that $\mathbf{D}_u^\xi=\overline{\partial}$, and have observed that its kernel is the $(2n-2)$-dimensional space of constant sections along $(\xi_-)_y$. Its index is $2n-2$, and it follows that it is surjective. Moreover, we this also follows from automatic transversality: If we write $(\xi_-)_y=\bigoplus_{i=1}^{n-1}L_i$, where $L_i=\langle v_i, Jv_i\rangle$ is a line bundle, and $v_1,Jv_1,\dots,v_n,Jv_n$ is any Hermitian basis of $(\xi_-)_y$, the first summand satisfies 
$$\mathbf{D}_u^\xi=\overline{\partial}=\bigoplus_{i=1}^{n-1} \overline{\partial}
$$
Here, the $i$-th summand $\overline{\partial}$ acts on $W^{1,2,\delta_0}(L_i)\oplus X_\Gamma^i$, where $X_\Gamma^i$ is the orthogonal projection to $L_i$ of $X_\Gamma$. Each of these is a Fredholm operator, and the Conley-Zehnder index at each puncture is
$$
\mu_{CZ}(-i\partial_t-\delta_0id)=\mbox{ind}(-\delta_0 id)-1=1
$$

Therefore, the Fedholm index of each summand $\overline{\partial}$ is $\mbox{ind}(\overline{\partial})=2-k+k=2$, and surjectivity of each summand follows either by automatic transversality (exactly as below), or by the fact that the kernel of each is the $2$-dimensional vector space of constant sections. We obtain again that $\mathbf{D}_u^\xi$ is surjective. 

The second summand satisfies
$$
\mathbf{D}_u^R=\overline{\partial}+s_H e_{11}
$$
In order to show that it is surjective, we use automatic transversality \cite{Wen1}. We need to check that $\mbox{ind }\mathbf{D}_u^R>c_1(\mathbf{D}_u^R)$, where $c_1(\mathbf{D}_u^R)$ denotes the \emph{adjusted first Chern number}, defined by 
$$
2c_1(\mathbf{D}_u^R)=\mbox{ind }\mathbf{D}_u^R-2+2g+\#\Gamma_{even},
$$
where $g$ is the genus, and $\Gamma_{even}$ is the set of punctures with even Conley-Zehnder index. 

Similarly as before, the Conley-Zehnder index of $\mathbf{D}_u^R$ at each puncture is
$$
\mu_{CZ}(-i\partial_t+2\pi \nabla_{-2\delta}^2H-\delta_0 id)= \mbox{ind}(2\pi \nabla_{-2\delta}^2H-\delta_0 id)-1=1
$$
And therefore 
$$
\mbox{ind}\;\mathbf{D}_u^R=2-k+k=2 
$$
On the other hand, the adjusted first Chern number is
$$
2c_1(\mathbf{D}_u^R)=\mbox{ind}\;\mathbf{D}_u^R -2 + \#\Gamma_{even}=0
$$
Then we have $$\mbox{ind}\;\mathbf{D}_u^R>c_1(\mathbf{D}_u^R),$$ and so $\mathbf{D}_u^R$ is surjective by automatic transversality. This finishes the proof of regularity.

\subsection{From Morse--Bott to Morse}
\label{MBtoM}

\begin{figure}[!ht] \centering
\includegraphics[width=0.40\linewidth]{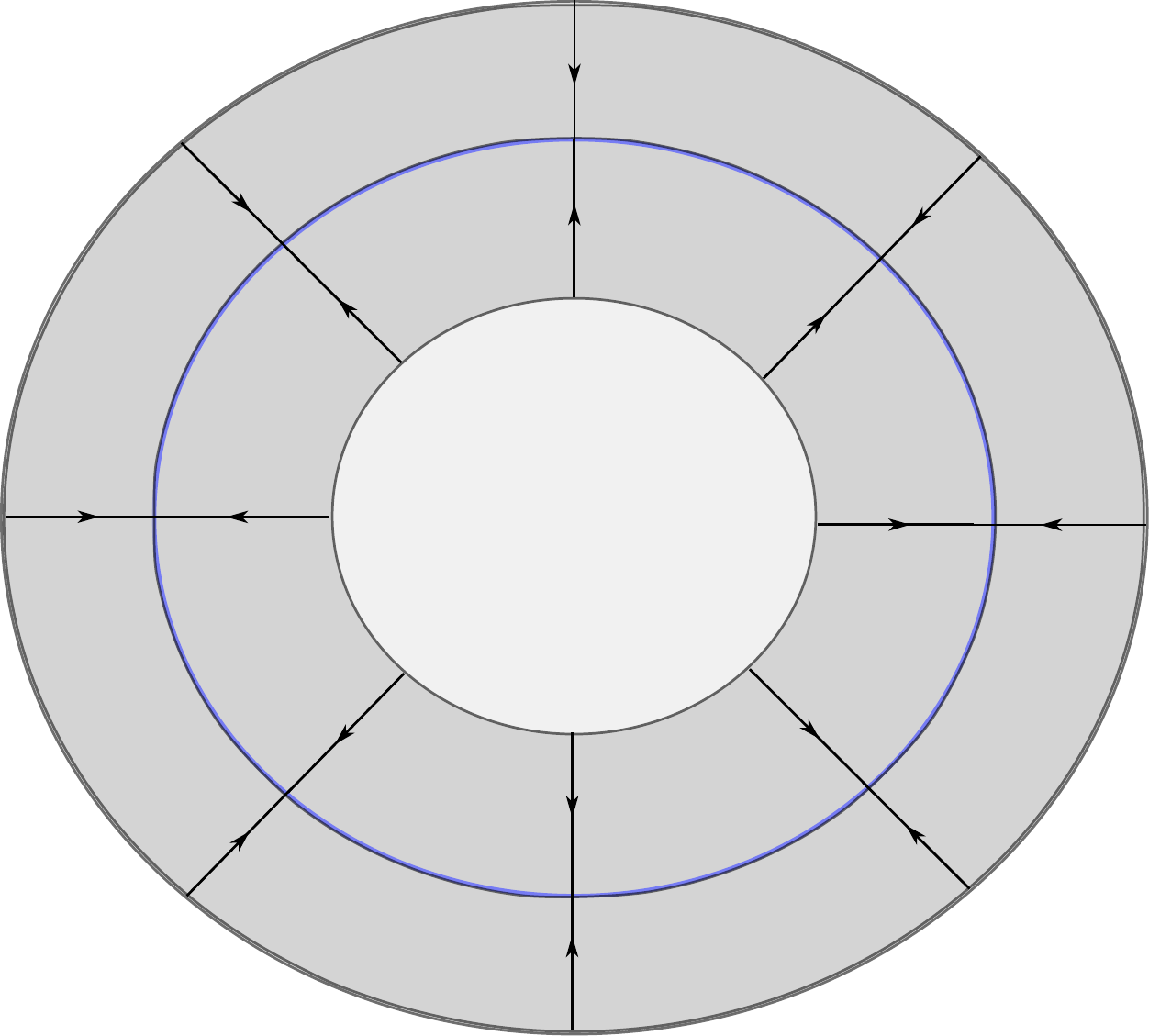}

\vspace{0.5cm}

\includegraphics[width=0.40\linewidth]{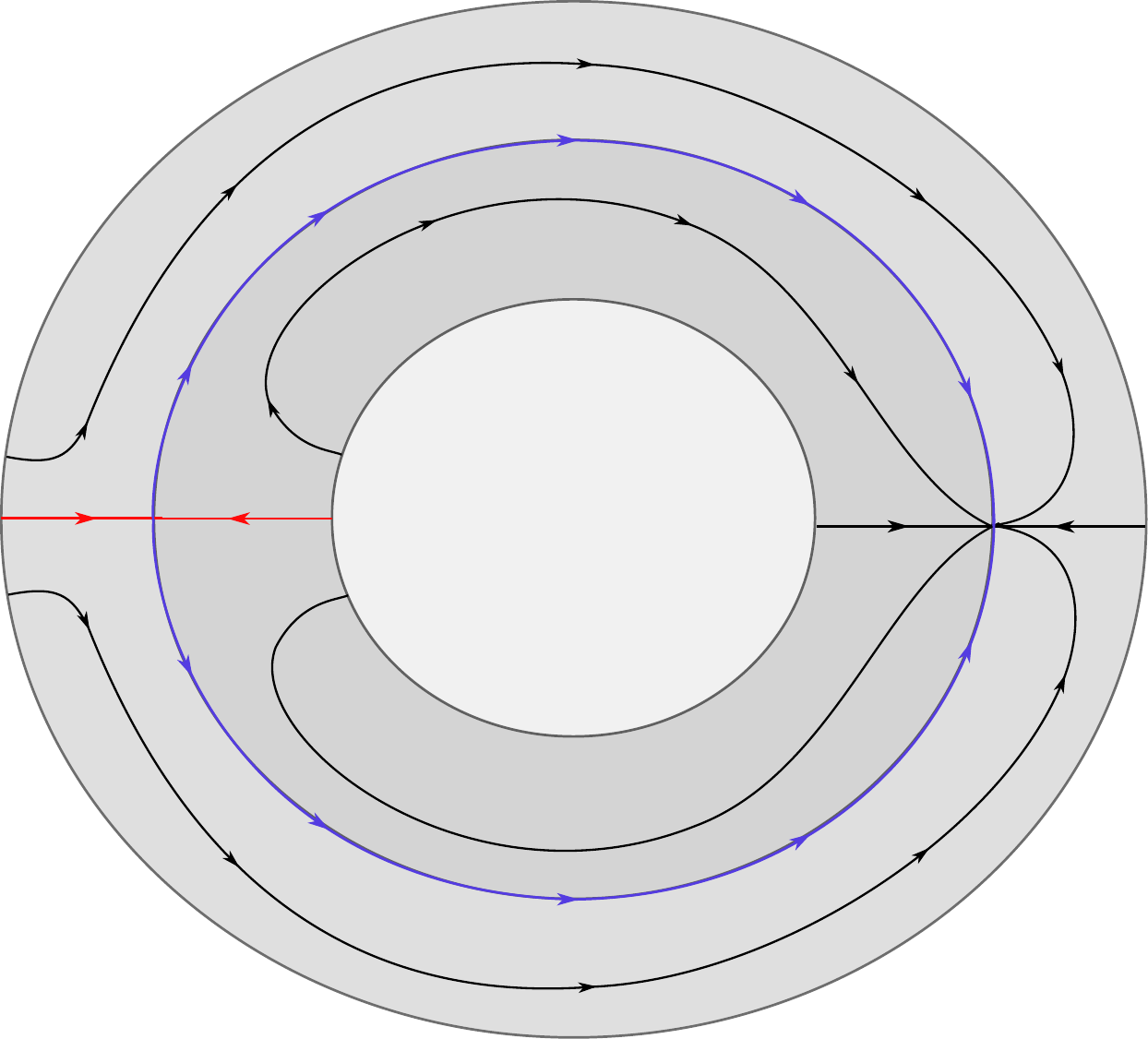}
\caption{\label{Morse} The pictures illustrate the Morse--Bott and the Morse scenarios, respectively, in the case where $Y=S^1$. The curves in the foliation are \emph{positive} Morse flow lines, which start being constant in $Y$ close to the boundary, where $H$ depends only on $r$. In the first picture, the central circle in blue is a critical manifold, whereas in the second we have two non-degenerate critical points. The ``evil twins'' are shown in blue (see Section \ref{Torsion}), going from the hyperbolic index 1 critical point $h=(min,0)$ to the maximum, whereas the unique flow lines approaching $h$ from each side are shown in red. These are the flow lines which correspond to rigid holomorphic curves that we shall count (except the red curve entering from the right, which corresponds to a non-regular positive genus curve). The generic behaviour is also seen in the picture, where a generic curve approaches the maximum at $+\infty$ and belongs to a family of flow lines of maximal dimension 2 ($2n$ in the higher dimensional setting).} 
\end{figure}

In this section, we clarify exactly how we are choosing our Morse function to do the Morse perturbation along the spine. We do this in two steps, because Morse--Bottness makes the proof of regularity manageable. First, choose $H$ to be given by $H(y,r)=f(r)$ where $f$ is a sufficiently small and positive Morse function on $I$ and has a unique critical point at 0 with Morse index 1 (which yields the Morse--Bott situation). And then, choose a sufficiently small positive Morse function $g$ on $Y$ and extend it to a neighbourhood of $Y\times \{0\}$ to make a further perturbation (obtaining the Morse case). 

The concatenation of the two steps is equivalent to having chosen $H$ to look like
$$H(y,r)=f(r)+\gamma(r)g(y),$$  
where $\gamma:I\rightarrow [0,1]$ is a smooth bump function satisfying $\gamma=0$ in the region $\{|r|>2\delta\}$ and $\gamma=1$ on $\{|r|\leq\delta\}$. 

We can also view the Morse case as a deformation of the Morse--Bott one. Indeed, we can define $H_t(y,r)=f(r)+t\gamma(r)g(y)$ for $t\in[0,1]$, a 1-parameter family of functions, which interpolates between the two scenarios. We obtain a corresponding 1-parameter families of SHS's and compatible almost complex structures $J^t$, which we can use for deformation arguments. In the case where $g$ is chosen small, then the Morse case can be viewed as a perturbation of the Morse-Bott case. Therefore, regularity for the curves $u^-_{y,a}$ in the Morse case, for \emph{sufficiently small $g$}, is deduced from that for the Morse--Bott one, by the implicit function theorem. Observe that this implies in particular that $ind(u^-_{y,a})\geq 1$, for every such curve. We state this in a theorem:

\begin{thm}[Fredholm regularity in the nearby Morse case]
For Morse data sufficiently close to Morse-Bott data, all the genus zero curves in the finite energy foliation are Fredholm regular. 
\end{thm}

In order to simplify the torsion computation of Section \ref{Torsion}, we will choose $g$ to have a unique maximum $max$ and minimum $min$. Both scenarios are depicted in Figure \ref{Morse} in the case of $Y=S^1$.

\subsection{Uniqueness of the curves in the Morse/Morse--Bott case}
\label{uniquenessMB}

In this section, we prove that the family of curves we constructed above are the unique curves (up to reparametrization and multiple covers) that asymptote Reeb orbits of the family $\{\gamma_p: p \in \mbox{crit}(H)\}$, and with positive asymptotics in different components of $M_Y$. We do this in both the Morse/Morse--Bott situations. We assume that $H$ has a unique critical point in the interval direction, at $r=0$ (and perhaps other critical points contained in $Y\times \{0\}$, in the Morse case). In the Morse--Bott case, we denote by $\gamma_y:=\gamma_{(y,0)}=\{(y,0)\}\times S^1$ the simply covered Reeb orbit corresponding to $y \in Y$.

\begin{lemma}\label{posvsneg}
Assume either the Morse or Morse--Bott cases. Let $N\in \mathbb{N}^{>0}$, and let $u$ be a $J$-holomorphic curve with positive asympotics of the form $\gamma_p=\{p\}\times S^1$, for which the number of positive punctures is bounded above by $N$. Then, we can find sufficiently small $\epsilon>0$ (depending only on $N$), such that $$\#\Gamma^-(u)\leq \#\Gamma^+(u),$$ where $\Gamma^\pm(u)$ denotes the set of positive/negative punctures of $u$. Here, we count punctures with the covering multiplicity of their corresponding asymptotic.  
\end{lemma}
\begin{proof}
By Remark \ref{remark11}, we have that the negative asymptotics of $u$ also correspond to critical points of $H$. If we denote by $p_z$ the critical point corresponding to the puncture $z$, the $\Omega_\epsilon$-energy of $u$ is then
$$
0\leq \mathbf{E}(u)=\int u^*\Omega_\epsilon=2\pi\left(\sum_{z\in \Gamma^+(u)}\kappa_z e^{\epsilon H(p_z)}-\sum_{z\in \Gamma^-(u)}\kappa_z e^{\epsilon H(p_z)}\right),
$$
where we denote by $\kappa_z\geq 1$ the covering multiplicity of the Reeb orbit $\gamma_{p_z}^{\kappa_z}$ corresponding to $z \in \Gamma^\pm(u)$.
In the Morse--Bott case, we have 
$$
\mathbf{E}(u)=2\pi e^{\epsilon H(0)}\left(\sum_{z\in \Gamma^+(u)}\kappa_z- \sum_{z\in \Gamma^-(u)}\kappa_z\right)= 2\pi e^{\epsilon H(0)}(\#\Gamma^+(u)-\#\Gamma^-(u)),
$$
and the claim follows.

In the Morse case, let us choose $\epsilon>0$ small enough so that 
$$
||e^{\epsilon H}||_{C^0}<1+\delta_0,
$$ 
for $\delta_0<1/N$.

If we assume by contradiction that $\#\Gamma^-(u)\geq \#\Gamma^+(u)+1$, then 
$$
\#\Gamma^+(u)+1\leq \#\Gamma^-(u)\leq \sum_{z\in \Gamma^-(u)}\kappa_z e^{\epsilon H(p_z)}\leq \sum_{z\in \Gamma^+(u)}\kappa_z e^{\epsilon H(p_z)}\leq (1+\delta_0)\#\Gamma^+(u),
$$
so that $N\delta_0\geq \#\Gamma^+(u)\delta_0\geq 1$, contradicting our choice of $\delta_0$. 
\end{proof}

\begin{thm}
\label{uniqueness}
Assume either the Morse or Morse--Bott scenario.

\vspace{0.5cm}

Let $u:\dot{\Sigma} \rightarrow \mathbb{R}\times M$ be a (not necessarily regular) $J$-holomorphic curve defined on some punctured Riemann surface $\dot{\Sigma}$ which is asymptotically cylindrical, and asymptotes simply covered Reeb orbits of the form $\gamma_p=\{p\}\times S^1$ at its positive ends. Assume that any two of the positive ends of $u$ lie in distinct components of $M_Y$. 

\vspace{0.5cm}

Then, for sufficiently small and uniform $\epsilon>0$, we have that, if $u$ is not a trivial cylinder over one of the $\gamma_p$'s, then $u$ is a curve of the form $u_{y,a}^\pm$ for some $y \in Y$, or a flow-line cylinder $u_\gamma$ (in the Morse scenario).
\end{thm}

\begin{proof}

We will consider two cases: either $u$ is completely contained in the region $\mathbb{R}\times M_Y$ (case A), or it is not (case B). 

\vspace{0.5cm}

\textbf{Case A.} This case is easily dealt with in the Morse--Bott scenario.  By assumption, we have that $u$ has a unique positive end. Since the $\gamma_p$'s are not contractible/nullhomologous inside $M_Y$, Lemma \ref{posvsneg} implies that $u$ is has one positive and one negative end, both simply covered, corresponding to Reeb orbits $\gamma_{p_\pm}$. But then the $\Omega_\epsilon$-energy of $u$ vanishes, and $u$ is necessarily is a cover of a trivial cylinder.

\vspace{0.5cm} 

In the Morse case, we show that $u$ is a flow line cylinder. Again, Lemma \ref{posvsneg} implies that $u$ is has one positive and one negative end, both simply covered, corresponding to Reeb orbits $\gamma_{p_\pm}$. Observe that, a priori, $u$ is not even necessarily a cylinder, since it may have positive genus. 

Let us first observe that in the degenerate case when $H\equiv 0$, we have that the projection $\pi_Y:\mathbb{R}\times M_Y \rightarrow Y \times I$ is holomorphic. Then, if $u$ is a holomorphic curve for this data, then so is $v:=\pi_Y \circ u$. Since the asymptotics of $u$ are covers of the $S^1$-fibers of $\pi_Y$, the map $v$ extends to a holomorphic curve in the closed surface $\Sigma$. But $Y\times I$ is exact, so that $v$ has to be constant by Stokes' theorem. This means that $u$ is necessarily a multiple cover of a trivial cylinder.

We then see that the space of \emph{stable holomorphic cascades} \cite{Bo} in $\mathbb{R}\times Y \times I \times S^1$, the objects one obtains as limits of honest curves when one turns off the function $H$, consists of finite collections of flow-line segments and covers of trivial cylinders. If we take $H_t=tH$, and we assume we have a sequence $\{u_n\}$ of $J_{t_n}$-holomorphic maps with $t_n\rightarrow 0$ (where $J^{t}$ is the almost complex structure corresponding to $H_{t}$), with one positive and negative simply covered orbits corresponding to critical points $p_\pm$, then we obtain a \emph{stable} holomorphic cascade $\mathbf{u}^H_\infty$ as a limiting object. Since the positive end of $u_n$ is simply covered, Lemma \ref{posvsneg} applied to $H=0$ implies that \emph{every} Reeb orbit appearing in $\mathbf{u}^H_\infty$ is simply covered, and therefore every of its holomorphic map components cannot be multiply covered. These can then only be trivial cylinders, but stability of the cascade means that it does not have trivial cylinder components. We conclude that the space of holomorphic cascades which glue to curves as in our hypothesis consists solely of flow-lines, which are regular by the Morse--Smale condition, and come in a $(ind_{p_+}(H)-ind_{p_-}(H)-1)$-dimensional family. The gluing results in \cite{Bo} give a 1-1 correspondence between index $d$ Morse families of curves for the non-degenerate perturbation and index $d$ regular holomorphic cascades, which in our situation is exactly what proves that our curve $u$ is a flow-line cylinder. This can be understood in terms of moduli spaces: consider the parametric moduli space
$$
\overline{\mathcal{M}}:=\overline{\mathcal{M}}(\mathbb{R}\times M_Y,\{J^t\}_{t\in [0,1]};\gamma_{p_+},\gamma_{p_-})
$$
$$
=\{(t,u): t \in [0,1], u \in \overline{\mathcal{M}}(\mathbb{R}\times M_Y,J^t;\gamma_{p_+},\gamma_{p_-})\} 
$$ 
where $\overline{\mathcal{M}}(\mathbb{R}\times M_Y,J^t;\gamma_{p_+},\gamma_{p_-})$ denotes the Gromov-compactified moduli space of $J^t$-holomorphic curves in $\mathbb{R}\times M_Y$ (of any genus) which have simply covered positive asymptotic $\gamma_{p_+}$, and negative asymptotic $\gamma_{p_-}$, and where the compactification at $t=0$ corresponds to stable holomorphic cascades. Since the Reeb orbits are simply covered, curves in the parametric moduli are somewhere injective. Then, by the Morse--Smale condition and the implicit function theorem, the parametric moduli space is, for sufficiently small $H$, a $(ind_{p_+}(H)-ind_{p_-}(H))$-dimensional compact manifold whose boundary contains $\overline{\mathcal{M}}(\mathbb{R}\times M_Y,J^0;\gamma_{p_+},\gamma_{p_-})$, only consisting of flow line cylinders. The flow-line parametric moduli space 
$$
\overline{\mathcal{M}}_{\mbox{flow-line}}:=\overline{\mathcal{M}}_{\mbox{flow-line}}(\mathbb{R}\times M_Y,\{J^t\}_{t\in [0,1]};\gamma_{p_+},\gamma_{p_-})
$$
$$
:=\{(t,u): t \in [0,1], u \mbox{ corresponds to a flow-line in } \overline{\mathcal{M}}(\mathbb{R}\times M_Y,J^t;\gamma_{p_+},\gamma_{p_-})\}
$$
is, a priori, a submanifold of $\overline{\mathcal{M}}(\mathbb{R}\times M_Y,\{J^t\}_{t\in [0,1]};\gamma_{p_+},\gamma_{p_-})$, which shares the boundary component $\overline{\mathcal{M}}(\mathbb{R}\times M_Y,J^0;\gamma_{p_+},\gamma_{p_-})$, and has its same dimension. By thinking of $\overline{\mathcal{M}}$ as a collar neighbourhood of this boundary component, we get that  
$$\overline{\mathcal{M}}=\overline{\mathcal{M}}_{\mbox{flow-line}},$$ from which our uniqueness follows. Observe that since there are finitely many critical points, we can take $H$ (or $\epsilon>0$) uniformly small.

\vspace{0.5cm}

We include two alternative proofs for cylinders, which are not formally needed, but may have independent interest:

\subparagraph*{Two extra proofs for cylinders in $\mathbb{R}\times M_Y$. 1st proof.}

Another way to prove uniqueness for \emph{cylinders} is by deforming the contact data along the spine to a SHS via a 1-parameter family $\mathcal{H}_\epsilon$ as we did before when proving regularity (see Section \ref{regcyls}), and using the fact from Floer theory that solutions to the Floer equation for a sufficiently $C^2$-small Hamiltonian are flow lines. This gives uniqueness for $\epsilon=0$. A similar moduli space argument recovers uniqueness for $\epsilon>0$ small. Indeed, for any sequence $u_{\epsilon_k}$ of $J_{\epsilon_k}$-holomorphic cylinders with $\epsilon_k\rightarrow 0$, Gromov compactness gives a building $\mathbf{u}_\infty$ as a limiting object, whose levels consist of flow line cylinders by uniqueness for $\epsilon=0$ (In the Morse--Bott case this argument also proves that no such $u$ can exist). 
If $\gamma_{p_\pm}$ are its positive/negative asymptotics, the implicit function theorem gives that the (Gromov compactified) parametric moduli space 
$\overline{\mathcal{M}_0}(\mathbb{R}\times M_Y,\{J_{\epsilon}\}_{\epsilon\in[0,\epsilon_0]},\gamma_{p_+};\gamma_{p_-})$ of holomorphic cylinders is, for sufficiently small $\epsilon_0>0$, a smooth manifold with boundary 
$$
\partial\overline{\mathcal{M}_0}(\mathbb{R}\times M_Y,\{J_{\epsilon}\}_{\epsilon\in[0,\epsilon_0]},\gamma_{p_+};\gamma_{p_-})
$$
$$
=-\overline{\mathcal{M}_0}(\mathbb{R}\times M_Y,J_0,\gamma_{p_+};\gamma_{p_-})\bigsqcup \overline{\mathcal{M}_0}(\mathbb{R}\times M_Y,J_{\epsilon_0},\gamma_{p_+};\gamma_{p_-}),
$$
and of dimension 
$$\dim \overline{\mathcal{M}_0}(\mathbb{R}\times M_Y,\{J_{\epsilon}\}_{\epsilon\in[0,\epsilon_0]},\gamma_{p_+};\gamma_{p_-})=\mbox{ind}_{p_+}(H)-\mbox{ind}_{p_-}(H)
$$  
Therefore, as before, there is a unique way of converging to $\mathbf{u}_\infty$ via non-broken simple curves, either through non-broken simple $J_0$-curves along $\mathcal{M}_0(\mathbb{R}\times M_Y,J_0;\gamma_{p_+},\gamma_{p_-})$ (non-compactified), which are necessarily flow-line cylinders, or through sequences of flow-line cylinders lying in the non-compactified moduli $\mathcal{M}_0(\mathbb{R}\times M_Y,J_{\epsilon_k},\gamma_{p_+};\gamma_{p_-})$ for $\epsilon_k\rightarrow 0$, $\epsilon_k\leq \epsilon_0$. We obtain that necessarily $u_{\epsilon_k}$ is a curve $u_\gamma$ with $\lim_{s\rightarrow \pm \infty}\gamma(s)=p_\pm$ for sufficiently large $k$. Again, $\epsilon_0>0$ can be taken uniformly small.

\subparagraph*{2nd proof, under the assumption that $ind(u)=1$, and $u$ is regular.} We exploit the $S^1$-symmetries that we have in $\mathbb{R}\times M_Y$. This approach is inspired by a similar proof in \cite{Smith}. Let $\widetilde{\mathcal{M}}_0:=\widetilde{\mathcal{M}}_0(\mathbb{R}\times M_Y;\gamma_{p_+},\gamma_{p_-})$ denote the moduli space of \emph{parametrized} holomorphic cylinders in $\mathbb{R}\times M_Y$, consisting of pairs $(j,u)$, where $j$ is a complex structure on the cylinder and $u: (\mathbb{R}\times S^1,j) \rightarrow (W,J)$ is a holomorphic curve which asymptotes the Reeb orbit $\gamma_{p_\pm}$ at $\pm \infty$. We have a right $S^1$-action on this moduli space, given by 
$$
(j,u).\alpha=(\varphi_\alpha^*j,u\circ \varphi_\alpha),
$$
where $\varphi_\alpha(s,t)=(s,t+\alpha)$ is the rotation by $\alpha \in S^1$. We also have a left $S^1$-action given by postcomposition with the natural $S^1$-action in $\mathbb{R}\times M_Y$, which preserves the almost complex structure $J$. Denote by $\mathfrak{g}=Lie(S^1)=i\mathbb{R}$ the Lie algebra of $S^1$.

We take an element $u \in \widetilde{\mathcal{M}}_0$, and we assume $ind(u)=1$ and $u$ regular, so that $\widetilde{\mathcal{M}}_0$ is a $1$-dimensional manifold, whose tangent space consists of a suitable collection of smooth vector fields along $u$. The subspace $u.\mathfrak{g}\subseteq T_u\widetilde{\mathcal{M}}_0$ generating the right $S^1$-action is nontrivial, since otherwise $u$ would be constant, being constant along a 1-dimensional submanifold and thus having non-isolated critical points. By dimensional reasons, then $u.\mathfrak{g}= T_u\widetilde{\mathcal{M}}_0$.  

For $z=(s,t)\in \mathbb{R}\times S^1$, we have an evaluation map 
$$
ev_z:T_u\widetilde{\mathcal{M}}_0\rightarrow T_{u(z)}(\mathbb{R}\times M_Y),
$$   
which evaluates an element of $T_u\widetilde{\mathcal{M}}_0$ at the point $z$. Since the Reeb orbits to which $u$ asymptotes are homogeneous for the action on $\mathbb{R}\times M_Y$, and surjectivity of a linear map is an open condition, any $z$ sufficiently close to $\pm \infty$ satisfies that the restriction of $ev_z$ to the subspace $\mathfrak{g}.u\subseteq T_u\widetilde{\mathcal{M}}_0$ which generates the left $S^1$-action maps surjectively to $T_{u(z)}L$, where $L$ is the image of the circle $\{s\}\times S^1$ containing $z$. $L$ cannot be a point, since $u$ would be constant, therefore $\mathfrak{g}.u$ cannot be a trivial subspace, so that $\mathfrak{g}.u= T_u\widetilde{\mathcal{M}}_0$, again by dimensional reasons. We conclude that 
$$
\mathfrak{g}.u=u.\mathfrak{g}=T_u\widetilde{\mathcal{M}}_0
$$
This means that for every $\xi \in \mathfrak{g}$, we can find a $\zeta \in \mathfrak{g}$ such that $u.\xi=\zeta.u$. Taking $\xi = i$, and exponentiating this equality we obtain
$$
u(e^{i \theta} z)=e^{\zeta \theta}.u(z),
$$
so that the image of $u$ is $S^1$-invariant. Since the flow-line cylinders also satisfy this, and they foliate $\mathbb{R}\times M_Y$, by choosing a point in the image of $u$ and moving it around with the action we see that $u$ must coincide with a flow-line cylinder along a non-isolated set of points, and by unique continuation, must coincide everywhere.

\vspace{0.5cm}

\textbf{Case B.} We first assume the Morse--Bott case, and deal with the Morse case via the gluing results in \cite{Bo}. The approach in this situation is to estimate the $\Omega_\epsilon$-energy of $u$, and to use a suitable branched cover argument. The details are as follows.  

\vspace{0.5cm}

Assume the Morse--Bott case. Since every positive puncture of $u$ corresponds to a critical point, Remark \ref{remark11} implies that so does every negative one. Let us denote by $\Gamma^\pm$ the set of positive/negative punctures of $u$, and for $z \in \Gamma^\pm$, let $\gamma^{\kappa_z}_{p_z}$ be the Reeb orbit corresponding to $z$, where $p_z\in \mbox{crit}(H)$, $\kappa_z\geq 1$, and $\kappa_z=1$ for $z \in \Gamma^+$. By assumption, we have that $\#\Gamma^+\leq k$, the number of components of $M_Y$. 

The $\Omega_\epsilon$-energy then has the following upper bound:

\begin{equation}
\label{firstineq}
\begin{split}
\mathbf{E}(u)=&\int_{\dot{\Sigma}}u^*\Omega_\epsilon\\
=&2\pi e^{\epsilon H(0)}\left(\#\Gamma^+-\sum_{z\in\Gamma^-}\kappa_z\right)\\
=&2\pi e^{\epsilon H(0)}(\#\Gamma^+-\#\Gamma^-)\\
\leq&2\pi\#\Gamma^+||e^{\epsilon H}||_{C^0}\\
\leq& 2\pi k ||e^{\epsilon H}||_{C^0}
\end{split}
\end{equation} 

By construction, we have that over the region $\Sigma_\pm\backslash\mathcal{N}(\partial \Sigma_\pm)$, the almost complex structure $J$ splits. This implies that the projection $$\pi_P^\pm:\mathbb{R}\times Y\times \Sigma_\pm\backslash\mathcal{N}(\partial \Sigma_\pm) \rightarrow\Sigma_\pm\backslash\mathcal{N}(\partial \Sigma_\pm)$$ is holomorphic.
Moreover, this is still true if we extend this region by adding a small collar $\{t=\mp \rho \in (-\delta, -\delta+\delta/9]\}$, as one can check. Denote by $B_\pm^\delta:=\Sigma_\pm\backslash\mathcal{N}(\partial \Sigma_\pm)\cup\{t \in(-\delta,-\delta+\delta/9]\}$.

For each $w \in B_\pm^\delta$, the hypersurface $E_w:=\mathbb{R}\times Y \times \{w\}$ is $J$-holomorphic, and $\pi_P^\pm$ has $E_w$ as fiber over $w$. Since the asymptotics of $u$ are away from $B_\pm^\delta$, the intersection of $u$ with any of the $E_w$ is necessarily a finite set of points, since they are restricted to lie in a compact part of the domain of $u$.

Assuming WLOG that $u$ indeed has a non-empty portion lying over the ``plus'' region $\Sigma_+\backslash{\mathcal{N}(\partial \Sigma_+)}$, by positivity of intersections we have that $u$ necessarily intersects every $E_w$ for $w\in B_\pm^\delta$. We may then find a $t_0\geq -\delta$, so that $\pi_P^+(u)$ is transverse to the circle $\{t_0\} \times S^1 \subseteq B_+^\delta$ (over all $k$ components of the collar), which exists by Sard's theorem, applied to $\pi_I \circ \pi_P^+$, where $\pi_I$ is the projection to the interval $I$. If we denote $B_+^{\delta;t_0}:= B_+^\delta\backslash\{t \in (t_0,-\delta+\delta/9]\}$, we have that $S_+:=(\pi_P^+\circ u)^{-1}(B_+^{\delta;t_0})\subseteq \dot{\Sigma}$ is a surface with boundary $\partial S_+=(\pi_P^+\circ u)^{-1}(\partial B_+^{\delta;t_0})$. The upshot of the discussion above is that the map
$$
F_u:=(\pi_P^+\circ u)\vert_{S_+}:S_+ \rightarrow B_+^{\delta;t_0}
$$ 
is a \emph{holomorphic branched cover}, having as degree the (positive) algebraic intersection number of $u$ with any of the $E_w$ (which is independent of $w \in B_+^{\delta;t_0}$). Call this degree $\deg^+(u):=\deg(F_u)$. We wish to show that $\deg^+(u)=1$, and so this map will be actually a biholomorphism.

\vspace{0.5cm}

Let us write $\partial S_+=\bigcup_{i=1}^l C_i$, where $C_i$ is a simple closed curve whose image under $F_u$ wraps around one of the circles $\{t_0\}\times S^1$ (a priori positively or negatively), with winding number $n_i\in \mathbb{Z}\backslash\{0\}$. Observe that necessarily one has that $l\geq k$, since $u$ intersects every $E_w$ at least once, and in particular for every $w$ on the $k$ circles $\{t_0\}\times S^1$. 

Both $B_+^{\delta;t_0}$ and $S_+$ carry the orientations given by $j_+$ and $J$, respectively, and $F_u$ preserves the orientation since it is holomorphic. Moreover, the boundary orientation on $\partial B_+^{\delta;t_0}$ induced from $j_+$ is given by $-\partial_\rho$, the outer normal, since $\partial_t=-\partial_\rho$ along this boundary, and $j_+(\partial_t)=K\partial_\theta$ by construction. Therefore the orientation on $\partial S_+$ is such that the image under $F_u$ of every $C_i$ wraps around $\{t_0\}\times S^1$ \emph{positively}, and then has $d\theta$-winding equal to $n_i>0$. By counting preimages under this projection of a point in each the $k$ circles $\{t_0\}\times S^1$, we obtain that

\begin{equation}\label{degreeplus}
\int_{\partial S_+} u^*d\theta = 2\pi\sum_{i=1}^l n_i=2\pi k \deg^+(u)
\end{equation}      

Using expression $\Omega_\epsilon=Kd\alpha_\pm+d\lambda_\pm$, equation (\ref{degreeplus}), the fact that $t_0\geq-\delta$, $u^*d\alpha_+\geq 0$, and Stokes' theorem, we have the following energy estimate:

\begin{equation}
\label{ineqs}
\begin{split}
\mathbf{E}(u)\geq& \int_{S_+}u^*\Omega_\epsilon\\
=&\int_{S_+}u^*\left(K d\alpha_++ d\lambda_+\right)\\
\geq&\int_{S_+}u^*d\lambda_+\\
=&\int_{\partial S_+}u^*(e^{t}d\theta)\\
=&2\pi k e^{t_0}\deg^+(u)\\
\geq & 2\pi k e^{-\delta}\deg^+(u)
\end{split}
\end{equation}

If we combine this with the inequality (\ref{firstineq}), we obtain 
\begin{equation}
\label{sss}
2\pi k ||e^{\epsilon H}||_{C^0} \geq 2\pi e^{\epsilon H(0)}\left(\#\Gamma^+-\#\Gamma^-\right)\geq 2\pi k e^{-\delta} \deg^+(u) \geq 2\pi k e^{-\delta}
\end{equation}
Now, since we have the freedom to choose $-\delta$ and $||\epsilon H||_{C^0}$ as close to zero as wished, so that the left and right hand sides of the above inequality are as close to each other as desired, one can easily see that $$\#\Gamma^+=k,\;\#\Gamma^-=0,\deg^+(u)=1$$ 

This proves that $u$ has no negative ends and precisely $k$ positive ends,  and that $F_u$ gives a biholomorphism between $S_+$ and $B_+^{\delta;t_0}$.

One can also see, even though we shall not need it in this proof, that the same reasoning and a string of inequalities completely analogous to (\ref{ineqs}), where we replace $S_+$ with $S_+\cup S_-$ for $S_-$ defined in the same way as $S_+$ (we can even take the same $t_0$) actually implies that $u$ has no portion lying over $\Sigma_-$; this would yield $2\pi k ||e^{\epsilon H}||_{C^0}\geq 4\pi k e^{-\delta}$, which is absurd if we choose $\delta$ and $\epsilon$ close enough to zero. 

\vspace{0.5cm}

It also follows from our energy estimates that $$\int_{S_+}u^*d\alpha_+=0$$ Since the integrand is non-negative, we get that $u^*d\alpha_+=0$, which means that the projection to $Y$ of $u(S_+)\subseteq u(\dot{\Sigma})$ lies in the image of some (not necessarily closed) $R_+$-Reeb orbit $\gamma$. In other words, $u(S_+)$ lies entirely in the almost complex 4-manifold 
$$
U:=\mathbb{R}\times \gamma \times\Sigma_+\subseteq \mathbb{R}\times Y \times \Sigma
$$

Choose now a point $y$ in the projection of $u(S_+)$ to $\gamma$, in this region. We shall prove that $u$ necessarily coincides with $u_{y,a}^+$, for some $a$, and hence there is only one such choice of $y$. We start by proving that the projection of $u(S_+)$ to $\gamma$ consists of just the point $y$, using the Morse--Bott assumption. 

Observe that using the Reeb flow along $\gamma$ we have a local coordinate $y=y(s)$, such that $y(0)=y$. Assume by contradiction that there is some $s\neq 0$ such that $y(s)$ belongs to the projection of $u(S_+)$ to $\gamma$, and consider the family of curves $\{u_{y(s),a}^+:a\in \mathbb{R}\}$. By choosing a suitable $a$, we obtain an intersection of $u_{y(s),a}^+$ and $u$, and observe that we may assume that $y(s)$ is such that $\gamma_{y(s)}$ is not an asymptotic orbit of $u$ (indeed, by positivity of intersections, $u$ intersects $u_{y(s^\prime),a^\prime}$ for $s^\prime$ and $a^\prime$ sufficiently close to $s$ and $a$, respectively, so that $y(s^\prime)$ also lies in the projection of $u$ to $Y$). Since now these two curves asymptote different Reeb orbits, they are not reparametrizations of each other, and hence intersections count positively, so that the total intersection in $U$ is at least 1.

On the other hand, the set $\mathbb{R} \times (u_M \cap (u_{y(s),a}^+)_M)$, where $u_M$ is the projection to $M$ of $u$, must be bounded in the $\mathbb{R}$-direction, since the corresponding Reeb orbits are bounded away from each other, and a sequence of intersections shooting off to infinity in this direction would yield an intersection of the orbits. This means that the standard intersection pairing is homotopy invariant (there are no contributions coming from infinity), but a priori the intersection points might move along with a given homotopy. If we choose to homotope by translating in the $\mathbb{R}$-direction, this can only happen if the projection of both curves to $M$ intersect the asymptotics of each other. However, the projection of the curve $u_{y(s^\prime),a^\prime}$ to $M$ does not intersect the asymptotics of $u$. But then one can homotope one curve away from the other by $\mathbb{R}$-translating. This would be absurd, precisely since intersection is invariant under homotopies. This proves the claim that $u$ is constant in $Y$.

We have obtained that the portion of $u$ which lies in the over $\Sigma_+\backslash \mathcal{N}(\partial \Sigma_+)$ is actually contained in the 3-manifold $M_y:=\mathbb{R}\times \{y\} \times \Sigma_+\backslash \mathcal{N}(\partial \Sigma_+)$, as is the corresponding portion of the curve $u_{y,a}^+$. But for dimensional reasons, this manifold has a unique 2-dimensional $J$-invariant distribution, given by $TM_y\cap JTM_y$. Therefore the tangent space to $u$ must coincide with the tangent space of some $u_{y,a}^+$ at every point in this region, both coinciding with this (integrable) distribution, which has unique maximal integral submanifolds. Alternatively, one can check that the holomorphic curve equation, in local holomorphic coordinates $(s,t)$, applied to $u=(a,y,z)\in \mathbb{R}\times Y \times \Sigma_+\backslash \mathcal{N}(\partial \Sigma_+)$ implies $\partial_s a=\Lambda_\epsilon(\partial_ty)=0$, $\partial_t a=-\Lambda_\epsilon(\partial_sy)=0,$ so that $a$ is indeed constant over this region once one has that $y$ is. The unique continuation theorem finally yields $u=u_{y,a}^+$. 

\vspace{0.5cm}

This proves the theorem in the Morse--Bott situation. The proof in the Morse case follows from uniqueness in the Morse--Bott one, and the gluing results in \cite{Bo}, since we have shown that the only possibilities for holomorphic cascades are built using curves $u_{y,a}^\pm$ (and flow segments of $H$).

\end{proof}

\vspace{0.5cm}

\subsection{From the SHS to a sufficiently non-degenerate contact structure}
\label{SHStoND}

For computations in SFT we need non-degenerate Reeb orbits and contact data. Therefore, we need to perturb the SHS $\mathcal{H}_\epsilon=(\Lambda_\epsilon,\Omega_\epsilon)$ to a nearby contact structure, which we do by replacing the stabilizing vector field $Z$ of the first section with another sufficiently nearby vector field. If we choose the Morse scenario, then all the Reeb orbits $\gamma_p$ will be non-degenerate. Since these are the only ones that one sees under a large action threshold $T$, these will be the ones whose nondegeneracy we actually care about. Of course, one could also assume that $(Y,\alpha_\pm)$ is nondegenerate after a Morse perturbation, but we will not need to do this at this stage.

\subparagraph*{Perturbation to contact data} Recall that we have defined an exact symplectic form $\omega_\sigma$ on the double completion, and we denote by $V_\sigma$ its associated Liouville vector field. We also have the ``vertical'' Liouville vector field $X$ associated to $d\lambda$, defined by expression (\ref{X}), and the stabilizing vector field $Z$ for the SHS $\mathcal{H}_\epsilon$, defined by (\ref{Z}).

For $s \in [0,1]$, define 
$$Z_s=\left\{\begin{array}{ll} V_\sigma+((1-s)\beta(t)+s)X, & \mbox{ in the region } E^{\infty,\infty}(t) \mbox{ where } t \mbox{ is defined}\\ 
\frac{\sigma}{\sigma+\sigma^\prime}\partial_r+sX_+,& \mbox{ in } E^{\infty,\infty}(t)^c \cap \{r > 1-\delta\}\\
-\frac{\sigma}{\sigma-\sigma^\prime}\partial_r+sX_-,& \mbox{ in } E^{\infty,\infty}(t)^c \cap \{r < -1+\delta\}
\end{array}\right. $$ 

We have that $Z_1=X_{\sigma}$ is the Liouville vector field associated to $\omega_{\sigma}$ and $Z_0=Z$. This yields a family of SHS's given by $$\mathcal{H}_{\epsilon,s}=(\Lambda_{\epsilon,s},\Omega_{\epsilon})_{s\in[0,1]}=(i_{Z_s}\omega_{\sigma}\vert_{TM^\epsilon},\omega_{\sigma}\vert_{TM^\epsilon})_{s\in[0,1]},$$ so that the 1-form can be written as 

\begin{equation}\label{modelA} \Lambda_{\epsilon,s}=\left\{\begin{array}{ll}
  e^{\epsilon H}(\epsilon\alpha + d\theta),& \mbox{ on } M_Y \\
 K\alpha_\pm + s\lambda_\pm, &\mbox{ on } M_P^\pm \\
 \sigma(F_\pm^\epsilon(\rho)) e^{\pm F^\epsilon_\pm(\rho)-1}\alpha_\pm+((1-s)\beta(G^\epsilon_\pm(\rho))+s)e^{G^\epsilon_\pm(\rho)}d\theta,& \mbox{ on } M_C^\pm\\
  \end{array}\right.
\end{equation}

Let us observe that $\Lambda_{\epsilon,s}$ is a contact form for $s>0$. Indeed, the 1-form $s\lambda_\pm$ is also Liouville in $\Sigma_\pm$, so that it may be used in place of $\lambda_\pm$ when using the Thurston trick in the construction of $\omega_{\sigma}$. Then, even though $\xi_{\epsilon,s}:=\ker \Lambda_{\epsilon,s}$ is not independent of the parameters, Gray's stability implies that, as long as $\epsilon$ and $s$ are positive and sufficiently close to zero, we may choose them at will, without changing the isotopy class of $\xi_{\epsilon,s}$.

\subparagraph*{Holomorphic curves for the contact data} Since we have shown that the genus zero holomorphic curves are regular for the SHS data, the implicit function theorem implies that they will survive a small perturbation to contact data, and will still be regular. 

In Section \ref{Torsion}, we will count holomorphic curves in this new family with asymptotics in a specific collection of Reeb orbits. While in general it might happen that an index 1 building for the SHS data might glue to an honest curved that needs to be counted, we will show that \emph{the index $1$ holomorphic buildings that may appear as limits of the curves that contributes to the count are actually perturbations of the curves in the original foliation}.

The proof of this is left for the end of Section \ref{Torsion}, when it will be clear what curves need to be counted, and the conditions that this imposes to the limit buildings (Proposition \ref{buildingsuni}).

\subsection{Proof of Theorem \ref{thm5d}}
\label{Torsion}

Once all the technical tools are in place, we are finally able to do the computation which is one of the goals of this document. We will prove Theorem \ref{thm5d}. We will assume that the parameters $\epsilon$ and $s$ have been fixed, so that we will be working in the SFT algebra $\mathcal{A}(\Lambda_{\epsilon,s})[[\hbar]]$ whose homology is $H_*^{SFT}(M,\xi_{\epsilon,s})$. Since for different choices of (sufficiently small) $\epsilon$ and $s$ all of the $\xi_{\epsilon,s}$'s are isotopic, by Gray's stability, all of these SFT's are isomorphic. We shall be taking coefficients in $R_{\Omega}=\mathbb{R}[H_2(M;\mathbb{R})/\ker \Omega]$, where $\Omega \in \Omega^2(M)$ defines an element in the annihilator $\mathcal{O}:=$Ann$(\bigoplus^k H_1(Y;\mathbb{R})\otimes H_1(S^1;\mathbb{R}))$. We view $\bigsqcup^k Y\times S^1=\bigsqcup^k Y\times S^1\times \{0\}$ as sitting in $M_Y\subseteq M$, and $H_1(Y;\mathbb{R})\otimes H_1(S^1;\mathbb{R})$ as a subspace of $H_2(M)$. We will show that $(M,\xi_{\epsilon,s})$ has $\Omega$-twisted $(k-1)$-torsion for every $[\Omega] \in \mathcal{O}$, so that in particular it has $(k-1)$-torsion for \emph{untwisted} version of the SFT algebra when the coefficient ring is $\mathbb{R}$ (or equivalently, $\Omega=0$).

Let us recall that for SFT to be defined, we need to introduce an \emph{abstract perturbation} of the Cauchy--Riemann equation. We shall be doing our computation \emph{prior} to introducing this perturbation, and prove by the end the section that this is a reasonable thing to do. See the end of this section for more details.

\subparagraph*{Computation of algebraic torsion} In order to be able to count in SFT, we need to be looking at index one curves, which after quotienting out the $\mathbb{R}$-action yield zero dimensional moduli spaces. The index formula (\ref{index}) tells us that this is the case for curves $u_{y,a}^-$ which asymptote index $2n$ critical points (maxima) at $k-1$ of the $k$ positive ends, and one index $1$ critical point at the remaining one. Indeed, if $u$ is such a curve, then:
$$
\mbox{ind}(u)= 2n(1-k) + 1+\sum_{i=1}^{k-1} 2n=2n(1-k)+2n(k-1)+1=1 
$$  

Given $e_1,\dots,e_{k-1}$ maxima and $h$ an index $1$ critical point, denote by $\gamma_{e_i}$ and $\gamma_h$ the corresponding non-degenerate Reeb orbits. Consider the moduli space $$\mathcal{M}=\mathcal{M}_{0}(W,J_{\epsilon,s};(\gamma_{e_1},\dots, \gamma_{e_{k-1}},\gamma_h),\emptyset)$$ of $\mathbb{R}$-translation classes of $k$-punctured genus zero $J_{\epsilon,s}$-holomorphic curves in $W=\mathbb{R}\times M$, which have no negative asymptotics, and have $\gamma_{e_1},\dots,\gamma_{e_{k-1}},\gamma_h$ as positive asymptotics. The uniqueness Theorem \ref{uniqueness} implies that every element in $\mathcal{M}$ is a genus zero curve in our foliation. Since every such curve is regular, after a choice of coherent orientations as in \cite{BoMon}, this moduli space is an oriented zero dimensional manifold, which can therefore be counted with appropiate signs. Now, a choice of coherent orientations for the moduli space of Morse flow lines induces a coherent orientation for the moduli space containing the curves $u^a_\gamma$, and we fix such a choice here onwards. We choose our function $$H:Y\times I \rightarrow \mathbb{R}^{\geq 0}$$ as made explicit in Section \ref{MBtoM}. The critical points of $H$ are then of the form $p=(y,0)$ for $y \in crit(g)$ , and we have $$ind_p(H)=ind_0(f)+ind_y(g)=1+ind_y(g)$$ In particular, there are no index zero critical points, and the only index 1 critical point is given by $(min,0)$, where $min$ is the unique minimum of $g$. We shall denote by $\mathcal{M}(H;p_-,p_+)$ the moduli space of \emph{positive} flow lines $\gamma$ connecting $p_-$ to $p_+$, that is, solutions to the ODE $\dot{\gamma}=\nabla H \circ \gamma$ with $\lim_{t\rightarrow \pm \infty}\gamma(t)=p_\pm$, quotiented out by the $\mathbb{R}$-action. Assuming the Morse--Smale condition for $H$, we have then that the zero dimensional moduli spaces correspond to critical points satisfying $ind_{p_+}(H)=ind_{p_-}(H)+1$. 

\vspace{0.5cm}

We then fix  $e_1,\dots,e_{k-1},h=(min,0)$ as above, and we let $q_{e_i}$ and $q_h$ be the generators in SFT corresponding to the Reeb orbits associated to these critical points. Let $$Q=q_{e_1}\dots q_{e_{k-1}}q_h,$$ which is an element of $\mathcal{A}(\Lambda_{\epsilon,s})[[\hbar]]$. In order to compute its differential, we need to count all of the rigid holomorphic curves which asymptote at its positive ends any of the Reeb orbits appearing in $Q$.

\vspace{0.5cm}

\textbf{Claim}. The holomorphic curves contributing to the differential of $Q$ are of either the two following types: 

\begin{itemize}
\item A holomorphic sphere $u_{y,a}^- \in \mathcal{M}$, which is in fact unique.

\item A holomorphic cylinder $u_\gamma$ inside $r=0$, connecting an index $2n-1$ critical point to one of the maxima $e_i$.
\end{itemize}

Indeed, using Theorem \ref{uniqueness}, we will prove below that only the somewhere injective curves in our foliation are involved in the computation of this differential, so that in particular none of their multiple covers contribute to the differential of $Q$. Assuming this, using that there are no index zero critical points, we see that there are only two possible ways to approach the critical point $(min,0)$, entering through the two different boundary components of $Y \times I$, and that there is a unique $\mathbb{R}$-class of the form $u_y^-$ which has $\gamma_h$ as a positive asymptotic (see Figure \ref{Morse}). Moreover, by observing that the generic behaviour is hitting a maxima, by a generic and small perturbation of the Morse function $H$ along different components of the spine, we can achieve what can be thought of as a ``relative version'' of the Morse--Smale condition for the different components of this function. That is, we can arrange that this translation class hits only maxima at all other ends, so that it actually defines an element of $\mathcal{M}$. The uniqueness Theorem \ref{uniqueness} above shows that in fact this is the only element in $\mathcal{M}$. Finally, ruling out the curves coming from the positive side (which have the wrong index and hence are not counted by SFT), we are left only with the curves listed in our claim.

\vspace{0.5cm}

In order to count the curves of type $u_\gamma$, we observe that positive flow lines going from an index $2n-1$ critical point $p$ to an index $2n$ critical point $q$ come in ``evil twins'' pairs $\gamma \leftrightarrow \overline{\gamma}$: by definition of the Morse index, we have only one positive eigendirection for the Hessian of the Morse function at $p$, and the flow lines approach this point on either side of this direction. The key observation is that, since $H_{2n}(Y\times I)=H_{2n}(Y)=0$ ($Y$ being a $(2n-1)$-manifold), we only have one generator of the top Morse homology chain group, which is necessarily closed under the Morse differential.  Therefore, after choosing a coherent orientation of the moduli spaces of curves, the evil twins come with different sign and cancel each other out, and hence $\# \mathcal{M}(H;p,q)=0$. 

\vspace{0.5cm}

Fix $\Omega$ a closed $2$-form so that $[\Omega] \in \mathcal{O}$. For $d \in H_2(M;\mathbb{R})$, we denote by $\overline{d} \in H_2(M;\mathbb{R})/\ker \Omega$ the class it induces, by $u_0$ the unique $\mathbb{R}$-translation class in $\mathcal{M}$, and, for a rigid holomorphic curve $v$, we denote by $\epsilon(v)$ the sign of $v$ assigned by our choice of coherent orientations. In particular, we know that $\epsilon(u_\gamma)=-\epsilon(u_{\overline{\gamma}})$. 

Observe that, for the $i$-th component of $M_Y$, Reeb orbits corresponding to critical points all define the same homology class $[S^1]_i \in H_1(M)$. We take as canonical representative of this class the $1$-cycle given by the Reeb orbit $\gamma_{e_i}$ over the unique maxima $e_i$. For every index $2n-1$ critical point $p$ lying in the $i$-th component, fix $\gamma$ an index $1$ flow line joining $p$ to the maxima $e_i$. Choose the spanning surface of $\gamma_p$ to be $F_p:= -\gamma \times S^1$, satisfying $\partial F_p = \gamma_p-\gamma_{e_i}$. Then, for this choice of spanning surfaces, the homology class associated to the holomorphic cylinder $u_\gamma$ is $[u_\gamma]=[F_p\cup u_\gamma]$. One thinks of $F_p$ as being attached to $u_\gamma$ at the negative end, corresponding to $\gamma_p$. Let $T_{\gamma,\overline{\gamma}}$ denote the $2$-torus $\overline{\gamma\cup \overline{\gamma}}\times S^1$. Observe that $$[u_\gamma]-[u_{\overline{\gamma}}]=[T_{\gamma,\overline{\gamma}}] \in H_1(Y)\otimes H_1(S^1)\subseteq \ker \Omega$$ 
Therefore, we have $\overline{[u_\gamma]}=\overline{[u_{\overline{\gamma}}]}$.

According to Proposition \ref{buildingsuni} below, the image of $Q$ under the differential $\mathbf{D}_{\epsilon,s}:\mathcal{A}(\Lambda_{\epsilon,s})[[\hbar]] \rightarrow \mathcal{A}(\Lambda_{\epsilon,s})[[\hbar]]$ is then given by 

\begin{equation}
\begin{split}
\mathbf{D}_{\epsilon,s}(Q)=&\epsilon(u_0)z^{\overline{[u_0]}}\hbar^{k-1} + \sum_{\substack{i=1,\dots,k-1\\ind_p(H)=2n-1}}\sum_{\gamma \in \mathcal{M}(H;p,e_i)}\epsilon(u_\gamma)z^{\overline{[u_\gamma]}}q_{p}\frac{\partial Q}{\partial q_{e_i}}\\
=&\epsilon(u_0)z^{\overline{[u_0]}}\hbar^{k-1} + \frac{1}{2}\sum_{\substack{i=1,\dots,k-1\\ind_p(H)=2n-1}}\sum_{\gamma \in \mathcal{M}(H;p,e_i)}(\epsilon(u_\gamma)+\epsilon(u_{\overline{\gamma}}))z^{\overline{[u_\gamma}]} q_{p}\frac{\partial Q}{\partial q_{e_i}}\\
=& \epsilon(u_0)z^{\overline{[u_0]}}\hbar^{k-1},
\end{split}
\end{equation}
which proves that our model has $(k-1)$-torsion with coefficients in $R_{\Omega}$. We have used that all orbits are simply covered, so no combinatorial factors appear, and we need not worry about asymptotic markers.

\subparagraph*{Why the computation works} We now justify the computation above. Let us recall first the fact that the abstract machinery of SFT requires the introduction of an abstract perturbation to the Cauchy--Riemann equation making every holomorphic curve of positive index regular. The basic facts about this perturbation scheme, which comes from the polyfold theory under development by Hofer--Wysocki--Zehnder, are that 

\begin{itemize}
\item Every Fredholm regular index 1 holomorphic curve gives rise to a unique solution to the perturbed problem, if the perturbation is sufficiently small.
\item If solutions to the perturbed problem with given asymptotic behaviour exist for small perturbations, then as the perturbation is switched off they give rise to a subsequence of curves which converge to a holomorphic building with the same asymptotic behaviour. 
\end{itemize}

Therefore, the index 1 curves in our foliation survive and are counted, but we need to make sure there are no extra curves which need to be taken into the count. That is, we need to study the space of holomorphic stable buildings of index 1 which may arise as limits of curves that contribute to the differential of $Q=q_{e_1}\dots q_{e_{k-1}}q_h$, which are exactly those which when the perturbation is turned on glue to index 1 curves which are seen by SFT. We also need to prove that we haven't added information when we perturbed from the SHS to the contact structure (which chronologically comes before the abstract perturbation), as was discussed at the end of Section \ref{SHStoND}.

In what follows, $J$ will denote the original $J_{\epsilon,0}$ we constructed in the Morse case. We have to prove the following:

\begin{prop}\label{buildingsuni} The space of connected $J$-holomorphic stable buildings of index 1 which may become curves contributing to the differential of $Q=q_{e_1}\dots q_{e_{k-1}}q_h$ after introducing an abstract perturbation, or after perturbing the SHS to a sufficiently nearby contact structure, have only one level, no nodes, and are somewhere injective and regular: they actually consist exactly of either an index 1 curve $u_\gamma$ for some positive flow line $\gamma$, or a curve $u^-_{y,a}$ for some $y \in Y$, $a \in \mathbb{R}$.
\end{prop}

\begin{proof}
Consider such a building $u$. Let us first observe that $u$ has in its top level positive asymptotics which are a \emph{subset} of $\{\gamma_{e_1},\dots,\gamma_{e_{k-1}},\gamma_h\}$. Therefore, it has \emph{at most $k$} of these, any two different positive punctures asymptote distinct simply covered Reeb orbits, and any two different of the latter asymptotics lie in \emph{distinct} components of $M_Y$. The uniqueness result \ref{uniqueness} implies that \emph{each component of $u$ is a simply covered curve in our foliation}, as follows by induction (by observing that flow-line cylinders connect Reeb orbits in the same component of $M_Y$).

Then the only possibilities of configurations for $u$ are shown in Figure \ref{possconfu}, depending on whether $u$ has a $u^-_{y,a}$ curve or not (the $u^+_{y,a}$ curves cannot possibly appear in $u$, since we are assuming that its components do not disappear after generic perturbations). 

\begin{figure}[t]\centering
\includegraphics[width=0.65\linewidth]{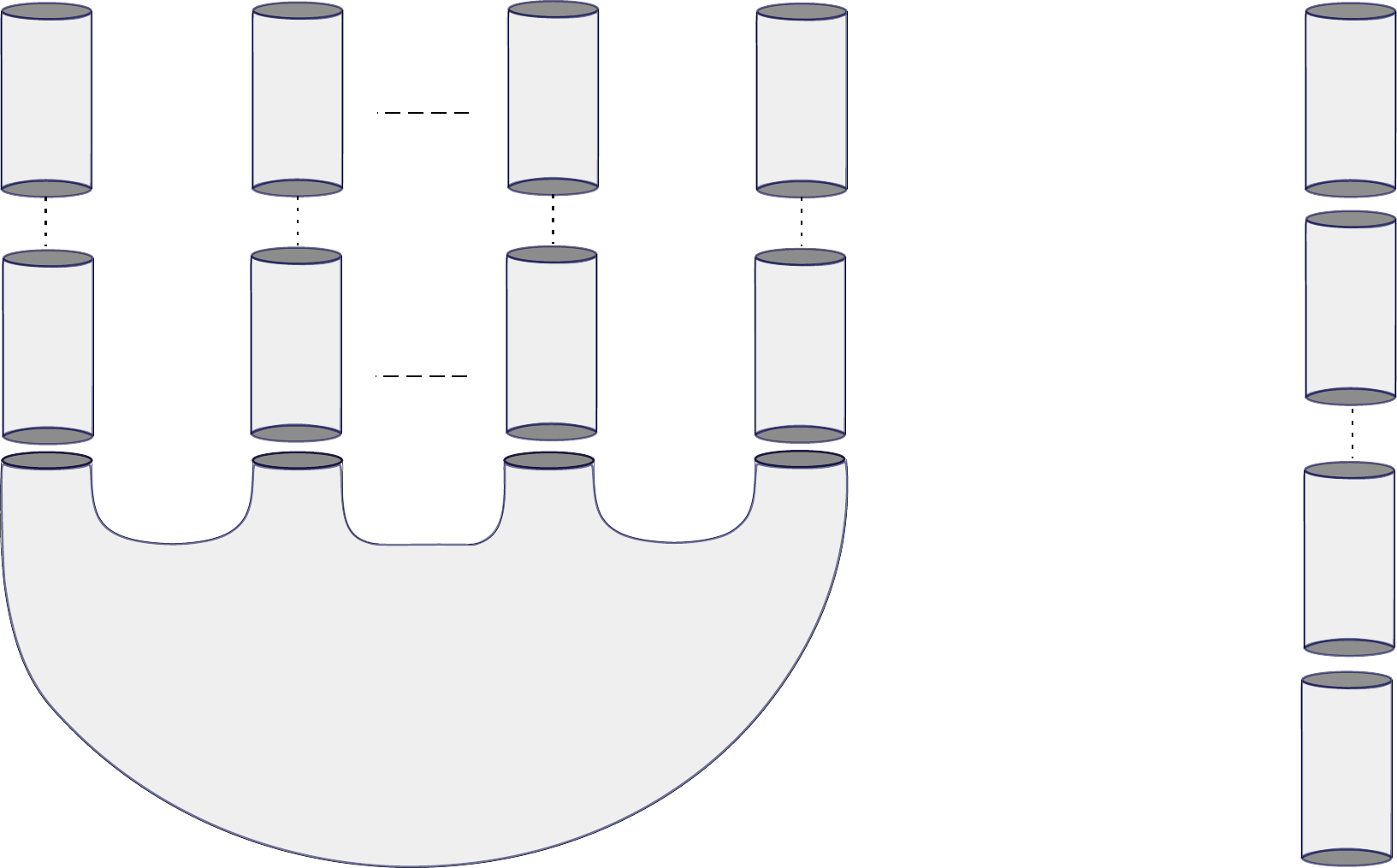}
\caption{\label{possconfu} the only two possible connected configurations one can build with curves in the foliation. The index restriction rules out the case of more than one level.}
\end{figure}

If a curve $u^-_{y,a}$ appears in $u$, since flow line cylinders have non-negative index, and by regularity, we have $1=\mbox{ind}(u)\geq \mbox{ind}(u^-_{y,a})\geq 1$. The stability condition, which means that there are no levels consisting solely of trivial cylinders, implies that $u=u^-_{y,a}$ is the unique element in $\mathcal{M}$. In the case where a curve $u^-_{y,a}$ does not appear in $u$, the same argument implies that $u$ corresponds to an index 1 flow line in $\mathcal{M}(H;p,e_i)$ for some $p$ with $\mbox{ind}_p(H)=2n-1$ and some $i\in\{1,\dots,k-1\}$. This concludes the proof.

\end{proof}

\section{Model B. Ruling out lower orders of torsion} 
\label{RulingOut}

In the second part of this document, we will construct another contact model for the same manifold $M=Y \times \Sigma$, model B, which is isotopic to the original model A. The role of spine and paper are swapped, so that these models may be considered dual in some sense. Whereas model A is suitable for computing that a certain order of torsion is achieved, model B was originally conceived for showing that no lower orders are also, with the aim that the order of torsion of the isotopy class both define is be precisely computed. Instead, we use it to prove Theorem \ref{nocob1}. The main features of model B is that closed Reeb orbits of low action correspond to pairs $(p,\gamma)$ of critical points $p$ of suitable Morse functions on $\Sigma_\pm$, and closed $R_\pm$-orbits $\gamma$ in $(Y,\alpha_\pm)$; and we will have a foliation of its symplectization by holomorphic hypersurfaces, for a suitable almost complex structure compatible with a SHS deforming the contact data, rather than holomorphic curves. These project to flow lines on the surface $\Sigma$, and the ones which do not cross the dividing set along which spine and paper glue together can all be identified with the symplectization of $(Y,\alpha_\pm)$ (the sign depending on the side of the dividing set). This identification can be done by projecting to the ``vertical'' hypersurfaces lying over any of the corresponding asymptotic critical points. The ones which do cross sides may be regarded as cylindrical completions of the Liouville domain $Y\times I$ (recall that we denote $I=[-1,1]$). Intersection theory (namely, Theorem \ref{posdegree} in Appendix \ref{AppSiefring}) may then be used to prove that any holomorphic curves necessarily lies in a hypersurface of this foliation, which is the key tool for restricting their behaviour. Once the setup is in place, ruling out torsion becomes a matter of understanding suitable index 1 buildings whose components are restricted to lie in these hypersurfaces. The reason why this is not a straightforward thing to do is that one needs a very good understanding of holomorphic curves in the completion of $Y\times I$, and in each symplectization of $(Y,\alpha_\pm)$, and how index of curves inside these relate to those of the same curves when considered as lying in the symplectization of $M$. After constructing the model, we will restrict our attention to a specific family of $5$-dimensional models for which we prove Theorem \ref{nocob1}.  
 
\subsection{Construction of model B} 

As in Section \ref{model}, consider $Y^{2n-1}$ such that $(Y\times I,d\alpha)$ is a cylindrical Liouville semi-filling, where $\alpha=\{\alpha_r\}_{r\in I}$ is a 1-parameter family of 1-forms such that $\alpha_r$ is contact for $r\neq 0$. We will denote by $\xi_\pm=\ker \alpha_\pm$ the contact structures on $Y_\pm:=Y\times \{\pm 1\}$, where $\alpha_\pm=\alpha_{\pm 1}=\alpha\vert_{Y_\pm}$, by $V$ the Liouville vector field associated to $d\alpha$, by $R_r$ the Reeb vector field associated to $\alpha_r$ for $r\neq 0$, and $R_\pm:=R_{\pm 1}$. Consider also $\Sigma=\Sigma_-\cup_\partial \Sigma_+$ a genus $g$ surface obtained by gluing a genus zero surface $\Sigma_-$ to a genus $g-k+1>0$ surface $\Sigma_+$, along $k$ boundary components in an orientation-reversing manner, and the manifold $M=Y\times \Sigma$.  

Take collar neighbourhoods of each boundary component of $\Sigma_\pm$ of the form $\mathcal{N}(\partial \Sigma_\pm):=[-1,0]\times S^1 $, with coordinates $(t_\pm,\theta) \in [-1,0] \times S^1$, such that $\partial\Sigma_\pm=\{t_\pm=0\}$. Glue these surfaces together along a cylindrical region $\mathcal{U}$, having $k$ components, each of which looks like $[-1,1]\times S^1 \ni (r,\theta)$, so that $r=\mp t_\pm$. We now use the fact that each surface $\Sigma_\pm$ carries a \emph{Stein} structure providing a Stein filling of its 1-dimensional convex contact boundary. That is, we may take Morse functions $h_\pm$ on $\Sigma_\pm$ which are plurisubharmonic with respect to a suitable complex structure $j_\pm$, which in the collar neighbourhoods look like $h_\pm(t_\pm,\theta)=e^{t_\pm}$ and $j_\pm(\partial_{t_\pm})=\partial_\theta$. Moreover, we have that $\lambda_\pm=d^{\mathbb{C}}h_\pm=-dh_\pm \circ j_\pm$ is a Liouville form, with Liouville vector field given by $\nabla h_\pm$, where the gradient is computed with respect to the metric $g_\pm=g_{d\lambda_\pm,j_\pm}=d\lambda_\pm(\cdot,j_\pm \cdot)$ (see Lemma 4.1 in \cite{LW} for details). Thus, in the collar neighbourhoods, we have $\lambda_\pm = e^{t_\pm}d\theta=e^{\mp r}d\theta$ and $\nabla h_\pm = \partial_{t_\pm}=\mp \partial_r$. Each regular level set of $h_\pm$, like the boundary $\partial\Sigma_\pm=h_\pm^{-1}(1)$, is a (rather uninteresting) 1-dimensional contact manifold with trivial contact structure $T\Sigma_\pm \cap j_\pm T\Sigma_\pm = 0$, and the Hamiltonian vector field $X_{h_\pm}=j_\pm \nabla h_\pm$, computed with respect to the symplectic form $d\lambda_\pm$ inducing the orientation in $\Sigma_\pm$, is tangent to these. We also have the basic formulas $$dh_\pm(\nabla h_\pm)= -dh_\pm(j_\pm X_{h_\pm})= \lambda_\pm(X_{h_\pm})=d\lambda_\pm(\nabla h_\pm,X_{h_\pm})=|\nabla h_\pm|^2$$$$ \lambda_\pm(\nabla h_\pm)=-dh_\pm(X_{h_\pm})=0$$ 

We shall be taking the orientation on $M= Y\times \Sigma$ to be the one given by the contact form $\alpha_-$ on the $Y$-factor, and in the $\Sigma$-factor we take the one inherited from $\Sigma_-$ and given by $d\lambda_-$ (which is opposite to the one in $\Sigma_+$ given by $d\lambda_+$). We will also take both $h_\pm$ to have a unique minimum in $\Sigma_\pm$, and we will make the convention of calling the minimum of $h_+$, the \emph{maximum}. This is because of the orientation considerations, and because it will indeed be the maximum of a suitably constructed Morse function (see Remark \ref{admissible} below).

We define, for $K_+>K_-\gg 0$ large constants, a function $g_0:[-1,1] \rightarrow \mathbb{R}$, such that it is constant equal to $K_-$ in $[-1,-1+\delta)$ for some small $\delta>0$, equal to $K_+$ in $(1-\delta,1]$, $g_0> 0$ everywhere, and $g_0^\prime\geq 0$ with strict inequality in the interval $(-1+\delta,1-\delta)$. The parameter $\delta>0$ is chosen so that the Liouville vector field in $Y \times [-1,1]$ is given by $V=\pm \partial_r$ in the corresponding components of $\{|r|>1-\delta\}$.  Define also a function $\gamma: [-1,1]\rightarrow \mathbb{R}$ so that $\gamma(r)=- e^{\mp r}$ in $\{|r|>1-\delta\}$, and $\gamma<0$. Take $$g_\epsilon(r)=g_0(r)+\epsilon^2\gamma(r),$$ where $\epsilon$ is chosen small enough so that $g_\epsilon>0$.

Now define 
$$h_\epsilon=\left\{\begin{array}{ll} K_\pm - \epsilon^2 h_\pm,& \mbox{ in } \Sigma\backslash \mathcal{U} \\
g_\epsilon,& \mbox { in } \mathcal{U} 
\end{array}\right.
$$

This is a (well-defined and smooth) function on $\Sigma$, and $\epsilon$ can be chosen small enough so that $h_\epsilon>0$.

\vspace{0.5cm}

Next, choose a smooth function $\psi:[-1,1]\rightarrow \mathbb{R}^+$ satisfying $\psi(r)=e^{\mp r}$ in $\{|r|>1-\delta/3\}$, $\psi\equiv 1$ on $[-1+\delta,1-\delta]$, $r\psi^\prime(r)\leq 0$. Define $$f_\epsilon(r)=\epsilon\psi(r),$$ so that $f_\epsilon>0$, $f_\epsilon^\prime\equiv 0$ along $[-1+\delta,1-\delta]$. 

With these choices, let $\Lambda_\epsilon$ be the 1-parameter family of 1-forms in $M=Y \times \Sigma$ given by 
\begin{equation}
\label{modelBform}
\Lambda_\epsilon=\left\{ \begin{array}{ll}\epsilon\lambda_\pm + h_\epsilon \alpha_\pm, &\mbox{ in } Y\times (\Sigma\backslash\mathcal{U})\\
f_\epsilon d\theta + h_\epsilon \alpha, &\mbox{ in } Y\times \mathcal{U}
\end{array}\right.
\end{equation}

We can slightly modify $\alpha$ so that the different expressions in $\Lambda_\epsilon$ glue together smoothly: we have $\alpha=e^{\pm r-1}\alpha_\pm$ in the corresponding components of the region $\{|r|\geq 1-\delta\}$), and we can replace $\alpha$ by a 1-form (of the same name) which looks like $\phi(r)\alpha_\pm$. Here, $\phi: \{|r|\geq 1-\delta\}\rightarrow \mathbb{R}$ is a smooth function which coincides with $e^{\pm r-1}$ near $|r|=1-\delta$, equals $1$ on $|r|=1$, and its derivative is non-negative/non-positive in the positive/negative components of $\{|r|\geq 1-\delta\}$.   

We shall refer to this family of contact forms as \emph{model B}, in contrast with \emph{model A}, given by the expression (\ref{modelA}).

\begin{lemma} \label{contact} $\;$
\begin{enumerate}
\item For a fixed $K_+\gg 0$, $\Lambda_\epsilon$ is a contact form for $K_-$ sufficiently close to $K_+$, and sufficiently small $\epsilon>0$.
\item Model A is isotopic to model B.

\end{enumerate}
\end{lemma}

We will prove this lemma by viewing model B as arising from the double completion construction as we did for model A, after suitably arranging our coordinates. This is something very useful which can be used to produce a suitable stable Hamiltonian structure deforming the contact data, as well as the required isotopy between both models. Morally, both models should be isotopic since they are ``supported'' by the same geometric decomposition underlying the dual SOBDs.  

\vspace{0.5cm} 

\begin{figure}[t]\centering
\includegraphics[width=0.65\linewidth]{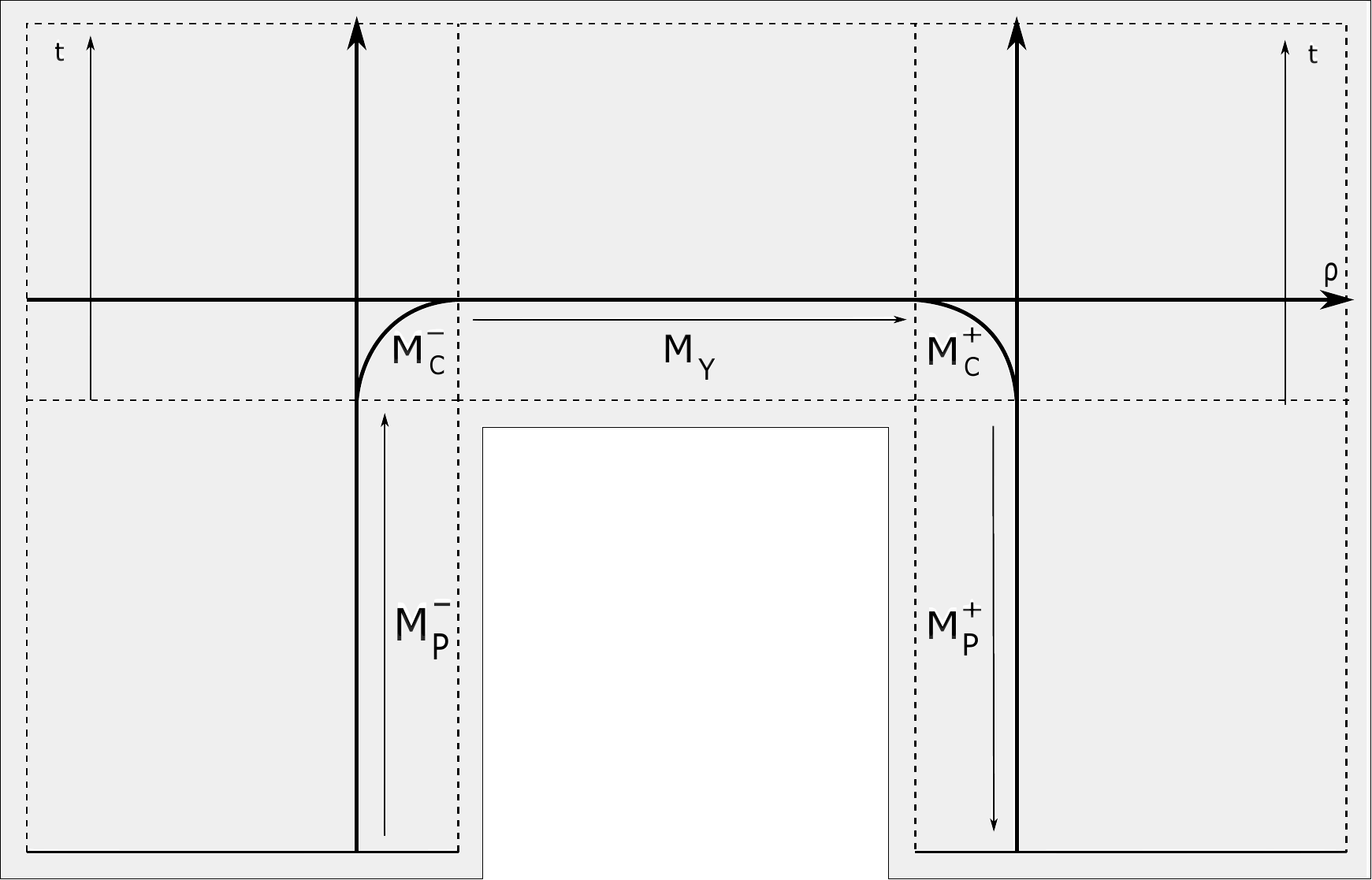}
\caption{\label{doubleB} The double completion.}
\end{figure}

We slightly change the coordinates in the definition of the double completion, in order to make it compatible with the model B construction. As in Section \ref{model}, denote by $\mathcal{N}^\infty(Y)=Y \times (1-\delta,+\infty]\times S^1$, and $\mathcal{N}^\infty(-Y)=Y \times (-\infty,-1+\delta)\times S^1$. Also, denote by $\Sigma_\pm^\infty$ the Liouville completion of $\Sigma_\pm$, obtained by attaching cylindrical ends of the form $(-\delta, + \infty)\times S^1$ at each boundary component. The coordinates $t_\pm$ and $\theta$ extend to the cylindrical ends in an obvious way, and we denote the latter by $\mathcal{N}(\partial \Sigma_\pm^\infty)$. We also have
$$
M_P^{\pm,\infty}= Y \times \Sigma_\pm^\infty
$$
$$
M^\infty_Y= Y \times \mathbb{R}\times S^1
$$

Take a variable $\rho$ parametrizing the $\mathbb{R}$-direction in the completion $Y\times \mathbb{R}$. We now define the double completion as the open manifold 
$$E^{\infty,\infty}=(-\infty,-1+\delta) \times M_P^{-,\infty} \;\bigsqcup\; (1-\delta,+\infty) \times M_
P^{+,\infty}\; \bigsqcup \; (-1,+\infty)\times M_Y^\infty /\sim,$$
where we identify $(\rho,y,t_-,\theta) \in (-\infty,-1+\delta) \times Y \times \mathcal{N}(\partial \Sigma^\infty_-)$ with $(t,y,\rho,\theta)\in (-1,+\infty) \times \mathcal{N}^\infty(-Y)$ if and only if $t=t_-$, and $(\rho,y,t_+,\theta) \in (1-\delta,+\infty) \times Y \times \mathcal{N}(\partial \Sigma^\infty_+)$ with $(t,y,\rho,\theta)\in (-1,+\infty) \times \mathcal{N}^\infty(Y)$ if and only if $t=t_+$.  We will denote by $E^{\infty,\infty}(t)$ the subset of $E^{\infty,\infty}$ where the variable $t$ is defined. See Figure \ref{doubleB}.

\begin{figure}[t]\centering
\includegraphics[width=0.65\linewidth]{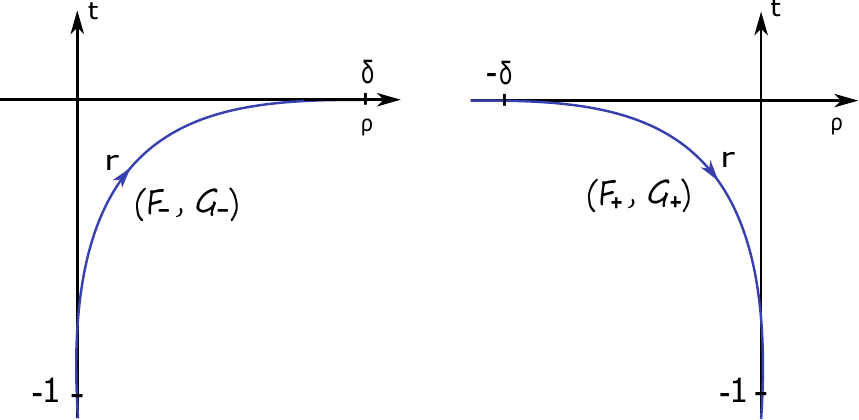}
\caption{\label{smoothB} The choice of smoothening.}
\end{figure}

Analogously as before, we choose smoothening functions $F_+,G_+:(1-\delta,1)\rightarrow \mathbb{R}$, and $F_-,G_-:(-1,-1+\delta)\rightarrow \mathbb{R}$ satisfying 

$$\left\{\begin{array}{ll}
 (F_+(r),G_+(r))=(r-1,0),&\mbox{ for } 1-\delta< r \leq 1-2\delta/3\\
 (F_-(r),G_-(r))=(0,r), &\mbox{ for } -1< r \leq -1+\delta/3\\
 G_+^\prime(r)<0,\;F_+^\prime(r)>0,& \mbox{ for } 1-2\delta/3 <r < 1-\delta/3\\
 F_-^\prime(r)>0,\; G_-^\prime(r)>0,& \mbox{ for } -1+\delta/3<r <-1+ 2\delta/3\\
(F_+(r),G_+(r))=(0,-r),& \mbox{ for } 1-\delta/3\leq r < 1\\
(F_-(r),G_-(r))=(r+1,0),& \mbox{ for } -1+2\delta/3\leq r < -1+\delta\\
\end{array}\right.$$

See Figure \ref{smoothB}.

For $L\geq 1$ and $Q\geq 0$, we define $M^{L,Q}$ as the codimension one submanifold of   $E^{\infty,\infty}$ given by 

$$ M^{L,Q}:= \partial E^{L,Q} \cap \Big(\{\rho \in [-L+\delta,L-\delta]\}\cup \{t\leq Q-1\} \cup \{t=Q\}\Big) \bigcup M_C^{\pm,L,Q},$$ 

where $E^{L,Q}=E^{\infty,\infty}\backslash(\{|\rho|>L\}\cup\{t>Q\}),$ and
\begin{equation}
\begin{split}
M_C^{\pm,L,Q}=&\{(\rho=F_\pm(r)\pm L,t=G_\pm(r)+Q,y,\theta):\\
&(r,y,\theta)\in \left(\mathbb{R}^\pm \cap \{1>|r|> 1-\delta\}\right)\times Y \times S^1\},
\end{split}
\end{equation}
are the smoothenings of the corners of $E^{L,Q}$.

We then have a decomposition $$M^{L,Q}=M_P^{\pm,L,Q}\bigcup M_C^{\pm,L,Q}\bigcup M_Y^{L,Q},$$ which we define exactly the same as for model A:
$$
M_P^{\pm,L,Q}:=(Y \times (\Sigma_\pm\backslash\mathcal{U}))\cap M^{L,Q}
$$
$$
M_Y^{L,Q}:= (Y \times [-1+\delta,1-\delta] \times S^1))\cap M^{L,Q}
$$
The region $M_P^{\pm,L,Q}$ should now be thought as the positive/negative spine, rather than the paper, but to avoid confusion we do not change the notation. We shall also drop the $L,Q$ from the notation when talking about $L=1,Q=0$.

If we define $$\widehat{E}^{\infty,\infty}=\bigcup_{L\geq 1}M^{L,L-1}\subset E^{\infty,\infty},$$ which is foliated by the disjoint submanifolds $M^{L,L-1}$, we may extend the function $h_\epsilon:M\rightarrow \mathbb{R}$ to $\widehat{E}^{\infty,\infty}$ by fixing suitable identifications of the $M^{L,L-1}$'s with the diffeomorphic manifold $M=M^{1,0}$. For instance, we may define
\begin{itemize}
\item $h_\epsilon(y,z,\rho):=h_\epsilon(z),$ 

for $(y,z,\rho) \in \widehat{E}^{\infty,\infty} \cap E^{\infty,\infty}(t)^c=Y \times \left(\Sigma_\pm\backslash \mathcal{N}(\partial \Sigma_\pm)\right) \times \{|\rho|\geq 1\}$
\item $h_\epsilon(F_\pm(r)\pm L,G_\pm(r)+L-1,y,\theta):=h_\epsilon(r)$, 

for $(F_\pm(r)\pm L,G_\pm(r)+L-1,y,\theta) \in M_C^{\pm,L,L-1}$.
\item $h_\epsilon(y,\rho,t,\theta):=h_\epsilon\left(\frac{1-\delta}{L-\delta}\rho\right)$, for $(y,\rho,t,\theta) \in Y \times [-L+\delta,L-\delta]\times S^1\subseteq M_Y^{L,L-1}$.
\item $h_\epsilon(y,\pm L,t,\theta)\equiv h_\epsilon(\pm 1),$ for $(y,\pm L,t,\theta) \in M_P^{\pm,L,L-1}\cap E^{\infty,\infty}(t)$.
\end{itemize} 

We also extend $\lambda_\pm$ and $\alpha$ to $E^{\infty,\infty}$ as in Section \ref{model}, so that they look like $\lambda=e^td\theta$ and $\alpha=e^{\mp \rho -1}\alpha_\pm$ over the cylindrical ends. We claim that the 1-form 
$$\Lambda_\epsilon:=\epsilon \lambda + h_\epsilon \alpha$$ can be made Liouville in $E^{\infty,\infty}$, after sufficiently shrinking the first derivative of $h_\epsilon$ (in the uniform norm), but keeping $h_\epsilon$ uniformly large. Again by the Thurston trick, we need only check that $h_\epsilon \alpha$ can be made Liouville in $Y \times \mathbb{R}\subset \{t=0\}$, over which $E^{\infty,\infty}$ is a ``fibration'' with symplectic fibers with cylindrical ends (which change topology). 

Indeed, we have
\begin{equation}
\label{Lioucillec}
d(h_\epsilon \alpha)^n =h_\epsilon^{n-1}(h_\epsilon -n h_\epsilon^\prime link(\alpha_\rho))dvol\wedge d\rho,
\end{equation}

where we write $d\alpha^n=dvol \wedge d\rho$ for some positive volume form $dvol$ in $Y$ with respect to the $\alpha_-$-orientation.
Observe that even though $$link(\alpha_\rho)=e^{n(\pm \rho -1)}link(\alpha_\pm)$$ over the cylindrical ends, and hence is not bounded as a function of $\rho$ on the whole real line, our choice of extension of $h_\epsilon$ implies that $h_\epsilon^\prime\equiv 0$ for $|\rho|\geq 2$. We may choose $h_\epsilon^\prime$ uniformly small enough on $[-2,2]$, without making $h_\epsilon$ small, by observing that $h_\epsilon^\prime=g_\epsilon^\prime=g_0^\prime + \epsilon \gamma^\prime$ can be made small by taking $\epsilon$ small and $K_-$ sufficiently close to $K_+$ (so that $g_0^\prime$ is small), and so the expression (\ref{Lioucillec}) is positive.

We may then compute the Liouville vector field associated to $h_\epsilon \alpha$, which is given by $$V_\epsilon:= \frac{h_\epsilon}{h_\epsilon+h_\epsilon^\prime dr(V)}V,$$ where $V$ is the Liouville vector field associated to $\alpha$ (cf.\ (\ref{Vsig})). Therefore it is a positive multiple of $V$, for suitably small $h_\epsilon^\prime$, for a multiple which is as arbitrarily close to $1$ as desired, and coincides with $\pm \partial_\rho$ over $\{|\rho|\geq 2\}$. It follows that the Liouville vector field associated to $\Lambda_\epsilon$ is $Z_\epsilon:=X +V_\epsilon$. Here, $X$ is the Liouville vector field of $\lambda$, defined by  
\begin{equation}
\label{X}
X=\left\{\begin{array}{ll}
  X_+, &\mbox{ on } \{r>1-\delta\} \\
  X_-, &\mbox{ on } \{r<-1+\delta\} \\
  \partial_t,& \mbox{ on } \{|r|<1\}\\
  \end{array}\right.
\end{equation} where $X_\pm$ is the Liouville vector field associated to $d\lambda_\pm$. Since $V_\epsilon$ can be made as arbitrarily close to $V$ as wished, $Z_\epsilon$ is transverse to $M^{L,Q}$ for every $L,Q$, and therefore, for the fixed choice $L=1$, the manifold $(M=M^{1,0},\Lambda_\epsilon\vert_M)$ is contact. This proves the first part of Lemma \ref{contact}. 

\vspace{0.5cm}

For the second part, we linearly homotope between the Liouville forms $\lambda_\sigma=\sigma \alpha + \lambda$ and $\Lambda_\epsilon=h_\epsilon \alpha + \epsilon\lambda$, which are the ones we used to construct models A and B, respectively. We consider
$$
\Lambda_t=g_t(\sigma,h_\epsilon)\alpha+g_t(1,\epsilon)\lambda,
$$
where $g_t(A,B)=(1-t)A+tB$ is the linear interpolation between the quantities $A$ and $B$. Again by Thurston's trick, to show that $\Lambda_t$ is Liouville for all $t$, we only need to check that $g_t(\sigma,h_\epsilon)\alpha$ and $g_t(1,\epsilon)\lambda$ are Liouville forms in the horizontal and vertical directions, respectively, and that the Liouville vector field $Z_t$ associated to $\Lambda_t$ is transverse to $M$ for every $t$. The latter 1-form is obviously Liouville and its Liouville vector field is clearly a multiple of that of $\lambda$. For the former, the computation is the same as for $h_\epsilon \alpha$ or $\sigma \alpha$, with the expression of the corresponding Liouville vector field being the same as for $V_\epsilon$ but with $h_\epsilon$ replaced by $g_t(\sigma,h_\epsilon)$. The path $\Lambda_t\vert_M$ is then an isotopy of contact forms interpolating between model A and model B. This proves Lemma \ref{contact}. $\square$

\vspace{0.5cm}

If we denote $D_\epsilon:= f_\epsilon g_\epsilon^\prime-f_\epsilon^\prime g_\epsilon$, the Reeb vector field is

$$
\arraycolsep=1.4pt\def\arraystretch{2.2}
R_\epsilon=\left\{\begin{array}{ll} \frac{1}{D_\epsilon + f_\epsilon g_\epsilon\frac{\partial \alpha_r}{\partial r}(R_r)}((g_\epsilon^\prime+g_\epsilon \frac{\partial \alpha_r}{\partial r}(R_r))\partial_\theta - f_\epsilon^\prime R_r),& \mbox{ in }M_C^\pm\\

\frac{1}{K_\pm + \epsilon^2(|\nabla h_\pm|^2-h_\pm)}(R_\pm + \epsilon X_{h_\pm}),& \mbox{ in } Y \times (\Sigma\backslash \mathcal{U})\\
\frac{\partial_\theta}{\epsilon}& \mbox{ in }M_Y,\\ 
\end{array}\right.
$$
Using that $\alpha=\phi(r)\alpha_\pm$, so that $\frac{\partial \alpha_r}{\partial r}(R_r)=\frac{\phi^\prime(r)}{\phi(r)}$, it can be checked that these expressions coincide with $R_\epsilon= \frac{1}{K_\pm}(R_\pm + \epsilon \partial_\theta)$ where $M_C^\pm$ glues together with $M_P^\pm$. We clearly see that $R_\epsilon$ is totally degenerate along $M_Y$. Recalling that $\widetilde{R}_r$ and $R_r$ are colinear (see Corollary \ref{compJ0}), we also see that $R_\epsilon$ is everywhere a linear combination of $\partial_\theta$ and $\widetilde{R}_r$. Using expression (\ref{somexpr}), we see that the coefficient of $R_\epsilon$ in front of $\partial_\theta$ can be written as $$\frac{g_\epsilon^\prime(\alpha(X_{g_\epsilon})+g_\epsilon)}{\alpha(X_{g_\epsilon})D_\epsilon+f_\epsilon g_\epsilon g_\epsilon^\prime},$$ and since we can make the derivative of $g_\epsilon$ as small as we like without making $g_\epsilon$ small, this expression can be made positive. This means that every closed Reeb orbit in $Y \times \mathcal{U}$ projects to $\mathcal{U}$ as a simple closed curve which is homologous to a positive multiple of a dividing circle.

Moreover, the expression of $R_\epsilon$ in the region $Y \times (\Sigma\backslash\mathcal{U})$ tells us that every pair $(p,\gamma)$, consisting of a critical point $p$ of $h_\pm$ and a closed $R_\pm$-orbit $\gamma$, gives rise to a closed Reeb orbit of $R_\epsilon$, of the form $\gamma_p:=(\gamma,p) \in Y \times \Sigma$. Also, observe that the closed Reeb orbits which do not correspond to critical points of $h_\pm$ (the ones that project to $\Sigma_\pm$ as closed Hamiltonian orbits, lying in the level sets of these functions) can be made to have arbitrarily large period by taking $\epsilon$ sufficiently close to zero. Therefore we can find an action threshold $T_\epsilon>0$, such that $\lim_{\epsilon \rightarrow 0} T_\epsilon = +\infty$, and such that every closed Reeb orbit of $R_\epsilon$ of action less than $T_\epsilon$ either lies in $Y \times \mathcal{U}$, or is a cover $\gamma_p^l$ of a simple closed Reeb orbit $\gamma_p$, for $l\leq N_\epsilon$, where $N_\epsilon:=\max\{n\in\mathbb{N}:\gamma^l \mbox{ has action less than } T_\epsilon\}$. This orbit has action given by $l h_\epsilon(p)\alpha_\pm(\dot{\gamma})= l(K_\pm - \epsilon^2h_\pm(p))\alpha_\pm(\dot{\gamma})$.

\subsection{Deformation to a SHS along the cylindrical region}

To obtain a SHS deforming the contact data, which we will need to construct holomorphic hypersurfaces, we homotope the Liouville vector field over $\{\rho \in [-1,1]\}$, as follows. Choose a bump function $\kappa:\mathbb{R}\rightarrow \mathbb{R}$ which equals $1$ outside of the unit interval, equals zero in $[-1+\delta,1-\delta]$, and sign$(\kappa^\prime(\rho))=$sign$(\rho)\neq 0$ for $\rho \in \{1-\delta<|\rho|<1\}$.

Define
$$
W^s_\epsilon:=\left\{\begin{array}{ll}Z_\epsilon, &\mbox{ in }E^{\infty,\infty}\backslash\left( \{|\rho|\leq 1\}\cap E^{\infty,\infty}(t)\right)\\
X +(s\kappa(\rho)+1-s)V_\epsilon, &\mbox{ in } \{|\rho|\leq 1\}\cap E^{\infty,\infty}(t) 
\end{array}\right.
$$

One can check that 
$$
\mathcal{H}_\epsilon^s:=(\Lambda_\epsilon^s,\Omega_\epsilon)=(i_{W_s^\epsilon}d\Lambda_\epsilon\vert_{TM}, d\Lambda_\epsilon\vert_{TM})_{s\in[0,1]}
$$
is a family of SHS's on $M$, which deform model B, and such that $\Lambda_\epsilon^s$ is contact for $s<1$. 

\vspace{0.5cm}

Focusing now in the case $L=1$, $Q=0$, we have $\Lambda_\epsilon^s=\epsilon \lambda+(s\kappa+1-s)h_\epsilon\alpha$, which we can explicit write as
$$
\Lambda_\epsilon^s=\left\{\begin{array}{ll} \epsilon \lambda_\pm + h_\epsilon\alpha_\pm,& \mbox{ in } M_P^\pm=Y\times (\Sigma\backslash \mathcal{U})\\
f_\epsilon(r)d\theta + (s\kappa(F_\pm(r)\pm 1)+1-s)h_\epsilon\alpha,& \mbox{ in } M_C^\pm=Y\times\left(\mathbb{R}^\pm\cap \{|r|>1-\delta\}\right)\times S^1\\
\epsilon d\theta+(1-s)h_\epsilon \alpha, &\mbox{ in } M_Y= Y \times [-1+\delta,1-\delta]\times S^1,
\end{array}\right.
$$
where we denote $f_\epsilon(r)=\epsilon e^{G_\pm(r)}$ (cf.\ (\ref{modelBform})).

\vspace{0.5cm}

One may also check that the Reeb vector field $R_\epsilon^s$ associated to $\mathcal{H}^s_\epsilon$ is given by
$$
R_\epsilon^s=\left\{\begin{array}{ll} R_\epsilon,& \mbox{ in } M_P^\pm\\ 
\frac{1}{D_\epsilon^s +f_\epsilon h_\epsilon \frac{\partial \alpha_r}{\partial r}(R_r)}((h_\epsilon^\prime + h_\epsilon \frac{\partial \alpha_r}{\partial r}(R_r))\partial_\theta-f_\epsilon^\prime R_r),& \mbox{ in } M_C^\pm\\
\partial_\theta/\epsilon, &\mbox{ in } M_Y
\end{array}\right. 
$$ 
where $D_\epsilon^s:= f_\epsilon h_\epsilon^\prime - (s\kappa(F_\pm \pm 1)+1-s)h_\epsilon f_\epsilon^\prime$. In particular, we still have that it is totally degenerate along $M_Y$, and we may still take the big action threshold $T_\epsilon>0$ we had before ($s$-independent), so that every closed $R_\epsilon^s$-orbit with action less than $T_\epsilon$ still corresponds to a critical point or lies in $Y \times \mathcal{U}$.

\begin{lemma}\label{sympcob}
The 2-form $\omega=d(e^\epsilon \Lambda_\epsilon^s)$ is symplectic in $(0,\epsilon_0]\times M$, for sufficiently small $\epsilon_0>0$ and for $0\leq s< 1$, and induces the orientation given by the natural product orientation (where the orientation on $M$ is as described before Lemma \ref{contact}).
\end{lemma}
\begin{proof}
We have 
$$ 
\omega = d(e^\epsilon \Lambda_\epsilon^s)= e^\epsilon(d\epsilon \wedge \Lambda_\epsilon^s + d\Lambda_\epsilon^s)
$$

Over the region $Y \times (\Sigma\backslash \mathcal{U})$, we get $\Lambda_\epsilon^s = \epsilon \lambda_\pm + h_\epsilon \alpha_\pm$ and $h_\epsilon=K_\pm -\epsilon^2 h_\pm$, and hence 
\begin{equation}
\begin{split}
d\Lambda_\epsilon^s&=\epsilon d\lambda_\pm + dh_\epsilon \wedge \alpha_\pm + h_\epsilon d\alpha_\pm\\
&= \epsilon d\lambda_\pm - \epsilon^2 dh_\pm \wedge \alpha_\pm + h_\epsilon d\alpha_\pm
\end{split}
\end{equation}

Therefore, using that $d\lambda_\pm \wedge dh_\pm = \lambda_\pm \wedge d\lambda_\pm=0$ (being pullbacks of 3-forms in $\Sigma$), and that the top power for $d\alpha_\pm$ is $d\alpha_\pm^{n-1}$, we obtain 
\begin{equation}
\begin{split}
\omega^{n+1}&=e^{(n+1)\epsilon}\binom{2}{n+1} h_\epsilon^{n-1}\left(\epsilon d\epsilon\wedge \Lambda^s_\epsilon\wedge d\lambda_\pm - \epsilon^2 d\epsilon \wedge \Lambda^s_\epsilon \wedge dh_\pm \wedge \alpha_\pm\right) \wedge  d\alpha_\pm^{n-1}\\
&=e^{(n+1)\epsilon}\binom{2}{n+1} h_\epsilon^{n-1}\left(\epsilon h_\epsilon d\epsilon\wedge \alpha_\pm \wedge d\lambda_\pm - \epsilon^3 d\epsilon \wedge \lambda_\pm \wedge dh_\pm \wedge \alpha_\pm\right) \wedge  d\alpha_\pm^{n-1}\\
&=e^{(n+1)\epsilon}\epsilon h_\epsilon^{n}\binom{2}{n+1}(d\epsilon \wedge \alpha_\pm \wedge d\alpha_\pm^{n-1}\wedge d\lambda_\pm) + O(\epsilon^3)
\end{split}
\end{equation}
 
Recalling that the orientation of $\Sigma$ coincides with the one for $\Sigma_-$ given by $d\lambda_-$, and is opposite to the one for $\Sigma_+$ given by $d\lambda_+$, and we chose to give $Y$ the $\alpha_-$-orientation, then the above is a positive volume form for $(0,\epsilon_0]\times Y \times \Sigma$, and $\epsilon \leq \epsilon_0$ sufficiently small.  

\vspace{0.5cm}

Over the region $Y \times \mathcal{U}$, we get $\Lambda_\epsilon^s = f_\epsilon d\theta + (s\kappa(F_\pm\pm 1)+1-s) h_\epsilon \alpha$, and therefore 
\begin{equation}
\begin{split}
d\Lambda_\epsilon^s=&f_\epsilon^\prime dr\wedge d\theta + \left(s h_\epsilon \kappa^\prime(F_\pm \pm 1)F_\pm^\prime+  h_\epsilon^\prime (s\kappa(F_\pm \pm 1)+ 1 - s)\right) dr \wedge \alpha\\
& + (s\kappa(F_\pm \pm 1)+ 1 - s) h_\epsilon d\alpha
\end{split}
\end{equation}

We can then compute that 

$$
\omega^{n+1}=e^{\epsilon(n+1)}h_\epsilon^{n-1}(s\kappa(F_\pm \pm 1)+ 1 - s)^{n-1}
$$
$$
.\left\{(s\kappa(F_\pm \pm 1)+ 1 - s)\left((n+1)h_\epsilon f_\epsilon + \binom{2}{n+1} link(\alpha_r)(f_\epsilon^\prime h_\epsilon -f_\epsilon h_\epsilon^\prime)\right)\right.
$$
$$
\left.-s\binom{2}{n+1} link(\alpha_r)f_\epsilon h_\epsilon \kappa^\prime(F_\pm \pm 1)F_\pm^\prime \right\}d\epsilon\wedge dvol \wedge dr \wedge d\theta,
$$

Since the sign of $link(\alpha_r)$ is the opposite to that of $\kappa^\prime(F_\pm \pm 1)F_\pm^\prime$ (whenever the latter expression is non-vanishing), the last summand is non-negative. Moreover, the sign of $link(\alpha_r)$ coincides with that of $f_\epsilon^\prime$ (again when the latter is non-vanishing) and $h_\epsilon^\prime$ can be taken small in comparison to $h_\epsilon$, so that the whole expression can be made positive. Observe we take $s$ strictly less than $1$, since for $s=1$ this expression vanishes identically whenever $\kappa=0$. 

\vspace{0.5cm}

Alternatively, we also have the freedom of taking the function $link(\alpha_r)$ to be as small as we want, by replacing $\alpha_\pm$ by $\alpha_{\pm\epsilon^\prime}$ for sufficiently small $\epsilon^\prime$, which also has the desired effect of making the above expression positive. 

\end{proof}

\subsection{Compatible almost complex structure} 
\label{compJmodelB}

We proceed now to define a $\mathcal{H}_\epsilon^s$-compatible (and non-generic) almost complex structure $J_\epsilon^s$ in $\mathbb{R}\times M$, analogously as we did in Section \ref{ACS}.

Recall that Proposition \ref{explicitJ} and Corollary \ref{compJ0} provide a distribution $\widetilde{\xi}_r$, and vector fields $\widetilde{V}$ and $\widetilde{R}_r$ on $Y\times I$, which respectively coincide with $\xi_r$, $\mp V$ and $\mp R_r$ away from an arbitrary small (but fixed) small neighbourhood of $r=0$, and such that $T(Y\times I)=\widetilde{\xi}_r\oplus \langle \widetilde{V},\widetilde{R}_r\rangle$ is a $d\alpha$-symplectic splitting. 

If we denote by $\xi_\epsilon^s=\ker \Lambda_\epsilon^s$, we may write 
\begin{equation}\label{xiB}
\xi^s_\epsilon=\left\{\begin{array}{ll} \xi_\pm \oplus L^s_\epsilon, &\mbox{ in } M_P^\pm \cup M_C^\pm\\
\widehat{\xi}^s_r \oplus L^s_\epsilon,& \mbox{ in } M_Y
\end{array}\right.
\end{equation}
where $\widehat{\xi}_r^s=\{v-(1-s)h_\epsilon\alpha(v)\partial_\theta:v \in \widetilde{\xi}_r\}$, and $L^s_\epsilon$ is a bundle of (real) rank 2. The latter may be written as 
\begin{equation}
L^s_\epsilon = \left\{\begin{array}{ll}
T\Sigma_\pm,& \mbox{ in crit}(h_\pm)\subseteq M_P^\pm\\ 
\left\langle \nabla h_\pm, h_\epsilon X_{h_\pm}-\epsilon |\nabla h_\pm|^2 R_\pm \right \rangle,& \mbox{ in } M_P^\pm\backslash \mbox{crit}(h_\pm)\\
\langle \widetilde{V}, \widetilde{W}\rangle,& \mbox{ in }  Y\times \mathcal{U},\\
\end{array}\right.
\end{equation}
where 
$$
\widetilde{W}:=f_\epsilon \widetilde{R}_r-(s\kappa(F_\pm\pm 1)+1-s)h_\epsilon \alpha(\widetilde{R}_r)\partial_\theta
$$

Over the region $M_P^\pm=Y\times (\Sigma_\pm\backslash\mathcal{U})$, we have that $L^s_\epsilon$ converges to $T\Sigma_\pm$ as $\epsilon \rightarrow 0$, and so the differential of the projection $\pi: Y \times (\Sigma_\pm\backslash\mathcal{U})\rightarrow \Sigma_\pm\backslash\mathcal{U}$ when restricted to $L^s_\epsilon$ converges to the identity (and is an isomorphism for every $\epsilon$). Therefore, we may define an $\mathbb{R}$-invariant almost complex structure $J^s_\epsilon$ in $\mathbb{R} \times Y \times (\Sigma_\pm\backslash\mathcal{U})$ by asking that it maps $\partial_a$ to $R^s_\epsilon$, and such that   
\begin{equation}\label{splitJ}
J^s_\epsilon\vert_{\xi^s_\epsilon}= J_\pm \oplus (\pi\vert_{L^s_\epsilon})^*j_\pm,
\end{equation} with respect to the splitting above, where $J_\pm$ is a $d\alpha_\pm$-compatible almost complex structure in $\xi_\pm$, and $j_\pm$ is a $d\lambda_\pm$-compatible almost complex structure on $\Sigma_\pm$ which satisfies $j_\pm(\partial_{t_\pm})=\partial_\theta$ on the collar neighbourhoods.   

\vspace{0.5cm} 

Choose now a $d\alpha$-compatible almost complex structure $J_0$ in the Liouville domain  $Y\times I$ as in Corollary \ref{compJ0}, such that its restriction to $\xi_r$ coincides with $J_\pm$ in $\{|r|>1-\delta\}$ (where $V=\pm \partial_r$, $\alpha_r=\phi(r)\alpha_\pm$, $R_r=\phi(r)^{-1}R_\pm$, $\xi_r=\xi_\pm$). 

\vspace{0.5cm} 

Over the region $Y \times \mathcal{U}=Y \times I \times S^1,$ observing that $\Lambda^s_\epsilon(\partial_\theta)=f_\epsilon\neq 0$, we have that the projection $\pi_Y: Y \times \mathcal{U}\rightarrow Y \times I$ gives an isomorphism $$d\pi_Y\vert_{\xi^s_\epsilon}:\xi_\epsilon^s \stackrel{\simeq}{\longrightarrow} T(Y \times I),$$ and thus $\pi_Y^* J_0$ is an almost complex structure on $\xi^s_\epsilon$. We will define the $\mathbb{R}$-invariant $J^s_\epsilon$ to coincide with this almost complex structure on $\xi^s_\epsilon$ only over the region $\mathbb{R}\times M_Y$ (and map $\partial_a$ to $R^s_\epsilon$).

\vspace{0.5cm}

In order to glue the two definitions along the region $\mathbb{R}\times M_C^\pm$, we compute, with the first definition (and recalling that $\alpha(\widetilde{R}_r)=\mp 1$ over this region), that close to $r=\pm1$ we have 
$$
J^s_\epsilon(\widetilde{W})=-(s\kappa(F_\pm \pm 1)+1-s)h_\epsilon \alpha(\widetilde{R}_r) d\pi^{-1}(j_\pm(\partial_\theta))
$$
$$
=(s\kappa(F_\pm \pm 1)+1-s)h_\epsilon \alpha(\widetilde{R}_r)\partial_{t_\pm}=(s\kappa(F_\pm \pm 1)+1-s)h_\epsilon \partial_r,
$$

With the second, and close to $r=\pm 1\mp \delta$, we obtain 
$$J^s_\epsilon(\widetilde{W})=f_\epsilon d\pi_Y^{-1}(J_0(\widetilde{R}_r))=-f_\epsilon d\pi_Y^{-1}(\widetilde{V})= f_\epsilon \partial_r$$

Therefore, in order to match up the two, we can take smooth functions $$\beta^s_{\epsilon,-}:[-1,-1+\delta) \rightarrow \mathbb{R}^+$$
$$\beta^s_{\epsilon,+}:(1-\delta,1] \rightarrow \mathbb{R}^+$$ such that $\beta^s_{\epsilon,\pm}= (s \kappa(F_\pm \pm 1)+1-s) h_\epsilon$ close to $\pm 1$, and $\beta^s_{\epsilon,\pm}= f_\epsilon$ close to $\pm 1 \mp \delta $, and define 
$$
J^s_\epsilon(\widetilde{W})=\beta^s_{\epsilon,\pm}(r) \partial_r
$$

This extends $J^s_\epsilon$ over $L^s_\epsilon$, and since $\pi_Y^*J_0\vert_{\xi_\pm}=J_\pm$, we get a well-defined almost complex structure $J^s_\epsilon$ over the whole model.

\paragraph*{Compatibility} One can check that $J_\epsilon^s$ is $\mathcal{H}_\epsilon^s$-compatible by explicit computation on the above splittings.

\begin{remark}
Observe that the splitting (\ref{xiB}) is holomorphic, and recalling from Section \ref{Liouvilledoms} that $c_1(\widetilde{\xi}_r)=0$, it follows that 
$$
c_1(\xi_\epsilon^1)=c_1(\widetilde{\xi}_r)+c_1(L_\epsilon^1)=0
$$
\end{remark}
  
\subsection{Foliation by holomorphic hypersurfaces} \label{holhyps}

The goal of this section is to construct a foliation of $\mathbb{R}\times M$ by holomorphic hypersurfaces (where the almost complex structure is the one compatible with the \emph{non-contact} SHS $\mathcal{H}_\epsilon^1$). 

\begin{remark} The proof below also gives a similar foliation of hypersurfaces if we take $J_\epsilon^s$ for $s<1$, but they will not be holomorphic along a region of the form $\mathbb{R}\times Y \times (-\delta^\prime,\delta^\prime)\times S^1$ of $M_Y$, where $\delta^\prime>0$ is such that $\widetilde{\xi}_r\neq \xi_r$ and $J_0$ is non-cylindrical along $Y \times (-\delta^\prime,\delta^\prime)$. This is the reason why we had to introduce a SHS deformation along $M_Y$.
\end{remark}  
 
\begin{prop}\label{curves} For $\epsilon \in (0,\epsilon_0]$, there exists a foliation $\mathcal{F}=\bigcup_{a,\nu} H_a^\nu$ of the symplectization of $\mathcal{H}_\epsilon^1$ by $J^1_\epsilon$-holomorphic hypersurfaces, which come in three types (see picture (\ref{fol2})):
\begin{enumerate}
\item (\textbf{cylindrical} hypersurfaces over critical points) For $p \in \mbox{crit}(h_\pm)$, we have a corresponding \emph{cylindrical} hypersurface of the form $\mathbb{R}\times H_p$, where $H_p:=Y \times \{p\}$, which is a copy of the symplectization of $(Y,\alpha_\pm)$.

\item (positive/negative \textbf{flow-line} hypersurfaces contained in one side of the dividing set $\{\rho=0\}$) Hypersurfaces of this type may be parametrized by $$H_a^\nu=\{(a(s),\nu(s)): s \in \mathbb{R}\} \times Y,$$ where $a:\mathbb {R}\rightarrow \mathbb{R}$ is a proper function, and $\nu:\mathbb{R}\rightarrow \mbox{int}(\Sigma_\pm)$ is a non-constant negative/positive reparametrization of a flow-line for $h_\pm$ (with respect to the metric induced by $d\lambda_\pm$ and $j_\pm$ inducing the positive orientation on $\Sigma_\pm$). 

In this case, if $\lim_{s \rightarrow \pm\infty}\nu(s)=p_\pm$, then $H_a^\nu$ is asymptotically cylindrical over $H_{p_\pm}$ (see Appendix \ref{AppSiefring} for a definition), and has exactly one positive and one negative end. If the side containing it is the positive one (and thus $\nu$ connects the maximum to a hyperbolic point), the positive end corresponds to the maximum. If it is $\Sigma_-$ (and thus $\nu$ connects the minimum to a hyperbolic point), it corresponds to the minimum.   

\item (\textbf{page-like} hypersurfaces crossing sides of the dividing set) Hypersurfaces of this type may be written as
$$H_a^\nu=H_a^{\nu,-}\bigcup H_a^{\nu,cst}\bigcup H_a^{\nu,+},$$
where 
$$
H_a^{\nu,cst}=\{a\}\times Y \times [-1,1] \subseteq \mathbb{R}\times M_Y,
$$
$$
H_a^{\nu,+}:=\{(a_+(s),\nu_+(s)):s\in[0,+\infty)\}\times Y\subseteq \mathbb{R}\times M_P^+
$$
$$
H_a^{\nu,-}:=\{(a_-(s),\nu_-(s)):s\in(-\infty,0]\}\times Y\subseteq \mathbb{R}\times M_P^-
$$ Here, $a_+:[0,+\infty)\rightarrow \mathbb{R}$ and $a_-:(-\infty,0]\rightarrow \mathbb{R}$ are proper functions, and $\nu= \nu_+ \sqcup\nu_-$, for $\nu_+: [0,+\infty)\rightarrow \Sigma_+$ a non-constant reparametrization of a negative Morse flow line of $h_+$, and $\nu_-: (-\infty,0]\rightarrow \Sigma_-$ a non-constant reparametrization of a positive Morse flow line of $h_-$, both satisfying $\nu_\pm(0) \in \partial \Sigma_\pm$. If $\lim_{s \rightarrow \pm \infty}\nu_\pm(s)=p_\pm$, the hypersurface $H_a^\nu$ is asymptotically cylindrical over $H_{p_\pm}$, having exactly two positive ends.    
\end{enumerate}

\begin{figure}[t]\centering
\includegraphics[width=0.65\linewidth]{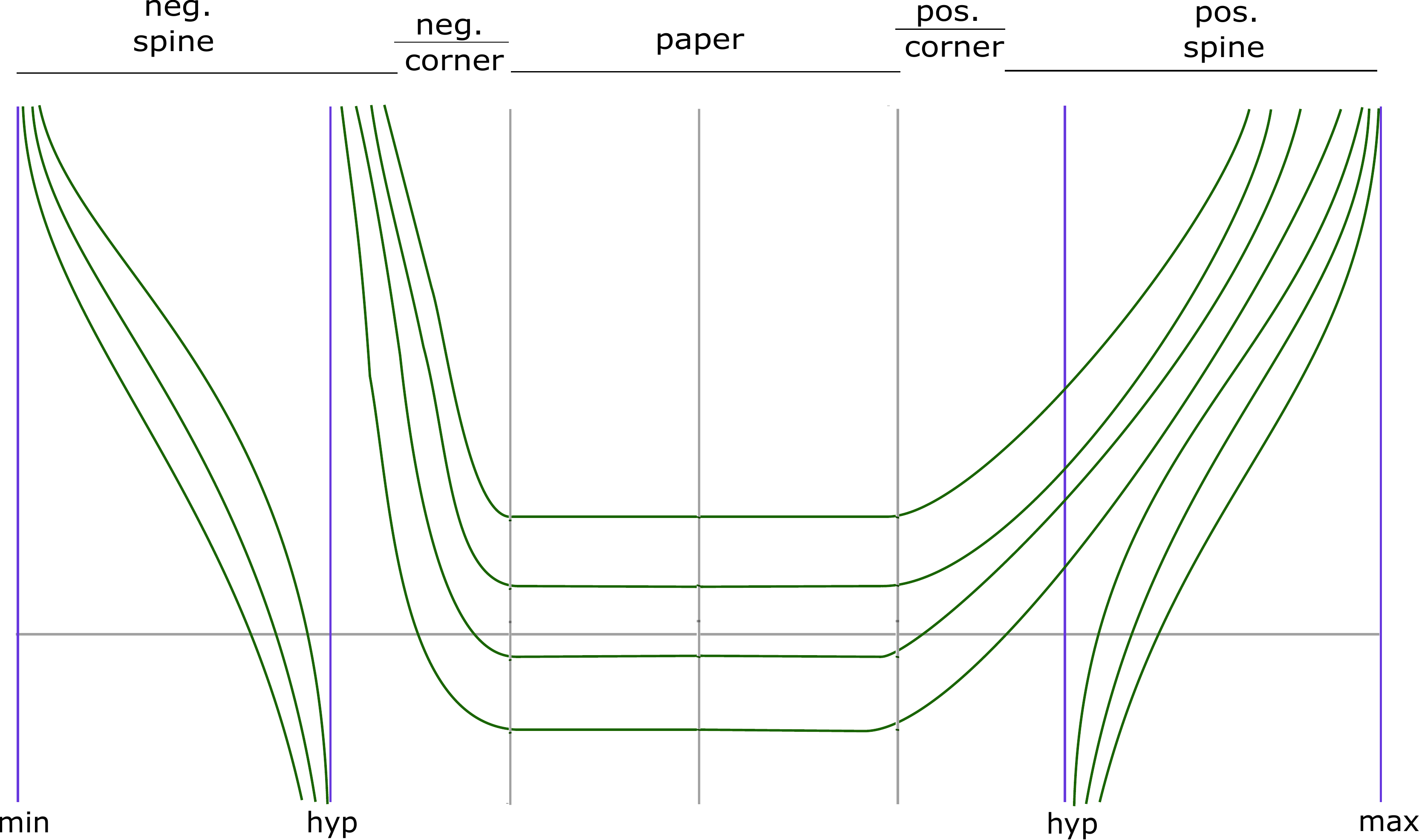}
\caption{\label{fol2} Sketch of the foliation by holomorphic hypersurfaces (with some hypersurfaces crossing sides left out). The vertical blue lines are the cylindrical hypersurfaces corresponding to critical points $min$ (minimum), $hyp$ (hyperbolic) or $max$ (maximum), and the green ones are the non-cylindrical $H_a^\nu$'s. The ones crossing sides are constant in the $\mathbb{R}$-direction along the paper $M_Y$. Compare this foliation to Figure \ref{fol}, its dual model A version.}
\end{figure}
\end{prop}

\begin{proof}
Denote by $\nabla^{g_\pm}$ the connection corresponding to the metric $g_\pm=d\lambda_\pm(\cdot,j_\pm \cdot)$. Observe that $\nabla ^{g_\pm} h_\pm$ takes values in $\xi^1_\epsilon$, and is non-vanishing away from the critical points of $h_\pm$.

\vspace{0.5cm}
 
Let us define a vector field $A$ on $M_P^\pm$ given by 
$$A(y,z)=J^1_\epsilon \nabla^{g_\pm} h_\pm(z)\in \xi^1_\epsilon, \mbox{ for } (y,z) \in M_P^\pm$$ 

Observe that we may write 
$$A = d\pi^{-1}j_\pm \nabla^{g_\pm} h_\pm =  X_{h_\pm} - \frac{\epsilon |\nabla h_\pm |^2}{h_\epsilon}R_\pm$$

We see that $A$ is a linear combination of the same vector fields as $R^1_\epsilon$, and given that the former is not colinear with $A$ since $A \in \xi^1_\epsilon$, we have that $$R_\pm \in \langle A, R^1_\epsilon \rangle$$ 

This implies that $$J^1_\epsilon R_\pm \in \langle \nabla^{g_\pm} h_\pm, \partial_a \rangle,$$ so that we may write  
$$
J_\epsilon^1 R_\pm=F\nabla^{g_\pm} h_\pm + G\partial_a,
$$
for smooth functions $F$ and $G$ on $M_P^\pm$. Indeed, one can compute that 
$$
F=\frac{\epsilon h_\epsilon}{h_\epsilon +\epsilon^2|\nabla h_\pm|^2},\; G=-h_\epsilon
$$
Observe that $F>0$, $G<0$.
Then
$$\langle R_\pm, J^1_\epsilon R_\pm\rangle \oplus \xi_\pm =\left\langle F\nabla^{g_\pm} h_\pm+G\partial_a\right\rangle \oplus TY$$ a $J^1_\epsilon$-invariant and integrable distribution over $M_P^\pm$. Integrating this distribution, one gets the hypersurfaces $H_a^\nu$ as in the statement. For instance, for the hypersurfaces of the second type, we take $a$ and $\nu$ to satisfy
$$\dot{\nu}(s)=\mp F(\nu(s))\nabla^{g_\pm} h_\pm(\nu(s))$$
$$\dot{a}(s)=\mp G(\nu(s))$$ 
For hypersurfaces of the third type, one can check that for $a_\pm$ and $\nu_\pm$ satisfying the same ODEs we can glue each piece as in the statement, to get a hypersurface $H_a^\nu$ whose tangent space is 
$$
TH_a^\nu=\left\{\begin{array}{ll} \langle R_\pm,J_\epsilon^1 R_\pm\rangle \oplus \xi_\pm,& \mbox{ over } M_P^\pm\cup  M_C^\pm\\
T(Y \times I), &\mbox{ over } M_Y
\end{array}\right.,
$$
which is clearly $J^1_\epsilon$-invariant.

It is clear that varying $a$ and $\nu$ we obtain a foliation of $\mathbb{R}\times M$, and the assertion about the asymptotic behaviour of the ends of $H_a^\nu$ is deduced by the fact that orientations change as we move from one side to the other.

\vspace{0.5cm} 

The condition that the holomorphic hypersurfaces $H_a^\nu$ are asymptotically cylindrical follows from the fact that they project to $\Sigma$ as the (reparametrized) Morse flow line $\nu$. This flow line may be written asymptotically as $\nu(s)=\mbox{exp}\;\eta_\pm(s)$ for $\eta_\pm(s)\in T_{z_\pm}\Sigma$ depending smoothly on $s$, where $z_\pm$ denotes a limiting critical point where $H_a^\nu$ has a positive/negative end, $\pm s>R$ for some $R>0$. The map $\eta_\pm$, and all of its derivatives, are exponentially decaying with respect to $s$. One can then extend $\eta_\pm$ as a section of $(\xi_\epsilon^1)^{\perp}_{H_{z_\pm}}$, the $d\Lambda_\epsilon^1$-symplectic normal bundle to $H_{z_\pm}$, which is the trivial fiber bundle with fiber $T_{z_\pm}\Sigma$, by being constant in the $Y$-direction. Thus, the exponential decay condition with respect to $Y$-derivatives is trivially satisfied. Since the function $a$ satisfies $a^\prime=G \circ \nu$, and $G$ is non-vanishing away from $\mathcal{U}$, $a$ has non-vanishing derivative over $M_P^\pm$, and thus it is invertible. We may then choose the reparametrization $\pm t=a(s)$, and define $\widetilde{\nu}(t)=\nu(a^{-1}(\pm t)),$ $\widetilde{\eta}_\pm(t)=\eta_\pm(a^{-1}(\pm t))$. We get 
\begin{equation}
\begin{split} H_a^\nu \cap (\{\pm a > R^\prime\}\times M)=&\bigcup_{(t,y) \in [R^\prime,+\infty)\times H_{z_\pm}}(\pm t,\mbox{exp}_y\widetilde{\eta}_\pm(t))\\=&\bigcup_{(t,y) \in [R^\prime,+\infty)\times H_{z_\pm}}\widetilde{\mbox{exp}}_{(\pm t,y)}\widetilde{\eta}_\pm(t,y),
\end{split}
\end{equation}
where $\widetilde{\mbox{exp}}$ denotes the exponential map $\widetilde{\mbox{exp}}_{(a,p)}(b,v)=(a+b,\mbox{exp}_pv)$. This, as defined in Appendix \ref{AppSiefring}, means that $H_a^\nu$ is positively/negatively asymptotically cylindrical over $H_{z_\pm}$, with asymptotic representatives $\widetilde{\eta}_\pm$.    
\end{proof}

\begin{remark}\label{symhyp}
All of the hypersurfaces $H_a^\nu$ for $\nu$ a flow line which does not change sides of the dividing set may be identified with the symplectization of $(Y,\alpha_\pm)$ (the sign depending on the side), and the ambient almost complex structure $J^s_\epsilon$, as well as the contact form $\Lambda^s_\epsilon$, are compatible with this identification. This is clear for a hypersurface of the form $\mathbb{R}\times H_p$, since $R^s_\epsilon\vert_{\mathbb{R}\times H_p}=\frac{1}{h_\epsilon(p)}R_\pm$ is a multiple of $R_\pm$, $\Lambda^s_\epsilon\vert_{\mathbb{R}\times H_p}=h_\epsilon(p)\alpha_\pm$, $\xi^s_\epsilon\vert_{\mathbb{R}\times H_p}=\xi_\pm$ and $J^s_\epsilon\vert_{\mathbb{R}\times H_p}$ preserves $\xi_\pm$ by construction. For a flow-line hypersurface $H_a^\nu$, and $p_\pm=\lim_{s\rightarrow \pm\infty}\nu(s)$ a critical point, the projection map  
$$P_a^\nu: H_a^\nu=\{(a(s),\nu(s)):s\in \mathbb{R}\}\times Y \rightarrow \mathbb{R}\times H_{p_\pm}$$
$$(a(s),\nu(s),y)\mapsto  (s,y,p_\pm)$$ 
gives an identification of $H_a^\nu$ with $\mathbb{R}\times H_{p_\pm}$, which, up to replacing $s$ by $-s$ over $M_P^-$, can be seen to preserve the data as well. In particular it preserves the almost complex structure and maps holomorphic curves in $H_a^\nu$ to holomorphic curves in $\mathbb{R}\times H_{p_\pm}$, and is identified with the symplectization of $(Y,\alpha_\pm)$. This can be easily seen: write 
$$
TH_a^\nu=\langle \partial_s=\pm(-F\nabla^{g_\pm}h_\pm - G\partial_a), R_\pm=J_\epsilon^s (\pm \partial_s) \rangle \oplus \xi_\pm 
$$
$$
T(\mathbb{R}\times H_{p_\pm})=\langle \partial_a, R_\pm\rangle \oplus \xi_\pm
$$

Recall that $J_\epsilon^s=J_\pm \oplus (\pi\vert_{L_\epsilon^s})^*j_\pm$ on $M_P^\pm$, for a line bundle $L_\epsilon^s$ isomorphic to $T\Sigma_\pm$. In particular, $J_\epsilon^s$ acts as $J_\pm$ in the $\xi_\pm$-summand of both $TH_a^\nu$ and $T(\mathbb{R}\times H_{p_\pm})$, and the differential $dP_a^\nu$ of $P_a^\nu$ is the identity on $TY$. Also, we have
$$
J_\epsilon^s(dP_a^\nu(\pm \partial_s))=J_\epsilon^s(\partial_a)= R_\pm = dP_a^\nu(R_\pm)=dP_a^\nu(J_\epsilon^s(\pm \partial_s))  
$$
So $P_a^\nu$ is holomorphic if we change the orientation of the $s$-parameter along $M_P^-$. 

Of course, one needs to be aware that the Conley--Zehnder indices of the orbits, and hence the index of curves, is not preserved under this identification (see Proposition \ref{indexthm}).

The hypersurfaces which do cross sides can be thought of as being cylindrical completions of the Liouville domain $Y \times I$, and the portions which lie in $M_P^\pm$, their cylindrical ends, can be identified with the symplectization of $(Y,\alpha_\pm)$ by projecting in the same way as above. Observe that the restriction of $J_\epsilon^1$ to the $H_a^{\nu,cst}$ piece of every such hypersurface is naturally identified with $J_0$ on $Y\times I$. 
\end{remark}

\begin{remark}\label{Weinstein} Assume that $(Y,\alpha_\pm)$ both satisfy the Weinstein conjecture, so that there are closed $\alpha_\pm$-Reeb orbits. Then, the hypersurfaces $H_a^\nu$ which stay on the same side of the dividing set contain finite energy $J^s_\epsilon$-holomorphic curves (for \emph{every} $s$). These are of the form
$$u^a_\gamma: \mathbb{R} \times S^1 \rightarrow \mathbb{R} \times M$$
$$u^a_\gamma(s,t)= (a(s),\nu(s),\eta(t)),$$ where $\nu$ and $a$ are as above, and $\eta$ is a \emph{closed} $\alpha_\pm$-Reeb orbit, satisfying $\lim_{s\rightarrow \pm\infty}\nu(s)=p_\pm$ for $p_\pm$ critical point. The asymptotic behaviour of $u_\gamma$ is thus the same as the hypersurface $H_a^\nu$ containing it. These cylinders correspond to trivial cylinders under the identification of the previous Remark \ref{symhyp}.

In general, there is no reason why holomorphic cylinders crossing sides should exist, and, in fact, there are examples for which they don't (see Section \ref{5dmodel}). 

\end{remark}

\begin{remark}\label{admissible}
In the language and notation of Section \ref{app} in Appendix \ref{AppSiefring}, the holomorphic foliation $\mathcal{F}$ is compatible with the $H_\epsilon$-admissible fibration $(\Sigma,\pi_\Sigma,Y,H_\epsilon)$, where $\pi_\Sigma:M=Y \times \Sigma \rightarrow \Sigma$ is the natural projection, and $H_\epsilon$ a suitable Morse function on $\Sigma$, defined as follows.

Following \cite{LW}, we take a function $G_0:[-1,1] \rightarrow \mathbb{R}$, such that it is constant equal to $-1$ in $[-1,-1+\delta),$ and equal to $1$ in $(1-\delta,1]$, $G_0(0)=0$, and $G_0^\prime\geq 0$ with strict inequality in the interval $(-1+\delta,1-\delta)$. Define also a function $\Gamma: [-1,1]\rightarrow \mathbb{R}$ so that $\Gamma(r)=\mp e^{\mp r}$ in $\{|r|>1-\delta-\nu\}$, for some small $\nu>0$ (so that $\Gamma^\prime$ does not vanish where $G_0^\prime=0$)), $\Gamma(0)=0$, and sign$(r\Gamma(r))< 0$ for $r\neq 0$. Take $$G_\epsilon(r)=G_0(r)+\epsilon\Gamma(r),$$ where $\epsilon$ is chosen small enough so that $G_\epsilon^\prime > 0$. 

Now define 
$$H_\epsilon=\left\{\begin{array}{ll} \pm 1 \mp \epsilon h_\pm,& \mbox{ in } \Sigma\backslash \mathcal{U} \\
G_\epsilon,& \mbox { in } \mathcal{U} 
\end{array}\right.
$$

Observe that this is a (well-defined and smooth) Morse function on $\Sigma$ whose critical points coincide precisely with those of $h_\pm$. 

We can easily construct a metric $g_\Sigma$ on $\Sigma$ which coincides with $g_\pm$ in $M_P^\pm$, and every positive flow line for $h_-$ is positive for $H_\epsilon$, as well as every negative flow line for $h_+$, so that every leaf in $\mathcal{F}$ projects to a positive flow line of $H_\epsilon$. 

The binding of $\mathcal{F}$ is $\bigcup_{p \in \crit(H_\epsilon)}H_p$, where each $H_p$ is a strong stable hypersurface. The sign function $\sign:\crit(H_\epsilon)\to\br{-1, 1}$ is defined by $\sign(p)=\pm 1$, if $p \in \crit(h_\pm)\subseteq \Sigma_\pm$, by observing that the restriction of $d\pi_\Sigma$ to the symplectic normal bundle to $H_p$ preserves the orientation over $\Sigma_-$, and reverses it over $\Sigma_+$ (with respect to the orientation on $\Sigma$).

We will use this fact in Section \ref{curvesvshyp}, together with Theorem \ref{posdegree} in Appendix \ref{AppSiefring}, to show that every holomorphic curve in $\mathbb{R}\times M$ has to lie in a leaf of $\mathcal{F}$.

\end{remark}

\subsection{Index Computations}

Let us begin by observing that for closed Reeb orbits $\gamma_p: S^1 \rightarrow M$ of the form $\gamma_p(t)=(\gamma(t),p)$, for $p \in $ crit$(h_\pm)$ and $\gamma$ a closed $\alpha_\pm$-Reeb orbit, we have a natural identification 
\begin{equation}\label{ident}
\Gamma(\gamma_p^*\xi^s_\epsilon)=\Gamma(\gamma^*\xi_\pm)\oplus \Gamma(L^s_\epsilon\vert_p) = \Gamma(\gamma^*\xi_\pm)\oplus \Gamma(T_p\Sigma_\pm),
\end{equation}
where the second summand denotes the space of sections of the trivial line bundle over $S^1$ with fiber $L^s_\epsilon\vert_p=T_p\Sigma_\pm$. Therefore any trivialization $\tau$ of $\gamma^*\xi_\pm$ extends naturally to a trivialization of $\gamma_p^*\xi^s_\epsilon$, which we shall also denote by $\tau$. 

\begin{prop}\label{indexthm} 

Consider a trivialization $\tau$ of $\gamma^*\xi_\pm$ over the simply covered $\alpha_\pm$-Reeb orbit $\gamma$, inducing a natural trivialization $\tau^l$ of $(\gamma^l)^*\xi_\pm$ for every $l$. Then for each sufficiently small $\epsilon>0$ there exists a covering threshold $N_\epsilon$, satisfying $\lim_{\epsilon \rightarrow 0}N_\epsilon=+\infty$, such that the Conley--Zehnder index of $\gamma^l_p$ for $l\leq N_\epsilon$ with respect to the induced trivialization of $(\gamma_p^l)^*\xi^s_\epsilon$ given by the identification (\ref{ident}) is 

$$
\mu_{CZ}^{\tau^l}(\gamma^l_p)=\left\{
\begin{array}{ll}
\mu_{CZ}^{\tau^l}(\gamma^l)+1,& \mbox{ if } p \mbox{ is the maximum or minimum}\\
\mu_{CZ}^{\tau^l}(\gamma^l),& \mbox{ if } p \mbox{ is hyperbolic}
\end{array} \right.
$$

Here, $\mu_{CZ}^\tau(\gamma)$ denotes the Conley--Zehnder index of $\gamma$ in the case this Reeb orbit is non-degenerate or Morse--Bott.

\begin{remark}
In the Morse--Bott case, one needs to add suitable weights adapted to the spectral gap of the corresponding operator to obtain $\mu_{CZ}^\tau(\gamma)$, as explained e.g.\ in \cite{Wen1}.
\end{remark}

\end{prop}

\begin{proof}
We need to compute the asymptotic operator associated to the Reeb orbit $\gamma_p$. This operator looks like 
$$
\mathbf{A}_{\gamma_p}:W^{1,2}(\gamma_p^*\xi^s_\epsilon)\rightarrow L^2(\gamma_p^*\xi^s_\epsilon)
$$
$$ \mathbf{A}_{\gamma_p}\eta=-J^s_\epsilon(\nabla_t \eta - T\nabla_\eta R_\epsilon),
$$
where $\nabla$ is any symmetric connection in $M$, and $T=\Lambda_\epsilon(\dot{\gamma}_p)=h_\epsilon(p) \alpha_\pm(\dot{\gamma})$ is the action of $\gamma_p$ (which we always assume to be parametrized so that $\dot{\gamma}_p=TR_\epsilon$).

If we let $H:\Sigma_\pm \backslash\mathcal{U} \rightarrow \mathbb{R}$ be given by 
$$
H:=\frac{1}{K_\pm+\epsilon^2(|\nabla h_\pm|^2-h_\pm)}
$$

then we may write the Reeb vector field as 
$$ R_\epsilon= H(R_\pm + \epsilon X_{h_\pm}),$$

and we have $H(p)=\frac{1}{h_\epsilon(p)}$, so that $$TH(p)=\alpha_\pm(\dot{\gamma})=:T_\gamma$$

We also have $d_pH=0$, since $d_ph_\pm=d_p|\nabla h_\pm|^2=0$ (the latter equality can be checked in local coordinates). Since we may always choose unitary trivializations of $(\gamma_p^*\xi^s_\epsilon,J^s_\epsilon,d\Lambda^s_\epsilon\vert_{\gamma_p^*\xi^s_\epsilon})$ making $J^s_\epsilon$ the standard almost complex structure, we may choose the connection $\nabla$ (at least locally) so that $J^s_\epsilon$ is parallel (i.e.\ we have $\nabla J^s_\epsilon=J^s_\epsilon \nabla$). We can then compute
\begin{equation}
\begin{split}
J^s_\epsilon \nabla R_\epsilon &= J^s_\epsilon (dH \otimes R_\pm + H \nabla R_\pm + \epsilon dH \otimes j_\pm \nabla h_\pm+ \epsilon H \nabla(j_\pm \nabla h_\pm))\\
&=dH \otimes J_\pm R_\pm + H J_\pm \nabla R_\pm -\epsilon dH \otimes \nabla h_\pm - \epsilon H \nabla^2h_\pm
\end{split}
\end{equation}

Observe that $dH\vert_{\Gamma(\gamma_p^*\xi^s_\epsilon)}=0$, and therefore, if $\eta = (\eta_1,\eta_2) \in W^{1,2}(\gamma_p^*\xi^s_\epsilon)=W^{1,2}(\gamma^*\xi_\pm)\oplus W^{1,2}(T_p\Sigma_\pm)$, we get 

\begin{equation}
\begin{split}
TJ^s_\epsilon \nabla_\eta R_\epsilon &= TH(p)J_\pm \nabla_{\eta_1}R_\pm - \epsilon TH(p)\nabla^2_ph_\pm\\
&= T_\gamma J_\pm \nabla_{\eta_1}R_\pm - \epsilon T_\gamma\nabla_p^2h_\pm
\end{split}
\end{equation}

Since $\nabla_t\eta=(\nabla_t \eta_1,\nabla_t \eta_2)$ and $J_\epsilon = J_\pm \oplus j_\pm$ preserves this splitting, we conclude that the asymptotic operator $\mathbf{A}_{\gamma_p}$ splits into a direct sum of the form 
\begin{equation}
\label{splitop}
\mathbf{A}_{\gamma_p}=\mathbf{A}_\gamma \oplus \mathbf{A}_{\epsilon,\gamma}^\pm,
\end{equation}
where 
$$\mathbf{A}_{\epsilon,\gamma}^\pm=-i \nabla_t - \epsilon T_\gamma\nabla^2_p h_\pm$$ Here, we have identified $j_\pm\vert_{T_p \Sigma_\pm}$ with the standard almost complex structure $i$. We shall refer to the operator $\mathbf{A}_{\epsilon,\gamma}^\pm$ as the \emph{normal asymptotic operator} associated to $\gamma$ (cf.\ Sec. \ref{curvesvshyp}, and App. \ref{AppSiefring}). 
 
Additivity of $\mu_{CZ}$ implies that
$$\mu_{CZ}^\tau(\gamma_p)=\mu_{CZ}^\tau(\gamma)+\mu_{CZ}(\mathbf{A}_{\epsilon,\gamma}^\pm)$$

Now, observe that all of the above is also valid if we replace $\gamma$ by $\gamma^l$. For $\epsilon>0$ small enough, the operator $\mathbf{A}_{\epsilon,\gamma^l}^\pm$ is non-degenerate, and the smaller $\epsilon$ is the larger we are allowed to take $T_{\gamma^l}$. Observing that we may take $T_\epsilon$ so that every orbit of action less than $T_\epsilon$ satisfies that this operator is non-degenerate, we have, as in Section \ref{indexcomp}, that
\begin{equation} \label{CZ1}
\mu_{CZ}(\mathbf{A}_{\epsilon,\gamma^l}^\pm)=1-ind_p(h_\pm) \in \{0,1\},
\end{equation}

and this computation is valid for $l\leq N_\epsilon$, where $$N_\epsilon:= \max\{n \in \mathbb{N}: \gamma^l \mbox{ has action less than } T_\epsilon\}$$
 
\end{proof}

\subsection{Holomorphic curves lie in hypersurfaces}
\label{curvesvshyp}

In this section, we shall make use of intersection theory for punctured holomorphic curves and holomorphic hypersurfaces, as outlined in Appendix \ref{AppSiefring}, to which we refer for the relevant definitions and notation. 

\vspace{0.5cm}

If we assume that $(Y,\alpha_\pm)$ are non-degenerate contact manifolds (which is enough if we only look at Reeb orbits over critical points of $h_\pm$, and we can always do---see Section \ref{pert}), then Assumptions \ref{standing-assumptions} are satisfied, where, using the nomenclature in the appendix, we set:

\begin{enumerate}\renewcommand{\theenumi}{\alph{enumi}}
\item $(M,\Lambda^1_\epsilon, \Omega_\epsilon)$ is our stable Hamiltonian manifold.
\item $V_\pm$ are any codimension-2 submanifold of the form $H_p:=Y \times \{p\}$ for $p$ a critical point of $h_\pm$, which are strong stable Hamiltonian hypersurfaces in $M$.
\item $\widetilde{J}= J^1_\epsilon$, which is compatible with every $H_p$.
\item $\widetilde{g}_J$ and $\widetilde{exp}$ are the associated metric and exponential map, as defined in the appendix.
\item $\widetilde{V}$ is any of the $H_a^\nu$'s, which we have shown are asymptotically cylindrical.
\end{enumerate}

Indeed, we explicitly check some of the conditions. The Reeb vector field is tangent to $H_p$, to which $\Lambda^1_\epsilon$ restricts as $\Lambda^1_\epsilon\vert_{H_p}=h_\epsilon(p)\alpha_\pm$, since $\lambda_\pm=d^{\mathbb{C}}h_\pm$ vanishes at $p$, and so $H_p$ is a strong contact hypersurface.

Moreover,  
\begin{equation} \label{split} 
\xi^1_\epsilon\vert_{H_p}=\xi_\pm \oplus (L^1_\epsilon)_p= \xi_\pm \oplus T_p\Sigma_\pm,
\end{equation} and $d\Lambda^1_\epsilon\vert_{H_p}=h_\epsilon(p)d\alpha_\pm + \epsilon d\lambda_\pm$. Therefore, the splitting (\ref{split}) yields a decomposition of $\xi^1_\epsilon\vert_{H_p}$ into $\xi_\pm=\xi^1_\epsilon\cap TH_p$, and the symplectic normal bundle $T_p\Sigma_\pm=\xi_\pm^{d\Lambda^1_\epsilon}$, both of which are preserved by $J^1_\epsilon=J\vert_{\xi_\pm}\oplus j_\pm$. This means that $J^1_\epsilon$ is $H_p$-compatible.

\vspace{0.5cm}

As we computed in Theorem \ref{indexthm}, we have that the splitting of the asymptotic operator $\mathbf{A}_{\gamma_p^l}$ into tangent and normal components is given by $$\mathbf{A}_{\gamma^l_p}=\mathbf{A}_{\gamma^l} \oplus \mathbf{A}_{\epsilon,\gamma_p^l}^\pm,$$
where the normal operator acting on $\Gamma(T_p\Sigma)$ is  
$$\mathbf{A}_{\epsilon,\gamma_p^l}^\pm=-J_0 \nabla_t - \epsilon T_{\gamma^l}\nabla^2_p h_\pm,$$ which has Conley--Zehnder index $\mu_{CZ}(\mathbf{A}_{\epsilon,\gamma_p^l}^\pm)=1-\mbox{ind}(h_\pm)$. Observe that, in order to be coherent with Section \ref{app} below, we may also write this as 
$$\mu_{CZ}(\mathbf{A}_{\epsilon,\gamma_p^l}^\pm)=\sign(p)(\mbox{ind}_p(H_\epsilon)-1),$$ where we use the Morse function $H_\epsilon$ and the sign function as defined in Remark \ref{admissible}, so that $\mathcal{F}$ is compatible with $(\Sigma,\pi_\Sigma,Y,H_\epsilon)$.

\begin{prop}
\label{uniqhyp}
Suppose $u:\dot{S}\rightarrow \mathbb{R}\times M$ is a finite energy $J^1_\epsilon$-holomorphic curve which has all of its positive ends asymptotic to Reeb orbits of the form $\gamma_p^l$, with $l\leq N_\epsilon$. Then the image of $u$ is contained in a  leaf of the foliation $\mathcal{F}$.
\end{prop}
\begin{proof}
Given such a curve $u$, by Stokes' theorem we get that its negative ends have action bounded by $T_\epsilon$, and so also correspond to critical points. Since all of its asymptotics project to $\Sigma$ as points, we may define a map $v:=\pi_\Sigma \circ u: S\rightarrow \Sigma$ between \emph{closed} surfaces. We are therefore in the situation of Theorem \ref{posdegree} of the Appendix \ref{AppSiefring}, which, since the map $\sign$ is surjective, implies the result.
\end{proof}

\begin{remark} Proposition \ref{uniqhyp} reduces the study of $J_\epsilon^1$-holomorphic curves and buildings inside $M$ to the study of $J_0$-curves/buildings inside the completion $W_0=\widehat{Y\times I}$ of the cylindrical Liouville semi-filling $Y\times I$. Here, $J_0$ is an almost complex structure as in Corollary \ref{compJ0}, and $W_0$ is obtained by attaching cylindrical ends to $\partial(Y\times I)=(Y_-,\alpha_-)\bigsqcup(Y_+,\alpha_+)$. Observe that holomorphic buildings inside $W_0$ are distributed along a main level, which can be identified with $Y\times I$ itself, and perhaps several upper levels, which come in two types, depending on whether they correspond to the symplectization of $(Y,\alpha_+)$ or $(Y,\alpha_-)$ (see Figure \ref{coolbuilding}). In further sections, we will refer to these upper levels as \emph{right} or \emph{left}, respectively. 

\begin{figure}[t]\centering
\includegraphics[width=0.60\linewidth]{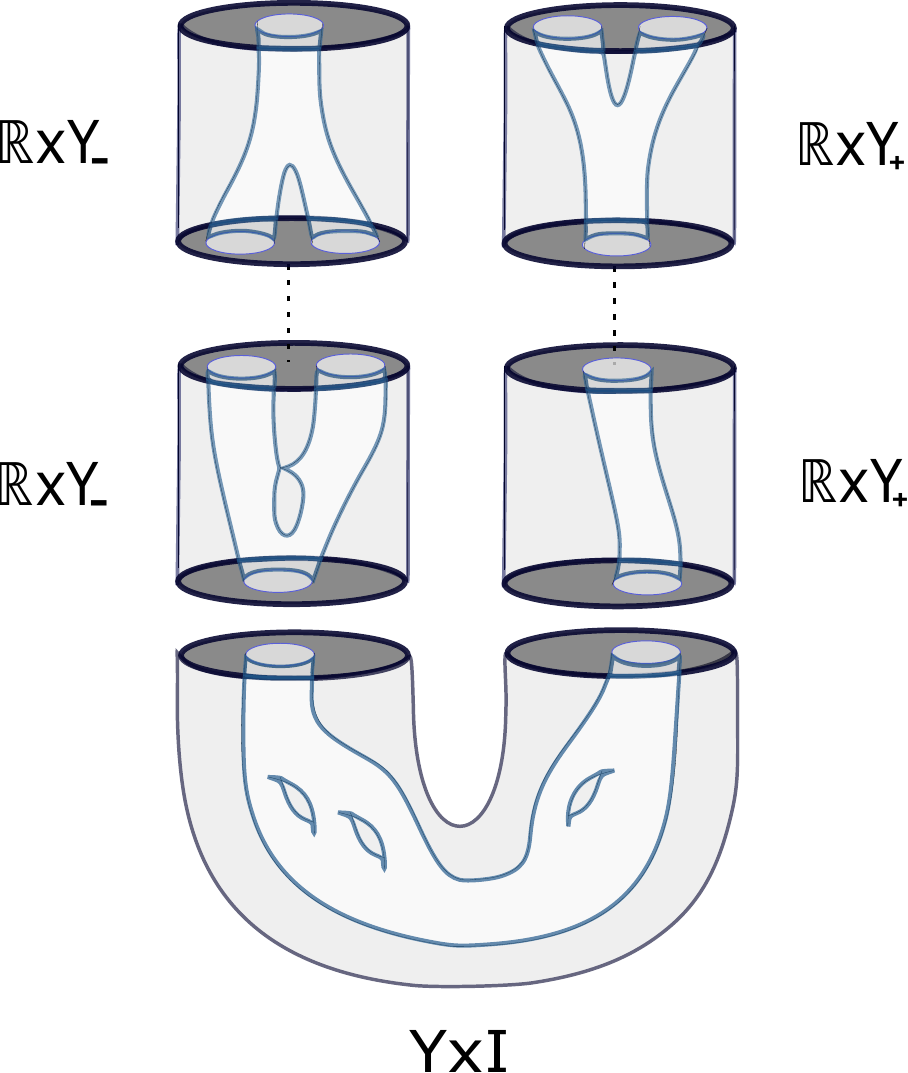}
\caption{\label{coolbuilding} A building in $W_0=\widehat{Y\times I}$, the completion of the cylindrical Liouville semi-filling $Y\times I$.}
\end{figure}

\end{remark}

\subsection{Regularity inside a hypersurface vs. regularity in the symplectization for curves of genus zero}

Consider $u:\dot{S}\rightarrow \mathbb{R}\times M$ a finite energy $J_\epsilon^1$-holomorphic curve, with asymptotics of the form $\gamma_p^l$, with $l\leq N_\epsilon$, and genus $g$. By Proposition \ref{uniqhyp}, $u(\dot{S})\subseteq H$, for some $H \in \mathcal{F}$. If $H$ is positively/negatively asymptotically cylindrical over $H_{p_\pm}$, then the positive/negative asymptotics of $u$ are of the form $\gamma^l_{p_\pm}$, for some simply covered orbit $\gamma\subseteq H_\pm$. Denote by $\Gamma=\Gamma^+\cup \Gamma^-$ the set of punctures of $u$ (where $\Gamma^\pm$ is the set of positive/negative punctures). Define a partition $\Gamma^\pm=\Gamma^\pm_0\cup \Gamma^\pm_1 \cup \Gamma^\pm_2$, where $\Gamma_i^\pm$ denotes the set of positive/negative punctures for which the corresponding asymptotic is of the form $\gamma_p^l$, with $\mbox{ind}_p(H_\epsilon)=i \in \{0,1,2\}$. 

 
\begin{lemma}\label{regvs}
Assume that $u\subseteq H$ has genus $g=0$, is Fredholm regular inside $H$, and that 
\begin{equation}\label{ATineq2}
\#\Gamma_0^-+\#\Gamma_1^-+\#\Gamma_1^++\#\Gamma_2^-<2
\end{equation}
Then $u$ is Fredholm regular in $\mathbb{R}\times M$.
\end{lemma}

\begin{remark}
Observe that, since all positive/negative asymptotics of $u$ correspond to the same critical point of $H_\epsilon$, one has that $\#\Gamma_i^\pm \neq 0$ implies $\#\Gamma_j^\pm =0$ for $j \neq i$. 
\end{remark}

\begin{proof}
Consider the complex splitting $N_u = N_u^H \oplus u^*N_H$, where $N_u$ is the (generalized) normal bundle to $u$, $N_u^H$ is the normal bundle to $u$ inside $H$, and $N_H$ is the normal bundle to $H$ inside $\mathbb{R}\times M$ (which is a complex line bundle). Then we get a matrix representation of the normal linearised Cauchy--Riemann operator $\mathbf{D}_u^N: W^{1,2,\delta_0}(N_u)\oplus X_\Gamma \rightarrow L^{2,\delta_0}(\overline{Hom}_\mathbb{C}(\dot{S}, N_u))$, given by 
$$\displaystyle
\mathbf{D}_u^N=\left(\begin{array}{cc} 
\mathbf{D}_{u;T_H}^N & \mathbf{D}_{u;NT}^N  \\
\mathbf{D}_{u;TN}^N& \mathbf{D}_{u;N_H}^N 
\end{array}\right),$$ where 
$$
\mathbf{D}_{u;T_H}^N:W^{1,2,\delta_0}(N_u^H) \oplus X_\Gamma\rightarrow L^{2,\delta_0}(\overline{Hom}_{\mathbb{C}}(\dot{S},N_u^H))
$$
and
$$
\mathbf{D}_{u;N_H}^N:W^{1,2,\delta_0}(u^*N_H) \rightarrow L^{2,\delta_0}(\overline{Hom}_{\mathbb{C}}(\dot{S},u^*N_H))
$$
are Cauchy--Riemann type operators, and the off-diagonal operators are tensorial. Here, $X_\Gamma$ is a finite-dimensional vector space of smooth sections which are supported near infinity and parallel to the Morse--Bott submanifolds containing the asymptotics of $u$. WLOG we can assume that they are tangent to $H$.

The assumption that $u$ is regular inside $H$ is equivalent to the surjectivity of $\mathbf{D}_{u;TH}^N$. Moreover, we claim that $\mathbf{D}_u^N$ is upper-triangular in this matrix representation, i.e.\ that 
$$
\mathbf{D}_{u;TN}^N=0
$$
This, which is basically a generalization of Lemma 3.8 in \cite{Wen1}, can be seen as follows: Take any metric making the splitting $T(\mathbb{R}\times M)\vert_H=T_H \oplus N_H$ orthogonal, and take $\nabla$ to be the associated Levi--Civita connection. This connection is symmetric and preserves this orthogonal splitting. If $\xi \in W^{1,2,\delta_0}(N_u^H)$, take a smooth path $t \rightarrow u_t$ of maps of class $W^{1,p,\delta}(\dot{S},H)$, with image inside $H$, such that $\partial_t\vert_{t=0}u_t=\xi$, $u_0=u$. Since $H$ is holomorphic, $\overline{\partial}u_t \in L^2(\overline{Hom}_{\mathbb{C}}(\dot{S},u_t^*TH))$ takes values in $u_t^*TH$. Then $\mathbf{D}_u^N\xi=\pi_N\nabla_t(\overline{\partial}u_t)\vert_{t=0}$ takes values in $N_u^H \subseteq u^*TH$, since $\nabla$ preserves the splitting, and the claim follows. 

Therefore, to show that $\mathbf{D}_u^N$ is surjective, it suffices to show that $\mathbf{D}_{u;N_H}^N$ is. For this, we can use the automatic transversality criterion in \cite{Wen1}, for the case of a line bundle. We need to check that 
\begin{equation}\label{ATineq}
\mbox{ind}(\mathbf{D}_{u;N_H}^N)>-2+2g+\#\Gamma_{even},
\end{equation}
where we denote by $\Gamma_{even}$ the set of punctures $z$ at which the Conley--Zehnder index of the asymptotic operator of $\mathbf{D}_{u;N_H}^N$ at $z$ is even. Observe that $\Gamma_{even}=\Gamma^+_1 \cup \Gamma^-_1$. Moreover, the operator $\mathbf{D}_{u;N_H}^N$ is asymptotic at each puncture $z$ of $u$ to the normal asymptotic operator associated to the corresponding Reeb orbit $\gamma_{p_z}^{l_z}$ at $z$, whose Conley--Zehnder index is $\mu_N(\gamma_{p_z}^{l_z})=|\mbox{ind}_{p_z}(H_\epsilon)-1|$. One can also show that $N_H$ is a trivial line bundle, by writing down an explicit trivialization, so that $c_1(N_H)=0$. Therefore, by the Riemann--Roch formula, the Fredholm index of $\mathbf{D}_{u;N_H}^N$ is
\begin{equation}\label{indN}
\begin{split}
\mbox{ind}(\mathbf{D}_{u;N_H}^N)=&\chi(\dot{S})+\sum_{z \in \Gamma^+}\mu_N(\gamma_{p_z}^{l_z})-\sum_{z \in \Gamma^-}\mu_N(\gamma_{p_z}^{l_z})\\
=&2-2g-\#\Gamma+\#\Gamma^+_0+\#\Gamma^+_2-\#\Gamma_0^--\#\Gamma_2^-\\
=&2-2g-2\#\Gamma_0^--2\#\Gamma_2^--\#\Gamma^+_1-\#\Gamma_1^-,
\end{split}
\end{equation}
where we used that $\#\Gamma=\sum_{i\in \{0,1,2\}}\#\Gamma^\pm_i$. Using that $\#\Gamma_{even}=\#\Gamma_1^++\#\Gamma_1^-$, one sees that the inequality (\ref{ATineq}) can only be satisfied if $g=0$, in which case it is equivalent to (\ref{ATineq2}). 
\end{proof}

\begin{remark}$\;$
\begin{itemize}[wide, labelwidth=!, labelindent=0pt]
\item Lemma \ref{regvs} implies that the holomorphic cylinders of Remark \ref{Weinstein} are regular index 1 cylinders in $\mathbb{R}\times M$. This generalizes the situation in \cite{LW}. 
\item Given the result of the previous section, combined with the fact that the symplectizations of $(Y_\pm,\alpha_\pm)$ lie naturally inside $M$ over the critical points of $h_\pm$, it is natural to ask whether properties for the SFT, or other Floer-type theories, of the former contact manifolds extend to properties for the corresponding theory for $M$. For instance, consider the following observations. 

If $u$ is a holomorphic curve inside $\mathbb{R} \times H_p=\mathbb{R}\times Y \times \{p\}$ for $p \in \mbox{crit}(h_\pm)$, and we denote by $\mbox{ind}_{Y_\pm}(u)$ its Fredholm index as a curve in the symplectization of $(Y_\pm,\alpha_\pm)$, we see that 
$$
\mbox{ind}(u)= \mbox{ind}_{Y_\pm}(u)+\chi(u)+(1-\mbox{ind}_p(h_\pm))(\#\Gamma^+(u)-\#\Gamma^-(u))
$$
This means that if $\mbox{ind}_{Y_\pm}(u)=1$ and $u$ is a cylinder with one positive and one negative end, then $\mbox{ind}(u)=1$. Moreover, if $p$ is the maximum or the minimum, it follows from Lemma \ref{regvs} that if $u$ is regular in $\mathbb{R}\times Y_\pm$, then it is also in $\mathbb{R}\times M$. In fact, this is also true for hyperbolic $p$, if the $\epsilon$-parameter is chosen small enough, and we consider curves with a fixed energy bound (Cor. \ref{regcylcyls}). Observe that non-degenerate and non-contractible $\alpha_\pm$-orbits embed as non-degenerate and non-contractible Reeb orbits in $\mathbb{R}\times M$. Also, given two $\alpha_\pm$-orbits with $\alpha_\pm$-action less than $T_\epsilon/h_\epsilon(p)$ and lying inside $\mathbb{R} \times H_p$, any holomorphic cylinder joining them has to lie inside $\mathbb{R}\times H_p$, by Proposition \ref{uniqhyp}. Therefore, if a primitive $l_\pm \in [S^1,Y]$ is such that $\alpha_\pm$ is $l_\pm$-admissible (in the sense of definition 10.16 in \cite{Wen3}), then our stable Hamiltonian structure will be $l_\pm$-admissible, after perturbing it to contact data. We obtain an embedding of the filtered cylindrical contact homologies
$$
HC_*^{l_\pm}(Y,\alpha_\pm;T_\epsilon/h_\epsilon(p))\hookrightarrow HC_*^{l_\pm}(M,\Lambda_\epsilon^1;T_\epsilon)
$$
Here, $HC_*^l(M,\Lambda;T)$ denotes the homology for a chain complex with generators the $\Lambda$-Reeb orbits in class $l$, and action less than $T$. This means that information coming from CCH still persists, and suggests that perhaps a direct limit, or a spectral sequence argument as in \cite{Jo} may recover an embedding of the full CCH.

\item We can do a the following sanity check: By Lemma \ref{regvs}, if a genus zero holomorphic curve $u$ lies in a hypersurface $H^\nu_a\simeq W_0$ corresponding to a flow-line $\nu$ joining the minimum to the maximum (which is the generic behaviour), and is regular in $W_0$, it is regular in $\mathbb{R}\times M$. Then, using that $c_1(\xi_\epsilon^1)=c_1(W_0)=0$, we see that 
\begin{equation}
\begin{split}
\mbox{ind}(u)=&(n-2)\chi(u)+\mu_{CZ}^{W_0}(u)+\#\Gamma\\
=&(n-3)\chi(u)+\mu_{CZ}^{W_0}(u)+2-\#\Gamma+\#\Gamma\\
=&\mbox{ind}_{W_0}(u)+2,
\end{split}
\end{equation}
where we denote by $\mbox{ind}_{W_0}(u)$ the index of $u$ inside $W_0$. Geometrically this makes sense, since the extra $2$ summand accounts for $\mathbb{R}$-translation, which is the same as reparametrization in the flow-line direction, and the fact that the unparametrized flow-line underlying $\nu$ has Morse index $1$. Similarly, one can analyze all the cases for which Lemma \ref{regvs} applies, and one finds that the index relations one gets are the ones one can deduce geometrically.
\end{itemize}
\end{remark}

\subsection{Obstruction bundles}\label{obbundles}

In this section we will deal with the existence of obstruction bundles over buildings of holomorphic curves in $\mathbb{R}\times M$. Given such a building $\mathbf{u}$, we wish to have a way of computing the number of gluings one obtains from $\mathbf{u}$ for generic $J$, in the case where its components fail to be transversely cut-out, but \emph{not too badly}. For this, we require that each component of $\mathbf{u}$ consists of either a regular curve, or a curve for which the dimension of the cokernel of the (normal) linearized Cauchy--Riemann operator is constant as the curve varies in its corresponding moduli space.  We will refer to these curves as \emph{not-too-bad}. If such is the case, one can construct an \emph{obstruction bundle} $\mathcal{O}_{\mathbf{u}}$, and orientable bundle whose fiber is the direct sum of the cokernels of these operators, where the sum varies over the components of $\mathbf{u}$ which are not regular. The base of this bundle is the domain of a \emph{pregluing map}. The idea is to preglue the components of $\mathbf{u}$, and impose that the resulting preglued curve be holomorphic. One can view the resulting equation as an obstruction to gluing in the form of a section of the obstruction bundle, whose algebraic count of zeroes is precisely the number of holomorphic gluings of $\mathbf{u}$ one obtains for sufficiently large gluing parameters. This technique, coming from algebraic geometry and Gromov--Witten theory, has been used e.g.\ in \cite{HT1} and \cite{HT2} in the context of ECH, to which we refer the reader for details (see also \cite{Hut} and \cite{HG} for more gentle introductions).

Since \cite{HT1,HT2} deal with very specific situations, we sketch the construction slightly more concretely in a more general case. We will discuss only the immersed and non-nodal case. The non-immersed case is slighlty more technical, and for simplicity we will omit it. One needs to replace the normal bundle with the restriction of the tangent bundle of $W$ to the curves in question, including a suitable Teichm\"uller slice parametrizing complex structures on the domain, and divide by the action of the automorphism group. The nodal case just adds and extra $S^1$ as a gluing parameter for each node to the base of the obstruction bundle.

Let $\mathbf{u}=(u_1^1,\dots,u_{r_1}^1,\dots, u^N_{1},\dots,u_{r_N}^N)$ be a building with $N\geq 1$ floors, $M \geq N$ components (where $M=\sum_{i=i}^N r_i$) denoted $u_j^i$, consisting of $\mathbb{R}$-translation classes of either regular curves or curves which are not-too-bad. We index the not-too-bad curves by $k \in \{1,\dots,r\}$. Assume all of the $u^i_j$ are immersed. Recall that trivial cylinder components are always regular, and glue uniquely if they appear as a component of a building. Therefore, we may assume that $\mathbf{u}$ contains no trivial cylinder components, since they do not affect the counts which we will consider. Then there is an orientable obstruction bundle 
$$
\mathcal{O}_{\mathbf{u}}\rightarrow \prod_{i=1}^{N-1}[R_i,+\infty)\times \prod_{i=1}^N\mathbb{R}^{r_i-1}\times \prod_{i=1}^N\left(\left(\prod_{j=1}^{r_i} \mathcal{M}^i_j\right)\Big/\mathbb{R}\right):=\mathcal{B}_{\mathbf{u}},
$$
where $R_i\gg 0$ is a gluing parameter for gluing the $i-th$ floor to the $i+1-th$ floor of $\mathbf{u}$, and $\mathcal{M}^i_j/\mathbb{R}$ is the moduli space of $\mathbb{R}$-classes of holomorphic curves containing $u^i_j$. The fiber of $\mathcal{O}_\mathbf{u}$ over an element in the base is $\bigoplus_{k=1}^r \mbox{coker}\;\mathbf{D}_k$, where $\mathbf{D}_k$ is the normal linearized Cauchy--Riemann operator of the $k$-th not-too-bad curve. 

The base of the obstruction bundle is that of the pregluing map. Observe that the case $r=0$ is when every component is regular, and one gets a unique gluing. The case $N=1$ is also considered, in which there is no gluing to do, but the components could be multiply covered, in which case we still wish to count their contributions in the form of a count of zeroes of a section of this bundle. Also, recall that a floor of a building is by definition an $\mathbb{R}$-translation class of its components, where the $\mathbb{R}$-shift is the same for all components in the same floor. This explains why we quotient the product of the moduli spaces by the diagonal $\mathbb{R}$-action. However, we still need to consider relative $\mathbb{R}$-shifts between different components of the same floor, since they give gluings which might not be related by a global $\mathbb{R}$-shift, and this is what the factor $\mathbb{R}^{r_i-1}$ keeps track of. In practice, given a floor, we choose a representative under the $\mathbb{R}$-action, and a arbitrary component of the floor which we consider to be at a fixed $\mathbb{R}$-level (say, at level $0$). Then we allow the other components to be shifted arbitrarily, according to a vector in $\mathbb{R}^{r_i-1}$. The set of potential gluings is independent of the chosen fixed component. 

\begin{remark}\label{orbifoldcount} In order to be in a setup which does not rely on the abstract perturbation scheme of polyfold theory, in practice, we should check that the moduli spaces $\mathcal{M}_j^i$ corresponding to not-too-bad curves are orbifolds, potentially with weights, as the worst case scenario. The obstruction bundle above should then be considered an orientable orbibundle, and the counts of zeroes of a generic section should be a weighted, rational count. The way to count, while independent of the abstract perturbation scheme, should be compatible with the counts as defined in polyfold theory (see e.g.\ \cite{HWZ10}, or \cite{Wen3} for a basic exposition on how to count zeroes of sections of orbibundles). In this document, the configurations that we will care about will satisfy this orbifold condition. Moreover, since we will not need to count zeroes explicitly, we will be content with the existence of an obstruction bundle, and the existence of a way of counting zeroes on an orbifold.
\end{remark}

One has a section $\mathfrak{s}$ of this bundle, the \emph{obstruction section}, defined roughly as follows. Take pregluing data $(\mathbf{r}, \mathbf{a},\mathbf{v}) \in \mathcal{B}_\mathbf{u},$ where
$$\mathbf{r}=(r_1,\dots,r_{N-1}),$$
$$\mathbf{a}=(a^1_1,\dots,a^1_{r_1-1},\dots,a^N_1,\dots,a^N_{r_N-1}),$$
$$\mathbf{v}=(v^1_1,\dots,v^1_{r_1},\dots,v^N_1,\dots,v^N_{r_N})$$
We identify the $\mathbb{R}$-translation classes $v^i_j$ with representatives, and we consider the components $v^i_{r_i}$ as the ones at a fixed $\mathbb{R}$-level, say. Consider a tuple $$\Psi=(\psi^1_1,\dots,\psi_{r_1}^1,\dots, \psi^N_1,\dots,\psi_{r_N}^N),$$ where $\psi^i_j$ is a section of the normal bundle of $v^i_j$ belonging to a suitable Hilbert completion of the space of smooth sections. Denote 
$$
\mathbf{a}.\mathbf{v}:=(v^1_1+a^1_1,\dots,v^1_{r_1-1}+a^1_{r_1-1},v^1_{r_1},\dots,v^N_1+ a^N_1,\dots,v^N_{r_N-1}+ a^N_{r_N-1},v^N_{r_N})
$$

One constructs a (only \emph{approximately} holomorphic) preglued curve $\oplus_{(\mathbf{r},\Psi)}\mathbf{a}.\mathbf{v}$ out of the building $\mathbf{v}$, with the property that it converges to $\mathbf{a}.\mathbf{v}$ as every component of $\mathbf{r}$ approaches $+\infty$. This is a standard construction, and is done by exponentiating $\psi^i_j+a^i_j$ along $v^i_j+a^i_j$ (set $a^i_{r_i}=0$), and by using suitable smooth cut-off functions $\beta^i_j$, which equal $1$ in the interior of $v^i_j+a^i_j$ and decay towards its cylindrical ends, translated in a way determined by the gluing parameters $\mathbf{r}$. 

If we impose that $\oplus_{(\mathbf{r},\Psi)}\mathbf{a}.\mathbf{v}$ be holomorphic, we get an equation of the form 
$$
\bigoplus_{i=1}^N\bigoplus_{j=1}^{r_i} \beta^i_j \Theta^i_j(\mathbf{r},\mathbf{a},\Psi)=0
,$$ 
where
$$
\Theta^i_j(\mathbf{r},\mathbf{a},\Psi)= \mathbf{D}^i_j\psi^i_j+\mathcal{F}^i_j(\mathbf{r},\mathbf{a},\Psi)
$$
Here, $\mathbf{D}^i_j$ is the linearized Cauchy--Riemann operator of $v^i_j$, and the second summand involves extra terms arising from the patching construction, mostly non-linear and depending only on the $\psi^m_n$'s for which $v^i_j$ is adjacent to $v^m_n$. For generic $J$, and fixed sufficiently large $R_i$, there exists a unique $\Psi=\Psi_{\mathbf{a},\mathbf{v}}$ such that 
$$
\psi^i_j \perp \ker\mathbf{D}^i_j, \mbox{ for every } i,j. 
$$
$$
\Theta^i_j(\mathbf{r},\mathbf{a},\Psi)=0, \mbox{ if } u^i_j \mbox{ (or equivalently } v^i_j \mbox{) is regular\ }.
$$
$$
\Theta^i_j(\mathbf{r},\mathbf{a},\Psi) \perp \mbox{im}\;\mathbf{D}^i_j, \mbox{ if } u^i_j \mbox{ (or equivalently } v^i_j \mbox{) is not-too-bad\ }.
$$
Then one can define 
$$
\mathfrak{s}(\mathbf{r},\mathbf{a},\mathbf{v})=\bigoplus_{k=1}^r \pi_k  \Theta_k(\mathbf{r},\mathbf{a},\Psi_{\mathbf{a},\mathbf{v}})=\bigoplus_{k=1}^r \pi_k  \mathcal{F}_k(\mathbf{r},\mathbf{a},\Psi_\mathbf{\mathbf{a},v}),
$$
where $\pi_k$ is the orthogonal projection to $\mbox{coker}\; \mathbf{D}_k$. Therefore, there is a 1-1 correspondence between the zeroes of $\mathfrak{s}(\cdot,\cdot,\mathbf{v})$ and the holomorphic gluings of $\mathbf{v}$. Moreover, $\mathfrak{s}(\cdot,\cdot,\mathbf{v})$ is transverse to the zero section for generic $J$. Therefore, in the case where the dimension of the base is the same as the rank of $\mathcal{O}_\mathbf{u}$, its algebraic count of zeroes, which is independent of the gluing parameters as long as they are sufficiently large, is the number of such holomorphic gluings. Observe that in that case the set of such zeroes is not an Euler class, but a \emph{relative} Euler class, by the non-compactness of the domain. In the case where the dimension of the base is strictly less than the rank, then generically there will be no gluing. If it is strictly greater, then gluings will come in (non-isolated) families.


\subparagraph*{Existence of obstruction bundles}\label{OBB} In the case where the hypothesis of the Lemma \ref{regvs} are not satisfied for a possibly multiply-covered curve $u:(\dot{S},j)\rightarrow H\subseteq W=\mathbb{R}\times M$, we wish to have a criterion for the existence of an obstruction bundle for building configurations containing $u$. We will assume that $u$ is is not-too-bad in $H$, which includes the case where $u$ is regular. We prove that $u$ is not-too-bad in $\mathbb{R}\times M$, which provides obstruction bundles in $\mathbb{R}\times M$.

As in the proof of Lemma \ref{regvs}, consider the splitting 
\begin{equation}
\label{splitop}
\displaystyle
\mathbf{D}_u^N=\left(\begin{array}{cc} 
\mathbf{D}_{u;T_H}^N & \mathbf{D}_{u;NT}^N  \\
0 & \mathbf{D}_{u;N_H}^N 
\end{array}\right)
\end{equation}
Assume first, for warmup, that $u$ is regular in $H$, which is equivalent to the surjectivity of the diagonal term $\mathbf{D}_{u;T_H}^N$. This implies that
$$
\mbox{coker} \;\mathbf{D}_u^N = \mbox{coker} \;\mathbf{D}_{u;N_H}^N
$$
In particular, one has 
\begin{equation}\label{warmup}
\dim\mbox{coker} \;\mathbf{D}_u^N  =\dim \mbox{coker} \;\mathbf{D}_{u;N_H}^N=\dim\ker \mathbf{D}_{u;N_H}^N-\mbox{ind}\;\mathbf{D}_{u;N_H}^N
\end{equation}
Therefore, if we assume that $\dim \ker \mathbf{D}_{u;N_H}^N$ is constant as $u$ moves around in its moduli space, the so will be $\dim\mbox{coker} \;\mathbf{D}_u^N$, since the other summand is. 

\vspace{0.5cm}

More generally, in the case where $u$ is not-too-bad inside $H$, we have that $\dim \ker \mathbf{D}_{u;T_H}^N$ only depends on the moduli space of $u$. We will assume this from now on.

The reason why we consider $\ker \mathbf{D}_{u;N_H}^N$ is because the automatic transverality criterion in \cite{Wen1} gives a bound on its dimension. We observe first that, whenever the hypothesis of Lemma \ref{regvs} fail, then we have 
$$
\mbox{ind}\;\mathbf{D}_{u;N_H}^N\leq 0
$$
Indeed, this follows from equation (\ref{indN}), if we assume that $g\geq 1$ or $\#\Gamma_0^-+\#\Gamma_1^-+\#\Gamma_1^++\#\Gamma_2^-\geq 2$. By \cite{Wen1}, we get a bound
\begin{equation}\label{boundker}
\dim \ker \mathbf{D}_{u;N_H}^N \leq K(c_1(N_H,u),\#\Gamma_{even}),
\end{equation}
where 
$$
2c_1(N_H,u)=\mbox{ind}\;\mathbf{D}_{u;N_H}^N-2+2g+\#\Gamma_{even}=-2(\#\Gamma_0^-+\#\Gamma_2^-)\leq 0
$$
and
$$
K(c,G)=\min\{k+l|k,l \in \mathbb{Z}^{\geq 0}, k\leq G, 2k+l>2c\}
$$

In particular, if either $\#\Gamma_0^-$ or $\#\Gamma_2^-$ are non-zero (which is the case only when $u$ lies in the cylindrical hypersurface over the minimum $min$ or maximum $max$), then $\mathbf{D}_{u;N_H}^N$ is injective.  In the case where $c_1(N_H,u)=0$, we have $K(0,G)=1$ for any $G$, and then we obtain 
$$
\dim \ker \mathbf{D}_{u;N_H}^N\leq 1
$$
Observe that in the case where $H$ is non-cylindrical, then $\dim \ker \mathbf{D}_{u;N_H}^N \geq 1$, since $\partial_a$ is then normal to $H$, and the almost complex structure is $\mathbb{R}$-invariant. We conclude that 
$$
\dim \ker \mathbf{D}_{u;N_H}^N=1,
$$
whenever $H$ is non-cylindrical, and this subspace is spanned by the $\mathbb{R}$-direction. Here we use the fact that $N_u$ needs to coincide with the contact structure at infinity, so that the section $\partial_a$ indeed decays as we approach the punctures, and so defines an element of the domain of $\mathbf{D}_{u;N_H}^N$. Observe that, by (\ref{warmup}), this concludes the existence of obstruction bundles in the case where $H$ is not $\mathbb{R}$-invariant, and $u$ is regular in $H$. For the not-too-bad case, we then observe that 
$$
\ker \mathbf{D}_u^N = \ker \mathbf{D}_{u;T_H}^N \oplus \langle \partial_a \rangle,
$$
and so its dimension is only dependent on the moduli containing $u$. Since the index is constant, then so is the dimension of the cokernel. This finishes the non-cylindrical case.

The only case where automatic transversality is non-conclusive is when $H$ is cylindrical over a hyperbolic critical point, in which case it only yields that $$\dim \ker \mathbf{D}_{u;N_H}^N \in \{0,1\}$$ Geometrically, it is clear that in the cylindrical case there are no holomorphic curves approximating $u$ in the normal directions to $H$, by non-degeneracy of the Morse functions on $\Sigma$, so this operator ``looks'' injective. We show this by a perturbation argument, when $H$ is cylindrical (not necessarily over a hyperbolic point), as follows. 

Observe that for $\epsilon=0$, the Morse function $H_\epsilon$ of Remark \ref{admissible} becomes constant, and, by (\ref{splitJ}), the almost complex structure $J_0^1=J_\pm \oplus j_\pm$ splits on $M_P^\pm$ over $\ker \Lambda_0^1=\xi_\pm \oplus T\Sigma_\pm$. In this case we get a Morse--Bott situation, and $\mathbb{R}\times M_P^\pm$ is foliated by cylindrical holomorphic hypersurfaces $\mathbb{R}\times Y \times \{p\}$ for \emph{any} $p\in M_P^\pm$. If $H$ is a cylindrical hypersurface corresponding to any critical point $p$ of $H_\epsilon$ for $\epsilon\geq 0$, and $u\subseteq H$, we have that $u^*N_H$ is the trivial line bundle over $\dot{S}$ with fiber $T_p\Sigma$. Moreover, we see that the Fredholm operator 
$$
\mathbf{D}_{u;N_H}^N: W^{1,2,\delta_0}(u^*N_H)\rightarrow L^{2,\delta_0}(\overline{Hom}_{\mathbb{C}}(\dot{S},u^*N_H))
$$ converges to the standard operator $\overline{\partial}$ as $\epsilon\rightarrow 0$ (and the domain is independent of $\epsilon$). Indeed, we can consider the extension of $\mathbf{D}_{u;N_H}^N$ to $W^{1,2,\delta_0}(u^*N_H)\oplus C^\infty(u^*N_H)$, which is its Morse--Bott version, including vector fields pointing along the Morse--Bott submanifolds, given by the same formula and same notation. If $x+iy$ is a local holomorphic coordinate for $\Sigma$ around $p$, then $\partial_x$ and $\partial_y$ are constant sections of $C^\infty(u^*N_H)$, and clearly this extended operator satisfies $\mathbf{D}_{u;N_H}^N\partial_x=\mathbf{D}_{u;N_H}^N\partial_y=0$ for $\epsilon=0$, since these vector fields generate holomorphic pushoffs of $u$ along $\Sigma$. The Leibnitz rule then gives the claim for the extended operator, and it follows that the same holds for the original one. Moreover, since the operator $\mathbf{D}_{u;NT}^N$ is linear, this shows that it vanishes.  So we have

\[
\displaystyle
\mathbf{D}_u^N=\mathbf{D}_{u,\epsilon}^N=\left(\begin{array}{cc} 
\mathbf{D}_{u;T_H}^N & \mathbf{D}_{u;NT}^N  \\
0 & \mathbf{D}_{u;N_H}^N 
\end{array}\right)\stackrel{\epsilon\rightarrow 0}{\longrightarrow}
\left(\begin{array}{cc} 
\mathbf{D}_{u;T_H}^N & 0  \\
0 & \overline{\partial} 
\end{array}\right)
\]

The operator $\overline{\partial}$ is injective, since the elements in its kernel are holomorphic sections which decay at the punctures. Since injectivity is an open condition in the usual Fredholm operator topology (it is equivalent to the existence of a left inverse), it follows that $\mathbf{D}_{u;N_H}^N$ is injective for sufficiently small $\epsilon>0$.
Therefore 
$$ 
\ker \mathbf{D}_u^N=\ker \mathbf{D}_{u;T_H}^N,
$$
which depends only on the moduli containing $u$, by assumption. This finishes all cases. 

However, there is a slight technical issue. The small $\epsilon$ needs to be depends a priori on the curve $u$. While this is perhaps just a technicality, as long as we consider curves $u$ with Morse--Bott asymptotics, of which the positive have bounded total action (for a fixed bound), we can assume that this operator is injective for such family of curves. Indeed, there is a finite number of the compactified moduli spaces containing those curves, since there is a finite number of possible asymptotic orbits for $u$. 

Observe that in the case of a regular cylinder inside a cylindrical hypersurface, in which case this operator has index zero, this implies that it is regular in $\mathbb{R}\times M$. This includes the case of a cylinder lying in a cylindrical hypersurface corresponding to a hyperbolic critical point, which is not covered by Lemma \ref{regvs}. The only other case not covered by that lemma is the case of a cylinder with two positive ends over a hyperbolic critical point, and lying in a non-cylindrical hypersurface. This case is not regular, but has an obstruction bundle of rank 1.    

We have proved the following two propositions:

\begin{cor}[cylinders]\label{regcylcyls}
Every cylinder with two positive ends over a hyperbolic critical point, and regular in a non-cylindrical hypersurface, has an obstruction bundle of rank 1. We can take $\epsilon>0$ sufficiently small so that cylinders which are regular inside their corresponding hypersurface are also regular in $\mathbb{R}\times M$, for every other case. For cylinders lying in cylindrical hypersurfaces corresponding to a hyperbolic critical point, we can ensure their regularity as long as we consider fixed action bounds on the positive asymptotics, or finite families of moduli spaces. 
\end{cor} 

More generally,

\begin{prop}\label{exisob} Assume that $u\subseteq H$ is not-too-bad inside $H$, and the rest of the hypothesis of Lemma \ref{regvs} fail (i.e.\ either $g>0$ or inequality (\ref{ATineq2}) fails). If $u$ does not lie in a cylindrical hypersurface corresponding to a hyperbolic critical point, then $u$ is not-too-bad in $\mathbb{R}\times M$. Moreover, for every $T>0$ action threshold, we can choose an $\epsilon>0$ sufficiently small, such that every curve $u$ which lies in such a hypersurface and is not-too-bad inside of it, and such the total action of its positive asymptotics is bounded by $T$, is not-too-bad in $M$.

In all above cases, this implies that there exists an obstruction bundle $\mathcal{O}_u$ for gluing $u$ to any building configuration which contains it.

\end{prop}

\begin{remark} It is not hard to show that, for small $\epsilon>0$, the rank of $\mathcal{O}_u$ is
\begin{equation}\label{rankO}
\mbox{rank}\;\mathcal{O}_u=\mbox{rank}\;\mathcal{O}^H_u -\mbox{ind}\;\mathbf{D}_{u;N_H}^N+\dim \ker \mathbf{D}_{u;N_H}^N, 
\end{equation}
where $\mathcal{O}^H_u$ is the obstruction bundle of $u$ inside $H$, and the second term is given by formula (\ref{indN}). In other words, the linear operator $\mathbf{D}_{u;NT}^N$ does not affect the cokernel. Recall that in the case for non-cylindrical hypersurfaces $H$, we have $\dim \ker \mathbf{D}_{u;N_H}^N=1$; and for cylindrical hypersurfaces, $\dim \ker \mathbf{D}_{u;N_H}^N=0$.

Indeed, from the splitting (\ref{splitop}), a priori we have that 
$$
\mbox{coker}\;\mathbf{D}_u^N=\mbox{coker}\left(\mathbf{D}_{u;T_H}^N\oplus \mathbf{D}_{u;NT}^N\right)\oplus \mbox{coker}\;\mathbf{D}_{u;N_H}^N
$$
We also have a short exact sequence 
$$
0\rightarrow \ker \pi_u \hookrightarrow \mbox{coker}\;\mathbf{D}_{u;T_H}^N \dhxrightarrow{\pi_u} \mbox{coker}\left(\mathbf{D}_{u;T_H}^N\oplus \mathbf{D}_{u;NT}^N\right) \rightarrow 0
$$

Therefore,
$$
\mbox{rank}\;\mathcal{O}_u^H=\dim\ker \pi_u + \dim \mbox{coker}\left(\mathbf{D}_{u;T_H}^N\oplus \mathbf{D}_{u;NT}^N\right)
$$

It follows that 
\begin{equation}\label{eqrk1}
\mbox{rank}\;\mathcal{O}_u = \mbox{rank}\;\mathcal{O}_u^H - \dim \ker \pi_u -\mbox{ind}\;\mathbf{D}_{u;N_H}^N + \dim\ker \mathbf{D}_{u;N_H}^N
\end{equation}
In general, we know from the theory of Fredholm operators that the dimension of the cokernel can jump down (but not up) under small perturbations. This formula illustrates that in the cylindrical case, since clearly (\ref{rankO}) holds for $\epsilon=0$, and how much the dimension jumps down is precisely $\dim \ker \pi_u$. However, we can show that $\ker \pi_u=0$, for small $\epsilon>0$. Indeed, on the other hand, we have
$$
\mbox{rank}\; \mathcal{O}_u= \dim \ker \mathbf{D}_u^N - \mbox{ind}\; \mathbf{D}_u^N
$$
Moreover, we have shown that $\ker \mathbf{D}_u^N = \ker \mathbf{D}_{u;T_H}^N \oplus \ker \mathbf{D}_{u;N_H}^N$, and so
\begin{equation}\label{eqrk2}
\mbox{rank}\; \mathcal{O}_u= \mbox{ind}\; \mathbf{D}_{u;T_H}^N+\mbox{rank}\; \mathcal{O}_u^H+ \dim \ker \mathbf{D}_{u;N_H}^N-\mbox{ind}\; \mathbf{D}_u^N
\end{equation}
Combining (\ref{eqrk1}) with (\ref{eqrk2}), we get
\begin{equation}\label{rkp}
\dim \ker \pi_u = - \mbox{ind}\; \mathbf{D}_{u;T_H}^N-\mbox{ind} \; \mathbf{D}_{u;N_H}^N + \mbox{ind}\;\mathbf{D}_u^N
\end{equation}
If $u=v \circ \varphi$ for $v$ somewhere injective, then 
\begin{equation}\label{11}
\mbox{ind}\; \mathbf{D}_{u;T_H}^N=(n-3)\chi(u)+(\mu_{CZ})^\tau_{W_0}-2Z(d\varphi),
\end{equation}
where $(\mu_{CZ})^\tau_{W_0}$ is the Conley--Zehnder index term inside $W_0$ with respect to $\tau$ (Here, we used that the $c_1$-term vanishes). Also,
\begin{equation}\label{12}
\mbox{ind}\; \mathbf{D}_{u}^N=(n-2)\chi(u)+(\mu_{CZ})_{W_0}^\tau+\#\Gamma^+_2+\#\Gamma^+_0-\#\Gamma^-_0-\#\Gamma^-_2-2Z(d\varphi)
\end{equation}

Putting together equations (\ref{indN}),(\ref{rkp}), (\ref{11}), and (\ref{12}), one readily obtains that $\ker \pi_u=0$. Therefore, $\pi_u$ is an isomorphism, and the claim follows. It also follows that
\begin{equation}\label{cokeq}
\mbox{coker}\;\mathbf{D}_u^N=\mbox{coker}\; \mathbf{D}_{u;T_H}^N \oplus \mbox{coker}\; \mathbf{D}_{u;N_H}^N
\end{equation}
  
\end{remark}

\subsection{Perturbing to non-degenerate contact data}
\label{pert}

In this section, we proceed analogously as in Section \ref{SHStoND}. From now on, we will assume that $(Y,\alpha_\pm)$ are non-degenerate contact manifolds. We can do this by taking Morse perturbations of the original ones of the form $(Y,e^{f_\pm}\alpha_\pm)$ where $f_\pm$ are small positive Morse functions supported around the Morse--Bott submanifolds of $\alpha_\pm$. We can then glue the exact symplectic cobordism 
$$
X_\pm:=\{(r,y): 0\leq \pm r - 1 \leq f_\pm(y)\}
$$
to the ends of $Y \times I$, to get a Liouville domain 
$$
X_-\cup (Y \times I) \cup X_+
$$
which is an exact symplectic cobordism with convex boundary $(Y,e^{f_+}\alpha_+)\bigsqcup (Y,e^{f_-}\alpha_-)$. We may then replace $\alpha_\pm$ with $e^{f_\pm}\alpha_\pm$ in our construction of model B.

\vspace{0.5cm}

Recall that we have constructed a $1$-parameter family of SHS's $\mathcal{H}_\epsilon^s$ which is contact for $s<1$,  so that we will take $s$ sufficiently close to $1$ as our contact data. Every $\mathcal{H}_\epsilon^s$-Reeb orbit with action bounded above by $T_\epsilon$ either lies in $Y\times \mathcal{U}$ or is a closed $R_\pm$-orbit $\gamma^l_p$ corresponding to a critical point of $h_\pm$, and $l\leq N_\epsilon$. Since the former are not necessarily non-degenerate, we will need a further perturbation along $Y\times \mathcal{U}$. For this, we take a generic $C^\infty$-small smooth function $F^s_\epsilon$ with compact support in $Y\times \mathcal{U}$, so that the contact form 
\begin{equation}\label{pertcon}
\widehat{\Lambda_\epsilon^s}:=e^{F_\epsilon^s}\Lambda_\epsilon^s\end{equation}
has, up to period $T_\epsilon$, only non-degenerate Reeb orbits along $Y\times \mathcal{U}$. We will also fix a choice of a $\widehat{\Lambda_\epsilon^s}$-compatible almost complex structure $\widehat{J_\epsilon^s}$ which is $C^\infty$-close to $J_\epsilon^s$.

By Lemma \ref{sympcob}, we have

\begin{lemma} \label{sympcob2} For $0<\epsilon_1\leq \epsilon_2\leq \epsilon_0$, with small $\epsilon_0>0$, the subdomain 
$$
X^s_{\epsilon_1,\epsilon_2}=\{(a,m)\in \mathbb{R}\times M: \epsilon_1 + F_{\epsilon_1}^s(m)\leq a \leq \epsilon_2+F_{\epsilon_2}^s(m)\}
$$
is an exact symplectic cobordism with convex boundary $(M,e^{\epsilon_2}\widehat{\Lambda^s_{\epsilon_2}})$ and concave boundary $(M,e^{\epsilon_1} \widehat{\Lambda_{\epsilon_1}^s})$.
\end{lemma} $\square$

\subsection{Bounding torsion from below in general}

In what follows, for a non-degenerate contact form $\Lambda$ and a constant $T > 0$, we will denote by $\mathcal{A}(\Lambda, T) \subset \mathcal{A}(\Lambda)$ the linear subspace generated by all the monomials of the form $q_{\gamma_1}\dots q_{\gamma_r}$
for which the total action is bounded by $T$, i.e.\
$$\sum_{j=1}^n \int_{\gamma_j}\Lambda < T$$

Since the energy of holomorphic curves contributing to $\mathbf{D}_{SFT}$ is nonnegative and given
by the action difference of the asymptotics, the operator $\mathbf{D}_{SFT}$ restricts to define a differential
$$\mathbf{D}_{SFT} : \mathcal{A}(\Lambda, T)[[\hbar]] \rightarrow \mathcal{A}(\Lambda, T)[[\hbar]]$$

The following generalizes lemma 4.15 in \cite{LW} to our more general situation:

\begin{lemma}\label{countiszero2}
Suppose $N$ is a nonnegative integer and $(M, \xi)$ is a closed contact
manifold admitting a contact form $\Lambda$, compatible almost complex structure $J$ and constant $T > 0$ with the following properties:

\begin{itemize} \item All Reeb orbits of $\Lambda$ with period less than $T$ are nondegenerate.

\item For every pair of integers $g \geq 0$ and $r \geq 1$ with $g + r \leq N + 1$, let $\overline{\mathcal{M}}^1_{g,r}(J;T)$ denote the space of all translation classes of index $1$ connected $J$-holomorphic buildings in $\mathbb{R}\times M$ with arithmetic genus $g$, no negative ends, and $r$ positive ends approaching orbits whose periods add up to less than $T$. Then $\overline{\mathcal{M}}^1_{g,r}(J;T)$ consists of buildings with possibly multiple floors, consisting of either regular or not-too-bad curves.

\item There is a choice of coherent orientations as in \cite{BoMon} for which, if one considers the set of curves obtained by gluing the buildings in $\overline{\mathcal{M}}^1_{g,r}(J;T)$ whose obstruction section has $0$-dimensional zero set, then, after making $J$ generic, their algebraic count is zero whenever $g+r \leq N + 1$. 

\end{itemize}

Then if $\mathbf{D}:\mathcal{A}(\Lambda, T)[[\hbar]] \rightarrow \mathcal{A}(\Lambda, T)[[\hbar]]$ is defined by counting solutions to a
sufficiently small abstract perturbation of the $J$-holomorphic curve equation, there is no $Q \in \mathcal{A}(\Lambda,T)[[\hbar]]$ such that $$\mathbf{D}(Q)=\hbar^N+\mathcal{O}(\hbar^{N+1})$$
\end{lemma}
\begin{proof}
We know that the buildings in $\overline{\mathcal{M}}^1_{g,r}(J;T)$ can be glued to honest regular curves, after making $J$ generic, and moreover they will survive after introducing an abstract perturbation. The number of such gluings is given by the zeroes of the obstruction sections of the associated obstruction bundles, in the case where they can actually be counted. The hypothesis implies that one can choose orientations in such a way that the sections $\mathfrak{s}$ for which the number is non-vanishing, one can find an ``evil twin'' section $\overline{\mathfrak{s}}$ whose count of zeroes cancels that of $\mathfrak{s}$. In the case where the building is an honest regular curve, then we know that this number is precisely one, so that this generalizes the situation in \cite{LW}. The argument now is the same as in that paper: we write $\mathbf{D}_{SFT}=\sum_{k \geq 1}D_k\hbar^{k-1}$ (for the perturbed data). Then for $Q\in \mathcal{A}(\Lambda,T)$, each term $D_k(Q)$ with $k\leq N$ contains a $q$-variable. So, if $Q \in \mathcal{A}(\Lambda,T)[[\hbar]]$ is arbitrary, we have
$$
\mathbf{D}_{SFT}(Q)=P + \mathcal{O}(\hbar^{N+1}),
$$
where $P$ is a polynomial in $\hbar$ of degree at most $N$ whose non-trivial terms contain at least a $q$-variable. This finishes the proof.   
\end{proof}


We include the proof of how to get that there is no $N$-torsion from the fact that there is no truncated $N$-torsion, in the presence of a suitable symplectic cobordism. This is slightly adapted from \cite{LW}.

\begin{lemma}\label{notorsion}
Assume that the hypothesis of Lemma \ref{countiszero2} are satisfied for $T=T_\epsilon$, $\Lambda=\Lambda_\epsilon$, $J=J_\epsilon$, and $N$ for all $\epsilon\leq \epsilon_0$. Assume also that for every $0<\epsilon_1\leq\epsilon_2\leq\epsilon_0$ there is an exact symplectic cobordism having $(M,e^{\epsilon_1} \Lambda_{\epsilon_1})$ as concave boundary, and $(M,e^{\epsilon_2}\Lambda_{\epsilon_2})$ as the convex one, and that $\lim_{\epsilon \rightarrow 0} T_\epsilon=+\infty$. Then $(M,\xi_\epsilon=\ker \Lambda_\epsilon)$ does not have algebraic $N$-torsion, for $\epsilon\leq \epsilon_0$.
\end{lemma}

\begin{proof}
Assume that the contact manifold $(M,\xi_\epsilon)$ has algebraic $N$-torsion, where $\epsilon\leq \epsilon_0$. Make all the necessary choices so that we have an $SFT$ differential $\mathbf{D}_\epsilon$ on $\mathcal{A}(e^{\epsilon}\Lambda_\epsilon)[[\hbar]]$ which computes $H_*^{SFT} (M,\xi_\epsilon)$, and behaves well under symplectic cobordisms. Then we can find some $Q \in \mathcal{A}(e^{\epsilon}\Lambda_{\epsilon})[[\hbar]]$ such that $\mathbf{D}_{\epsilon}(Q)=\hbar^{N}$. Writing $Q=Q_0+\mathcal{O}(\hbar^{N+1})$, we have that 
$$
\mathbf{D}_{\epsilon}(Q_0)=\hbar^N+\mathcal{O}(\hbar^{N+1}),
$$ 
since $\mathbf{D}_\epsilon$ never lowers the degree in $\hbar$. Since $Q_0$ is a polynomial in the $q$-variables, there is some $T>0$ so that $Q_0$ lies in $\mathcal{A}(\Lambda_{\epsilon},T)[[\hbar]]$.

Using that $T_\epsilon \rightarrow +\infty$ as $\epsilon \rightarrow 0$, we can find some $\epsilon_1 \in (0,\epsilon]$ sufficiently small so that $e^{\epsilon_1} T_{\epsilon_1}>T$. Using the exact symplectic cobordism between $e^{\epsilon_1}\Lambda_{\epsilon_1}$ and $e^{\epsilon}\Lambda_{\epsilon}$, we obtain a chain map $$\Phi: \mathcal{A}(e^{\epsilon}\Lambda_{\epsilon},T)[[\hbar]]\rightarrow \mathcal{A}(e^{\epsilon_1}\Lambda_{\epsilon_1},T)[[\hbar]]=\mathcal{A}(\Lambda_{\epsilon_1},e^{-\epsilon_1}T)[[\hbar]]$$ Since $\mathcal{A}(\Lambda_{\epsilon_1},e^{-\epsilon_1}T)[[\hbar]]$ embeds into $\mathcal{A}(\Lambda_{\epsilon_1},T_{\epsilon_1})[[\hbar]],$ we may compose $\Phi$ with this inclusion map to get a new chain map 
$$
\Psi:\mathcal{A}(e^{\epsilon}\Lambda_{\epsilon},T)[[\hbar]] \rightarrow \mathcal{A}(\Lambda_{\epsilon_1},T_{\epsilon_1})[[\hbar]]
$$

But then $$\mathbf{D}_{\epsilon_1} \Psi(Q_0)= \Psi\mathbf{D}_{\epsilon}(Q_0)=\hbar^{N}+\mathcal{O}(\hbar^{N+1}),$$ which contradicts Lemma \ref{countiszero2}.
\end{proof} 

\begin{remark}
Lemmas \ref{sympcob2} and \ref{notorsion} are the reason why model B, instead of model A, is suitable for ruling out values of torsion.
\end{remark}

\section{Giroux torsion implies algebraic 1-torsion in higher dimensions}
\label{GirouxTorsion1-tors}

In this section, we address a conjecture in the paper \cite{MNW} (conjecture 4.14). The philosophy behind the result is that if a contact manifold contains special kinds of Giroux domains as subdomains, then this can be detected algebraically by the SFT machinery. The proof follows by the same sort of recipe that we have developed so far, so that not much more effort is needed; the crucial higher-dimensional input is the uniqueness Theorem \ref{uniqueness}, which we need to slightly adapt.   

\vspace{0.5cm}

Let us begin by introducing the main definitions that we shall use.

\begin{defn}[Giroux]\label{idLD} Let $\Sigma$ be a compact $2n$-manifold with boundary, $\omega$ a symplectic form on the interior $\mathring{\Sigma}$ of $\Sigma$, and $\xi_0$ a contact structure on $\partial \Sigma$. The triple $(\Sigma,\omega,\xi_0)$ is an \emph{ideal Liouville domain} if there exists an auxiliary 1-form $\beta$ on $\mathring{\Sigma}$ such that
\begin{itemize}
\item $d\beta = \omega$ on $\mathring{\Sigma}$.
\item For any smooth function $f:\Sigma \rightarrow [0,\infty)$ with regular level set $\partial\Sigma=f^{-1}(0)$, the 1-form $\beta_0=f\beta$ extends smoothly to $\partial \Sigma$ such that its restriction to $\partial \Sigma$ is a contact form for $\xi_0$.
\end{itemize} 
The 1-form $\beta$ is called a \emph{Liouville form} for $(\Sigma,\omega,\xi_0)$.
\end{defn}

In \cite{MNW}-terminology, we say that an oriented hypersurface $H$ in a contact manifold $(M,\xi)$ is a \emph{$\xi$-round hypersurface} modeled on some closed contact manifold $(Y, \xi_0)$ if it is transverse to $\xi$ and admits an orientation preserving identification with $S^1\times Y$ such that $\xi \cap TH = TS^1 \oplus \xi_0$. 

Given an ideal Liouville domain $(\Sigma,\omega,\xi_0)$, the \emph{Giroux domain} associated to it is the contact manifold $\Sigma \times S^1$ endowed with the contact structure $\xi=\ker(f(\beta + d\theta))$, where $f$ and $\beta$ are as before, and $\theta$ is the $S^1$-coordinate. Away from $V(f)=\partial \Sigma$, the vanishing locus of $f$, this contact structure coincides with the \emph{contactization} $\ker(\beta + d\theta)$. Over $V(f)$ it is just given by $\xi_0 \oplus TS^1$, so that $V(f)\times S^1$ is a $\xi_0$-round hypersurface modeled on $(\partial \Sigma,\xi_0)$.

It is not hard to see that one has a neighbourhood theorem for any $\xi$-round hypersurface $H$ in the boundary of a contact manifold $M$. One may find a collar neighbourhood of the form $[0,1)\times H\subseteq \Sigma \times S^1$, on which $\xi$ is given by the kernel of a contact form $\beta_0+sd\theta$, where $s$ is the coordinate on the interval, where $H \subseteq \partial M$ corresponds to $s=0$, $\theta$ the coordinate in $S^1$, and $\beta_0$ a contact form for $\xi_0$ (\cite{MNW}, lemma 4.1). Using these collar neighbourhoods, one has a well-defined notion of gluing of two Giroux domains along boundary components modeled on isomorphic contact manifolds (see Section \ref{GirouxSOBDs} below). 

We also have a blow-down operation for round hypersurfaces lying in the boundary. If $H$ is a $\xi$-round boundary component of $(M, \xi)$, with orientation opposite
the boundary orientation, consider the collar neighborhood $[0,1)\times H$ as before. Let $\mathbb{D}$ be the disk of radius $\sqrt{\epsilon}$ in $\mathbb{R}^2$. The map $\Psi: (re^{i\theta},y) \mapsto (r^2,\theta, y)$ is a diffeomorphism from $(\mathbb{D} \backslash \{0\}) \times Y$ to $(0, \epsilon) \times S^1 \times Y$ which pulls back $\beta_0 + s dt$ to the contact
form $\beta_0 + r^2d\theta$. Thus we can glue $\mathbb{D} \times Y$ to $M \backslash H$ to get a new contact manifold in which $H$ has been replaced by $Y$, and the $S^1$-component of $H$ has been capped off.

Given a contact embedding of the interior of a Giroux domain $G_\Sigma
:=\Sigma\times S^1$ inside a contact manifold $(M,\xi)$, we shall denote by $\mathcal{O}(G_\Sigma)\subseteq H^2(M;\mathbb{R})$ the annihilator of $H_1(\Sigma;\mathbb{R})\otimes H_1(S^1;\mathbb{R})$, when the latter is viewed as a subspace of $H_2(M;\mathbb{R})$. If $N\subseteq (M,\xi)$ is a subdomain resulting from gluing together a collection of Giroux domains $G_{\Sigma_1},\dots, G_{\Sigma_k}$, we shall denote $\mathcal{O}(N)=\mathcal{O}(G_{\Sigma_1})\cap\dots\cap\mathcal{O}(G_{\Sigma_k})\subseteq H_2(M;\mathbb{R})$.

The following theorem implies Theorem \ref{AlgvsGiroux} (see Example \ref{GTDom} below): 

\begin{thm}\label{Girouxthen1torsion}
If a contact manifold $(M^{2n+1},\xi)$ admits a contact embedding of a subdomain $N$ obtained by gluing two Giroux domains $G_{\Sigma_-}$ and $G_{\Sigma_+}$, such that $\Sigma_+$ has a boundary component not touching $\Sigma_-$, then $(M^{2n+1},\xi)$ has $\Omega$-twisted algebraic 1-torsion, for every $\Omega \in \mathcal{O}(N)$. Moreover, it is also algebraically overtwisted if $N$
contains any blown down boundary components.
\end{thm}

The motivating example of an explicit model of such subdomain is the following:

\begin{example}\label{GTDom}
Let $\alpha_\pm$ be the 1-forms in $\mathbb{R}^{2n-1}$ defined by
$$
\alpha_\pm = \pm e^{\sum_{i=1}^{n-1} t_i} d\theta_0 + e^{-t_1}d\theta_1+\dots + e^{-t_{n-1}}d\theta_{n-1},
$$
with coordinates $(\theta_0,\dots, \theta_{n-1},t_1,\dots, t_{n-1}) \in \mathbb{R}^{2n-1}$. It was proven in \cite{MNW} that there exists manifolds $Y$ which are compact quotients of $\mathbb{R}^{2n+1}$, which are $\mathbb{T}^{n}$-bundles over $\mathbb{T}^{n-1}$, on which $\alpha_\pm$ are also defined, and constitute a Liouville pair. 

Let us then consider $(Y,\alpha_+,\alpha_-)$ a Liouville pair on one of these particular manifolds $Y^{2n-1}$. As in the introduction, consider the \emph{Giroux $2\pi$-torsion domain} modeled on $(Y,\alpha_+,\alpha_-)$, given by the contact manifold $(GT,\xi_{GT}):=(Y\times [0,2\pi]\times S^1,\ker\lambda_{GT})$, where 
\begin{equation}\label{lambdaGT}
\begin{split}
\lambda_{GT}&=\frac{1+\cos(r)}{2}\alpha_++\frac{1-\cos(r)}{2}\alpha_-+\sin (r)d\theta\\
\end{split}
\end{equation} and the coordinates are $(r,\theta)\in [0,2\pi] \times S^1$.

Observe that we may write $\lambda_{GT}=f(\alpha+d\theta),$ where $$f=f(r)=\sin (r),$$$$\alpha=\frac{1}{2}(e^{u(r)}\alpha_++e^{-u(r)}\alpha_-),$$$$u(r)=\log \frac{1+\cos(r)}{\sin (r)},$$ Here, $u:(0,\pi)\rightarrow \mathbb{R}$ is an orientation reversing diffeomorphism, whereas $u:(\pi,2\pi)\rightarrow \mathbb{R}$ is an orientation-preserving one, which is used to pull-back the Liouville form $\frac{1}{2}(e^u\alpha++e^{-u}\alpha_+)$ defined on $Y \times \mathbb{R}$ via a map of the form $v=id\times u$. 

This means that we may view $(GT,\xi_{GT})$ as being obtained by gluing two Giroux domains of the form $GT_-=Y \times [0,\pi] \times S^1$ and $GT_+=Y \times [\pi,2\pi] \times S^1$, along their common boundary, an $\xi_{GT}$-round hypersurface modeled on $(Y,\alpha_-)$. 

\vspace{0.5cm}

Therefore, from Theorem \ref{Girouxthen1torsion} we obtain Theorem \ref{AlgvsGiroux} as a corollary.

\end{example}

\subsection{Giroux SOBDs}
\label{GirouxSOBDs}

We shall be interested in a specially simple kind of SOBDs, which arise on manifolds which have been obtained by gluing a family of Giroux domains along a collection of common boundary components, each a round hypersurface modeled on some contact manifold. Such is the case of the Giroux $2\pi$-torsion domains $GT$. Basically, these SOBDs are obtained by declaring suitable collar neighbourhoods of each gluing hypersurface to be the paper components, whereas the spine components are the complement of these neighbourhoods. They have the desirable features that the fibrations are trivial, and that they have 2-dimensional pages (see Appendix \ref{SOBDap} for definitions). 

\vspace{0.5cm}  

Let $$(G_\pm=\Sigma_\pm\times S^1, \xi_\pm=\ker(f_\pm(\beta_\pm+d\theta)))$$ be two Giroux domains which one wishes to glue along a round-hypersurface $H=Y\times S^1$ modeled on some contact manifold $(Y,\xi_0)$ and lying in their common boundaries. Fix choices of collar neighbourhoods $\mathcal{N}_\pm(H)$ of $H$ inside $G_\pm$, of the form $Y \times S^1 \times [0,1]$. Take coordinates $s_\pm \in [0,1]$ so that $H=\{s_\pm=0\}$, and $\theta \in S^1$, such that the contact structures $\xi_\pm$ are given by the kernel of the contact forms $\beta_0+s_\pm d\theta$. Here, $\beta_0$ is a contact form for $\xi_0$. This means that, in these coordinates, one can take $\beta_\pm=\beta_0/s_\pm$ and $f=f(s_\pm)=s_\pm$, and the corresponding \emph{ideal} Liouville vector fields are $V_\pm=-s_\pm \partial_{s_\pm}$. We glue $\mathcal{N}_\pm(H)$ together in the natural way, by taking a coordinate $s\in [-1,1]$, so that $s=-s_-$ and $s=s_+$, and $s=0$ corresponds to $H$. We thus obtain a collar neighbourhood $\mathcal{N}(H)=\mathcal{N}_+(H) \bigcup_\Phi \mathcal{N}_-(H)\simeq H\times [-1,1]$, where we denote by $\Phi$ the resulting gluing map. Doing this for each of the boundary components that we glue together, we obtain a decomposition for $$M:=G_+\bigcup_{\Phi} G_-=\left(\Sigma_+\bigcup_{\Phi} \Sigma_-\right)\times S^1,$$ given by $$M=M_P \cup M_\Sigma,$$ where $M_P$ is the disjoint union of the collar neighbourhoods of the form $\mathcal{N}(H)$ for each of the gluing round-hypersurfaces $H$, and $M_\Sigma$ is the closure of its complement in $M$. We have a natural fibration structure $\pi_P$ on $M_P$, which is the trivial fibration over the disjoint union of the hypersurfaces $H$, so that the pages are identified with the annuli $P:=S^1\times [-1,1]$. We also have an $S^1$-fibration $\pi_\Sigma$ on $M_\Sigma$, which is also trivial. It has as base the disjoint union of $\Sigma_+$ and $\Sigma_-$ minus the collar neighbourhoods, which we denote by $\Sigma$. Let us assign a \emph{sign} $\epsilon(G_\pm):=\pm 1$ to each of the Giroux subdomains $G_\pm \subseteq M$.

\vspace{0.5cm}

We shall refer to the SOBDs obtained by the procedure described above as \emph{Giroux} SOBDs, where are also allowing the number of \emph{Giroux subdomains} involved (which is the same as the number of spine components) to be arbitrary. 

\vspace{0.5cm}

We will also fix collar neighbourhoods of the components of $\partial M$, which correspond to the non-glued hypersurfaces, in the very same way as we did before for the glued ones. Their components look like $\mathcal{N}_\Sigma(H):=H\times [0,1] = Y\times S^1\times [0,1]$, for some non-glued hypersurface $H$. There is a coordinate $s\in [0,1]$, such that $\partial M=\{s=0\}$, and such that the contact structure on the corresponding Giroux domain is given by $\ker(\beta_0+s d\theta)$, for some contact form $\beta_0$ for $(Y,\xi_0)$. Denote by $\mathcal{N}_\Sigma$ the disjoint union, over all unglued $H$'s, of all of the $\mathcal{N}_\Sigma(H)$'s. 

\vspace{0.5cm}

We may glue together the contact structures $\xi_\pm$, and present the resulting contact structure as the kernel of a glued 1-form, as follows. If we denote by $\lambda_\pm= f_\pm(\beta_\pm + d\theta)$, we can extend the expression $\lambda_+=\beta_0+sd\theta$, a priori valid on $\mathcal{N}_+(H)$, to $\mathcal{N}(H)$, by the same formula. Observe that, by choice of our coordinates, the resulting 1-form glues smoothly to $\overline{\lambda}_-:=f_-(\beta_--d\theta)$. Therefore, we may still think of the contact structure $\overline{\xi}_-=\ker \overline{\lambda}_-$ as a \emph{contactization} contact structure over the region $(\Sigma_-\times S^1)\backslash\mathcal{N}_-(H)$, with the caveat that we need to switch the orientation in the $S^1$-direction. Alternatively, we may think of it as obtained from $\lambda_-$ by replacing $f_-$ by $-f_-$ and $\beta_-$ by $-\beta_-$.

\begin{figure}[t]
\centering
\includegraphics[width=0.68\linewidth]{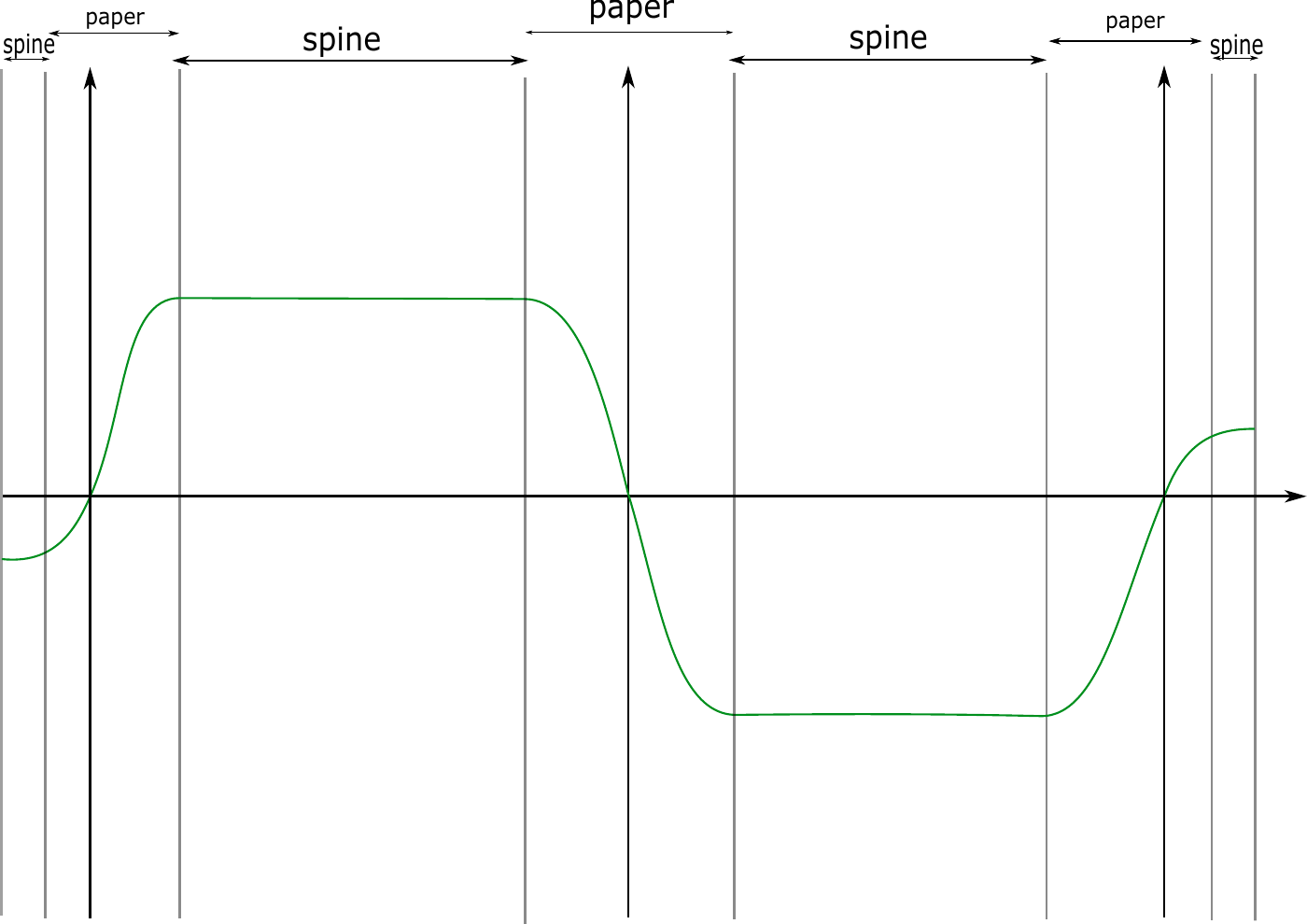}
\caption{\label{fGiroux} The isotoped function $f$ in the extended Giroux torsion domain $GT^\epsilon$.}
\end{figure}

Denote the resulting contact form by $\lambda:=\lambda_+\cup_\Phi \overline{\lambda}_-$. We may globally write it as $\lambda=f(\beta+d\theta).$ Here, $f$ is a function which is either strictly positive or strictly negative over the interior of the Giroux domains, and vanishes precisely at the (glued and unglued) boundary components. The $1$-form $\beta$ coincides with the Liouville forms $\pm \beta_\pm$ where these are defined, and is undefined along said boundary components. From this construction, in $M$, each of the subdomains $G\subset M$ used in the gluing procedure carries a \emph{sign} $\epsilon(G)$, which we define as the sign of the function $f\vert_{\mathring{G}}$. Observe that $f$ coincides with $\epsilon(G)$ along every $H \times \{s_\pm=1\}$, the inner boundary components of the collar neighbourhoods. Therefore, we may isotope it relative every collar neighbourhood to a smooth function which is constant equal to $\epsilon(G)$ along $G\backslash \left((M^\delta_P \cap G) \cup \mathcal{N}^\delta_\partial\right)$ (see Figure \ref{fGiroux}). Here, $M^\delta_P$ and $\mathcal{N}^\delta_\partial$ denote small $\delta$-extensions in the interval direction of $M_P$ and $\mathcal{N}_\Sigma$, respectively, so that now $|s| \in [0,1+\delta]$. The regions $M_P^\delta\backslash M_P$ play the role of the smoothened corners, as in our model A construction. Observe that isotoping $f$ as we just did does not change the isotopy class of $\lambda$, by Gray stability (version with boundary), and has the effect of transforming the Reeb vector field of $\lambda$ into $\epsilon(G)\partial_\theta$ along $M_Y$. Observe also that along the paper, we have $\lambda=\beta_0+\gamma$, where $\gamma=sd\theta$ is a Liouville form for $S^1\times [-1,1]$. Since this manifold is trivially a cylindrical Liouville semi-filling, we may view $\lambda$ as a model A contact form. Compare this with (\ref{modelA}): the difference lies in the parameters and the Morse/Morse--Bott function $H$, which do not change the isotopy class of $\lambda$. 

We can also blow-down boundary components in the Giroux subdomains, and in this case we may extend our SOBD by declaring the $\mathbb{D}\times Y$ we glued in the blow-down operation to be part of the paper $M_P$ (so that its pages are disks), and the contact form $\lambda$ also extends in a natural way. We will include in the definition of Giroux SOBD the assumption that any blown-down component is part of the paper.  

It will be convenient to consider an $\epsilon$-extension of our SOBD, which we call $M^\epsilon$, by gluing small collars to the boundary. To each component $\mathcal{N}_\Sigma(H)$ for $H$ lying in $\partial G$ for a subdomain $G$ of $M$, we glue a collar neighbourhood of the form $H \times [-\epsilon,0]$ for some small $\epsilon>0$. We extend our function $f$ so that the extended version coincides with $\epsilon(G)id$ near $s=0$, $f^\prime(-\epsilon)=0$, and $sgn(f^\prime(s))=\epsilon(G)\neq 0$ for $s>-\epsilon$ (see Figure \ref{fGiroux}). We will define the paper $M^\epsilon_P$ to be the union of $M_P$ and the region $H \times [-\epsilon/2,1+\delta]$, and the spine $M_\Sigma^\epsilon$ to be the union of $M_\Sigma$ with $H \times [-\epsilon,-\epsilon/2]$. The point of this extension is that now the Reeb vector field of the extended $\lambda$ coincides with $-\epsilon(G)\partial_\theta$ along the boundary. 

\vspace{0.5cm}

\begin{remark}
Because of the above discussion on orientations, from which we gathered that the $S^1$-orientation that we need depends on the sign of the Giroux domain, we will not allow our manifold $M$ to have ``Giroux cycles'' of odd order. That is, we want to rule out the case where we have a sequence $(G_0,\dots,G_{N-1})$ of Giroux subdomains of $M$, where $G_i$ has been glued to $G_{i+1}$ (modulo $N$) along some collection of boundary components, and $N>1$ is odd. This condition is to be taken as part of the definition of a Giroux SOBD.

\end{remark}

\subsection{Prequantization SOBDs (non-trivial case)} \label{PrequantSOBDs}

We can generalize the previous construction to the case where the $S^1$-bundles are not necessarily globally trivial, but trivial on the boundary components which are glued together.

\begin{defn}\cite{DiGe}\label{IdealLDomDIGE} Let $\Sigma$ be a compact $2n$-manifold with boundary, $\omega_{int}$ a symplectic form on $\mathring{\Sigma}$, $\xi_0$ a contact structure on $\partial \Sigma$, and $c$ a cohomology class in $H^2
(\Sigma; \mathbb{Z})$. The tuple $(\Sigma, \omega_{int}, \xi_0, c)$ is called an \emph{ideal Liouville
domain} if for some (and hence any) closed 2-form $\omega$ on $\Sigma$ with $-[\omega/2\pi] = c \otimes \mathbb{R} \in H^2(\Sigma;\mathbb{R})$ there exists a 1-form $\beta$ on $\mathring{\Sigma}$ such that 
\begin{itemize}
\item $d\beta = \omega_{int} - \omega$ on $\mathring{\Sigma}$.
\item For any smooth function $f:\Sigma \rightarrow [0,\infty)$ with regular level set $\partial\Sigma=f^{-1}(0)$, the 1-form $\beta_0=f\beta$ extends smoothly to $\partial \Sigma$ such that its restriction to $\partial \Sigma$ is a contact form for $\xi_0$.
\end{itemize} 
\end{defn}
 
An ideal Liouville domain $(\Sigma, \omega_{int}, \xi_0, 0)$ is an ideal Liouville domain in the sense of Giroux, where we may take $\omega = 0$. In the case where $c\vert_{\partial \Sigma}=0$, so that we may take $\omega=0$ in any collar neighbourhood  $[0,1)\times \partial \Sigma$ of the boundary, we will call $(\Sigma, \omega_{int}, \xi_0, c)$ an \emph{ideal strong symplectic filling}.

\vspace{0.5cm}

Now let $\pi: M\rightarrow \Sigma$ be the principal $S^1
$-bundle over $\Sigma$ of (integral) Euler class $c$.
Choose a connection 1-form $\psi$ with curvature form $\omega$ on this bundle. As in the trivial case $c=0$ (where we may take $\psi=d\theta$), the 1-form $f(\psi + \beta)=f\psi + \beta_0$ defines a contact structure $\xi$ on $M$, and we call $(M, \xi)$ the \emph{contactization} of $(\Sigma, \omega_{int}, \xi_0, c)$. Observe that $d(\psi + \beta)=\pi^*(\omega + \omega_{int}-\omega)=\pi^*\omega_{int}$, where we identify $\beta$ with $\pi^*\beta$, so that on $\mathring{\Sigma}$, $\xi$ may be regarded as the prequantization contact structure corresponding to $\omega_{int}$. 

\begin{defn}\cite{DiGe} Let $B$ be a closed, oriented manifold of dimension $2n$. An ideal Liouville splitting of class $c \in H^2(B; \mathbb{Z})$ is a decomposition $B = B_+ \bigcup_\Gamma B_-$ along a two-sided (but not necessarily connected) hypersurface $\Gamma$, oriented as the boundary of $B_+$, together with a contact structure $\xi_0$ on $\Gamma$ and symplectic forms $\omega_{\pm}$ on $\pm\mathring{B}_\pm$, such that $(\pm B_\pm, \omega_\pm, \xi_0, \pm c\vert_{B_\pm})$
are ideal Liouville domains.
\end{defn}

Proposition 3.6 in \cite{DiGe} tells us that an $S^1$-invariant
contact structure $\xi = \ker(\beta_0 + f\psi)$ on the principal $S^1$-bundle $\pi: M \rightarrow B$ defined by $c \in H^2(B; \mathbb{Z})$ leads to an ideal Liouville splitting of $B$ of class $c$, along the dividing set $$\Gamma=\{b \in B: \partial_\theta\vert_b \in \xi\}=f^{-1}(0)$$ 

If we require that $c \in H^2(B,\Gamma;\mathbb{Z})$, so that $c\vert_{\Gamma} =0$, then the $S^1$-principal bundle $\pi$ is trivial in any collar neighbourhoods $[-1,1]\times \Gamma$ of $\Gamma$, along which we may take $\psi=d\theta$, $\omega=0$. In this situation, we may regard the total space $M$ as carrying what we will call a \emph{prequantization SOBD}, which is obtained in the same way as we defined a Giroux SOBD. That is, we declare suitable collar neighbourhoods of $\Gamma$ to be paper components, which we can do by the triviality assumption on $c$ along $\Gamma$. The only difference now is that the spine is no longer globally a trivial $S^1$-bundle. Reciprocally, one can glue prequantization spaces over ideal strong symplectic fillings to obtain a manifold $M$ (possibly with non-empty boundary) carrying an $S^1$-invariant contact structure $\xi=\ker(\beta_0+f\psi)$, and a prequantization SOBD. The gluing construction is the same as for Giroux SOBDs, and is a particular case of the gluing construction of Thm.\ 4.1 in \cite{DiGe}.  

Observe that, using the coordinates of the previous section, the contact $1$-form $\alpha=\beta_0+f\psi$ defining $\xi$ satisfies that $d\alpha=ds \wedge d\theta>0$ is positive along the pages $S^1\times [-1,1]$, and, after isotoping $f$ in the same way as for Giroux SOBDs, its Reeb vector field is tangent to the $S^1$-fibers of $\pi_\Sigma$ along $M_\Sigma$. Moreover, along the boundary, it induces the distribution $\ker \alpha \cap T(\partial M)=\xi_0 \oplus TS^1$, with characteristic foliation given by the $S^1$-direciton, and its Reeb vector field coincides with $R_0$, the Reeb vector field of $\beta_0$. From this, and having the definition in \cite{LVHMW}, and the standard one due to Giroux, both in mind, one can define

\begin{defn}
A \emph{Giroux form} for a prequantization SOBD (in particular, for a Giroux SOBD) is any contact form inducing a contact structure which is isotopic to the $S^1$-invariant contact structure $\xi=\ker(\beta_0+f\psi)$.
\end{defn} 

With this definition, the contact structure $\lambda=\lambda_+\cup_{\Phi}\overline{\lambda}_-$ (as well as the model A form $\Lambda$ of the proof below) and $\lambda_{GT}$ defined in (\ref{lambdaGT}) are Giroux forms for the Giroux SOBD we have defined on $GT$ via the previous section. One can similarly define a Giroux form for our model $M=Y\times \Sigma$ of the main sections to be any contact form inducing a contact structure isotopic to the contact structure induced by any model A or model B contact form. For prequantization SOBDs, observe that any $S^1$-invariant contact structure in a principal bundle can always be written as $\alpha= \beta_0 + f\psi$ for a horizontal $f$ and a connection form $\psi$. Therefore, we have the following tautology, which we spell out: 

\vspace{0.5cm}

\emph{Any $S^1$-invariant contact form $\alpha$ inducing an $S^1$-invariant contact structure in a principal bundle defined by a cohomology class $c$ satisfying $c\vert_{\Gamma}=0$, where $\Gamma$ is a set of dividing hypersurfaces for $\alpha$, is a Giroux form for the induced prequantization SOBD}.

\subsection{Proof of Theorem \ref{Girouxthen1torsion}}

The proof of this theorem is simply a reinterpretation of what we did in the construction of model A. 

\begin{proof}[Thm. (\ref{Girouxthen1torsion})]

Let $N$ be a subdomain as in the hypothesis, carrying a Giroux form $\lambda=f(\beta+d\theta)$ obtained by gluing. Since the contact embedding condition is open, and we are assuming that there are boundary components of $\Sigma_+$ not touching $\Sigma_-$, we can find a small $\epsilon>0$ such that $M$ admits a contact embedding of the $\epsilon$-extension $N^\epsilon$. Endow this extension with a Giroux SOBD $N^\epsilon=M^\epsilon_P \cup M^\epsilon_\Sigma$ as in the previous section. On this decomposition, add corners where spine and paper glue together as we explained before. Also choose a small Morse/Morse--Bott function $H$ in $M_\Sigma^\epsilon$, which lies in the isotopy class of $f$, and vanishes as we get close to the paper. With this data, we then may construct a Morse/Morse--Bott model A contact form $\Lambda$ on $N$ which lies in the isotopy class of $\lambda$, along with a SHS deforming it. 

Since in our situation we are allowing $\Sigma_\pm$ to be more general than a semi-filling $Y \times I$, we need to specify what we mean by the Morse--Bott situation. We will take our Morse--Bott function $H$ so that it depends only on the interval parameter along the collar neighbourhoods $\mathcal{N}_\Sigma$ close to the boundary and along a slightly bigger copy of $M_P^\epsilon$, and matches the function $f$ close to $\partial N$. In particular, $\partial N$ is a Morse--Bott submanifold. We also impose that $H$ is Morse in the interior of the components of $M_\Sigma$ which are away from the boundary. For simplicity, we will assume that, besides the boundary, $H$ only has exactly one Morse--Bott submanifold close to each boundary component of $M_\Sigma$ which is glued to a paper component, of the form $Y\times \{t\}$ for some $t$ (see Figure \ref{Girouxd}). The Morse situation is then obtained by a perturbation of this situation obtained by choosing Morse functions along the Morse--Bott submanifolds. Observe that we may always choose the interior Morse-Bott submanifolds so that they lie in the cylindrical ends of the Liouville domains $\Sigma_\pm$.

\begin{figure}[t]
\centering
\includegraphics[width=0.80\linewidth]{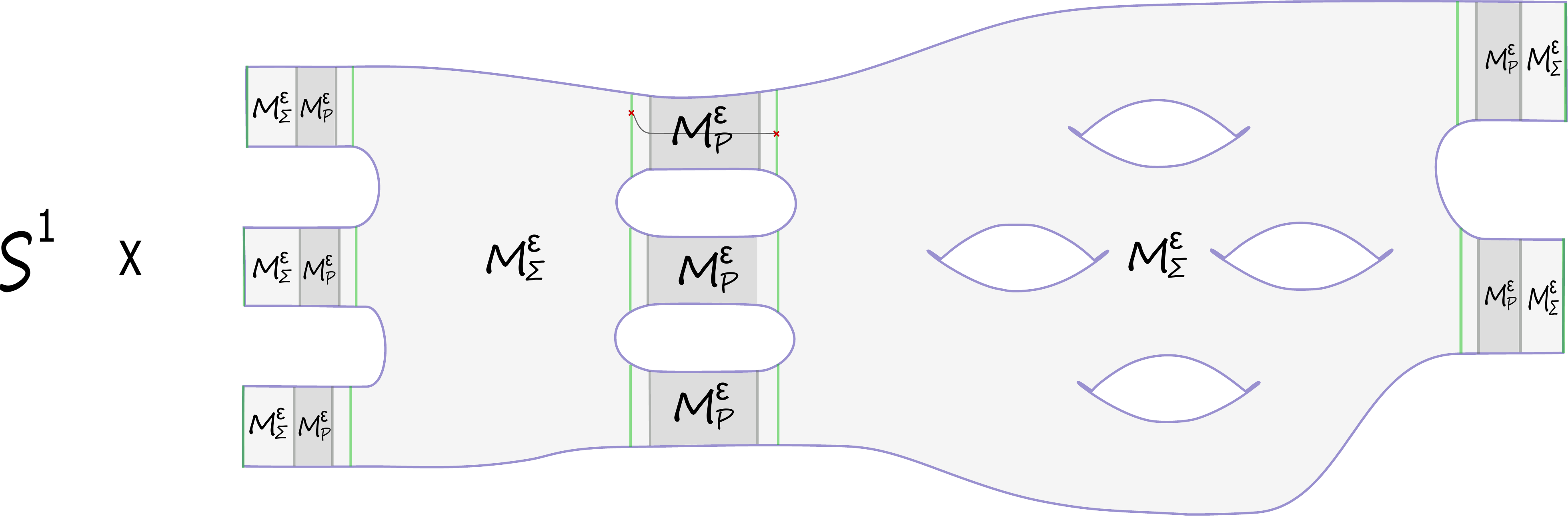}
\caption{\label{Girouxd} Sketch of an example of (the $\epsilon$-extension of) a domain $N$ consisting of two Giroux domains glued along three boundary components. The paper components of $M_P$ all are of the form $Y \times I\times S^1$ for some contact manifold $Y$. We draw the Morse--Bott submanifolds of $H$ in green, inside the spine but close to the paper. These are a ``barrier'' for page-like holomorphic curves, whose role is to prevent these from venturing inside the interior of the spine components, so that all have different asymptotics (this is what we need to adapt Theorem \ref{uniqueness}). We also draw an index $1$ holomorphic cylinder (of the first kind) connecting a Reeb orbit corresponding to a maximum to one corresponding to an index 1 critical point (in the Morse case), which are the ones we will count.}
\end{figure} 

We have a corresponding almost complex structure $J^\prime$ compatible with the SHS, for which we get a stable finite energy foliation by holomorphic cylinders, which come in two types: either they are obtained by gluing constant lifts of the cylindrical pages and flow-line cylinders over $M^\epsilon_\Sigma$ along the corners (first kind); or they are flow-line cylinders completely contained in $M_\Sigma^\epsilon$ (second kind). Both kinds have as asymptotics simply covered Reeb orbits $\gamma_p$ corresponding to critical points $p$ of $H$, either along the boundary, or at interior points of $\Sigma$ (see Figure \ref{foliationGirouxTorsion}). These cylinders have two positive ends if its corresponding critical points lie in Giroux subdomains with different sign, and one positive and one negative end if these signs agree. In this particular case, where there are only two Giroux subdomains involved, this happens only for cylinders of the second kind. The Fredholm index of any of them is given by the Conley--Zehnder index term, since the $c_1$-term in the index formula vanishes. The Conley--Zehnder index of $\gamma_p$, with respect to the trivialization $\tau$ of the contact structure along $\gamma_p$ coming from $S^1$-invariance, is given by $\mu_{CZ}^\tau(\gamma_p)=\mbox{ind}_p(H)-n$. 

We can adapt the proof of Fredholm regularity in the Morse-Bott model A case to this particular situation, to obtain regularity for cylinders of both kinds. There is no change in the proof for cylinders of the second kind, which corresponds to our Floer theory approach of Section \ref{regcyls}. For regularity for cylinders of the first kind, we adapt the proof of Section \ref{reggenuszero}. For this, we observe that our choices of interior Morse-Bott submanifolds imply that the cylindrical ends of the cylinders of the first kind lie completely over cylindrical ends of the Liouville domains $\Sigma_\pm$ (see Figure \ref{Girouxd}). This means that we may assume that the obvious analogous version of Assumptions \ref{Assumptions} hold, and the rest of the proof is a routine adaptation. Fredholm regularity in the sufficiently nearby Morse case follows from the implicit function theorem.   

\begin{figure}[t]
\centering
\includegraphics[width=0.70\linewidth]{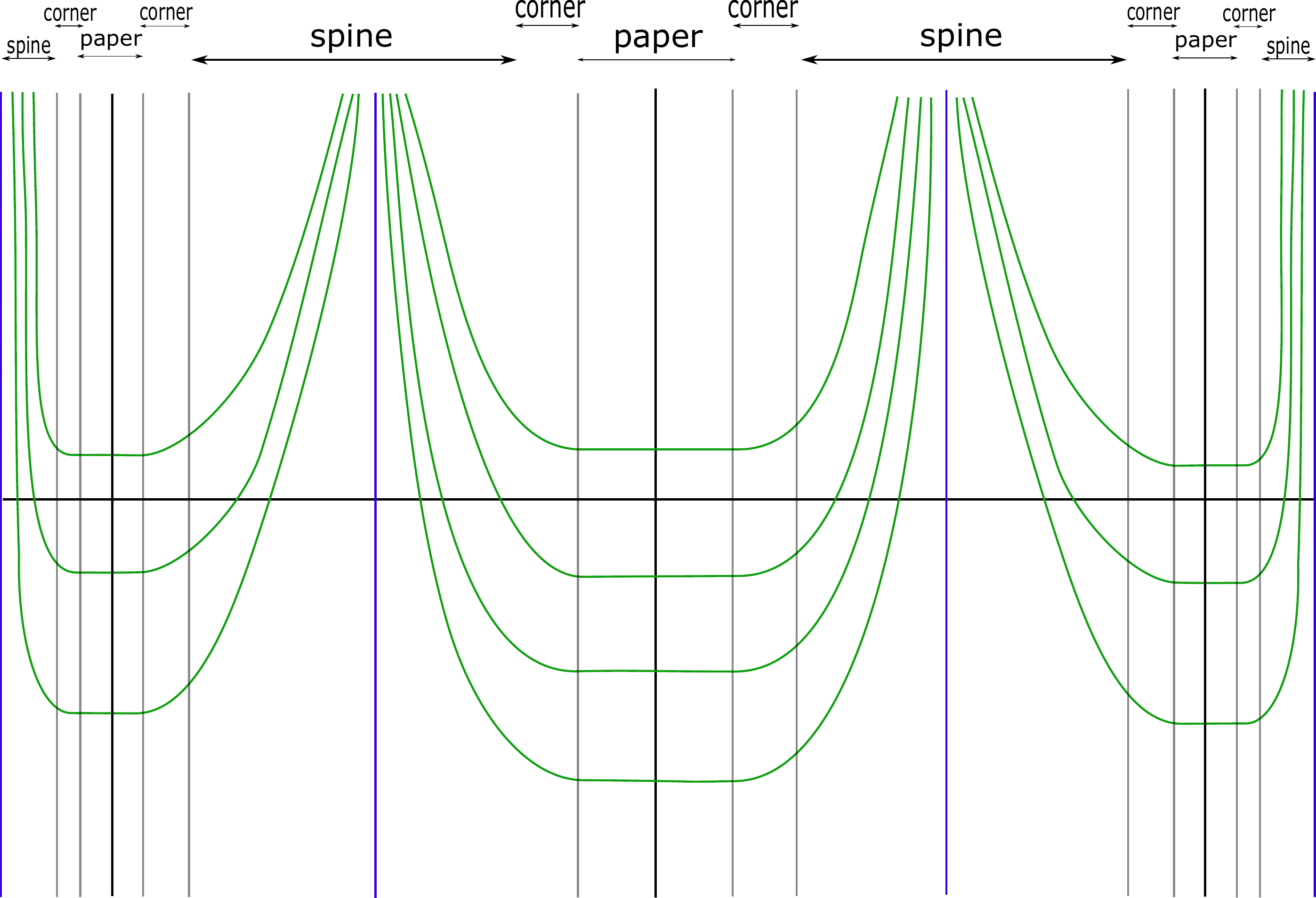}
\caption{\label{foliationGirouxTorsion} The foliation by holomorphic cylinders (with respect to SHS data) of the symplectization of $GT_\epsilon$.}
\end{figure}

We also have a version of the uniqueness Theorem \ref{uniqueness}, adapted to this situation. We have a few small changes, as follows. In the hypothesis we require that the positive asymptotics are simply covered; correspond to critical points all of which, in the Morse--Bott case, belong to a Morse--Bott submanifold of the type already described; and any two lie in distinct components of $M_\Sigma$, both of which are away from the boundary, separated by a paper component. Since the number of these components of $M_\Sigma$ is two, we have that the number of positive punctures is bounded by 2, which incidentally agrees with the number of boundary components of each paper component. We may then replace $k$ by $2$ in the proof. We have more than just one Morse--Bott submanifold now, but, for case B where $u$ lies over some paper component $P=Y \times I \times S^1$, we can still get a similar upper bound on the energy:

\begin{equation}
\begin{split}
\mathbf{E}(u)\leq & 2\pi\left(\sum_{z \in \Gamma^+}e^{\epsilon H(p_z)}- \sum_{z \in \Gamma^-}T_z\right)\\
\leq &2\pi \sum_{z \in \Gamma^+}e^{\epsilon H(p_z)} \\
\leq &4\pi ||e^{\epsilon H}||_{C^0},
\end{split}
\end{equation}
where $T_z>0$ is the action of the Reeb orbit corresponding to $z \in \Gamma^-$, which a priori might even lie in $M\backslash N$ (but not a posteriori). Here, the energy is computed with respect to a $2$-form which is the derivative of a contact form for $\xi$ outside of $N$, and which restricts to the $2$-form of the fixed model A SHS on $N$. Observe that we have chosen the model A SHS so that it is contact near $\partial N$, so we can take such a $2$-form. 

Using that the number of boundary components of each paper component is 2, the very same arguments we used before give a lower bound on the energy of the form
$$
\mathbf{E}(u)\geq 4 \pi e^{-\delta}\mbox{deg}_P(u),
$$
for some small $\delta>0$. Here, we denote by $\mbox{deg}_P(u)$ the covering degree of the restriction to $u$ of the projection to $I \times S^1$, over the component $P$. The rest of the proof is the same. Observe that since we are assuming that the critical points appearing in the positive asymptotics are away from the boundary of $N$, and the only way of venturing into $M\backslash N$ is to escape through the latter, any holomorphic curve with these positive asymptotics which leaves $N$ would necessarily need to go through a paper component, so is dealt with by case B in the proof of \ref{uniqueness} (which deliberately does not need to assume that the whole curve stays over $N$).

\vspace{0.5cm}

Make all the necessary choices to have an SFT differential 
$$
\mathbf{D}_{SFT}:\mathcal{A}(\Lambda)[[\hbar]]\rightarrow \mathcal{A}(\Lambda)[[\hbar]],
$$
which computes $H_*^{SFT}(N,\ker \Lambda;\Omega\vert_N)$, where $\Omega \in \mathcal{O}(N)$. These includes a choice of coherent orientations and spanning surfaces (which will choose more suitably below), and an abstract perturbation, for which we still have a perturbed foliation by holomorphic cylinders. Extend these choices to $M$ so as to be able to compute $H_*^{SFT}(M,\xi;\Omega)$.

Let $e$ be a maximum (index $2n$) of $H$ in $M_\Sigma \cap (\Sigma_-\times S^1)$, and let $h$ be an index $1$ critical point of $H$ in $M_\Sigma \cap (\Sigma_+\times S^1)$. We take both to lie in the Morse--Bott manifolds of the Morse--Bott case, before the Morse perturbation. Denote by $q_e$, $q_h$ the corresponding SFT generators. Define 
$$
Q=q_hq_e,
$$
an element of $\mathcal{A}(\Lambda)$. There is a unique (perturbed) cylinder $u$ of the first kind which has $\gamma_e$ and $\gamma_h$ as positive asymptotics, and its index is given by 
$$
\mbox{ind}(u)=\mu_{CZ}^\tau(u)=\mbox{ind}_e(H)-n+\mbox{ind}_h(H)-n=1
$$
If we choose $H$ so that it does not have any minimum (which we can since $\Sigma$ has non-empty boundary), uniqueness tells us that any other holomorphic curve over $M$ which may contribute to the differential of $Q$ is a flow-line cylinder completely contained in $M_\Sigma$, connecting an index $2n-1$ critical point $p$ with $e$. 

Denote by $\overline{d}$ the element in $H_2(M;\mathbb{R})/\ker \Omega$ defined by any $d \in H_2(M;\mathbb{R})$, by $\mathcal{M}(H;p,e)$ the space of positive flow-lines connecting $p$ with $e$, $u_\gamma$ the flow-line cylinder corresponding to a flow-line $\gamma$, and $\epsilon(u)$ the sign of the holomorphic curve $u$ assigned by a choice of coherent orientations. As before, elements of $\mathcal{M}(H;p,e)$ come in evil twins pairs $\gamma \leftrightarrow \overline{\gamma}$, and we have that $\epsilon(u_\gamma)=-\epsilon(u_{\overline{\gamma}})$. This is because $H_{2n}(\Sigma_\pm)=0$, being $2n$-manifolds with non-empty boundary, and we may take coherent choices compatible with orientations of Morse flow-lines moduli spaces. 
As explained in Section \ref{Torsion}, once can choose suitable spanning surfaces such that $\overline{[u_\gamma]}=\overline{[u_{\overline{\gamma}}]}$. 

Then, we obtain that 
\begin{equation}
\begin{split}
\mathbf{D}_{SFT}Q=&z^{\overline{[u]}}\hbar+\sum_{\substack{\gamma \in \mathcal{M}(H;p,e)\\ \mbox{ind}_p(H)=2n-1}} \epsilon(u_\gamma)z^{\overline{[u_\gamma]}}q_pq_h\\
=& z^{\overline{[u]}}\hbar+\frac{1}{2}\sum_{\substack{\gamma \in \mathcal{M}(H;p,e)\\ \mbox{ind}_p(H)=2n-1}} \left(\epsilon(u_\gamma)+\epsilon(u_{\overline{\gamma}})\right) z^{\overline{[u_\gamma]}}q_pq_h\\
=&z^{\overline{[u]}}\hbar,\\
\end{split}
\end{equation}  
which proves that $(N,\ker \Lambda)$ has $\Omega\vert_N$-twisted 1-torsion. This implies that $(M,\xi)$ has $\Omega$-twisted 1-torsion, since our uniqueness theorem gives that there are no holomorphic curves with asymptotics in $N$ which venture into $M\backslash N$, and therefore $H_*^{SFT}(N,\ker \Lambda; \Omega\vert_N)$ embeds into $H_*^{SFT}(M,\xi; \Omega)$. Here, we use that $\ker \Lambda$ is isotopic to $\ker \lambda=\xi\vert_N$.

\vspace{0.5cm}

For the second statement, assume that we have a blown-down boundary component, so that the corresponding $\mathbb{D}\times Y$ is a paper component with disk pages. As in model A, we have that the disk pages lift as finite-energy holomorphic planes with a single positive asymptotic $\gamma_p$, corresponding to a critical point $p$ in $M_\Sigma$. If $u$ is such a plane, its index is 
$$
\mbox{ind}(u)=(n-2)\chi(u)+2\chi(u)+ind_p(H)-n=ind_p(H)
$$
Take $p$ so that $ind_p(H)=1$, and let $P=q_{\gamma_p}$. Since there is no minima for $H$, by our uniqueness theorem we know that
$$
\mathbf{D}_{SFT}\left(z^{-\overline{[u]}}P\right)=1
$$  
Since the choice of coefficients is arbitrary, $(N,\ker \Lambda)$ is algebraically overtwisted, which implies that $(M,\xi)$ is.

\end{proof}

In view of the definition of a prequantization SOBD, and the fact that we may still construct a foliation by flow-line holomorphic cylinders over a non-trivial prequantization space (Section \ref{prequantcurves} below), we have the following:

\begin{conjeorem}
Consider a principal $S^1$-bundle obtained by gluing a collection of prequantization spaces over ideal strong symplectic fillings, such that there are two of them, one with a boundary component not touching the other after the gluing. Then the resulting $S^1$-invariant contact structure has (untwisted) algebraic 1-torsion. 
\end{conjeorem}

The only technical ingredient that is missing is the regularity of the holomorphic flow-line cylinders for sufficiently small $H$, in dimension at least $5$. In dimension $3$ one can use automatic transversality, and the conjeorem is actually a theorem. One may choose $H$ as in the proof above, and so page-like curves are regular by our proof for model A. However, the Floer theory approach for regularity is not readily generalizable to the non-trivial case---basically because of curvature---so that further work/a different approach is needed.

\section{Non-fillable 5-dimensional model with no Giroux torsion}\label{5dmodel}

For this section, we fix $Y=ST^*X$ to be the unit cotangent bundle of a hyperbolic surface $X$ with respect to a choice of hyperbolic metric, and $\pi: Y \rightarrow X$ the natural projection. The 1-form $\alpha_-$ is the standard Liouville form, and $\alpha_+$ is a prequantization contact form. This means that $\alpha_+$ is a connection form with curvature $d\alpha_+=\pi^* \omega$, where $\omega$ is a symplectic form on $X$ representing $c_1(Y)=-2+2g(X)>0$ in $H^2(X;\mathbb{R})\simeq\mathbb{R}$ when $Y$ is viewed as a complex line bundle, or, in other words, area$(\omega)=\int_X \omega = -\chi(X)$. This example of $Y \times I$ was originally constructed in \cite{McD}. The following theorem is a reformulation of Corollary \ref{corGT} of the introduction. We shall prove:

\begin{thm}\label{5d}
Let $Y$ be the unit cotangent bundle of a hyperbolic surface, which fits into a cylindrical Liouville semi-filling $Y\times I$ with $\partial (Y\times I)=(Y,\alpha_-)\bigsqcup (Y,\alpha_+)$, for $\alpha_-$ the standard Liouville form, and $\alpha_+$ a prequantization contact form. Then, the resulting isotopy class of contact structures $\xi$ induced by models A and B on the 5-manifold $M=Y \times \Sigma$, with dividing set consisting of $k\geq 3$ circles (recall Figure \ref{themanifold}), has no Giroux torsion. 
\end{thm}

This follows as a corollary of Theorem \ref{nocob1}, proven below in Section \ref{proofofnocob}. We will also dig into the SFT algebra of this example. We shall investigate whether these examples have 1-torsion (for any $k$, not just $k=1,2$ in which case we already know they do). We will classify all possible building configurations of Lemma \ref{countiszero2} that can contribute to 1-torsion in the whole SFT algebra. It turns out that there are 35 possibilities, and out of those, only one (!) has a nontrivial contribution. And indeed, we have fairly compelling geometric reasons to believe that we get 1-torsion by counting curves in this particular configuration (see Section \ref{esp} below). We give details of this classification in Section \ref{proooof}.

\begin{remark}
Originally, we aimed to show that these models do not achieve 1-torsion, which justifies the fairly involved pursuit of this classification, after which we were aware of this quaint configuration.  
\end{remark}   

\subsection{Curves on symplectization of prequantization spaces: Existence and uniqueness}
\label{prequantcurves}

Consider $(Y,\xi=\ker e^H \lambda)$ a prequantization space over an integral symplectic base $(X,\omega)$, where $H$ is a Hamiltonian on $X$, and $d\lambda=\pi^*\omega$. It is a reasonably standard construction that any choice of $\omega$-compatible complex structure on $X$ naturally lifts to a compatible $J$ on the symplectization of $Y$, for which there exists a finite energy foliation $\mathcal{F}_H$ by $J$-holomorphic cylinders. These project to $X$ as flow lines of $H$, and their asymptotics correspond to circle fibers over the critical points of $H$ (e.g.\ see \cite{Sie}, or our construction of flow-line cylinders in model A). In the case where $Y$ is not the trivial $S^1$-bundle, one can denote $(\theta,y)\mapsto \theta^*y$ the global $S^1$-action on the $S^1$-principal bundle $Y$. Then, these $J$-holomorphic cylinders can be parametrized by
$$
u_\gamma:\mathbb{R}\times S^1 \rightarrow \mathbb{R}\times Y
$$ 
\begin{equation}\label{nontrivcyl}
u_\gamma(s,t)= (a(s),t^*\gamma(s)),
\end{equation}
where $\dot{a}(s)=e^{H (\gamma(s))}$, $\gamma$ solves the ODE  $\dot{\gamma}=\widehat{\nabla H}(\gamma(s))$, and $\widehat{\nabla H}$ denotes the lift of $\nabla H$ to $\xi$.

The case where the base is Liouville, so that the symplectic form $\omega=d\alpha$ is exact, corresponds to the case of the contactization $(Y=X\times S^1, e^{H}(\alpha + d\theta))$, and the parametrization (\ref{nontrivcyl}) coincides with (\ref{trivcyl}) (for $\epsilon=1$) up to reparametrization in the fiber direction. In this case, as in our model A, one can use holomorphic cascades to show that, for sufficiently $C^2$-small $H$, these flow line holomorphic cylinders are the unique somewhere injective holomorphic \emph{curves} having asymptotics corresponding to critical points (see the proof of Theorem \ref{uniqueness}, case A). In the simpler case where one assumes the curve in question is a cylinder, one can use Floer theory as in \cite{Wen3}, or Section \ref{regcyls}.

In the non-exact case, one could try the holomorphic cascade argument we used, which goes as follows. In the degenerate case $H=0$, the projection $\mathbb{R} \times Y \rightarrow X$ is holomorphic, and every Reeb orbit in $Y$ is a multiple of the $S^1$-fiber, lying in a Morse--Bott family which is naturally identified with $X$. The symplectization of $Y$ is foliated by trivial cylinders, projecting to a point in $X$. Any punctured holomorphic curve in $\mathbb{R}\times Y$ asymptotes $S^1$-fibers, and therefore, when projected to $X$, extends to a closed holomorphic curve in $X$. But now it is not necessarily true that this curve will be constant, since $X$ is no longer exact. However, if we ask that $X$ is a surface, one can look at the degree of this map. By holomorphicity, it has non-negative degree, and has zero degree if and only if it is constant. The price to pay, in order to show that this degree is zero, is a restriction on the genus.

\vspace{0.5cm}

Let $X$ be a closed hyperbolic surface (and $Y$ necessarily not a trivial $S^1$-fibration). We will show the following uniqueness result:

\begin{lemma}[Uniqueness]
\label{uniqpre}
Let $X$ be a closed surface of genus $g\geq 1$, and let $\omega$ be an integral area form on $X$. Let $\pi:Y\rightarrow X$ denote the prequantization space over $X$, with a connection (contact) form $\lambda$ whose curvature form satisfies $d\lambda=\pi^*\omega$. Then, for any choice of $\omega$-compatible complex structure on $X$ lifting to a compatible almost complex structure $J$ on $\mathbb{R}\times Y$,  any action threshold $T>0$, and any Morse Hamiltonian $H:X\rightarrow \mathbb{R}$, we can find sufficiently small $\epsilon>0$, such that any genus  $g^\prime<g$ holomorphic curve on $(\mathbb{R}\times Y,e^{\epsilon H}\lambda)$ whose positive asymptotics all correspond to critical points of $H$, and have total action bounded by $T$, is a multiple cover of a flow line cylinder in the finite energy foliation $\mathcal{F}_H$ induced by $H$.
\end{lemma}

\begin{proof} 
\textbf{First proof: holomorphic cascades.}
We follow the argument we started before the statement of the lemma. Recall that positive degree maps between equidimensional closed manifolds induce injections in rational cohomology, as follows easily from Poincar\'e duality (i.e.\ the nondegeneracy of the cup product). In the degenerate case $H=0$, if we apply this to the closed curve in $X$, and the map it induces in the degree 1 part of the cohomology, this implies that $\mathbb{Q}^{2g}$ injects into $\mathbb{Q}^{2g^\prime},$ so that $g\leq g^\prime$. We conclude that for $g^\prime<g$, the projection to $X$ of a curve as in the statement is constant. This means that it is a cover of one of the trivial cylinders in the foliation. 

In the non-degenerate case where $H$ is small but non-identically zero Morse perturbation, one can use holomorphic cascades as in the proof of Theorem \ref{uniqueness}, case A. 

\vspace{0.5cm}

\textbf{Second proof: Siefring intersection theory.} The proof is similar to the proof of Proposition \ref{uniqhyp}, only, in dimension 4. 

The $S^1$-invariance of the contact structure $\xi$ induced by $\lambda$ induces a natural trivialization $\tau$ of $\xi$ along the Reeb orbits $\gamma_p$ corresponding to critical points $p$ of $H$. For any action threshold $T>0$, a standard computation implies that one can find sufficiently small $H$ such that the asymptotic operators corresponding to Reeb orbits $\gamma^l_p$ with actions bounded by $T$ are nondegenerate. Their Conley--Zehnder indices with respect to $\tau$ are given by $$\mu_{CZ}^\tau(\gamma^l_p)=ind_p(H)-1,$$ where $ind_p(H)$ is the Morse index of $H$ at $p$. 

Moreover, the codimension-2 holomorphic foliation $\mathcal{F}_H$ (by flow-line cylinders with respect to $H$) is compatible with the $H$-admissible fibration $(X,\pi,S^1,H)$ for $Y$, in the sense of Section \ref{app} (second paragraph), and the sign function $\sign: \crit(H)\rightarrow \{-1,1\}$ in this case is constant equal to $1$. Given any punctured holomorphic curve $u$ as in the statement, Stokes' theorem implies its negative ends also correspond to critical points of $H$, and then, if we consider its projection to $X$, it extends to a closed surface. Theorem \ref{posdegree} then implies that this map has nonnegative degree and has zero degree precisely when $u$ is a cover of a leaf of $\mathcal{F}_H$. But again Poincar\'e duality implies that it has to have zero degree, and this finishes the proof. 
\end{proof}

\subsection{Obstruction bundles for pair of pants, and super-rigidity}\label{obpants}

In this section, we consider a $J$-holomorphic curve $u: \dot{S}\rightarrow \mathbb{R}\times Y$, where $\dot{S}$ is a 3-punctured sphere, with one positive end and two negative ones, which is a multiple cover of a flow line cylinder $v:=u_\gamma$. We assume that the positive end of $u$ has action bounded by a fixed $T>0$, and we have chosen $H$ so that the Conley--Zehnder indices of $u$ with respect to a natural trivialization (defined below) correspond to the Morse index of the underlying critical points. We show that $u$ admits an obstruction bundle. This is used in Section \ref{proooof}.

In fact, we will show this via Proposition \ref{unb} below, which is more general, and may be of independent interest. 

In the current situation, let $\varphi: \dot{S}\rightarrow \mathbb{R}\times S^1$ be the corresponding branched cover, so that $u=u_\gamma \circ \varphi$, and has branching $B=Z(d\varphi)=1$. Its (generalized) normal bundle is then $N_u=\varphi^*N_v$. Let $p,q \in \mbox{crit}(H)$ be the critical points corresponding to the positive and negative ends of $v$, respectively. We have a natural trivialization $\tau$, which extends globally, for which $c_1^\tau(N_u)=\varphi^*c_1^\tau(N_v)=0$, and $\mu_{CZ}^\tau(\gamma_p)= \mbox{ind}_p(H)-1$.

We have an associated normal Cauchy--Riemann operator
$$
\mathbf{D}_u^N: W^{1,2,\delta_0}(N_u) \rightarrow L^{2,\delta_0}(\overline{Hom}_\mathbb{C}(\dot{S},N_u)),
$$
whose index is given by
\begin{equation}
\begin{split}
\mbox{ind}\;\mathbf{D}_u^N=&\mbox{ind}(u)-2B\\
=&-\chi(\dot{S})+\mbox{ind}_p(H)-1-2(\mbox{ind}_q(H)-1)-2\\
=&\mbox{ind}_p(H)-2\;\mbox{ind}_q(H)\\
\end{split}
\end{equation}

The adjusted first Chern number is 
$$
2c_1(N_u,u)= \mbox{ind}\;\mathbf{D}_u^N-2+\#\Gamma_{even}=\mbox{ind}_p(H)-2\;\mbox{ind}_q(H)-2+\#\Gamma_{even}
$$

\begin{enumerate}[wide, labelwidth=!, labelindent=0pt]
\item \textbf{Case} $(\mbox{ind}_p(H),\mbox{ind}_q(H))\in \{(2,0),(1,0),(0,0)\}$.

We have $(\mbox{ind}\;\mathbf{D}_u^N,c_1(N_u,u))\in \{(2,0),(1,0),(0,-1)\}$, respectively, so that $u$ is regular. Observe that in all case we have $$\mbox{ind}(u)=\mbox{ind}\;\mathbf{D}_u^N+2=\mbox{ind}(v)+2Z(d\varphi),$$ so that $u$ satisfies the obvious necessary condition for regularity to be achieved (see \cite{Wen6}, remark 1.5). 

\item \textbf{Case} $(\mbox{ind}_p(H),\mbox{ind}_q(H))=(2,2)$.

Then, $(\mbox{ind}\;\mathbf{D}_u^N,c_1(N_u,u))=(-2,-2)$ so that $\mathbf{D}_u^N$ is injective.  

\item \textbf{Case} $(\mbox{ind}_p(H),\mbox{ind}_q(H))=(2,1)$.

We get $(\mbox{ind}\;\mathbf{D}_u^N,c_1(N_u,u))=(0,0)$. By Prop.\ 2.2 in \cite{Wen1} we know that $\dim \ker \mathbf{D}_u^N \leq 1$. Observe that there is a section $\eta \in \ker \mathbf{D}^N_v$ which generates the 1-dimensional kernel of the operator of $v$ (which is tangent to $X$ along $\gamma$ and collinear with $X_H$). The pullback $\pi^*\eta$ then lies in $\ker \mathbf{D}_u^N$, and we conclude that $\dim \ker \mathbf{D}_u^N=1$.

\item \textbf{Case} $(\mbox{ind}_p(H),\mbox{ind}_q(H))=(1,1)$. 

We have $(\mbox{ind}\;\mathbf{D}_u^N,c_1(N_u,u))=(-1,0)$. By Prop.\ 2.2 in \cite{Wen1} we know that $\dim \ker \mathbf{D}_u^N \leq 1$, but it is non-conclusive in this case. However, it follows from Proposition \ref{unb} below that $\dim \ker \mathbf{D}_u^N =0$.
\end{enumerate}

We will deal with the last case by an adaptation of Hutchings' \emph{magic trick} (see \cite{Hut2}), as used in proposition 7.2 of \cite{Wen6} in the case of closed curves. We redo Wendl's computation in a punctured setup.

For a holomorphic degree $d$ branched cover $\varphi: (\dot{S},j)\rightarrow (\dot{\Sigma},i)$ between punctured Riemann surfaces $\dot{S}=S\backslash\Gamma(S)$ and $\dot{\Sigma}=\Sigma\backslash\Gamma(\Sigma)$, define by $$\kappa(\varphi)=\sum_{w\in\Gamma(S)}\kappa_w\geq \#\Gamma(S)$$ the \emph{total multiplicity} of $\varphi$, where $\kappa_w=Z(d_w\varphi)+1$ is the multiplicity of $\varphi$ at the puncture $w$ (where $Z(d_w\varphi)$ is the vanishing order of $d\varphi$ at $w$). Observe that 
$$
\kappa(\varphi)=\sum_{z \in \Gamma(\Sigma)}\left(\sum_{w \in \varphi^{-1}(z)}\kappa_z\right)=d\#\Gamma(\Sigma)
$$  
Therefore $d\#\Gamma(\Sigma)\geq \#\Gamma(S)$, with equality if and only if there is no branching at the punctures. We will refer to the sum of the branching of $\varphi$ at interior points of $\dot{S}$ with the branching at the punctures as the \emph{total branching} of $\varphi$.

\begin{prop}\label{unb} Let $v:(\dot{\Sigma},i)\rightarrow (W,J)$ be a somewhere injective possibly punctured holomorphic curve in a 4-dimensional symplectic cobordism, with index zero and $\chi(v)=0$. Assume that there exist a trivialization $\tau$ for which the asympotics of $v$ have vanishing Conley--Zehnder index. Assume also that $v$ satisfies $\ker \mathbf{D}_v^N = 0$.  Then $\mathbf{D}_u^N$ is also injective for any multiple cover $u=v \circ \varphi$ of $v$, such that $\varphi$ has strictly positive total branching. 
\end{prop}

\begin{remark} The proposition does not apply (a priori) only to totally unbranched curves, and the unpunctured case recovers Wendl's results for the torus. It follows from the proof that one can also assume $\chi(v)\geq 0$. However, there are no holomorphic planes as in the hypothesis, since then $0=\mbox{ind}\; \mathbf{D}_v^N=\chi(v)+2c_1^\tau(N_v)=1+2c_1^\tau(N_v)$, but $c_1^\tau(N_v)$ is an integer. One would then only obtain a result for spheres, which is already shown in \cite{Wen6} by a much shorter proof, and where the injectivity of $\mathbf{D}_v^N$ is automatic, and the conclusion holds for every multiple cover. Observe also that we have not assumed that $J$ is generic. For generic $J$, one can use Theorem \ref{PM} below to obtain the conclusion for \emph{every} multiple cover. By a small modification of the proof, one can drop the index zero condition as well as the condition on $\chi(v)$, and impose instead that $c_1^\tau(N_v)=0$ for the given trivialization $\tau$. However, then $\mbox{ind}\;\mathbf{D}_v^N=\chi(v)$, and if this is negative, then generically $v$ does not exist. One then only obtains a potentially interesting result for certain index $1$ planes, or index $2$ spheres.
\end{remark}

In the proof of the proposition, we will use Siefring's intersection pairing, and the adjunction formula for punctured curves. We refer to \cite{Wen8} (especially theorem 4.4) for the statement, and details on the necessary definitions. In particular, the symbol $\overline{\sigma}(u)$ denotes the \emph{spectral covering number} of $u$. This is defined as the sum, over all punctures, of the covering multiplicity of any eigenfunction, corresponding to the eigenspaces of the extremal eigenvalues, of the associated asymptotic operators. Since the multiplicity of the eigenfunction divides that of the asymptotic orbit, we have $\overline{\sigma}(u)=\#\Gamma(u)$ whenever the asymptotics of $u$ are simply covered. Also, the symbol $\delta(u)$ denotes the algebraic count of double points and critical points of $u$. It vanishes if and only if $u$ is embedded. The symbol $\delta_\infty(u)$ is the algebraic count of ``hidden double points at infinity''. If different punctures of $u$ asymptote different orbits, and all of them are simply covered, then this number vanishes (see Thm. 4.17 in \cite{Wen 8} for the general case).

\begin{proof}(Prop. \ref{unb}) As explained in appendix B of \cite{Wen6}, one can construct an almost complex structure $J_N$ on the total space of the normal bundle $\pi: N_v \rightarrow \dot{\Sigma}$, such that there is a 1-1 correspondence between holomorphic curves $u_{\eta}:(\dot{S},j) \rightarrow (N_v,J_N)$ and sections $\eta \in \ker \mathbf{D}_u$ along holomorphic branched covers $\varphi=\pi \circ u_\eta: (\dot{S},j)\rightarrow (\dot{\Sigma},i)$. 

Let $u=v\circ\varphi$ be a branched cover of $v$. We will assume nothing on $\chi(v)$ until the very end. Let $d:=\mbox{deg}(\varphi)\geq 1$. We proceed by induction on $d$. The base case $d=1$ holds by assumption, and so take $d > 1$. Assume by contradiction that there exists a nonzero $\eta \in \ker \mathbf{D}_u^N$, and take the corresponding $u_\eta$. We can view $v$ as a $J_N$-holomorphic embedding inside $N_v$ (as the zero section), and we have that $u_\eta$ is homologous to $d[v]$, where $[v]$ is a \emph{relative} homology class. By the induction hypothesis, we can assume that $u_\eta$ is somewhere injective. Observe that $v$, as a holomorphic map to $N_v$, is embedded, and each of its asymptotics is distinct and simply covered, so that we get $\delta(v)=\delta_\infty(v)=0$, and $\overline{\sigma}(v)=\#\Gamma(v)$.

Since the asymptotics of $v$ have vanishing Conley--Zehnder index, their extremal winding numbers also vanish. Therefore the extremal eigenfunctions are constant, so that their covering multiplicity coincides with that of the corresponding asymptotic orbit. If we denote by $\kappa(u_\eta)$ the sum of all the multiplicities of the asymptotics of $u_\eta$ (which coincides with $\kappa(u):=\kappa(\varphi)$), we obtain that the spectral covering number $\overline{\sigma}(u_\eta)=\kappa(u_\eta)=\kappa(u)=d\#\Gamma(v)$.

We can then compute its Siefring self-intersection number:
\begin{equation}
\label{uno}
\begin{split}
u_\eta*u_\eta&=2(\delta(u_\eta)+\delta_\infty(u_\eta))+c_1^\tau(u_\eta^*TN_v)-\chi(u_\eta)+\overline{\sigma}(u_\eta)-\#\Gamma(u_\eta)\\&=2(\delta(u_\eta)+\delta_\infty(u_\eta))+dc_1^\tau(v^*TN_v)-\chi(u_\eta)+d\#\Gamma(v)-\#\Gamma(u)
\\
\end{split}
\end{equation}

Moreover, Riemann--Hurwitz yields
\begin{equation}\label{cinco}
\chi(u_\eta)=\chi(u)=d\chi(v)-Z(d\varphi)
\end{equation}
Identifying the normal bundle to $v$ inside $N_v$ with $N_v$ itself, the adjunction formula then produces
\begin{equation}\label{tres}
c_1^\tau(v^*TN_v)=\chi(v)+c_1^\tau(N_v)
\end{equation}
Moreover, we have that
\begin{equation}\label{indzeroN}
0=\mbox{ind}\; \mathbf{D}_v^N = \chi(v) + 2c_1^\tau(N_v)
\end{equation}
Using (\ref{indzeroN}) and the vanishing of the Conley-Zehnder indices, we obtain 
\begin{equation}\label{dos}
\begin{split}
u_\eta*u_\eta&=d^2v*v\\
&=d^2(c_1^\tau(v^*TN_v)-\chi(v))\\
&=d^2 c_1^\tau(N_v)\\
&=-\frac{d^2}{2}\chi(v)
\end{split}
\end{equation}

Combining (\ref{uno}), (\ref{cinco}), (\ref{tres}), (\ref{indzeroN}) and (\ref{dos}), we get
\begin{equation}
\begin{split}
0\leq & 2(\delta(u_\eta)+\delta_\infty(u_\delta))\\
=&-\frac{d^2}{2}\chi(v)-d(\chi(v)+c_1^\tau(N_v))+d\chi(v)-Z(d\varphi)+\#\Gamma(u)-d\#\Gamma(v)\\
=&\frac{d(1-d)}{2}\chi(v)-Z(d\varphi)+\#\Gamma(u)-d\#\Gamma(v)
\end{split}
\end{equation}
If we assume that $\chi(v)=0$ (or $\chi(v)\geq 0$), the above is negative, and we obtain a contradiction. 
\end{proof}

\begin{remark} We could have obtained obstruction bundles for case 4 by using the same deformation argument as in Section \ref{OBB}, i.e by taking $\epsilon H$, and $\epsilon \rightarrow 0$. The operator $\mathbf{D}_u^N$ becomes $\overline{\partial}$, which is injective, and this is an open condition. We could also have appealed to \cite{HT1}, where the existence of obstruction bundles for covers of trivial cylinders is shown. We chose this approach because it is more general and has independent interest. Conversely, replacing the normal bundle $N_v$ with $v^*N_H$ in the above proof, this approach reduces the existence of obstruction bundles for curves inside $\mathbb{R}\times H_{hyp}$ (the only case where automatic transversality is non-conclusive) to the somewhere injective or totally unbranched case. So the approach we took in Section \ref{OBB} works better for our particular case, since it proves it for all cases at once. \end{remark}

In order to obtain Theorem \ref{superig}, we need another result. Recall that a closed hyperbolic Reeb orbit is one whose linearized return map does not have eigenvalues in the unit circle. In dimension $3$, if their Conley--Zehnder index is even (in one and hence every trivialization), then they are hyperbolic.

\begin{thm}\cite{Wen7}\label{PM} For generic $J$, all index $0$ punctured holomorphic curves in a $4$-dimensional symplectic cobordism, which are unbranched multiple covers, which have zero Euler charactersitic, and such that all of their asymptotics are hyperbolic, are Fredholm regular. 
\end{thm}

This is proved by the same arguments as in \cite{Wen6} (which only deals with closed curves), and we will not include a proof. Although this paper has recently been withdrawn due to a gap in a lemma, which was crucial for the super-rigidity results there claimed, the hypothesis included in Theorem \ref{PM} bypasses that gap. Indeed, the problem in the paper has to do with the failure in general of \emph{Petri's condition} (see the recent blog post \cite{Wen9} explaining this). In the current situation, for a curve $u$ as in Theorem \ref{PM}, the normal Cauchy-Riemann operator $\mathbf{D}_u^N$ is defined on a line bundle $N_u$ whose adjusted first Chern number is non-positive: we have $2-2g-\#\Gamma=0$, so that $2-2g-\#\Gamma_{odd}=\#\Gamma_{even}$, and therefore
$$
2c_1(N_u,u)=\mbox{ind}(u)-2+2g+\#\Gamma_{even}=-\#\Gamma_{odd}\leq 0 
$$  
This means that either $\ker \mathbf{D}_u^N=0$, or every element in $\ker \mathbf{D}_u^N$ is nowhere vanishing. Therefore Petri's condition is vacuously satisfied (this is also pointed out in the aforementioned blog post). 

With this said, theorem D in that paper needs to be generalized to the setting of asymptotically cylindrical curves in cobordisms. This is mostly a matter of putting it in the proper functional-analytic setup for punctured curves, and almost nothing else changes. Since all orbits are hyperbolic, Fredholm indices of multiple covers are related to indices of the underlying somewhere injective curves via the same multiplicative relations (involving the Riemann--Hurwitz formula) as in the closed case. In fact, the multiplicativity of the Conley--Zehnder indices is all what is needed. With that in mind, the same dimension counting argument as in Theorem B of said paper works in the case of hyperbolic punctures.

\begin{remark}[conjecture] Even though we shall not need it in this document, we also expect the following to be provable from arguments in \cite{Wen6}: \emph{For generic $J$, all punctured holomorphic curves in a $2n$-dimensional symplectic cobordism, which are unbranched multiple covers, and such that all of their asymptotics are hyperbolic, are Fredholm regular.} 
\end{remark}

\begin{proof} (Super-rigidity, Theorem \ref{superig}) In the case where $v$ is as in the hypothesis of Theorem \ref{superig}, generically we can assume it is immersed and regular. Since $\mbox{ind}\;\mathbf{D}_v^N=0$, then $\ker \mathbf{D}_v^N=0$ by regularity, and so the hypothesis of Proposition \ref{unb} are satisfied for generic $J$. If we also assume that its totally unbranched covers satisfy the conclusion of the proposition (which is also a generic condition by Theorem \ref{PM}), the result follows.
\end{proof}

\subsection{Some remarks and useful facts}

Holomorphic curves in the symplectization of $(Y=ST^*X,\alpha_-)$, the unit cotangent bundle of a hyperbolic surface $X$, where $\alpha_-$ is the standard Liouville form, where studied e.g.\ in \cite{CL09,Luo,Sch}. We shall use the following facts about closed $\alpha_-$-orbits:
\begin{itemize}
\item They project to $X$ as closed hyperbolic geodesics.
\item They are non-degenerate.
\item Their Conley--Zehnder index with respect to a natural trivialization of $\xi_-$ coincides with the Morse index (in the calculus of variations sense) of the corresponding geodesic. Therefore, they vanish, since negative curvature implies that no conjugate points exist. 
\end{itemize}

In what follows, we will denote by $H_X$ a choice of a Morse function satisfying the Morse--Smale condition on $X$. Choose a (non-generic) cylindrical almost complex structure $J_+$ on the symplectization of $(Y,\alpha_+=e^{\epsilon^\prime H_X}\lambda)$, coming from a lift of a complex structure on $X$ (as described in Section \ref{prequantcurves}). As before, this yields a foliation $\mathcal{F}_{X}$ of $\mathbb{R}\times Y$ by $H_X$-flow-line cylinders, as explained in Section \ref{prequantcurves}, for $\epsilon^\prime>0$ chosen small enough so that Lemma \ref{uniqpre} holds for $T=T_\epsilon$. We will call $\gamma_{p;q}$ the simply covered Reeb orbit in $M$ corresponding to $(p,q)\in \mbox{crit}(h_+)\times \mbox{crit}(H_X)$. Take $J_0$ a cylindrical almost complex structure in $W_0=\widehat{Y\times I}$, as given by Corollary \ref{compJ0}. Recall that in Section \ref{compJmodelB} we have defined an almost complex structure $J_\epsilon^1$ in $\mathbb{R}\times M$, which is compatible with the SHS $\mathcal{H}_\epsilon^1$, and for which we have a foliation by holomorphic hypersurfaces $\mathcal{F}$. By a generic perturbation of $J_0$ (and hence of $J_\epsilon^1$ along the hypersurfaces in $\mathcal{F}$), we can assume that:

\begin{itemize}
\item[\textbf{(A)}] Every somewhere injective holomorphic curve $u$ in the completion $W_0=\widehat{Y\times I}$ which intersects the main level is regular, and therefore satisfies ind$_{W_0}(u)\geq 0$ (where this denotes the index computed in $W_0$). 
\end{itemize}

The tangent space to the $S^1$-fiber in $Y$ lies in the contact structure $\xi_-$, which then can be globally trivialized by $\partial_\varphi$, the generating vector field of the $S^1$-action. Therefore $\{\partial_\varphi, J_- \partial_\varphi\}$ is a global frame for $\xi_-$, so that $c_1(\xi_-)=0$. 
We also have that $c_1(\xi_+)=0$. This can be seen directly, since we have $c_1(\xi_+)=\pi^*c_1(TX)=\pi^*[\omega]=[d\alpha_+]=0$. We may also appeal to the fact that we have a homotopy of hyperplane fields given by $\xi_r=\ker \alpha_r$, $r \in I$, between $\xi_-$ and $\xi_+$, so that one can push-forward the orientation induced by $\alpha_-$ on $\xi_-$ under this homotopy (or the other way round), which gives the opposite orientation on $\xi_+$ as induced by $\alpha_+$. From this we obtain 

\begin{itemize}
\item[\textbf{(B)}] $c_1(\xi_+)=-c_1(\xi_-)=0$. 
\end{itemize}

Moreover, in the notation of Example \ref{sl2ex}, we have identifications $\alpha_+=k^*$, $\alpha_-=h^*$, $\xi_-=\langle k=\partial_\theta, l\rangle$, $\xi_+=\langle l,h\rangle$, $R_+=k$, $R_-=h$. In particular, we can use our construction of Proposition \ref{explicitJ} and Corollary \ref{compJ0}, so that:
$$c_1(Y\times I)=0$$
   
Therefore, using expression (\ref{xiB}), we obtain
\begin{itemize} 
\item[\textbf{(C)}] $c_1(\xi_\epsilon^1)=0$.
\end{itemize} 

Observe that, since closed Reeb orbits for $\alpha_-$ project down to closed hyperbolic geodesics, they are non-contractible, and Reeb orbits for $\alpha_+$ project to points. This implies the following fact:

\begin{itemize}
\item[\textbf{(D)}] There are no holomorphic cylinders in our model crossing sides of the dividing set.
\end{itemize}

Observe that this does not happen for the $3$-dimensional cases in \cite{LW}, where in fact the hypersurfaces are cylinders. The authors then use it to prove that there is $1$-torsion whenever $\Sigma_-$ or $\Sigma_\pm$ is disconnected (\cite{LW}, Thm. 3).

\vspace{0.5cm}

We can easily see that every multiple of the $S^1$-fibers is not contractible in $Y$. Indeed, the homotopy exact sequence for the fibration $S^1$-fibration $\pi:Y=ST^*X \rightarrow X$ gives $$0\rightarrow \pi_1(S^1)=\mathbb{Z}\rightarrow \pi_1(Y)\rightarrow \pi_1(X)\rightarrow 0,$$ so that the free group generated by the fiber injects into $\pi_1(Y)$. In fact, one can compute via Van-Kampen's theorem, that 
$$
\pi_1(Y)=\left\langle t,a_i,b_i \;(i=1,\dots, g(X)): \prod [a_i,b_i]=t^{-\chi(X)},[a_i,t]=[b_i,t]=1\right\rangle,
$$
where the generator $t$ is represented by the $S^1$-fiber, and the other generators correspond to generators of $\pi_1(X)$. 


\vspace{0.5cm}

Observe that we are allowed to take a \emph{generic} almost complex structure $J_-$ on $(\mathbb{R}\times Y_-,d(e^a\alpha_-))$ and use it in the construction of our model. Since the index of every cylinder is necessarily zero, and every multiply-covered cylinder is unbranched, and multple covers of trivial cylinders are again trivial, we obtain:

\begin{itemize}
\item[\textbf{(E)}] Every holomorphic cylinder in $(\mathbb{R}\times Y, d(e^a\alpha_-))$ is necessarily trivial. And every holomorphic curve of index zero is necessarily a trivial cylinder.
\end{itemize}

Recall that for a given free homotopy class of simple closed curves in a hyperbolic surface $X$, there is a unique geodesic representative. This implies that if two Reeb orbits in $(Y,\alpha_-)$ are joined by a holomorphic cylinder in $W_0$ with two positive ends (which by (\textbf{E}) is necessarily contained in the main level $Y\times I$), they correspond to the same geodesic, but with different orientation. In particular, they are \emph{disjoint} Reeb orbits. We conclude:

\begin{itemize}
\item[\textbf{(F)}] In $W_0$, every holomorphic cylinder with two positive ends on a left upper level joins disjoint Reeb orbits.
\end{itemize}

We make a couple of further general observations, which hold for any curve $u$ in general, and follow from Riemann--Hurwitz:

\begin{itemize}
\item[\textbf{(G)}] If $u$ is a non-constant holomorphic 1-punctured torus, which is a multiple cover of a somewhere injective curve $v$, then $v$ is either a holomorphic plane, or $u=v$ is somewhere injective. 
\end{itemize}

Indeed, By Riemann--Hurwitz we have 
$$-1=\chi(u)=d\chi(v)-B=d(2-2g(v)-1)-B=d(1-2g(v))-B,$$ where $d\geq 1$ is the covering degree and $B\geq 0$ is the total ramification. Then, $$d(1-2g(v))=B-1\geq -1$$ If $g(v)=0$, then $v$ is a plane. If $g(v)\geq 1$, then $-d\geq -1$, so that $d=1$, and $u=v$.

\vspace{0.5cm}
Similarly, in general, we have 

\begin{itemize}
\item[\textbf{(H)}] If $u$ is a non-constant holomorphic 3-punctured sphere, which is a multiple cover of a somewhere injective curve $v$, then $v$ is either a holomorphic cylinder, a holomorphic plane, or $u=v$ is somehwere injective. 
\end{itemize}


Given a holomorphic curve (or building) $u$ with asymptotics corresponding to critical points of $h_\pm$, and then lying in a hypersurface of the foliation $\mathcal{F}$, we will view it as a punctured curve in $W_0$, when convenient. We will denote by $\mbox{ind}_M(u)$ (and $\mbox{ind}_{W_0}(u)$) its Fredholm index when viewed as a curve in $\mathbb{R}\times M$ (resp. in $W_0$). By \textbf{(C)}, we have $$\mbox{ind}_M(u)=(\mu^\tau_{CZ})_M(u)$$$$ \mbox{ind}_{W_0}(u)=-\chi(u)+(\mu^\tau_{CZ})_{W_0}(u)$$ By Proposition \ref{indexthm}, given $l\leq N_\epsilon$, then for orbits $\gamma^l_{p;q}$ for hyperbolic $p$, we have $(\mu^\tau_{CZ})_M(\gamma^l_{p;q})= (\mu^\tau_{CZ})_{W_0}(\gamma^l_{p;q})=\mbox{ind}(q)-1$; and for orbits $\gamma^l_{p;q}$, for $p$ the maximum or minimum, we have $(\mu^\tau_{CZ})_M(\gamma^l_{p;q})= (\mu^\tau_{CZ})_{W_0}(\gamma^l_{p;q})+1=\mbox{ind}(q)$. 

\subsection{Investigating 1-torsion}\label{proooof}

In this section, we will study the portion of the SFT differential of the contact manifold $(M=ST^*X\times Y,\xi)$ which has the potential to yield $1$-torsion, where we denote by $\xi$ the isotopy class defined by models A and B on $M$. Recall from the introduction that we have a power series expansion $\mathbf{D}_{SFT}=\sum_{k\geq 1}D_k \hbar^{k-1}$ of said differential, where $D_k$ is a differential operator of order $\leq k$. This operator is defined by counting holomorphic curves in $\mathbb{R}\times M$ with $|\Gamma^+|+g=k$, where $\Gamma^+$ is the set of positive punctures, and $g$ is the genus. We consider the untwisted version of the SFT, where we do not keep track of homology classes. Basically, we will compute the projection of the operators $D_1$ and $D_2$ to the base field $\mathbb{R}$, but considering their actions on Reeb orbits only up to a large action threshold. In other words, we forget the $q$-variables from these operators (or rather, their restriction of these operators to the subalgebra generated by Reeb orbits with action bounded by a suitably fixed bound), so that we only consider curves with no negative ends. This is obviously motivated by Lemma \ref{countiszero2} (for $N=1$) and Lemma \ref{notorsion}. While the projection of $D_1$ to $\mathbb{R}$ vanishes identically (there are no holomorphic disks), the computation for $D_2$ is rather involved. The result is rather unexpected: among precisely $35$ possibilities, there is basically only one way of obtaining $1$-torsion. While we cannot prove rigorously that $1$-torsion indeed arises, in the next section we provide a heuristic argument as to why we expect this to be true.

Observe that $\xi$ is also the isotopy class given by $\ker \widehat{\Lambda^s_\epsilon}$ for $\epsilon>0$ small and $0\leq s < 1$, where $\widehat{\Lambda^s_\epsilon}$ is the perturbed version of $\Lambda_\epsilon^s$ as defined in equation (\ref{pertcon}), for which we have a compatible $\widehat{J}_\epsilon^s$. Observe that for $k\leq 2$, we already know that our model has $1$-torsion, so we set $k\geq 3$. In the context of lemmas \ref{sympcob2}, \ref{countiszero2}, and \ref{notorsion}, we will set $T:=T_\epsilon, \Lambda:=\widehat{\Lambda^s_\epsilon}, J:=\widehat{J^s_\epsilon}, N=1$. Also, we have $J_\epsilon^1$, the almost complex structure compatible with the non-contact and degenerate SHS $\mathcal{H}_\epsilon^1$, for which we have a foliation of holomorphic hypersurfaces. Recall that $\overline{\mathcal{M}}^1_{g,r}(J_\epsilon^1;T_\epsilon)$ denotes the moduli space of all translation classes of index $1$ connected $J_\epsilon^1$-holomorphic buildings in $\mathbb{R}\times M$ with arithmetic genus $g$, no negative ends, and $r$ positive ends approaching orbits whose periods add up to less than $T_\epsilon$. We will prove the existence of a choice of coherent orientations such that elements in $\overline{\mathcal{M}}^1_{g,r}(J_\epsilon^1;T_\epsilon)$, whenever $g+r\leq 2$, cancel in pairs (after introducing a SHS-to-contact perturbation and a generic perturbation), \emph{except for a single sporadic building configuration}. This unusual cofiguration consists of a single-level punctured torus, and we will discuss it in the next section. This means that the only way of obtaining 1-torsion is to differentiate the positive asymptotic of one of these \emph{sporadic} configurations. In particular, configurations corresponding to cylinders come in cancelling pairs, and this is what we will use to show that there is no Giroux torsion in our model. It follows that, for $0\leq s <1$, $F_\epsilon^s$ sufficiently small, and $\epsilon\leq \epsilon_0$, the moduli space $\overline{\mathcal{M}}^1_{g,r}(\widehat{J_\epsilon^s};T_\epsilon)$ satisfies the same conditions.

\vspace{0.5cm}
 
If $\mathbf{u}$ is an element of $\overline{\mathcal{M}}^1_{g,r}(J_\epsilon^1;T_\epsilon)$ with $g+r\leq 2$, since $r\geq 1$ by exactness, then it can only glue to one of the following: 

\begin{itemize}
\item \textbf{Case 1.} A disk with one positive end.

\item \textbf{Case 2.} A cylinder with two positive ends.

\item \textbf{Case 3.} A 1-punctured torus with one positive end. 

\end{itemize}

\begin{figure}[t]\centering
\includegraphics[width=0.65\linewidth]{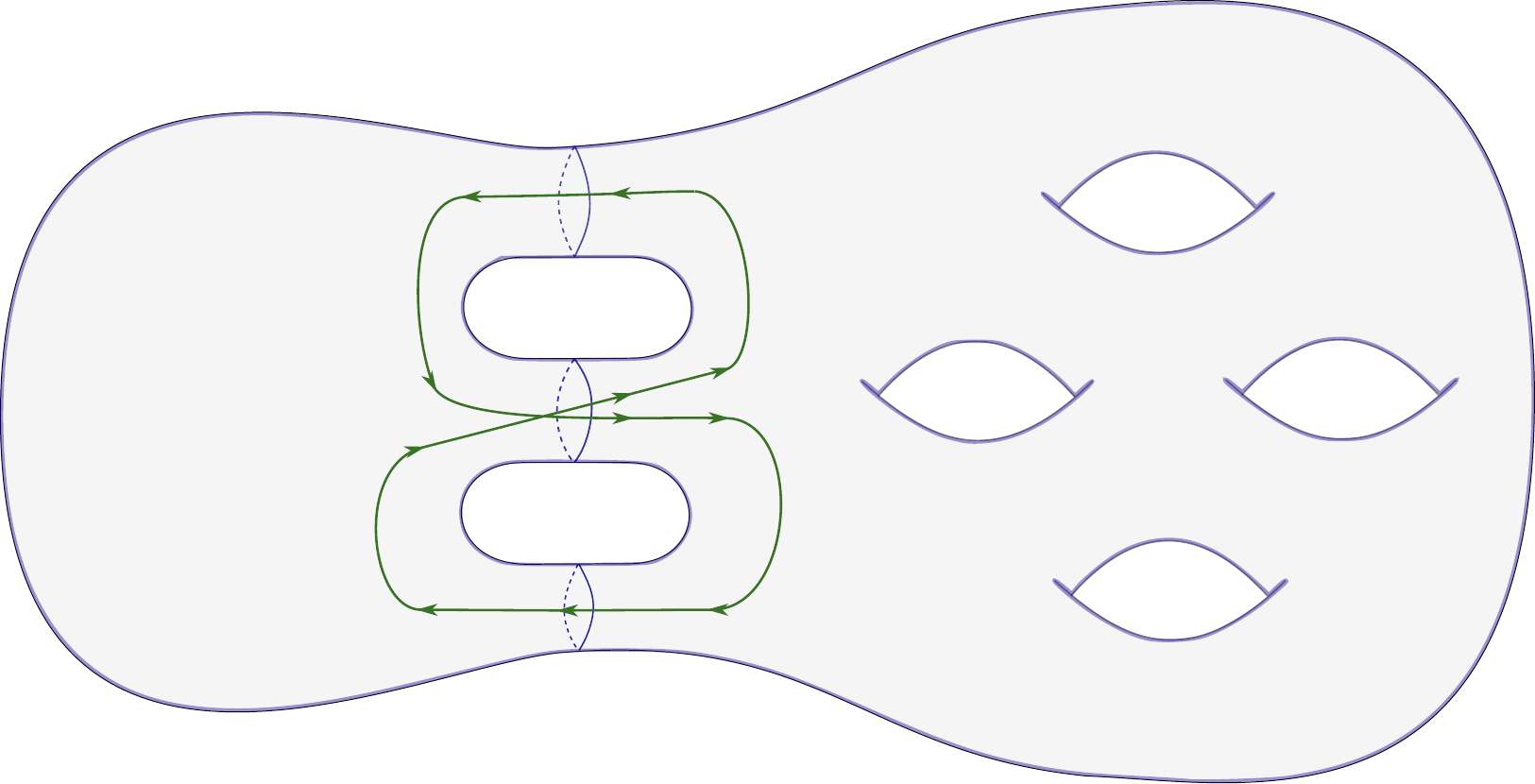}
\caption{\label{curveint} A simple closed curve in $\Sigma$ with non-zero intersection number with the projection of the positive asymptotics which lie in $Y \times \mathcal{U}$, of a hypothetical building $\mathbf{u} \in \overline{\mathcal{M}}^1_{g,r}(J_\epsilon^1;T_\epsilon)$ with $g+r\leq 2$, containing such. It follows that no such $\mathbf{u}$ exists.}
\end{figure}

Case 1 is ruled out, since Reeb orbits in $Y$ are non-contractible (for both the prequantization form and the standard Liouville form), as are the Reeb orbits in $Y \times \mathcal{U}$. The other orbits, which project to $\Sigma$ as closed Hamiltonian orbits and to $Y$ as closed $\alpha_\pm$-orbits may be contractible, but they are ruled out by the definition of $T_\epsilon$ and $\overline{\mathcal{M}}^1_{g,r}(J_\epsilon^1;T_\epsilon)$. This means that there is no $0$-torsion, so this already shows our model is tight (we already knew this from model A, since both $(Y,\alpha_\pm)$ are hypertight).

\vspace{0.5cm}

So we deal with the other two cases. Case 2 is the one which is relevant for Theorem \ref{nocob1} and Corollary \ref{corGT}. Let's first observe that the Reeb orbits in the asymptotics of every floor of $\mathbf{u}$ are of the form $\gamma_{p;q}^l$, for some $l\leq N_\epsilon$. This is certainly the case for asymptotics outside of the cylindrical region, by the action restriction given by $T_\epsilon$. Since Reeb orbits like these project to $\Sigma$ as points, the projection of $\mathbf{u}$ to $\Sigma$ is a cobordism between the sum of the positive asymptotics of the top floor of $\mathbf{u}$ lying in $Y \times\mathcal{U}$ (which are positive multiples of the dividing circles) and zero. On the other hand, since $\mathbf{u}$ has at most two positive asymptotics, there is at least one component of $Y \times \mathcal{U}$ which does not contain asymptotics of $\mathbf{u}$ (here we use $k\geq 3$). But then, one can find a simple closed curve which passes through this component and has non-zero intersection number with the sum of these asymptotics, which is absurd (see Figure \ref{curveint}). By induction, removing the top floor at each step, we obtain the claim.

We can then appeal to Proposition \ref{uniqhyp}, which yields the fact that each component of $\mathbf{u}$ lies in one of the holomorphic hypersurfaces of the foliation, so that we may view them as a building inside $W_0=\widehat{Y\times I}$. 

\vspace{0.5cm} 

\textbf{Case 2. A cylinder with two positive ends.} 

\begin{figure}[t]\centering
\includegraphics[width=0.15\linewidth]{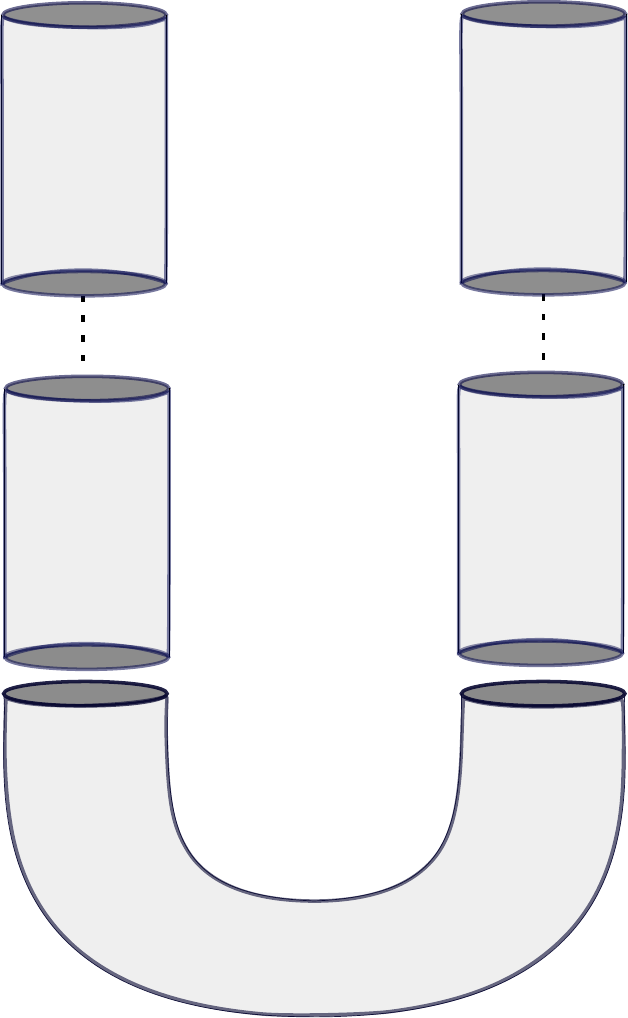}
\caption{\label{poss} The only possible a priori combinatorics for $u$  in the case $g=0$, $r=2$.}
\end{figure}

We claim that since $\mathbf{u}$ glues to a cylinder with no negative ends, and since the Reeb orbits $\gamma_{p;q}$ are non-contractible then $\mathbf{u}$ can only consist of a cylinder with two positive ends in the bottom floor, together with two chains of cylinders on top of its ends (see Figure \ref{poss}). Indeed, no capping off is possible, and no components in consecutive floors can glue along more than a single Reeb orbit (otherwise genus would be positive). So, each component has at most 2 positive ends, since $\mathbf{u}$ does, and at most one negative end, while no pair of pants component appears. There also cannot be any nodes, since no bubbles can appear and the orbits are non-contractible. This shows the claim. 

By \textbf{(D)}, no component of $\mathbf{u}$, all cylinders, can have asymptotics in different sides. This means that if it has a component in a left upper level, it does not have a component in a right upper one (and viceversa). By adding trivial cylinders if necessary, we assume that the two strings of cylinders in the upper levels have the same length. Also, recall that we may assume that $\mathbf{u}$ is a stable building, so that it does not have levels consisting only of trivial cylinders.

We readily see that $\mathbf{u}$ cannot consist solely of right upper level components (corresponding to the prequantization space), since its bottom level would be a cylinder with two positive ends and then cannot be a cover of a flow line cylinder, which contradicts Lemma \ref{uniqpre}. We also see that it cannot consist solely of left upper level components, since its bottom level cannot be a trivial cylinder (it has no negative ends), which contradicts observation \textbf{(E)}.

Then $\mathbf{u}$ has a nontrivial component in the main level, which we call $u_0$.  

\vspace{0.5cm}

\begin{figure}[t]\centering
\includegraphics[width=0.50\linewidth]{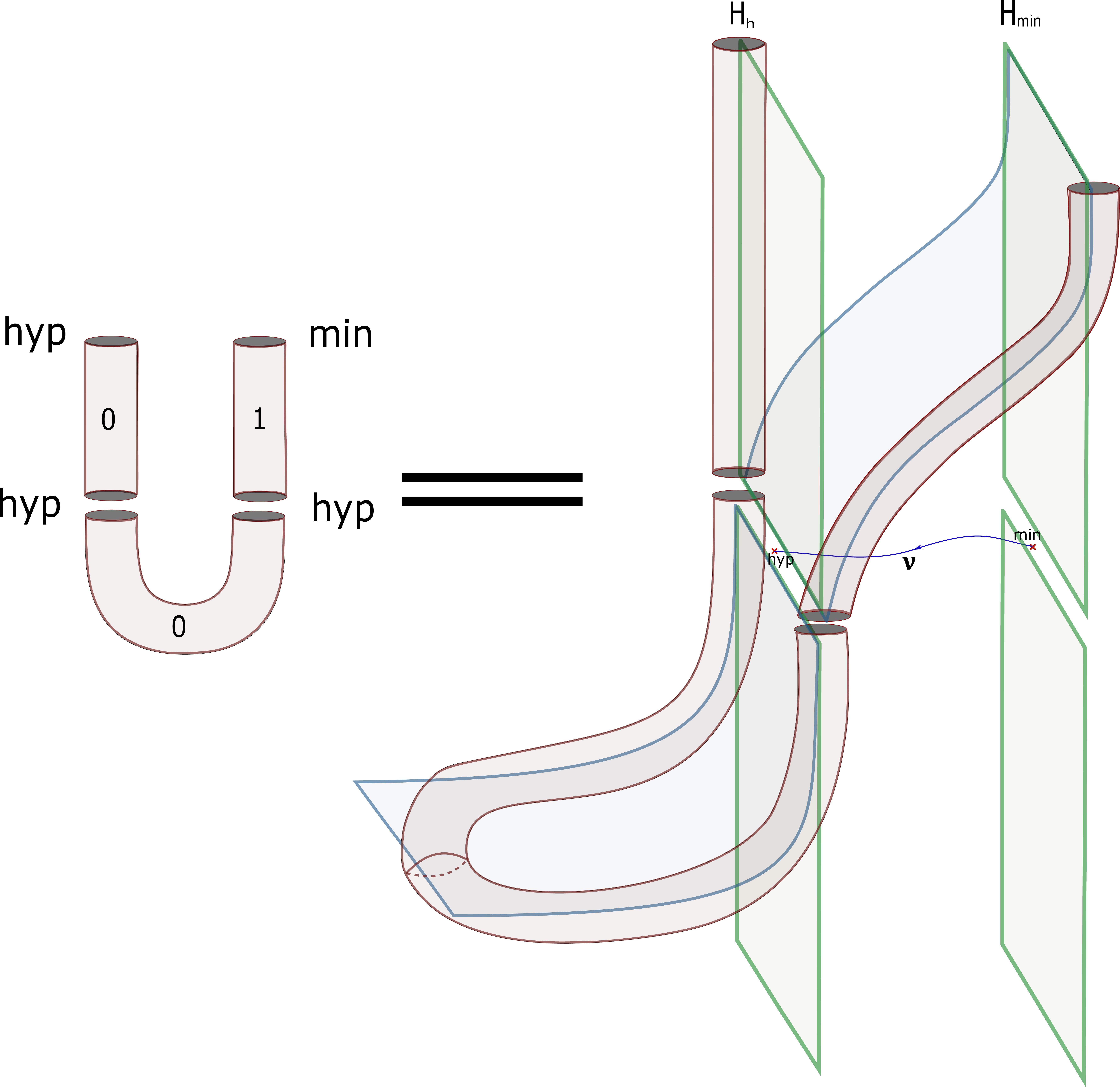}
\caption{\label{excurve} An example of a possible configuration for the index $1$ building $u$ for the case $g=0$, $r=2$, with a non trivial main level. The picture on the right depicts the situation. The bottom floor lies in a hypersurface corresponding to a flow line which crosses sides, and has a nontrivial intersection with $\mathbb{R}\times M_Y$, and the upper floor consists of cylinders. One lies along the hypersurface corresponding to the flow line $\nu$ joining a hyperbolic point $hyp$ with the minimum $min$, the other lies inside the symplectization hypersurface $\mathbb{R}\times H_{hyp}$. The picture on the left is a more diagramatic representation, and includes the Fredholm index of each component. This configuration has an ``evil twin'', which cancels it out when counted.}
\end{figure}

\textbf{Case 2.A. Both asymptotics of $u_0$ lie on a left level.} In this case, every component in the upper levels of $\mathbf{u}$ is either contained in a holomorphic hypersurface $H_a^\nu$ corresponding to an index 1 Morse flow line $\nu$ connecting a hyperbolic point $hyp \in \Sigma_-$ to the minimum $min \in \Sigma_-$, or is contained in either a hypersurface $\mathbb{R}\times H_{hyp}$ or $\mathbb{R}\times H_{min}$ lying over $hyp$, or $min$. Label by $\bold{hyp}$ or by $\bold{min}$ the Reeb orbits appearing in $\mathbf{u}$, according to whether they lie over a hyperbolic point, or the minimum. Since $u_0$ lies in a hypersurface, its two positive asymptotics necessarily have the same label. Thus, to each of its two strings of cylinders, we have associated an ordered sequence of labels, which, by the asymptotic behaviour of the hypersurfaces (depicted in Figure \ref{fol2}), from bottom to top, can only look like $(\bold{hyp},\dots,\bold{hyp},\bold{min},\dots,\bold{min})$. The number of $\bold{hyp}$'s or $\bold{min}$'s may be zero, but not both (see Figure \ref{excurve}). Observation \textbf{(E)} implies that all the upper components correspond to trivial cylinders (under Remark \ref{Weinstein}). Therefore their index in $M$ is either $0$ (in which case both asymptotics have the same label) or $1$ (in which case they differ). We then see that the only possibility for $\mathbf{u}$ is the one depicted in Figure \ref{excurve}, since all others will have index different from $1$.

Observe that $u_0$ has only hyperbolic asymptotics, and, if it is multiply covered, it is necessarily unbranched. Since $J_0$ is generic, Theorem \ref{PM} implies that $u_0$ is regular in $W_0$, and the same is true for its upper components, which are trivial. We therefore have a finite number of these configurations (up to $\mathbb{R}$-translation). While this is a case where Lemma \ref{regvs} does not apply, we still can apply Proposition \ref{regcylcyls}, and we obtain an obstruction bundle of rank $1$ for this configuration. 

And here we make a crucial observation: this configuration has an ``evil twin'', obtained by replacing the index $1$ upper component, which lies over an index $1$ Morse flow-line $\gamma$ connecting $hyp$ to $min$, with the index $1$ cylinder lying over the flow-line $\overline{\gamma}$ (the unique other flow-line connecting $hyp$ to $min$). If we take our coherent orientation to be compatible with the Morse orientation, these twin configurations will cancel out after introducing a generic perturbation. 

\vspace{0.5cm}

We may also write down the associated obstruction bundle for the above configuration. Observe that the moduli space of $u_0$ is just a copy of $\mathbb{R}$, corresponding to its $\mathbb{R}$-translations, since it is otherwise geometrically isolated (in both the hypersurface and the $\Sigma$ directions). In the case where $u_0$ is multiply covered, it is an unbranched cover, with degree $d$, say, and then we would perhaps need to associate a suitable weight to this moduli defined in terms of $d$, if we actually wanted to count the gluings (see Remark \ref{orbifoldcount}). The other components are a trivial cylinder, which is not included in the obstruction bundle, and a regular index $1$ cylinder, which is rigid after dividing by $\mathbb{R}$. If $u_0$ is multiply covered with degree $d$, then so is the index $1$ component. We then see that its obstruction bundle is the trivial bundle
$$
\mathcal{O}_\mathbf{u}=\mathbb{R}\times [R,+\infty) \rightarrow [R,+\infty),
$$ 
where $R$ is a gluing parameter. Observe that the rank and the dimension of the base coincide, so the configuration might indeed glue a finite number of times.

\vspace{0.5cm}

\textbf{Case 2.B. Both asymptotics of $u_0$ lie on a right level.}

Every component in the upper levels of $\mathbf{u}$ is either contained in a holomorphic hypersurface $H_a^\nu$ corresponding to an index 1 Morse flow line $\nu$ connecting a hyperbolic point $hyp \in \Sigma_+$ to the maximum $max \in \Sigma_+$, or is contained in either a hypersurface $\mathbb{R}\times H_{hyp}$ or $\mathbb{R}\times H_{max}$ lying over $hyp$, or $max$. Moreover, Lemma \ref{uniqpre} implies that, when $\mathbf{u}$ is viewed as lying inside $W_0$, they are (necessarily unbranched) covers of flow-line cylinders. Denote by $q_1,q_2 \in \mbox{crit}(H_X)\subseteq X$ the critical points corresponding to the two positive asymptotics of $u_0$. Label by $(\bold{hyp};\mbox{ind}_q(H_X))$ or by $(\bold{max};\mbox{ind}_q(H_X))$ the Reeb orbits appearing in $u$, according to whether they lie over a hyperbolic point, or the maximum, and where $q \in \crit(H_X)$ is the corresponding critical point. Again, since $u_0$ lies in a hypersurface, the first component of the labels of its two positive asymptotics necessarily agree.

Observe that, if $v$ denotes a right upper level component of $u$, one has 
$$1=\mbox{ind}_M(u)=\mbox{ind}_M(u_0)+\sum_{v}\mbox{ind}_M(v)\geq \mbox{ind}_M(v)\geq \mbox{ind}_{W_0}(v)\geq 0$$
Denote by $u_0^\prime$ the somewhere injective curve underlying $u_0$, which, since Reeb orbits are non-contractible, is a cylinder over which $u_0$ is unbranched, and satisfies $\mbox{ind}_{W_0}(u_0^\prime)=\mbox{ind}_{W_0}(u_0)$. By \textbf{(A)}, we have that 
$$
0\leq \mbox{ind}_{W_0}(u_0^\prime)=\mbox{ind}_{W_0}(u_0)=\mbox{ind}_p(H_X)-1+\mbox{ind}_q(H_X)-1\leq\mbox{ind}_M(u)=1,
$$
so that 
$$3\geq \mbox{ind}_p(H_X)+\mbox{ind}_q(H_X)\geq 2$$

If the two labels of the positive ends of $u_0$ are $\bold{max}$, then its index in $M$ is 
$$
\mbox{ind}_M(u_0)= \mbox{ind}_p(H_X)+\mbox{ind}_q(H_X)\geq 2
$$
Since $1=\mbox{ind}_M(u) \geq \mbox{ind}_M(u_0)\geq 2$, we get a contradiction. Then both asymptotics of $u_0$ have $\bold{hyp}$ as the first component of their label, which in particular implies $\mbox{ind}_{W_0}(u_0)=\mbox{ind}_M(u_0)$. 

\begin{figure}[t]\centering
\includegraphics[width=0.99\linewidth]{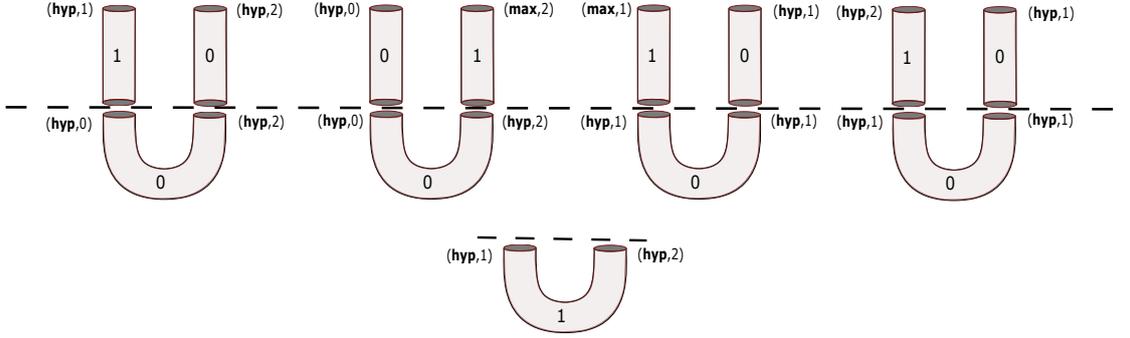}
\caption{\label{cylconf} All 5 possible configurations with positive ends on a right upper level (up to obvious symmetries). The dotted lines separates the main level from the upper levels. The upper level components correspond to unbranched covers of flow-line cylinders over $X$, which are regular in $W_0$. All configurations have ``evil twins'', obtained by replacing the index $1$ component with its mirror. These pairs cancel each other when counted using an obstruction bundle after perturbing.}
\end{figure}

Assume $\mbox{ind}_p(H_X)+\mbox{ind}_q(H_X)=2$, so that $u_0$ has index zero in $W_0$. In the case that $(\mbox{ind}_p(H_X),\mbox{ind}_q(H_X))=(0,2)$, the automatic transversality criterion in \cite{Wen1} implies that it is regular in $W_0$. If $(\mbox{ind}_p(H_X),\mbox{ind}_q(H_X))=(1,1)$, then we use Theorem \ref{PM} to conclude again that $u_0$ is regular in $W_0$. The only possibilities for $\mathbf{u}$ are shown in Figure \ref{cylconf}, which are finite up to $\mathbb{R}$-translation, all of them have an ``evil twin'' in the Morse theory sense, and an associated obstruction bundle of rank 1.

We have only one case left: $(\mbox{ind}_p(H_X),\mbox{ind}_q(H_X))=(1,2)$ (and both labels $\bold{hyp}$). Then $\mbox{ind}_{W_0}(u)=\mbox{ind}_M(u_0)=1$, and $u=u_0$ is non-broken, i.e.\ has only one level, with a non-trivial main level. The bottom component $u_0$ is regular inside its hypersurface, as can be checked via automatic transversality. The resulting configuration, depicted in Figure \ref{cylconf}, also has a rank 1 obstruction bundle and an evil twin, since it lies over an index $1$ flow-line connecting a hyperbolic point in $\Sigma_+$ to the minimum in $\Sigma_-$. 

We conclude that, after perturbing to generic $J$, the count of gluings of the configurations considered is zero, and finishes case 2.

\vspace{0.5cm}

We may also study the obstruction bundles for all the configurations in Figure \ref{cylconf}, and see that we are in the situation explained in Remark \ref{orbifoldcount}. For the ones with two floors, the obstruction bundle is exactly the same as that of the configuration in Figure \ref{excurve}, i.e.\ the trivial $\mathbb{R}$ bundle over $[R,+\infty)$, for large gluing parameter $R\gg 0$ (with potentially some associated weights, in the multiply covered case). For the unique configuration with a single floor, let us assume, as the worst case scenario, that $u_0$ is an unbranched cover of degree $d$. Then we see that its moduli space is $2$-dimensional, of the form $\mathcal{M}_0\cong \mathbb{R}\times \mathcal{M}_H$, where the $\mathbb{R}$-component corresponds to its $\mathbb{R}$-shifts, and $\mathcal{M}_H$ denotes its moduli space inside the hypersurface $H$ containing $u_0$. Therefore, $\mathcal{M}_0\backslash \mathbb{R}\cong \mathcal{M}_H$, and this is, as worst, a $1$-dimensional orbifold. These come only in two types: a circle, or an interval. Observe that in the latter case, its ends necessarily correspond to configurations also listed in Figure \ref{cylconf}, those with positive ends with labels $\{(\bold{hyp},1), (\bold{hyp},2)\}$. We conclude that $\mathcal{O}_\mathbf{u}$ is (at worst) an $\mathbb{R}$-bundle over a circle, or an interval, with perhaps some associated weight in order to count the zeroes of the obstruction section properly.

\vspace{0.5cm} 

\textbf{Case 3. 1-punctured torus with one positive end.} 

\begin{figure}[t]\centering
\includegraphics[width=0.13\linewidth]{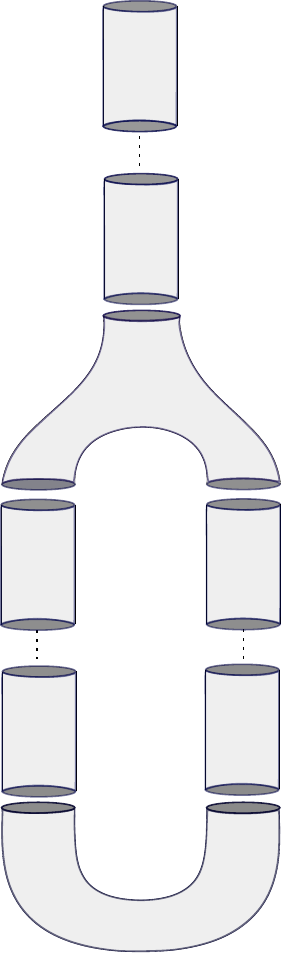}
\caption{\label{buildings2} A priori ``most broken'' possible configuration for the building $u$ for the case $g=1$, $r=1$. All others are obtained by removing breaking from this one.}
\end{figure} 

This case will not influence the proof of Theorem \ref{nocob}. We will discuss it in order to illustrate our techniques, and introduce the sporadic configuration.

\vspace{0.5cm} 

The a priori combinatorial description of $u$ is shown in Figure \ref{buildings2}, which follows from the non-contractibility of the Reeb orbits and the absence of bubbling (in particular there are no nodes). Each of its components is contained in a holomorphic hypersurface of the foliation $\mathcal{F}$.  

\vspace{0.5cm}

\textbf{Case 3.A: The positive end of $u$ lies in a right upper level.}
 
We see that $u$ cannot consist solely of right upper level components. Indeed, by Lemma \ref{uniqpre}, the bottom component of this building is a cover of a flow line cylinder, but this cannot happen since it does not have negative ends. Then $u$ has a nontrivial main level component, which we call $u_0$. 

\textbf{Case 3.A.1: $u_0$ is a cylinder with two positive ends.}

Let $p,q \in \mbox{crit}(H_X)$ be the critical points corresponding to the positive ends of $u_0$, and $u_0^\prime$ be the underlying somewhere injective curve (which is a cylinder over which $u_0$ is unbranched). 

By \textbf{(A)}, we have 
\begin{equation}\label{constr}
0\leq \mbox{ind}_{W_0}(u_0^\prime)=\mbox{ind}_p(H_X)+\mbox{ind}_q(H_X)-2
\end{equation}

\begin{figure}[]\centering
\includegraphics[width=0.9\linewidth]{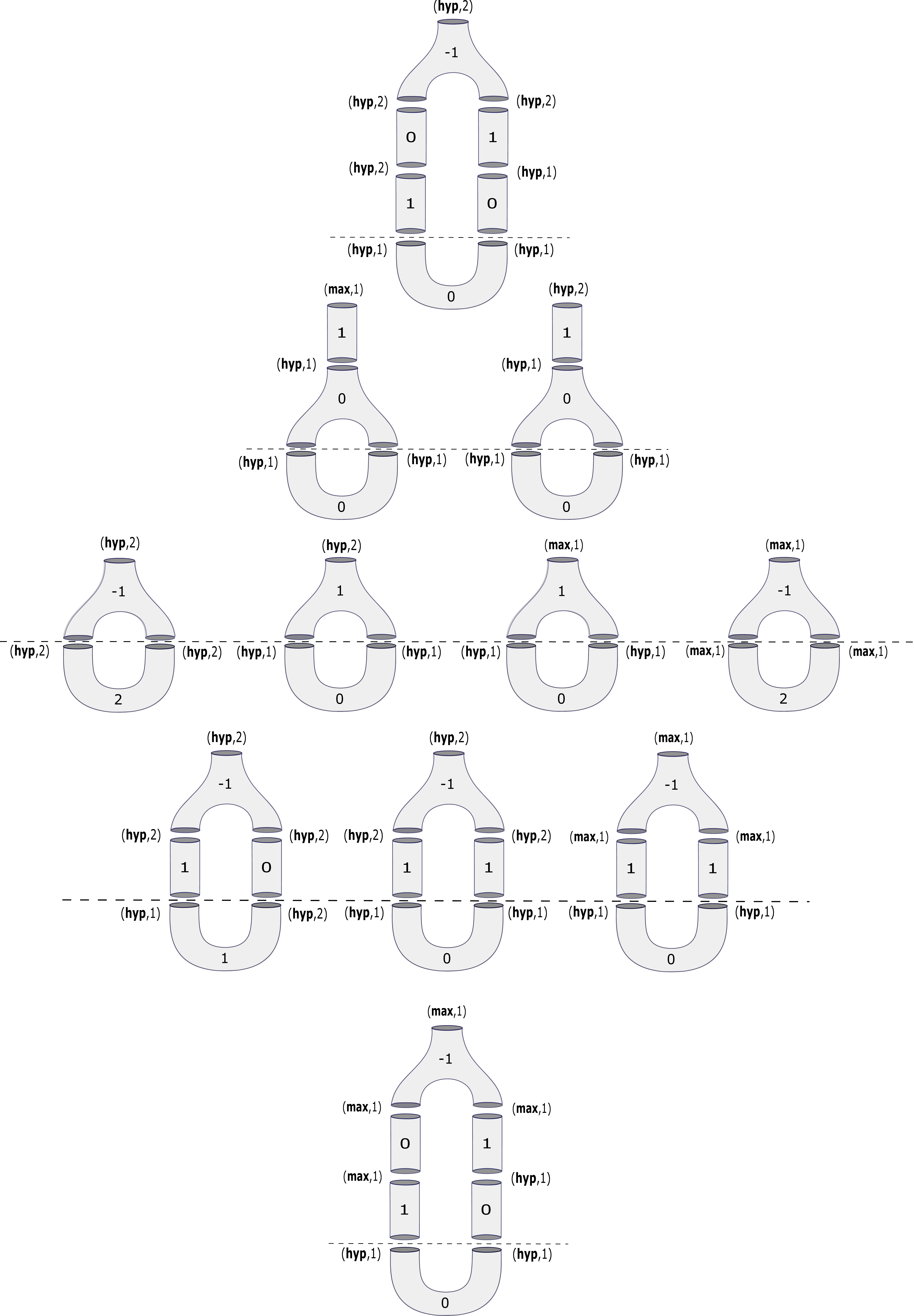}
\caption{\label{posstori} All 11 possibilities in the case where $u_0$ is a cylinder with two positive ends in a right level, ordered by breaking configuration.}
\end{figure}

Using this condition, the fact that each upper component is a cover of a cylinder (by Lemma \ref{uniqpre}) and therefore has non-negative index, and careful index computations, one can classify every possibility. These are listed in Figure \ref{posstori}. We will not include the details, since they are tedious and straightforward. One also uses that the pair of pants components has to have the same label on both negative ends, since it is a cover of a flow-line cylinder, and lies in a hypersurface. Also, that the Morse index (for both Morse functions $H_X$ in $X$ and $H_\epsilon$ in $\Sigma$) goes up as we read from bottom to top. 

Proposition \ref{unb} and Proposition \ref{exisob} give obstruction bundles for the pair of pants component in each configuration. And all the other components are regular inside their corresponding hypersurfaces. Indeed, $u_0$ is always regular by automatic transversality, except in the case where its two positive ends are hyperbolic. In the latter case, we use that $J$ is generic, and $u_0$ is unbranched, and appeal to Theorem \ref{PM}. All the other cylinders are unbranched covers of flow line cylinders, and so regular in the hypersurface. Then each configuration has an obstruction bundle. 

\begin{figure}[]\centering
\includegraphics[width=0.5\linewidth]{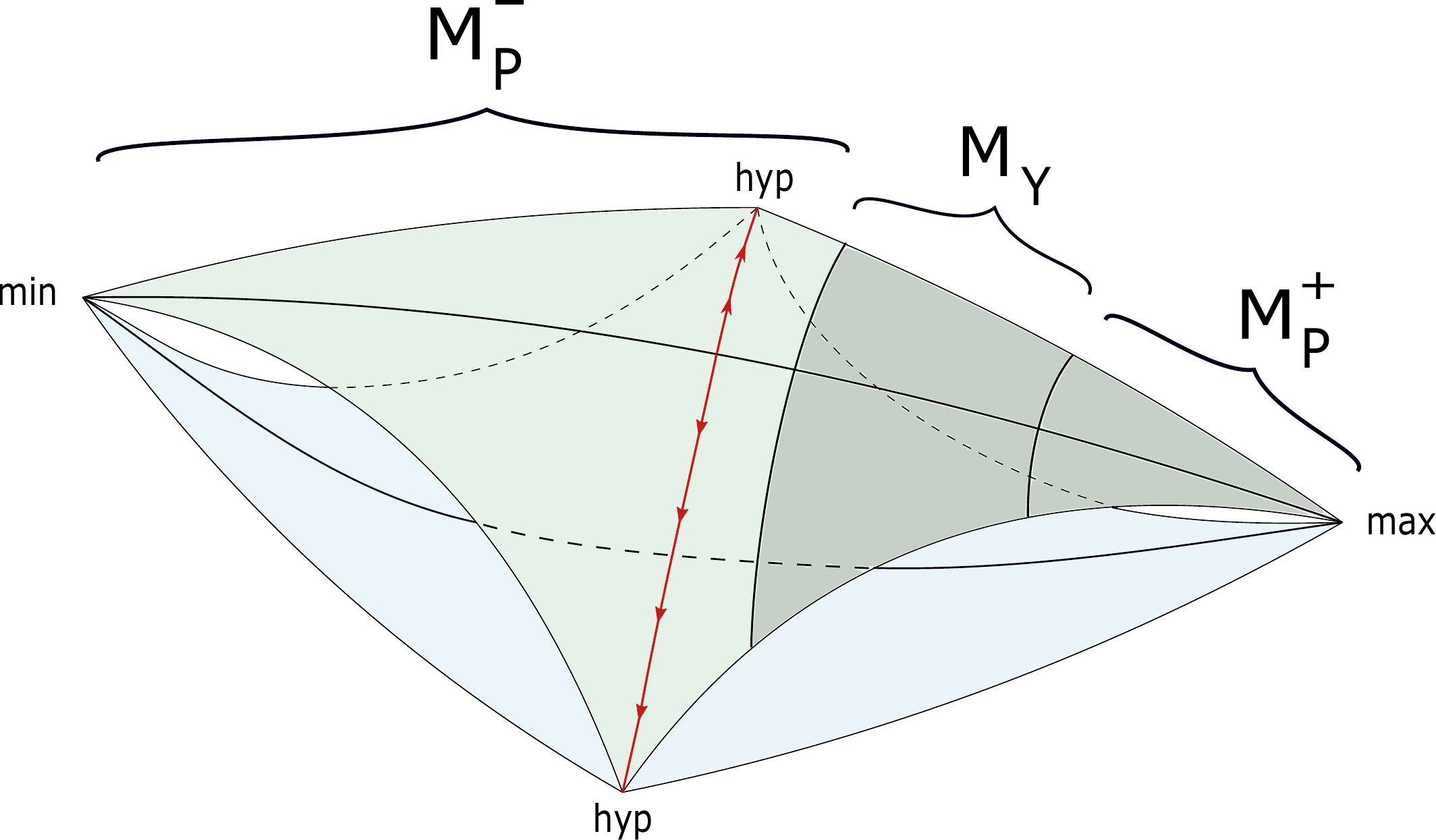}
\caption{\label{diamonds} The twin Morse diamonds. They are glued by the vertices, and give twin moduli spaces, for curves lying over the flow-lines of the diamond. For each moduli space of curves which do not venture into $M_P^-$, there is a 2-dimensional submoduli space parametrized by a subregion of the shaded region, which after quotienting by the $\mathbb{R}$-action can be identified with the red line. The rest of the moduli space is seen inside the corresponding hypersurface. Observe that the twin flow lines induce different signs on the curves over them, and so the twin red lines have different orientation as submoduli. Of course, there is a symmetric picture for curves which do not venture into $M_P^+$. Observe that the twin phenomena also holds for curves which do venture to both sides.}
\end{figure}

Observe that in the case where the first component of the labels of $u_0$ is $\bold{hyp}$, then it has a twin, as we observed before. For the case where it is $\bold{max}$, we claim that the moduli space of $u_0$ also has a twin \emph{moduli space}. Indeed, in this case $u_0$ has index 2, and lies either over an index 2 flow-line, or an index 1 flow-line between the maximum and a hyperbolic point in $\Sigma_-$. Since $u_0$ cannot venture into the latter region (by the maximum principle), its moduli space is parametrized by a 2-dimensional subregion of a ``Morse diamond''. The twin moduli space lies over the twin diamond, and the corresponding sections of the twin obstruction bundles also cancel out (see Figure \ref{diamonds}).

\vspace{0.5cm}

\textbf{Case 3.A.2: $u_0$ is a 1-punctured torus. }

\begin{figure}[]\centering
\includegraphics[width=0.7\linewidth]{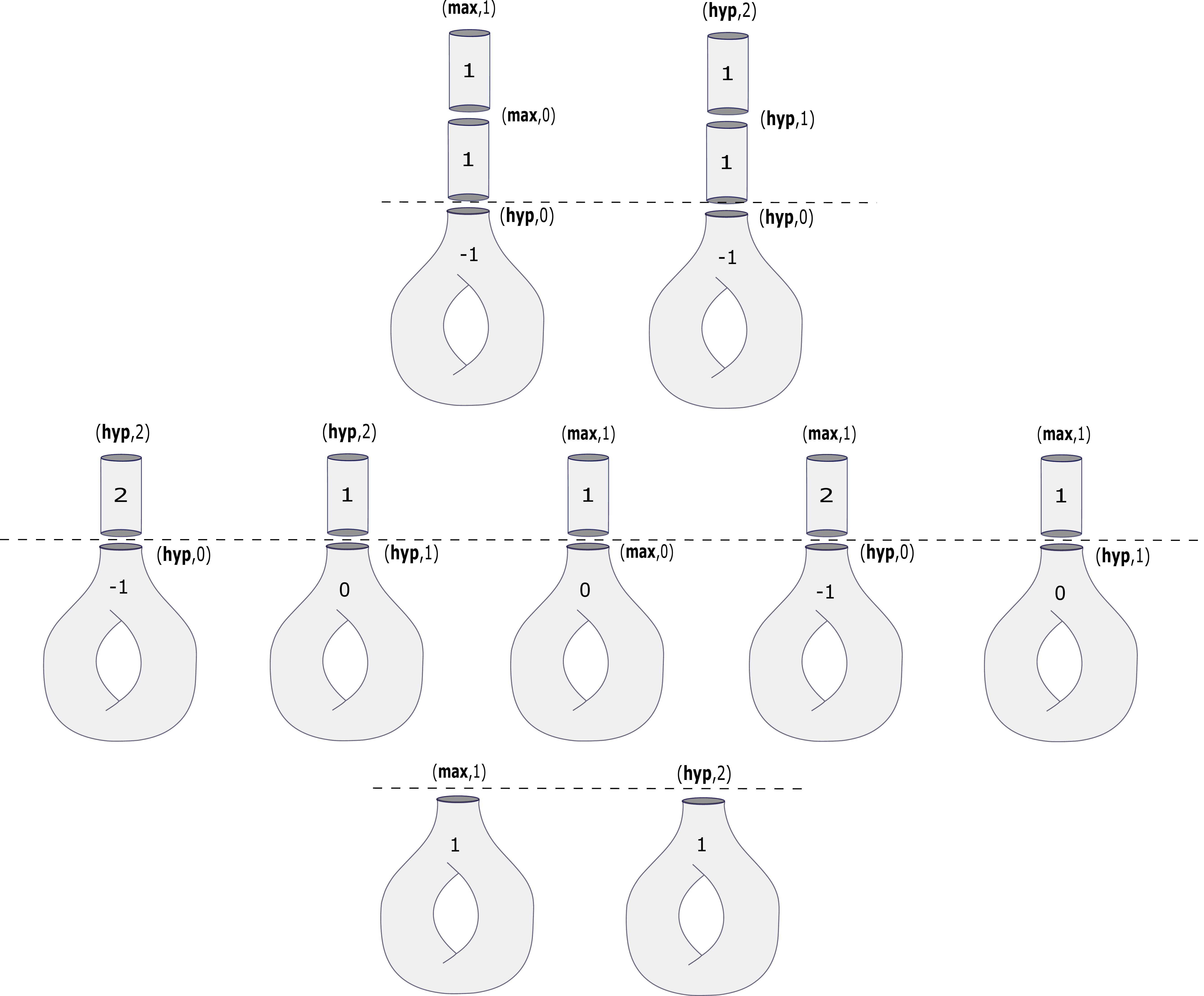}
\caption{\label{moreconfigs} All 9 possibilities in the case where $u_0$ is a 1-punctured torus with the asymptotic in a right level, ordered by breaking configuration.}
\end{figure}

All possibilities are classified in Figure \ref{moreconfigs}. Since Reeb orbits are non-contractible, observation \textbf{(G)} implies that $u_0$ is somewhere injective, and so it is regular inside its hypersurface. By Proposition \ref{exisob}, each configuration then has an obstruction bundle, and at least one twin (for the configuration with one level and label $(\bold{max},1)$, we use the Morse diamond argument of Figure \ref{diamonds}).

\vspace{0.5cm}

\textbf{Case 3.B} \textbf{The positive end lies in a left upper level}. 

Denote by $u_0$ the bottom level of $\mathbf{u}$.

\vspace{0.5cm}

\textbf{Case 3.B.1} $u_0$ \textbf{is a cylinder with two positive ends}. 

\begin{figure}[t]\centering
\includegraphics[width=0.35\linewidth]{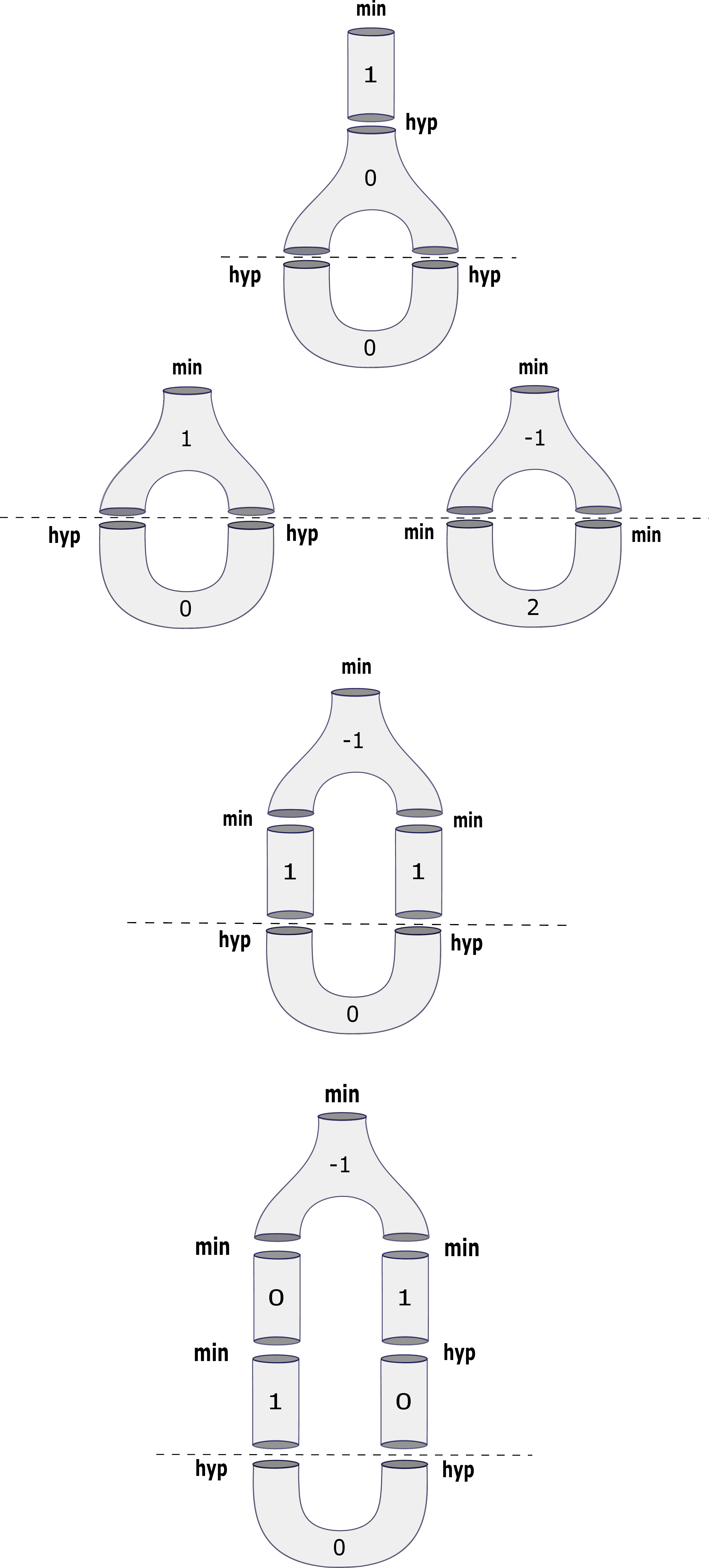}
\caption{\label{configs3} All 5 possibilities in the case where $u_0$ is a cylinder with asymptotics in a left upper level.}
\end{figure}

The curve $u_0$ cannot lie in an upper level, by observation \textbf{(E)}. Then we can classify all configurations, depicted in Figure \ref{configs3}. The pair of pants components cannot be a cover of a (trivial) cylinder, by \textbf{(F)} applied to the cylinder components. By \textbf{(H)}, it follows that it is somewhere injective, and therefore regular in its hypersurface. Then every configuration has an obstruction bundle, and also an evil twin.

\vspace{0.5cm}

\textbf{Case 3.B.2} $u_0$ \textbf{is a 1-punctured torus.} 

\begin{figure}[t]\centering
\includegraphics[width=0.7\linewidth]{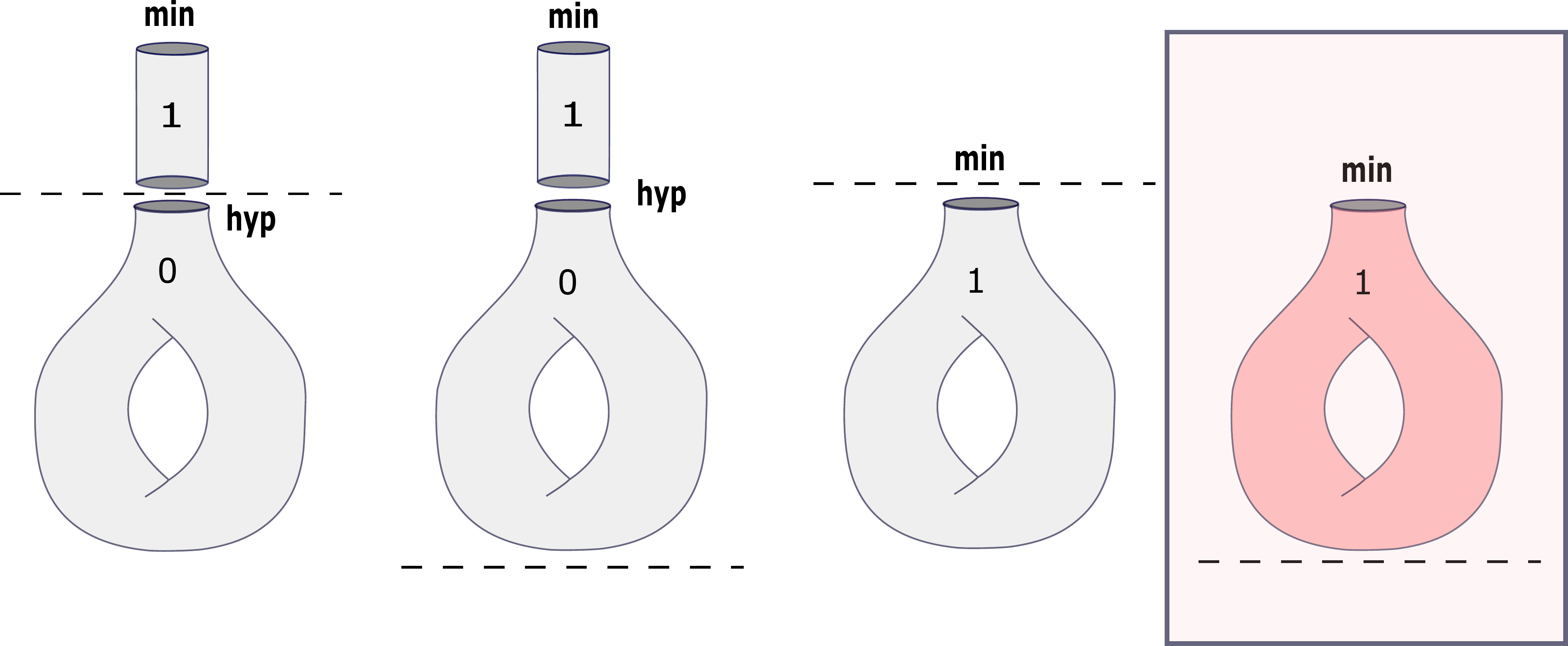}
\caption{\label{moreconfigs2} All 4 possibilities in the case where $u_0$ is a 1-punctured torus with the asymptotic in a left level. The sporadic configuration is depicted in red.}
\end{figure}

As in case 2, we associate labels of the form $(\bold{hyp},\dots,\bold{hyp},\bold{min},\dots,\bold{min})$. By the stability condition and \textbf{(E)}, there can't be two consecutive labels on which are equal, since then we would have a trivial cylinder in $M$. Therefore the only possibilities are the ones listed in Figure \ref{moreconfigs2}. By \textbf{(G)}, and since orbits are non-contractible, $u_0$ is somewhere injective, and so it is regular in its hypersurface. So are the cylinder components, which correspond to trivial cylinders under the identification of Remark \ref{symhyp}. So each configuration has an obstruction bundle. All of them have evil twins, obtained as before, \emph{except the last configuration (!)}.

\subsection{The sporadic configuration}\label{esp}

In this section, we address Conjecture \ref{1tc}. Let $u$ be the sporadic configuration depicted in Figure \ref{moreconfigs2}. Let us first observe that $u$ is in fact regular in $\mathbb{R}\times M$. Indeed, by \textbf{(G)} it is somewhere injective, and so regular in the hypersurface $H:=\mathbb{R}\times H_{min}$. Equation (\ref{indN}) implies that $\mbox{ind}\;\mathbf{D}_{u;N_H}^N=0$, and since $\mathbf{D}_{u;N_H}^N$ is injective ($H$ is cylindrical), we have that it is surjective. By (\ref{cokeq}), we deduce the claim.

\subparagraph*{Some string topology.} Let us recall the definition of Goldman--Turaev string bracket and co-bracket operations, in the case of a hyperbolic surface. This is the starting point of Chas--Sullivan's string topology \cite{CS}. We follow the exposition of \cite{Sch}, which builds on \cite{CL09}.

Let $X$ be a closed surface of genus at least $2$. The set of non-trivial free homotopy classes of loops (i.e.\ a conjugacy class of the fundamental group) is countable. We choose some ordering, label by $[i]$ the $i$-th one, and make the convention that $[-i]$ denotes a change in orientation of the corresponding geodesic. Let $V$ denote the vector space generated by the free homotopy classes of non-contractible loops on $X$. We define the string co-bracket
$$
\Delta: V \rightarrow V \otimes V
$$
as follows. Represent $[i]$ by a string $s$ in general position, so that it has finite self-intersections $x_1,\dots,x_k$. At any $x_j$ there are two directions, which we orient by the orientation on $X$. We can resolve the string $s$ at each self-intersection $x_j$, obtaining two new strings $s_1^j$ and $s_2^j$. For instance, $s_1^j$ is obtained by following the first direction out of $x_j$ and taking the piece of $s$ connecting $x_j$ back to itself, and similarly for $s_2^j$. Define
$$
\Delta([i])=\sum_{j=1}^k [s_1^j]\otimes [s_2^j]-[s_2^j]\otimes [s_1^j]
$$
This co-bracket is a Lie co-bracket, i.e.\ it is bilinear, co-anti-symmetric and satisfies the co-Jacobi identity. Similarly, one define the string bracket
$$
\nabla: V\times V \rightarrow V,
$$
which is a Lie bracket (see \cite{Sch}).

In this notation, we have the coefficients $A^i_{j,k}$ for $\Delta$, defined by
$$
\Delta([i])=\sum_{j,k}A^i_{j,k}[j]\otimes [k]
$$
On the SFT side, denote by $d^i$ the SFT-count of sporadic configurations inside $\mathbb{R}\times ST^*X$, whose positive asymptotic corresponds to $[i]$. Cieliebak--Latschev show, by looking at the boundary of index 1 moduli spaces of holomorphic curves with Lagrangian boundary components in $T^*X$, the existence of linear relationships between the count of index 1 curves in $\mathbb{R}\times ST^*X$, and the coefficients of the above string operations. This allows them to compute the full SFT Hamiltonians for $T^*X$ and $\mathbb{R}\times ST^*X$. In the case of 1-punctured tori, one has the following (\cite{Sch},thm. 2.2.1): 
\begin{equation}\label{dArel}
d^i=\sum_{j} A^i_{j,-j}
\end{equation}

\subparagraph*{Relationship to 1-torsion} Figure 17 in \cite{CL09} gives an example of a closed geodesic $[i]$ for which $d^i$ is non-zero. Denote by $q_i$ the corresponding SFT generator in $\mathbb{R}\times M$. Since we have shown that these sporadic configuration is regular in the latter, we have that
$$
\mathbf{D}\left(\frac{q_i}{d^i}\right)=\hbar + \mathcal{O}(\hbar^2)
$$
Observe that we have used that there are no holomorphic disks.

Now, all the other contributions to the above differential come from potential index 1 building nodal configurations $\mathbf{u}$ with only one positive end approaching $[i]$, which glue to honest curves after introducing an abstract perturbation. If such building is entirely contained in $\mathbb{R}\times H_{min}$, then 
$$
1=\mbox{ind}_M\;(\mathbf{u})=1-\#\Gamma^-(\mathbf{u}),
$$ from which we obtain that $\#\Gamma^-(\mathbf{u})=0$. And if $\mathbf{u}$ is not contained in this cylindrical hypersurface, then from the discussion of Section \ref{proooof} we know that it has a twin, obtained by exploiting the symmetries of the Morse function on $\Sigma$. However, we do not know whether the latter building configurations have associated obstruction bundles, and indeed there are a priori many of them to classify (they can have any number of negative punctures, as well as any positive genus). In any case, this is strong geometric evidence for the following:

\vspace{0.5cm}

\textbf{Conjecture} Any building configuration that may contribute to the differential of $q_i$, and does not lie in $\mathbb{R}\times H_{min}$, cancels its twin when counted in SFT.

\vspace{0.5cm}

If the conjecture were true, then we have that the term $\mathcal{O}(\hbar^2)$ \emph{does not have any $q$ variables}. Since each term of the form $1+\mathcal{O}(\hbar)$ not involving $q$ variables can be inverted as a power series in $\hbar$, we obtain
$$
\mathbf{D}\left(\frac{q_i}{d^i}(1+\mathcal{O}(\hbar))^{-1}\right)=\hbar,
$$
and so our $5$-dimensional model would have 1-torsion (which would prove Conjecture \ref{1tc}).

\vspace{0.5cm}

We observe that this absence of $q$-variables is a purely $5$-dimensional phenomenon, in the following sense: Since $\mbox{ind}_{W_0}(u)=-\chi(u)$ for a curve $u$ inside $H=\mathbb{R}\times H_{min}$, the only possible contributions for the differential of $q_i$ \emph{inside} $H$, come from our sporadic configuration, as well as 3-punctured spheres with 1-positive end at $[i]$, and two negative ones, at $[j]$, and $[k]$, say. However, the latter configurations have index $-1$ in $\mathbb{R}\times M$, so do not contribute to $\mathbf{D}q_i$ in the ambient manifold. Moreover, from \cite{CL09} and \cite{Sch} we have that the SFT count of such curves inside $H$, which we denote $a^i_{j,k}$, coincides with the string coefficients $a^i_{j,k}=A^i_{j,k}$. Observe that we know from (\ref{dArel}) that there is at least one $j$ for which $A_{j,-j}^{i}$ is non-zero. So indeed we can find 3-punctured spheres in $H$ which contribute non-trivially to the SFT count in $H$, but which do not in $\mathbb{R}\times M$. 

\subsection{Non-SFT proof of Theorem \ref{nocob1}}\label{proofofnocob}

We are now ready to obtain an application from our knowledge of the SFT algebra of $ST^*X\times \Sigma$.

\begin{proof}[Thm. \ref{nocob1}] 

Let $(M_0,\xi_0)$ be a $5$-dimensional contact manifold with Giroux torsion, and let $i:(GT,\xi_{GT})\hookrightarrow (M_0,\xi_0)$ be a contact embedding. Consider the $5$-dimensional model $(M=ST^*X\times \Sigma,\xi)$, with $k\geq 3$, where we view $\xi=\ker \Lambda_\epsilon$ as the contact structure induced by a model B contact form $\Lambda_\epsilon$. The $\epsilon$-parameter does not change the isotopy class, and we will choose it suitably below. We may take a contact form $\Lambda_0$ for $\xi_0$, such that it coincides with a model A contact form on $i(GT)$ (as considered in the proof of Theorem \ref{Girouxthen1torsion}). Assume that $(W,d\lambda)$ is an exact cobordism with convex end $(M_0,\xi_0)$ and concave end $(M,\xi)$. Let $J_1$ be a model B almost complex structure on the negative half-symplectization of $(M,\xi)$ (which is compatible with $\Lambda_\epsilon$). Let also $J_0$ be a $\Lambda_0$-compatible almost complex structure on the positive half-symplectization of $(M_0,\xi_0)$, which coincides with a model A almost complex structure on $i(GT)$. 

By attaching small trivial symplectic cobordisms to $W$, we can assume that $\lambda\vert_{M_0}$ and $\lambda\vert_{M}$ are positive and large constant multiples of $\Lambda_0$ and $\Lambda_\epsilon$, respectively, where the constant can be chosen $\epsilon$-independent. Consider the Liouville completion $\widehat{W}$, obtained from $W$ by attaching one positive and one negative cylindrical ends, and the natural extension of $\lambda$ to $\widehat{W}$. Choose a $\lambda$-compatible almost complex structure $J$, such that $J$ coincides with $J_0$ and $J_1$ on their respective ends, and such that $J$ is generic along the main level. Consider also the compactification $\overline{W}^\infty$, obtained by adding $+\infty$ and $-\infty$ to $\widehat{W}$, on which we have extensions of both the symplectic form $d\lambda$ and $J$ (see \cite{Wen4} for more details on this construction). 

By the proof of Theorem \ref{Girouxthen1torsion}, there exists a distinguished ($\mathbb{R}$-translation class of a) regular index 1 $J_0$-holomorphic cylinder $u_0$, lying in the upper levels of $\widehat{W}$. Let $T$ be the total action of its two positive ends. Recall that by the model B construction, there is an action threshold $T_\epsilon>0$, such that $\lim_{\epsilon\rightarrow 0}T_\epsilon=+\infty$ and such that every $\Lambda_\epsilon$-closed orbit with action less than $T_\epsilon$ corresponds to a critical point in $M$, or lies along the spine $Y \times \mathcal{U}$. Observe that this construction also holds if we multiply the contact form by a constant positive number which is $\epsilon$-independent. Take $\epsilon>0$ small enough so that $T_\epsilon>T$. By Lemma \ref{sympcob}, by having attached a suitable exact symplectic cobordism to the concave end of $W$ prior to taking the completion, we can assume that $\lambda\vert_{M}$ is a constant ($\epsilon$-independent) multiple of the model B contact form $\Lambda_\epsilon$ corresponding to this particular $\epsilon>0$. By Stokes' theorem, we get that the action of the Reeb orbits appearing in any building configuration arising as a limit of a curve in the moduli space of $u$ is bounded by $T_\epsilon$. Therefore, in what follows, we can make use of the discussion in Section \ref{proooof}, from which we gather that, whenever this cylinder breaks in lower levels, the Reeb orbits we obtain all correspond to critical points, and none lies over the spine.

Denote by $\overline{\mathcal{M}}$ the connected component containing $u_0$ of the nodal compactification of the space of finite energy $J$-holomorphic curves in $\overline{W}^\infty$. It is a 1-dimensional compact manifold with boundary, having one boundary component corresponding to the projection of $u_0$ to $M_0$.

We observe first that no element of $\overline{\mathcal{M}}$ can break in an upper level. Indeed, by the uniqueness Theorem \ref{uniqueness} (more precisely, its adaptation as in the proof of Theorem \ref{Girouxthen1torsion}) all of the upper components of a stable building $u \in \partial\overline{\mathcal{M}}$ lying in $\mathbb{R}\times M_0$ have to correspond to flow-line cylinders of a Morse function (which have non-negative index). But then $u$ can only consist of an index $1$ cylinder in the bottom end, and a chain of index $0$ cylinders on top of it, which are necessarily trivial. Since $u$ is stable, this configuration is not allowed unless $u$ coincides with its bottom level, but then $u$ is non-broken.  

Similarly, if an element in $\overline{\mathcal{M}}$ breaks with a non-trivial main level component, and no lower level ones, then uniqueness gives that there is only one possible breaking configuration. It consists of a single upper level, having an index $1$ somewhere injective flow line cylinder, together with a trivial (index zero) cylinder; and a nontrivial somewhere injective index $0$ main level component consisting of a cylinder with two positive ends, glued along simply covered Reeb orbits corresponding to critical points of Morse index $1$ and $3$. This configuration has an evil twin, since the index $1$ component does, and it also consists of somewhere injective components. Therefore, we may glue this configuration to a unique honest cylinder. Then we obtain two $1$-dimensional moduli spaces, each having an open end on each twin configuration. We may then canonically identify them along these open ends, and obtain a new moduli space, which we still call $\overline{\mathcal{M}}$, for which we repeat this process.

Assume that an element in $\overline{\mathcal{M}}$ breaks in a lower level of $\overline{W}^\infty$, corresponding to the negative symplectization of $M$. Since the cobordism is exact, there are no holomorphic caps, so that all components of all possible breaking configurations are still cylinders. Then we can argue similarly as above, but using an obstruction bundle. From the discussion in Section \ref{proooof}, we get an obstruction bundle and a twin. Recall also that we have also made sure that we are in the situation of Remark \ref{orbifoldcount}. We may choose finite and sufficiently large gluing parameters, such that the twin configuration glues to honest holomorphic cylinders, whose number is the algebraic count of zeroes of the section of the corresponding obstruction bundle. This number, which is independent of the gluing parameters as long as they are large enough, is the opposite of the original configuration. Since we know by construction that the latter glues at least once, its evil twin will also. And we can proceed as before. Observe that we haven't needed an abstract perturbation.

After these identifications, we will have constructed a 1-dimensional moduli space which has only one boundary component, which is a contradiction. This finishes the proof.
\end{proof} 

\section{Addressing Conjecture \ref{conj}}\label{conj14}

We adapt Lemma \ref{countiszero2} to the case of a holomorphic foliation compatible with admissible fibration, which in particular applies for the examples of Conjecture \ref{conj}. In order to find lower bounds on the algebraic torsion for these fibrations, one needs to find obstruction bundles inside the leaves, and mirrors inside the symplectization of the binding. Observe that this is precisely the condition that failed in the $5$-dimensional models discussed before: the sporadic configurations do not have twins.

Let $(M,\mathcal{H}=(\lambda,\omega))$ be a stable Hamiltonian manifold, and assume that its symplectization is equipped with a compatible $J$, and a codimension-2 $J$-holomorphic foliation $\mathcal{F}$ with binding $V$, which is compatible with an $f$-admissible fibration $(\Sigma,p,Y,f)$ (in the sense of Section \ref{app}). Assume that $\lambda$ is contact around the critical points of $f$, so that the fibers $Y_w$ over $w \in \mbox{crit}(f)$ are contact submanifolds of $M$, and, without loss of generality, that $f>0$. By the symplectic neighbourhood theorem, in a small neighbourhood $\mathcal{U}\subset \Sigma$ of a critical point $w$, we have a local model neighbourhood $\mathbb{R}\times\mathcal{V}\subset \mathbb{R}\times M$ for the codimension-2 symplectic submanifold $\mathbb{R}\times Y_w\subset \mathbb{R}\times M$, where $$(\mathcal{V}=Y \times \mathcal{U}, \lambda=f\alpha+\epsilon \beta)$$ Here, $\epsilon>0$, which we may change without changing the isotopy class of $\lambda$, $\beta$ is a Liouville form on $\mathcal{U}$ which vanishes at $w$, and $\alpha=(\lambda\vert_{Y_w})/f$ is a contact form for $Y$ (compare with our model B contact forms (\ref{modelBform})). As follows from our computations for model B, for any action threshold $T>0$, we can find $\epsilon>0$ small enough so that if $\gamma$ is a closed Reeb orbit inside $V$ with action less than $T$, then its normal Conley-Zehnder index with respect to $V$ satisfies the hypothesis of Theorem \ref{posdegree} below.  

Assume also that the associated sign function is surjective. This means that there is a dividing set of curves in $\Sigma$, which do not contain any critical points of $f$, such that they partition $\Sigma$ into two (possibly disconnected and non-empty) positive and negative regions $\Sigma_\pm$, satisfying $\mbox{sign}^{-1}(\pm 1)\subset \Sigma_\pm$. The asymptotic behaviour of the leaves in the foliation $\mathcal{F}$ is dictated by the endpoints of their corresponding flow line.

In this situation, one can show that curves which are not-too-bad inside a leaf of $\mathcal{F}$, are not-too-bad in $\mathbb{R}\times M$, by the same proof as in Section \ref{obbundles} (using the above local models for the case of a hyperbolic critical point). By Theorem \ref{posdegree}, holomorphic curves with asymptotics in $V$ lie in the leaves of $\mathcal{F}$.

\begin{lemma}\label{countiszero3}
In the above situation, suppose $N$ is a nonnegative integer and assume that:

\begin{itemize} \item All closed Reeb orbits of $\mathcal{H}$ inside $V$ with period less than $T$ are nondegenerate in $V$.

\item For every pair of integers $g \geq 0$ and $r \geq 1$ with $g + r \leq N + 1$, let $\overline{\mathcal{M}}^1_{g,r}(J;T)$ denote the space of all translation classes of index $1$ connected $J$-holomorphic buildings in $\mathbb{R}\times M$ with arithmetic genus $g$, no negative ends, and $r$ positive ends, containing only orbits which lie in $V$ and whose periods add up to less than $T$. Then  $\overline{\mathcal{M}}^1_{g,r}(J;T)$ consists of buildings with possibly multiple floors, consisting of curves which are not-too-bad in the leaves of $\mathcal{F}$ .

\item Denote by $\overline{\mathcal{M}}^1_{g,r}(J;T,V)$ the elements of $\overline{\mathcal{M}}^1_{g,r}(J;T)$ with every component lying inside $V$. Then there is a choice of coherent orientations for $V$ for which the algebraic count of curves obtained by gluing the buildings in $\overline{\mathcal{M}}^1_{g,r}(J;T,V)$, which are curves in $V$, is zero whenever $g+r \leq N + 1$ (after introducing a generic perturbation, and if the zero set of the obstruction section is actually $0$-dimensional). 

\end{itemize}

Then if $\mathbf{D}:\mathcal{A}(\Lambda, T)[[\hbar]] \rightarrow \mathcal{A}(\Lambda, T)[[\hbar]]$ is defined by counting solutions to a
sufficiently small abstract perturbation of the $J$-holomorphic curve equation, there is no $Q \in \mathcal{A}(\Lambda,T)[[\hbar]]$ such that $$\mathbf{D}(Q)=\hbar^N+\mathcal{O}(\hbar^{N+1})$$

\end{lemma}

\begin{remark} Consider $(M=Y \times \Sigma,\xi=\ker \Lambda)$ as in Theorem \ref{thm5d}, and $\Lambda$ a model B contact form, for which we have a foliation $\mathcal{F}$ as above. Then to prove Conjecture \ref{conj}, it is enough to check the second and third hypothesis of the above lemma (with $N=k-2$) for some $Y$, since we can appeal to Lemmas \ref{sympcob2} and \ref{notorsion}. Surely, this is quite a technical challenge, and will be pursued in future work.
\end{remark}

\begin{proof} This basically follows from the symmetries of the Morse function $f$ in $\Sigma$, which we already observed in Section \ref{proooof} (i.e.\ index 1 flow-lines come in evil twin pairs, whereas index 2 flow-lines lie in evil twin Morse diamonds). Combining these symmetries compatibly with a choice of leaf-wise coherent orientations, which coincides with the given one over the critical points, we are in the hypothesis of Lemma \ref{countiszero2}. 

\end{proof}

\appendix
\section{CCH and Lutz--Mori twists}
\label{LMtwists}

This appendix is devoted to describe the Lutz--Mori twists as defined in \cite{MNW}, and to use them to our model A. In our particular case, one can use CCH (cylindrical contact homology) to distinguish the twisted model A contact structures for different amounts of twisting, something which has been done already in \cite{MNW} (Theorem 8.13) for very similar models, so we will not include any details.  

\subsection{The twisted contact structures}

Take $(Y^{2n-1},\alpha_-,\alpha_+)$ a Liouville pair, and denote by $\xi_\pm =\ker \alpha_\pm$. Consider the \emph{Giroux $2\pi l$-torsion domain} modeled on $(Y,\alpha_-,\alpha_+)$, which is the contact manifold $$GT_l^+:=(Y \times [0,2\pi l] \times S^1,\xi_{GT}=\ker\lambda_{GT}),$$ where $$\lambda_{GT}=\frac{1+\cos (r)}{2}\alpha_++\frac{1-\cos (r)}{2}\alpha_-+\sin (r) d\theta,$$ both of whose boundary components are $\xi_{GT}$-round hypersurfaces modeled on $(Y,\xi_+)$. We can also consider the modified version given by $$GT_l^-:=(Y \times [\pi/2,2\pi l+\pi/2] \times S^1,\xi_{GT}=\ker\lambda_{GT}),$$ whose boundaries are now modeled on $(Y,\xi_-)$. 

\vspace{0.5cm}

Observe that our model A SHS $(\Lambda_\epsilon,\Omega_\epsilon)$ has also ($k$ copies of) round hypersurfaces $H_\pm$ modeled on $(Y,\xi_\pm)$, corresponding respectively to the boundary components of the regions $Y \times (\Sigma_\pm \backslash \mathcal{N}(\partial\Sigma_\pm))$. Indeed, recall that $\xi_\epsilon=\xi_\pm \oplus T\Sigma_\pm$ over these particular regions. We may therefore perform an $l$-fold \emph{Lutz--Mori twist} along the $H_\pm$. This can be done in two equivalent ways. The first consists in cutting our model along each of the $k$ copies of $H_-$ and gluing a copy of $GT_l^-$, and the second in doing the same for $H_+$, along which we glue $GT_l^+$. Since the cylindrical region $M_Y$ gluing $Y \times (\Sigma_\pm \backslash \mathcal{N}(\partial\Sigma_\pm))$ together consists of $k$ copies $Y \times I \times S^1$, all of which have a contact structure which is isotopic to a Giroux $2\pi$-torsion domain, the result we get is the same for both gluings. This yields a manifold which is diffeomorphic to $M$, with a stable Hamiltonian structure which coincides with the original outside of $M_Y$, and which we may deform exactly as in Section \ref{MBtoM} to obtain a contact structure on $M$ which we denote $\xi_{l}$. This contact structure is homotopic as almost contact structures to the contact structure $\xi_{\epsilon,s}$ which we constructed before. Since $\xi_{l}$ has Giroux torsion, Section \ref{GirouxTorsion1-tors} implies that it has $\Omega$-twisted algebraic $1$-torsion (for $l>0$), so that they are indistinguishable from this invariant alone. Corollary \ref{nocob} follows. 

On the other hand, recall that the \emph{cylindrical contact homology} (CCH) of $(M,\xi)$ with respect to a homotopy class $h$ is the homology of a chain complex Reeb orbits in the class $h$, and the differential counts holomorphic cylinders with one positive and one negative end between them. One can compute, in the same way that it is done in \cite{MNW}, or, in a more detailed fashion, in \cite{Wen3} (for the case of $\mathbb{T}^3$), that the CCH of $(M,\xi_l)$ with respect to the free homotopy class of loops $h$ corresponding to $t \mapsto (const,const,t) \in Y \times I \times S^1 \subseteq M$ (any of the $k$ different choices) is given by 
\begin{equation}
HC_*^h(M,\xi_{l})=\left\{\begin{array}{ll} H_{*+n}(Y)^{\oplus_{l}} &, \mbox { if } k \neq 2\\
H_{*+n}(Y)^{\oplus_{2l}} &, \mbox { if } k = 2 
\end{array}\right.
\end{equation}
The case $k=2$ is different since then the homotopy classes for both components of $M_Y$ agree. This implies that $\xi_{l}$ is not isotopic to $\xi_{l^\prime}$ if $l \neq l^\prime$, since an isotopy acts trivially in $\pi_1(M)$.

Moreover, if we assume that $(Y_\pm, \xi_\pm)$ has no contractible Reeb orbits (which can be achieved in the examples of \cite{MNW}), so that any other homotopy class $h^\prime$ of the subgroup $G:=\pi_1(\Sigma)\times \{1\}$ different from $h$ is not achieved by any Reeb orbit in $M$, and hence $HC_*^{h^\prime}(M,\xi_{l})=0$, we have that they are not contactomorphic by a map preserving $G$. If we further assume that $\pi_1(Y)$ has trivial center and is solvable, then, by lemma 8.15 in the aforementioned paper, they are not contactomorphic at all.   

\section{Spinal Open Book Decompositions (SOBD's)}\label{SOBDap}

We have decided to include a definition of a higher-dimensional SOBD, since it has already proved useful for us to use related language. This comes with the disclaimer that it is the version suited for the purposes that we have pursued in this document, and can definitely be relaxed/simplified/modified (for instance, one might need to if one wanted to prove a higher-dimensional version of the Giroux correspondence). The treatment will be therefore rather speculative. The original 3-dimensional notion was defined in \cite{LVHMW}.

\begin{defn}\label{SOBDdef}
A spinal open book decomposition (SOBD) on a compact oriented $(2n+1)$-manifold $M$, possibly with boundary, is a decomposition
$$M = M_S \cup M_P,$$ where the pieces $M_S$ and $M_P$ (called the \emph{spine} and \emph{paper} respectively) are smooth compact $(2n+1)$-dimensional submanifolds with disjoint interiors, which are glued along common boundary components, and which carry the following additional structure:

\begin{itemize}

\item \textbf{(spine fiber bundle)} A smooth fiber bundle $\pi_S: M_S \rightarrow \Sigma$ with typical fiber a closed, connected and oriented contact manifold $C^{2(n-k)+1}$, all of which are either disjoint from $\partial M_S$, or contained in it (in which case they are either contained in a component of $\partial M_P$ or of $\partial M$). 

Here, $\Sigma$ is a (possibly disconnected) $2k$-dimensional compact oriented manifold (with $k\leq n$) which has the structure of a strong symplectic filling, so that in particular its connected components (called \emph{vertebrae}) all have nonempty boundary.

\item \textbf{(paper fiber bundle)} A smooth fiber bundle $\pi_P : M_P \rightarrow Y$ over a (possibly disconnected) closed $(2k-1)$-dimensional contact manifold $Y$, with typical fiber a (possibly disconnected) oriented manifold $P^{2(n-k+1)}$ which can be given the structure of a Liouville domain, so that the connected components of these fibers (called \emph{pages}) all have nonempty boundary. The pages meet $\partial M_P$ transversely, either intersecting $\partial M_S$ in a disjoint union of fiber of $\pi_S$, or intersecting $\partial M$. 

\item \textbf{(compatibility along the interface)} One can require the following compatibility condition where the spine and paper fiber bundles glue together: 

There exists an \emph{interface region} of the form 
$$
M_I:=(\partial M_S \cap \partial M_P) \times [-1,1]=C \times Y \times [-1,1],
$$
which is the region along which $M_S$ and $M_P$ glue together, such that the spine fiber bundle 
$$
\pi_{S}\vert_{M_I}:M_I \rightarrow Y \times [-1,1] 
$$
is the trivial one, and if $\pi_0: M_I \rightarrow Y$ denotes the trivial fibration, then 
$$
\pi_P\vert_{M_I}=f^* \pi_0,
$$
is the pullback bundle with respect to some contactomorphism $f: Y \rightarrow Y$.   

In particular, this implies that the connected components of $\partial \Sigma$, which correspond to components of $\partial M_S$ which meet $M_P$, are all diffeomorphic to a component of $Y$.

\end{itemize}
\end{defn} 

The most interesting special case is when $k=n$, which means that $C=S^1$, so that the spine fiber bundle can be endowed with the structure of a principal $S^1$-bundle, and the pages are surfaces. The 3-dimensional case is $k=n=1$, on which also $Y=S^1$. 

\begin{remark}
According to the definition above, we are allowing both the spine and the paper to touch $\partial M$ (in contrast with the definition in \cite{LVHMW}, in which only the paper is allowed to touch $\partial M$). Since in practice sometimes one may wish to allow spine components touching the boundary (like the proof of Theorem \ref{Girouxthen1torsion}), we have decided on this definition. 

One perhaps could drop the compatibility along the interface, to get something more general. The map $f$ of the definition is in place since in general one might wish to allow different boundary components of a page to touch the same spine component (it follows by definition all of these components are diffeomorphic to $C$), and $f$ keeps track of that. A special case is when $f=id$, in which the paper fiber bundle is trivial along the interface. In the $3$-dimensional case, a motivating example is $f:S^1\rightarrow S^1$ mapping $z \mapsto z^n$. One could attempt, more generally, to replace $f$ with a smooth map $F:C \times Y \rightarrow C \times Y$, such that $\pi_P\vert_{M_I}=\pi_0 \circ (F \times Id_{[-1,1]})$, satisfying that $F$ is transverse to every $C \times \{pt\}$, and no connected components of the manifold $F^{-1}(C\times \{pt\})$ is regularly homotopic to $\{pt\} \times Y$. Morally, the fibers of $\pi_P\vert_{M_I}$ are not ``meridional''. This would be an attempt at generalizing the notion of a rational open book as defined in \cite{BEVHM}, but it is not entirely clear how to include the contact/symplectic-topological data so as to provide a proper generalization. 

If we  let the spine and paper bundles be globally trivial, then their total spaces are interchangeable, as well as the roles of $Y$ and $C$, and we obtain a dual SOBD (something we used in the duality between models A and B of the main sections).
  
We have chosen to define the spine bundle as having a base which is not necessarily a Liouville domain, since there are interesting examples on whose base is not necessarily so. For instance the prequantization SOBDs as in Section \ref{PrequantSOBDs}, which arise naturally on principal $S^1$-bundles defined by an Euler class vanishing along a dividing set of hypersurfaces.

This definition is not suited for proving that SOBDs support contact structures in the sense of Giroux, which one could attempt in a different way that we have (in the particular case of prequantization SOBDs, and models A and B). If one wished to pursue a purpose like that, one probably would need extra data (such as a choice of foliation along the boundary or a choice of a characteristic distribution, whatever is convenient). An ideal approach would be to define a suitable notion of Giroux form by a easy-to-check set of conditions mimicking the standard $3$-dimensional situation (which contains as examples the contact forms of Sections \ref{GirouxSOBDs} and \ref{PrequantSOBDs} obtained by gluing). We will not pursue that here.   

\end{remark}

\section{Intersection theory of punctured curves and pseudoholomorphic hypersurfaces (with Richard Siefring)}\label{AppSiefring}

\subsection{Basic theory}
In this Appendix we summarize some facts that will be proved in the in-preparation work \cite{Sie2}, which extends some of the intersection theory for punctured pseudoholomorphic curves from \cite{Sie08, Sie11} to higher dimensions,
and we use these facts to prove a result
(Theorem \ref{posdegree}) needed for proving Proposition \ref{uniqhyp} and the uniqueness Lemma \ref{uniqpre}.

\vspace{0.5cm}

We will consider a closed, orientable
manifold $M^{2n+1}$ equipped with a nondegenerate
stable Hamiltonian structure $(\lambda, \omega)$
and denote the associated hyperplane distribution and Reeb vector field by
$\xi=\ker\lambda$ and $R$ respectively.
A codimension-$2$ submanifold $V\subset M$ is said to be a
\emph{stable Hamiltonian hypersurface} of $M$ if
the pair $\lambda'$, $\omega'$ defined by
\[
\lambda':=i^{*}\lambda \qquad \omega':=i^{*}\omega
\]
where $i:V\hookrightarrow M$ is the inclusion map is a stable Hamiltonian structure on $V$.
We say a stable Hamiltonian hypersurface $V\subset M$ is a
\emph{strong stable Hamiltonian hypersurface}
if, in addition, $V$ is invariant under the flow of $R$.
We note that this condition is equivalent to requiring that
the Reeb vector field $R'$ of the stable Hamiltonian structure $\lambda'$, $\omega'$ on $V$
to equal $R$ at all points in $V$.
Indeed, using naturality of contractions
we have that
\begin{equation}\label{e:restricted-forms}
i_{v}\lambda'=i^{*}\bp{i_{v}\lambda} \quad\text{ and }\quad i_{v}\omega'=i^{*}\bp{i_{v}\omega}
\quad\text{ for $v\in TV\subset TM|_{V}$.}
\end{equation}
If $V$ is invariant under the flow of $R$, then $R|_{V}$ is a section of $TV$ so the above shows that
along $V$ we will have 
\[
i_{R}\lambda'=i^{*}\bp{i_{R}\lambda}=1 \qquad\text{ and }\quad i_{R}\omega'=i^{*}\bp{i_{R}\omega}=0
\]
so
$R|_{V}$ is the Reeb vector field of $\lambda'$.
On the other hand,
if $R'(p)=R(p)$ for all $p\in V$, then it's clear that $V$ is invariant under the flow of $R$.

Given a
stable Hamiltonian hypersurface $V$ of $M$,
the first part of \eqref{e:restricted-forms} implies that a vector field $v$ in $TV$ satisfies
$i_{v}\lambda'=0$ precisely when $i_{v}\lambda=0$.  Therefore, we have that the
hyperplane distribution
$\xi':=\ker\lambda'$ is given by
\[
\xi'=\xi|_{V}\cap TV.
\]
Letting $\xi_{V}^{\perp}$ denote the symplectic complement
\[
\xi_{V}^{\perp}=\br{v\in \xi|_{V}\,\middle|\, \omega(v, w)=0 \quad\forall w\in \xi'},
\]
to $\xi'$
we obtain a splitting 
\[
\xi|_{V}=\xi'\oplus\xi_{V}^{\perp}
\]
of $\xi$ along $V$.
We note that since $\xi'$ is a symplectic subbundle of
$(\xi|_{V}, \omega)$
its symplectic complement $\xi_{V}^{\perp}$ is a symplectic subbundle as well.

If $V$ is a
strong stable Hamiltonian hypersurface
then using that
\[
TV=\R R'\oplus\xi'=\R R\oplus\xi'
\]
we have that
\[
\xi_{V}^{\perp}=TV^{\perp}:=\br{v\in \xi|_{V}\,\middle|\, \omega(v, w)=0 \quad\forall w\in TV}.
\]
We thus have a splitting
\begin{equation}\label{e:splitting}
TM|_{V}=\R R\oplus(\xi', \omega')\oplus (\xi_{V}^{\perp}, \omega)
\end{equation}
so that the first two summands give $TV$,
and the linearized flow of $R$ along $V$ preserves this splitting along with the symplectic structure on the
second two summands.

Let $\gamma:S^{1}\approx\R/\Z\to M$ be a $T$-periodic orbit of $R$, i.e.\ $\gamma$ satisfies the equation
\[
\dot\gamma(t)=T\cdot R(\gamma(t))
\]
for all $t\in S^{1}$.
Assuming that 
$\gamma(S^{1})\subset V$, 
we can choose a symplectic trivialization of the
hyperplane distribution
$\xi=\xi'\oplus \xi_{V}^{\perp}|_{\gamma(S^{1})}$ along $\gamma$ which respects the splitting
\eqref{e:splitting}, i.e.\ one of the form
\[
\Phi=\Phi_{T}\oplus\Phi_{N}:\xi'\oplus \xi_{V}^{\perp}|_{\gamma(S^{1})} \to
S^{1}\times (\R^{2n-2}, \omega_{0})\oplus(\R^{2}, \omega_{0})
\]
with $\Phi_{T}$ and $\Phi_{N}$ symplectic trivialization of $\xi'|_{\gamma(S^{1})}$ and
$\xi_{V}^{\perp}|_{\gamma(S^{1})}$ respectively.
Given such a trivialization we can define the Conley--Zehnder index of the orbit
$\gamma$ viewed as an orbit in $M$ as usual by
\[
\mu^{\Phi}(\gamma)
:=\mu_{CZ}(\Phi(\gamma(t))\circ d\psi_{Tt}(\gamma(0))\circ\Phi(\gamma(0))^{-1})
\]
where $\psi:\R\times M\to M$ is the flow generated by $R$ and where
$\mu_{CZ}$ on a path of symplectic matrices starting at the identity and ending at a matrix without $1$
in the spectrum is as defined in \cite[Theorem 3.1]{HWZ}.  But since, as observed above,
$d\psi_{t}$ preserves the splitting \eqref{e:splitting}, we can also consider the
Conley--Zenhder indices that arise from the restrictions of $d\psi_{t}$ to $\xi'|_{\gamma}$ and
$\xi_{V}^{\perp}|_{\gamma}$.
In particular we define
\[
\mu^{\Phi_{T}}_{V}(\gamma)
:=\mu_{CZ}(\Phi_{T}(\gamma(t))\circ d\psi_{Tt}(\gamma(0))|_{\xi'}\circ\Phi_{T}(\gamma(0))^{-1})
\]
which is the Conley--Zehnder index of $\gamma$ viewed as a periodic orbit lying in $V$, and
\[
\mu_{N}^{\Phi_{N}}(\gamma)
:=\mu_{CZ}(\Phi_{N}(\gamma(t))\circ d\psi_{Tt}(\gamma(0))|_{\xi_{V}^{\perp}}\circ\Phi_{N}(\gamma(0))^{-1})
\]
which we will call the \emph{normal Conley--Zehnder index of $\gamma$ relative to $\Phi_{N}$}.
We note that basic properties of Conley--Zehnder indices which can be found, e.g.\ in
\cite[Theorem 3.1]{HWZ} show that these quantites are related by
\[
\mu^{\Phi}(\gamma)=\mu_{V}^{\Phi_{T}}(\gamma)+\mu_{N}^{\Phi_{N}}(\gamma).
\]

We now recall that a complex structure $J$ on $\xi$ is said to be \emph{compatible} with the
stable Hamiltonian structure $(\lambda, \omega)$
if the bilinear form on $\xi$ defined by
$\omega(\cdot, J\cdot)|_{\xi\times\xi}$
is symmetric and positive definite.
We will denote the set of compatible complex structures on $\xi$ by $\J(M, \xi)$.
Given a
strong stable Hamiltonian hypersurface
$V\subset M$
and a choice of $J\in\J(M, \xi)$ we will say that $J$ is \emph{$V$-compatible} if $J$ fixes the
hyperplane distribution
$\xi'=TV\cap\xi$ along $V$.
We will denote the set of such complex structures by
$\J(M, V, \xi)$.
We claim that a $J\in\J(M, V, \xi)$ necessarily also fixes
$\xi_{V}^{\perp}$.  Indeed, compatibility of $J$ with
$\omega$
implies that
$\omega(J\cdot, J\cdot)|_{\xi\times\xi}=\omega|_{\xi\times \xi}$
since we can compute
\begin{align*}
\omega(Jv, Jw)&
=\omega(w, J(Jv))
& \text{symmetry of $\omega(\cdot, J\cdot)|_{\xi\times\xi}$} \\
&=\omega(w, -v)
& J^{2}=-\mathds{1} \\
&=\omega(v, w)
\end{align*}
for any sections $v$, $w$ of $\xi$.
Thus if $v\in\xi|_{V}$ is
$\omega$-orthogonal
to every vector in $\xi'=\xi\cap TV$, then $Jv$ is also
$\omega$-orthogonal
to every vector in $\xi'$ provided that $J$ fixes $\xi'$.
Thus $J$ fixes $\xi_{V}^{\perp}$ as claimed.
We note that since both the linearized flow $d\psi_{t}$ of $R$
and a compatible $J\in\J(M, V, \xi)$ preserve the splitting
\eqref{e:splitting}, the asymptotic operator
\[
\A_\gamma h(t)=-J\left.\frac{d}{ds}\right|_{s=0}d\psi_{-Ts}h(t+s)
\]
of a periodic orbit $\gamma$ lying in $V$ also preserves the splitting.  We will write
\[
\A_\gamma=\A_\gamma^{T}\oplus\A_\gamma^{N}:W^{1,2}(\xi')\oplus W^{1,2}(\xi_{V}^{\perp})
\to L^{2}(\xi')\oplus L^{2}(\xi_{V}^{\perp})
\]
to indicate the resulting splitting of the operator.

Continuing to assume that $V\subset M$ is a strong
stable Hamiltonian hypersurface and
$J\in\J(M, V, \xi)$ is a $V$-compatible complex structure, we extend 
$J$ to an $\R$-invariant almost complex structure $\tilde J$ on $\R\times M$
in the usual way, i.e.\ so that $\tilde J$ satisfies
\begin{equation}\label{e:R-invariant-extension}
\tilde J\partial_{a}=R \qquad\text{ and }\qquad \tilde J|_{\pi^{*}\xi}=\pi^{*}J
\end{equation}
where $a$ is the coordinate along $\R$ and $\pi:\R\times M\to M$ is the canonical projection.
We note that for such an almost complex structure, the submanifold
$\R\times V$ of $\R\times M$ is $\tilde J$-holomorphic since
we assume that $R$ is tangent to $V$ and that $J$ fixes
$\xi'=\xi\cap TV$.

As with holomorphic curves, one can consider $\tilde J$-holomorphic hypersurfaces which
are asymptotic to cylindrindical $\tilde J$-holomorphic hypersurfaces of the form
$\R\times V$ with $V$ a
strong stable Hamiltonian hypersurface.
Before giving a more precise definition,
we introduce some more geometric data on our manifold.
Given a $J\in\J(M, \xi)$ we can define a Riemannian metric
\begin{equation}\label{e:metric-M}
g_{J}(v, w)=\lambda(v)\lambda(w)+\omega(\pi_{\xi}v, J\pi_{\xi}w)
\end{equation}
where
$\pi_{\xi}:TM\approx\R R\oplus\xi\to\xi$
is the projection onto $\xi$ along $R$.
We can extend $g_{J}$ to a metric $\tilde g_{J}$ on $\R\times M$
by forming the product metric with the standard metric on $\R$,
i.e.\ by defining
\begin{equation}\label{e:metric-RxM}
\tilde g_{J}:=da\otimes da+\pi^{*}g_{J}.
\end{equation}
We will denote the exponential maps of
$g_{J}$ and $\tilde g_{J}$ by $\exp$ and
$\widetilde\exp$ respectively and note that these are related by
\[
\widetilde\exp_{(a, p)}(b, v)=(a+b, \exp_{p}v).
\]
We note that if $J$ is $V$-compatible for some strong
stable Hamiltonian hypersurface $V\subset M$,
then the symplectic normal bundle $\xi_{V}^{\perp}$ is the
$g_{J}$-orthogonal complement of $TV$ in $TM|_{V}$, and that
$\pi^{*}\xi_{V}^{\perp}$ is the $\tilde g_{J}$-orthogonal complement of $T(\R\times V)$ in
$T(\R\times M)|_{\R\times V}$.
We further note that, since $V$ is assumed to be compact, the restrictions of
$\exp$ and $\widetilde\exp$ to $\xi_{V}^{\perp}$ and $\pi^{*}\xi_{V}^{\perp}$ respectively
are embeddings on some neighborhood of the zero sections.

Now consider a pair $V_{+}$, $V_{-}$ of
strong stable Hamiltonian hypersurfaces
and assume that
$V:=V_{+}\cup V_{-}$ is also a
strong stable Hamiltonian hypersurface.
We let $J\in\J(M, V, \xi)$ be a $V$-compatible $J$ with associated $\R$-invariant almost complex structure
$\tilde J$ on $\R\times M$.
We are interested in $\tilde J$-holomorphic submanifolds which outside of a compact set can be described
by exponentially decaying sections of the normal bundles to $V_{+}$ and $V_{-}$,  More precisely,
we say that a $\tilde J$-holomorphic submanifold $\tilde V\subset\R\times M$ is
\emph{positively asymptotically cylindrical over $V_{+}$} and
\emph{negatively asymptotically cylindrical over $V_{-}$} if there exists an $R>0$ and sections
\begin{gather*}
\eta_{+}:[R, +\infty)\to C^\infty(\xi_{V_{+}}^{\perp}) \\
\eta_{-}:(-\infty, -R] \to C^\infty(\xi_{V_{-}}^{\perp}) \\
\end{gather*}
so that
\begin{gather*}
\tilde V\cap\left([R, +\infty)\times M\right)=\bigcup_{(a, p)\in [R, +\infty)\times V_{+}} \widetilde\exp_{(a, p)}\eta_{+}(a, p) \\
\tilde V\cap\left((-\infty, -R]\times M\right)=\bigcup_{(a, p)\in (-\infty, -R]\times V_{-}} \widetilde\exp_{(a, p)}\eta_{-}(a, p)
\end{gather*}
and so that there exist constants $M_{i}>0$, $d>0$ satisfying
\[
\abs{\widetilde\nabla^{i}\eta_{\pm}(a, p)}\le M_{i}e^{-d \abs{a}}
\]
for all $i\in\N$ and $\pm a\in [R, +\infty)$, where
$\tilde\nabla$ is the extension of a connection $\nabla$ on  $\xi_{V}^{\perp}$
to a connection $\tilde \nabla$ on $\pi^{*}\xi_{V}^{\perp}$ defined by requiring
$\tilde\nabla_{\partial_{a}}\eta(a, p)=\partial_{a}\eta(a, p)$.
We will refer to the sections $\eta_{+}$ and $\eta_{-}$ respectively as \emph{positive}
and \emph{negative asymptotic representatives of $\tilde V.$}

Our main goal here is to understand the intersection properties of punctured
pseudoholomorphic curves with
asymptotically cylindrical pseudoholomorphic hypersurfaces.  The main difficulty 
arises from the noncompactness of the manifolds in question.
Indeed a punctured pseudoholomorphic curve whose image is not contained in the
$\tilde J$-holomorphic hypersurface $\tilde V$ may have punctures limiting to a periodic orbits lying in
$V_{+}$ or $V_{-}$.
In this case, it's not a priori clear that the intersection number between the curve and the hypersurface is finite.
Even assuming this intersection number is finite, it is not homotopy invariant as intersections can be lost or created
at infinity.
We will see below that these difficulties can be dealt with via higher-dimensional analogs of
techniques developed in \cite{Sie11}.

Before presenting the relevant results it will be convenient
to establish some standard assumptions and notations for the next several
definitions and results.

\begin{assumptions}\label{standing-assumptions}
We assume that:
\begin{enumerate}\renewcommand{\theenumi}{\alph{enumi}}
\item $(M, \lambda, \omega)$ is a closed, manifold with nondegenerate stable Hamiltonian structure
$(\lambda, \omega)$ with $\xi=\ker\lambda$ and $R$ the associated Reeb vector field,

\item $V_{+}\subset M$, $V_{-}\subset M$, and $V=V_{+}\cup V_{-}$ are
strong stable Hamiltonian hypersurfaces
of $M$,

\item $J\in\J(M, V, \xi)$ is a $V$-compatible complex structure on $\xi$ and
$\tilde J$ is the $\R$-invariant almost complex structure on $\R\times M$ associated to $J$ (defined by
\eqref{e:R-invariant-extension},

\item $\tilde g_{J}$ is the Riemannian metric on $\R\times M$ defined by \eqref{e:metric-RxM} and
$\widetilde\exp$ is the associated exponential map,

\item $\tilde V\subset\R\times M$ is a $\tilde J$-holomorphic hypersurface which is positively asymptotically cylindrical over $V_{+}$ and negatively asymptotically cylindrical over $V_{-}$,

\item $\eta_{+}:[R, +\infty)\to C^\infty (\xi_{V_{+}}^{\perp})$ and
$\eta_{-}:(-\infty, -R]\to C^{\infty}(\xi_{V_{-}}^{\perp})$ are, respectively, positive and negative asymptotic representatives of $\tilde V$,

\item $C=[\Sigma, j, \Gamma=\Gamma^{+}\cup\Gamma^{-}, \widetilde{u}=(a, u)]$ is a finite-energy $\tilde J$-holomorphic curve, and
at $z\in\Gamma$, $C$ is asymptotic to $\gamma_{z}^{m_{z}}$
(with $\gamma_{z}^{m_{z}}$ indicating the $m_{z}$-fold covering of a simple periodic orbit $\gamma_{z}$),

\item $\Phi$ is a trivialization of $\xi_{V}^{\perp}$ along every periodic orbit lying in $V$
which occurs as an asymptotic limit of $C$. 

\end{enumerate}
\end{assumptions}

The following theorem can be seen as a generalization of
Theorem 2.2 in \cite{Sie08}.  Proof will be given in \cite{Sie}.

\begin{thm}\label{t:normal-asymptotics}
Assume \ref{standing-assumptions} and 
assume that at $z\in\Gamma^{+}$, $C$ is asymptotic to a periodic orbit
$\gamma_{z}^{m_{z}}\subset V_{+}$.
Then there exists an $R'\in\R$, a smooth map
\[
u_{T}:[R', \infty)\times S^{1}\to [R, \infty)\times V_{+}
\]
and a smooth section
\[
u_{N}:[R', \infty)\times S^{1}\to u_{T}^{*}\pi^{*}\xi_{V_{+}}^{\perp}
\]
so that the map
\begin{equation}
(s, t)\mapsto\widetilde\exp_{u_{T}(s, t)}u_{N}(s, t)
\end{equation}
parametrizes $C$ near $z$.
Moreover,
if we assume that the image of $C$ is not a subset of the asymptotically cylindrical hypersurface $\tilde V$,
then
\begin{equation}\label{e:normal-asymptotics}
u_{N}(s, t)-\eta_{+}(u_{T}(s, t))=e^{\mu s}[e(t)+r(s, t)]
\end{equation}
for all $(s, t)\in[R', \infty)\times S^{1}$
where:
\begin{itemize}
\item $\mu<0$ is a negative eigenvalue of the normal asymptotic operator
$\A_{\gamma_{z}^{m_{z}}}^{N}$,
\item $e\in\ker(\mathbf{A}_{\gamma_{z}^{m_{z}}}^{N}-\mu)\setminus\br{0}$ is an eigenvector with eigenvalue $\mu$, and
\item $r:[R', \infty)\times S^1 \to u_{T}^{*}\pi^{*}\xi_{V}^{\perp}$ is a smooth section satisfying
exponential decay estimates of the form
\begin{equation}\label{e:normal-asymptotics-remainder}
\abs{\tilde\nabla^{i}_{s}\tilde\nabla^{j}_{t} r(s, t)}\le M_{ij}e^{-d\abs{s}}
\end{equation}
for some positive contants $M_{ij}$, $d$ and all $(i, j)\in\N^{2}$
\end{itemize}

Similarly, if we 
assume that at $z\in\Gamma^{-}$, $C$ is asymptotic to a periodic orbit
$\gamma_{z}^{m_{z}}\subset V_{-}$., then there exists an $R'\in\R$, a smooth map
\[
u_{T}:(-\infty, R'] \times S^{1}\to (-\infty, -R]\times V_{-}
\]
and a smooth section
\[
u_{N}:(-\infty, R'] \times S^{1}\to u_{T}^{*}\pi^{*}\xi_{V_{-}}^{\perp}
\]
so that the map
\[
(s, t)\mapsto\widetilde\exp_{u_{T}(s, t)}u_{N}(s, t)
\]
parametrizes $C$ near $z$.
Moreover, if the image of $C$ is not contained in $\tilde V$, then 
$u_{N}(s, t)-\eta_{-}(u_{T}(s, t))$ satisfies a formula of the form \eqref{e:normal-asymptotics}
for all $(s, t)\in(-\infty, R']\times S^{1}$, where now:
\begin{itemize}
\item $\mu>0$ is a positive eigenvalue of the normal asymptotic operator
$\A_{\gamma_{z}^{m_{z}}}^{N}$,
\item $e\in\ker(\A_{\gamma_{z}^{m_{z}}}^{N}-\mu)\setminus\br{0}$, as before,
is an eigenvector with eigenvalue $\mu$, and
\item $r:(-\infty, R']\to u_{T}^{*}\pi^{*}\xi_{V}^{\perp}$ is a smooth section satisfying
exponential decay estimates of the form \eqref{e:normal-asymptotics-remainder}
for some positive contants $M_{ij}$, $d$.
\end{itemize}
\end{thm}

We note that the bundles of the form
$u_{T}^{*}\pi^{*}\xi_{V}^{\perp}$
ocurring in the statement of this theorem are trivializable
since they are complex line bundles over a space which retracts onto $S^{1}$.
In any trivialization the eigenvector $e$ from formula \eqref{e:normal-asymptotics}
satisfies a linear, nonsingular ODE,
and thus is nowhere vanishing since we assume it is not identically zero.
Since the ``remainder term'' $r$ in the formula \eqref{e:normal-asymptotics}
converges to zero, we conclude that the functions
$u_{N}(s, t)-\eta_{\pm}(u_{T}(s, t))$ are nonvanishing for sufficiently large $\abs{s}$.
However, since zeroes of this function can be seen to correspond to intersections between
the curve $C$ and the hypersurface $\tilde V$ occuring sufficiently close to the punctures of $C$, we conclude
that all intersections between $C$ and $\tilde V$ are contained in a compact set.
Moreover, since intersections between $C$ and $\tilde V$ can be shown to be isolated and of positive local order
(see e.g.\  \cite[Lemma 3.4]{IP}),
we conclude that the algebraic intersection number between $C$ and $\tilde V$ is finite:

\begin{cor}\label{c:finite-intersections}
Assume \ref{standing-assumptions} and assume that
no component of  the curve $C$
has image contained in 
the $\tilde J$-holomorphic hypersurface $\tilde V$.  Then the algebraic intersection number
$C\cdot \tilde V$, defined by summing local intersection indices, is finite and nonnegative, and
$C\cdot \tilde V=0$ precisely when $C$ and $\tilde V$ do not intersect.
\end{cor}

This corollary deals with the first difficulty in understanding intersections between punctured curves and asymptotically cylindrical hypersurfaces described above, namely the finiteness of the intersection number. A second consequence of the asymptotic formula from Theorem \ref{t:normal-asymptotics}, again stemming
from the fact that the quantities $u_{N}(s, t)-\eta_{\pm}(u_{T}(s, t))$ are nonzero for sufficiently large
$\abs{s}$, is that the normal approach of the curve $C$ has a well-defined winding number
relative to a trivialization $\Phi$ of $\xi_{V}^{\perp}|_{\gamma_{z}}$.  This winding will be given by
the winding of the eigenvector from formula \eqref{e:normal-asymptotics} relative to $\Phi$, and for a given 
puncture $z$ of $C$, we will denote this quantity by
\[
\wind_{rel}^\Phi((C; z), \tilde V)=\wind(e).
\]
Combining this observation with the characterization of the Conley--Zehnder
index in terms of the asymptotic operator from \cite[Definition 3.9/Theorem 3.10]{HWZ} leads
to the following corollary.

\begin{cor}\label{c:normal-rel-winding}
Assume \ref{standing-assumptions}, and assume that
no component of the curve
$C=[\Sigma, j, \Gamma_{+}\cup\Gamma_{-}, \tilde u=(a, u)]$
has image contained in
the holomorphic hypersurface $\tilde V$.  Then:
\begin{itemize}
\item If $z\in\Gamma$ is a positive puncture at which $\tilde u$ limits to $\gamma_z^{m_z}\subset V_{+}$ then
\begin{equation}\label{e:wind-rel-pos}
\wind_{rel}^{\Phi}((C; z), \tilde V)\le \fl{\mu_{N}^{\Phi}(\gamma_z^{m_z})/2}=:\alpha_N^{\Phi;-}(\gamma_z^{m_z}).
\end{equation}

\item If $z\in\Gamma$ is a negative puncture at which $\tilde u$ limits to $\gamma_z^{m_z}\subset V_{-}$ then
\begin{equation}\label{e:wind-rel-neg}
\wind_{rel}^{\Phi}((C; z), \tilde V)\ge \ceil{\mu_{N}^{\Phi}(\gamma_z^{m_z})/2}=:\alpha_N^{\Phi;+}(\gamma_z^{m_z}).
\end{equation}
\end{itemize}
\end{cor}

The numbers $\alpha_N^{\Phi;-}(\gamma)$ and $\alpha_N^{\Phi;+}(\gamma)$ are, respectively, the biggest/smallest winding number achieved by an eigenfunction of the normal asymptotic operator of any orbit $\gamma$ corresponding to a negative/positive eigenvalue. Observe that we have the formulas 

\begin{equation}
\label{alphas}
\mu_{N}^{\Phi}(\gamma)= 2\alpha_N^{\Phi;-}(\gamma)+p_N(\gamma)=2\alpha_N^{\Phi;+}(\gamma)-p_N(\gamma),
\end{equation}
where $p_N(\gamma) \in \{0,1\}$ is the \emph{normal parity} of the orbit $\gamma$ (which is independent of the trivialization $\Phi$).

We will see in a moment that this corollary
can be used to 
deal with the second difficulty in understanding intersections between 
punctured curves and asymptotically cylindrical hypersurfaces described above, namely, the fact
that the algebraic intersection number may not be invariant under homotopies.
We first introduce some terminology.
Assuming again \ref{standing-assumptions} and that
no component of
$C$ is a subset of $\tilde V$, we define the asymptotic
intersection at the punctures of $C$ in the following way:
\begin{itemize}
\item
If for the positive puncture $z\in\Gamma_{+}$, $\gamma^{m_z}_{z}\subset V_{+}$, we define the 
\emph{asymptotic intersection number $\delta_{\infty}((C; z); \tilde V)$ of $C$ at $z$ with $\tilde V$} by
\begin{equation}\label{e:local-asymp-inum-pos}
\delta_{\infty}((C; z), \tilde V)=\fl{\mu_{N}^{\Phi}(\gamma_z^{m_z})/2}-\wind_{rel}^{\Phi}((C; z), \tilde V).
\end{equation}

\item 
If for the negative puncture $z\in\Gamma_{-}$, $\gamma^{m_z}_{z}\subset V_{-}$, we define the 
\emph{asymptotic intersection number $\delta_{\infty}((C; z); \tilde V)$ of $C$ at $z$ with $\tilde V$} by
\begin{equation}\label{e:local-asymp-inum-neg}
\delta_{\infty}((C; z), \tilde V)=\wind_{rel}^{\Phi}((C; z), \tilde V)-\ceil{\mu_{N}^{\Phi}(\gamma^{m_z}_{z})/2}.
\end{equation}

\item For all other punctures $z\in\Gamma_{\pm}$
(i.e.\ those for which $\gamma_{z}$ is not contained in $V_{\pm}$),
we define
\begin{equation}\label{e:local-asymp-inum-zero}
\delta_\infty((C; z), \tilde V)=0.
\end{equation}
\end{itemize}
We then define the \emph{total asymptotic intersection number of $C$ with $\tilde V$} by
\begin{equation}\label{e:global-asymp-inum}
\delta_{\infty}(C, \tilde V)=\sum_{z\in\Gamma}\delta_{\infty}((C; z), \tilde V).
\end{equation}

We observe that as a result of Corollary \ref{c:normal-rel-winding} the local and total asymptotic intersection numbers are always nonnegative.

\vspace{0.5cm}

Continuing to assume \ref{standing-assumptions} (but no longer necessarily that no component of
the curve $C$ has image contained in $\tilde V$),
we can use the trivialization $\Phi$ of $\xi_{V}^{\perp}$ along the asymptotic periodic orbits of $C$ lying
in $V$ to construct a perturbation $C_{\Phi}$ of $C$ in the following way.
For each puncture $z\in\Gamma$ for which we the asymptotic limit $\gamma_{z}^{m_{z}}$ lies in $V$,
we first extend $\Phi$ to a trivialization
$\Phi:\xi_{V}^{\perp}|_{U_{z}}\to U_{z}\times\R^{2}$ on some open neighborhood 
$U_{z}\subset V$ of the asymptotic limit $\gamma_{z}$.  Then we consider the asymptotic parametrization
\[
(s, t)\mapsto \widetilde\exp_{u_{T}(s, t)}u_{N}(s, t)
\]
from Theorem \ref{t:normal-asymptotics} above
for $(s, t)\in [R,+\infty)\times S^{1}$ or $(-\infty, -R]\times S^{1}$ as appropriate,
where $R>0$ is chosen large enough so that $u_{T}$ has image contained in the neighborhood $U_{z}$ of
$\gamma_{z}$ on which the trivialization $\Phi$ has been extended.
We then perturb the map by replacing the above parametrization of $C$ near $z$ by the map
\[
(s, t)\mapsto \widetilde\exp_{u_{T}(s, t)}\bp{u_{N}(s, t)+\beta(\abs{s})\Phi(u_{T}(s, t))^{-1}\varepsilon}
\]
where $\beta:[0, \infty)\to[0, 1]$ is a smooth cut-off function equal to $0$ for $s<\abs{R}+1$ and equal to
$1$ for $\abs{s}>\abs{R}+2$, and $\varepsilon \ne 0$ is thought of as a number in $\C\approx\R^{2}$.
Given this, we can then define the 
\emph{relative intersection number $i^{\Phi}(C, V)$
of $C$ and $\tilde V$ relative to the the trivialization $\Phi$}
by
\[
i^{\Phi}(C, V):=C_{\Phi}\cdot V.
\]
It can be shown that this number is independent of choices made in the construction of $C_{\Phi}$ provided
the perturbations are sufficiently small.

Again continuing to assume \ref{standing-assumptions},
we now define the
\emph{holomorphic intersection product of $C$ and $\tilde V$} by
\[
C*V:=i^{\Phi}(C, V)
+\sum_{\stackrel{z\in\Gamma_{+}}{\gamma_{z}\subset V_{+}}}\fl{\mu_{N}^{\Phi}(\gamma_{z}^{m_{z}})/2}
-\sum_{\stackrel{z\in\Gamma_{-}}{\gamma_{z}\subset V_{-}}}\ceil{\mu_{N}^{\Phi}(\gamma_{z}^{m_{z}})/2}
\]
The key facts about the holomorphic intersection product are now given in the following theorem which generalizes
\cite[Theorem 2.2/4.4]{Sie11}.
\begin{thm}[Generalized positivity of intersections]
With $\tilde V$ and $C$ as in \ref{standing-assumptions}, assume that $C$ is not contained in $\tilde V$,
the holomorphic intersection product $C*\tilde V$ depends only on the relative homotopy classes of $C$ and $\tilde V$.
Moreover, if the image of $C$ is not contained in $\tilde V$, then
\[
C*\tilde V=C\cdot \tilde V+\delta_{\infty}(C, \tilde V)\ge 0
\]
where $C\cdot \tilde V$ is the algebraic intersection number, defined by summing local intersection indices,
and $\delta_{\infty}(C, \tilde V)$ is the total asymptotic intersection number, defined by
\eqref{e:local-asymp-inum-pos}--\eqref{e:global-asymp-inum}.
In particular, $C*\tilde V\ge 0$ and equals zero if and only if $C$ and $\tilde V$ don't intersect and
all asymptotic intersection numbers are zero.
\end{thm}
The proof of this theorem follows along very similar lines to Theorem 2.2/4.4 in \cite{Sie11}.
The essential point is that the relative intersection number can be shown to be given by the formula
\[
i^{\Phi}(C, V)=
C\cdot \tilde V
-\sum_{\stackrel{z\in\Gamma_{+}}{\gamma_{z}\subset V_{+}}}\wind^\Phi_{rel}((C; z), \tilde V)
+\sum_{\stackrel{z\in\Gamma_{-}}{\gamma_{z}\subset V_{-}}}\wind^\Phi_{rel}((C; z), \tilde V).
\]
The result will then follow from Corollary \ref{c:normal-rel-winding} above. Detailed proof will be given in \cite{Sie}.

Analogous to the case in four dimensions studied in \cite{Sie11}, the
$\R$-invariance of the almost complex structure in the set-up here allows one to compute
the holomorphic intersection number and (in some cases) the algebraic intersection
number of a holomorphic curve and a holomorphic hypersurface with respect
to asymptotic winding numbers and intersections of each object with the asymptotic limits
of the other.

Before stating the relevant results we will first make some additional assumptions. We will henceforth assume that:
\begin{assumptions}\label{assumptions-2}$\;$
\begin{enumerate}\renewcommand{\theenumi}{\alph{enumi}}
\item $V_{+}$ and $V_{-}$ are disjoint,
\item $\xi_{V}^{\perp}$ of $V=V_{+}\cup V_{-}$ is trivializable, 
\item $\Phi:\xi_{V}^{\perp}\to V\times \R^{2}$ is a global trivialization,
\item $\tilde V$ is connected, and
\item the projection $\pi(\tilde V)$ of $\tilde V$ to $M$ is an embedded codimension-$1$ submanifold of
$M\setminus V$.
\end{enumerate}
\end{assumptions}
Under these assumptions, $\tilde V$ has a well-defined winding
$\wind_{\infty}^{\Phi}(\tilde V, \gamma)$ relative to $\Phi$
around any orbit $\gamma\subset V=V_{+}\cup V_{-}$
which can be defined by considering the asymptotic representatives $\eta_{+}$ or $\eta_{-}$ as appropriate
and computing
\[
\wind_{\infty}^{\Phi}(\tilde V, \gamma)=\lim_{\abs{s}\to\infty}\wind \Phi^{-1}\eta_{\pm}(s, \gamma(\cdot)),
\]
or, equivalently by
\[
\wind_{\infty}^{\Phi}(\tilde V, \gamma)=\wind_{rel}^{\Phi}((\R\times\gamma; \pm\infty), \tilde V).
\]
As in Corollary \ref{c:normal-rel-winding} above, it follows from the
asymptotic formula from Theorem \ref{t:normal-asymptotics} above and the characterization of the
Conley--Zehnder index from \cite{HWZ} that
\[
\wind_{\infty}^{\Phi}(\tilde V, \gamma)\le \fl{\mu_{N}^{\Phi}(\gamma)/2}=\alpha_N^{\Phi;-}(\gamma)
\]
if $\gamma\subset V_{+}$ and
\[
\wind_{\infty}^{\Phi}(\tilde V, \gamma)\ge \ceil{\mu_{N}^{\Phi}(\gamma)/2}=\alpha_N^{\Phi;+}(\gamma)
\]
if $\gamma\subset V_{-}$.

The following theorem, which can be seen as the higher dimensional version of
\cite[Corollary 5.11]{Sie11}, gives a computation of the algebraic intersection number of
$\tilde V$ and $\R$-shifts of the curve $C$ in terms of the asymptotic data, and the intersections of
each object with the asymptotic limits of the other.

\begin{thm}\label{t:intersection-computation}
Assume \ref{standing-assumptions} and \ref{assumptions-2} and that the curve $C$ is connected and
not equal to an orbit cylinder and not contained in $\R\times V$.
For $c\in\R$, denote by $C_{c}$ the curve obtained from translating $C$ in the $\R$-coordinate by $c$.
Then for all but a finite number of value of $c\in\R$, the algebraic intersection number
$C_{c}\cdot \tilde V$ is given by the formulas:

\begin{equation*}
\begin{split}
C_{c}\cdot \tilde V
&=C\cdot (\R\times V_{+}) \\
&\hskip.25in
+\sum_{\stackrel{z\in\Gamma_{+}}{\gamma_{z}\in V_{+}}} 
\bp{\max\br{m_{z}\wind_{\infty}^{\Phi}(\tilde V, \gamma_{z}), \wind^{\Phi}_{rel}((C;z), \R\times V_{+})}-\wind^{\Phi}_{rel}((C;z), \R\times V_{+})} \\
&\hskip.25in
+ \sum_{z\in\Gamma_{-}}m_{z}(\R\times \gamma_{z})\cdot \tilde V  \\
&\hskip.25in
+\sum_{\stackrel{z\in\Gamma_{-}}{\gamma_{z}\in V_{-}}}
\bp{m_{z}\wind_{\infty}^{\Phi}(\tilde V, \gamma_{z})
-
\min\br{m_{z}\wind_{\infty}^{\Phi}(\tilde V, \gamma_{z}), \wind^{\Phi}_{rel}((C;z), \R\times V_{-})}} \\
&\hskip.25in
+\sum_{\stackrel{z\in\Gamma_{-}}{\gamma_{z}\in V_{+}}}
\bp{\wind_{rel}^{\Phi}((C;z), \R\times V_{+})-m_{z}\wind_{rel}^{\Phi}(\tilde V, \gamma_{z})} \\
\end{split}
\end{equation*}
\begin{equation*}
\begin{split}
&=\sum_{z\in\Gamma_{+}}m_{z}(\R\times \gamma_{z})\cdot \tilde V  \\
&\hskip.25in
+\sum_{\stackrel{z\in\Gamma_{+}}{\gamma_{z}\in V_{+}}}
\bp{\max\br{m_{z}\wind_{\infty}^{\Phi}(\tilde V, \gamma_{z}), \wind_{rel}^{\Phi}((C; z), \R\times V_{+})}-m_{z}\wind_{\infty}^{\Phi}(\tilde V, \gamma_{z})} \\
&\hskip.25in
+C\cdot (\R\times V_{-}) \\
&\hskip.25in
+\sum_{\stackrel{z\in\Gamma_{-}}{\gamma_{z}\in V_{-}}}
\bp{\wind_{rel}^{\Phi}((C;z), \R\times V_{-})-\min\br{m_{z}\wind_{\infty}^{\Phi}(\tilde V, \gamma_{z}),\wind_{rel}^{\Phi}((C;z), \R\times V_{-})}} \\
&\hskip.25in
+\sum_{\stackrel{z\in\Gamma_{+}}{\gamma_{z}\in V_{-}}}
\bp{m_{z}\wind_{\infty}^{\Phi}(\tilde V, \gamma_{z})-\wind_{rel}^{\Phi}((C;z), \R\times V_{-})}
\end{split}
\end{equation*}

with each of the grouped terms always nonnegative.
\end{thm}

The nonnegativity of the terms in the above formulas allows us to establish a
convenient set of conditions which will guarantee that the  projections 
$\pi(C)$ and $\pi(\tilde V)$ of the curve and hypersurface to $M$ do not intersect.
Indeed, if $\pi(C)$ and $\pi(\tilde V)$ are disjoint then
$C_{c}$ and $\tilde V$ are disjoint for all values of $c\in\R$ and hence
the algebraic intersection number of $C_{c}$ and $\tilde V$ is zero for all values of $c\in\R$.
Since the formulas from Theorem \ref{t:intersection-computation} compute this number (for all but a finite
number of values of $c\in\R$) in terms of nonnegative quantities, we can conclude that all terms in the above formulas vanish.
We thus obtain the following corollary, which generalizes \cite[Theorem 2.4/5.12]{Sie11}.

\begin{cor}\label{c:nonintersection-conditions}
Assume that all of the hypotheses of Theorem \ref{t:intersection-computation}
hold and that $\pi(C)$ is not contained in $\pi(\tilde V)$.
Then the following are equivalent:

\begin{enumerate}
\item $\pi(C)$ and $\pi(\tilde V)$ are disjoint.

\item  All of the following hold:
\begin{enumerate}
\item
None of the asymptotic limits of $C$ intersect $\pi(\tilde V)$.

\item
$\pi(C)$ does not intersect $V=V_{+}\cup V_{-}$.

\item
For any puncture $z$ at which $C$ has asymptotic limit $\gamma^{m}$ lying in $V=V_{+}\cup V_{-}$,
$\wind^{\Phi}_{rel}((C; z), V)=m \wind_{\infty}^{\Phi}(\tilde V, \gamma)$.
\end{enumerate}
\end{enumerate}
\end{cor}

\subsection{An application to holomorphic foliations}\label{app}
Here we apply the above theory to prove a 
result about the behavior of pseudoholomorphic curves in the presence
of certain foliations by asymptotically cylindrical holomorphic hypersurfaces.
This result is used in the proofs of Proposition \ref{uniqhyp} and the uniqueness Lemma \ref{uniqpre}. 

We consider a manifold $M^{2n+1}$
equipped with a nondegenerate stable Hamiltonian structure $(\lambda, \omega)$,
and we will denote by the triple $(\Sigma, p, Y^{2n-1})$
a smooth fibration $p:M\to\Sigma$
over a closed, oriented surface $\Sigma$ with fiber diffeomorphic to $Y$.
Given a Morse function $f$ on $\Sigma$ we say the fibration
$(\Sigma, p, Y)$ is \emph{$f$-admissible}
if for each critical point
$w\in\crit(f)$ of the function, the fiber
$Y_{w}=p^{-1}(w)$ over $w$
is a stable Hamiltonian hypersurface.
We will denote $f$-admissible fibrations by quadruples
$(\Sigma, p, Y, f)$.

Given an $f$-admissible fibration $(\Sigma, p, Y, f)$ for $M$ and a critical point
$w\in\crit(f)$, we note that at any point $y\in Y_{w}$,
the derivative
$p_{*}$ at $y$ determines a linear isomorphism
\[
p_{*}(y):(\xi_{Y_{w}}^{\perp})_{y}\to T\Sigma_{w}
\]
from the symplectic normal bundle of $Y_{w}$ in $M$ at $y\in Y_{w}$
to the tangent space of $\Sigma$ at $w$.
This map will either be orientation preserving at every point in $Y_{w}$ or orientation 
reversing at every point in $Y_{w}$.
This allows us to define a sign function
\[
\sign:\crit(f)\to\br{-1, 1}
\]
on the set of critical points of $f$ by requiring
\[
\sign(w)=
\begin{cases}
1 & \text{ if $p_{*}|_{\xi_{Y_{w}}^{\perp}}$ is everywhere orientation preserving} \\
-1 &\text{ if $p_{*}|_{\xi_{Y_{w}}^{\perp}}$ is everywhere orientation reversing}.
\end{cases}
\]
Choosing at each $w\in\crit(f)$ with $\sign(w)=1$ an 
orientation preserving linear map $\Phi_{w}:T\Sigma_{w}\to (\R^{2}, \omega_{0})$
and 
at each $w\in\crit(f)$ with $\sign(w)=-1$ an 
orientation reversing linear map $\Phi_{w}:T\Sigma_{w}\to (\R^{2}, \omega_{0})$
we obtain a global orientation-preserving trivialization
\[
\Phi_{w}\circ p_{*}|_{\xi_{Y_{w}}^{\perp}}\to(\R^{2}, \omega_{0})
\]
of $\xi_{Y_{w}}^{\perp}$ which can be homotoped to a symplectic trivialization.
Thus in an $f$-admissible fibration there is a preferred homotopy class of
symplectic trivialization of the symplectic normal bundle
to $p^{-1}(\crit (f))$.

Assuming still that $M$ is a stable Hamiltonian manifold, 
we assume $\tilde J$ is a compatible almost complex structure
on $\R\times M$, and we let $\pi:\R\times M\to M$ denote the canonical projection.
We say that a codimension-$2$ foliation $\F$ of $\R\times M$ is an
\emph{$\R$-invariant,  asymptotically cylindrical, $\tilde J$-holomorphic foliation} if:
\begin{itemize}
\item $\mathcal{F}$ is invariant under translations in the $\R$-coordinate,
\end{itemize}
and if there exists a strong stable Hamiltonian hypersurface $V\subset M$
so that:
\begin{itemize}
\item $\R\times V$ has $\tilde J$-invariant tangent space,

\item the union of leaves of $\mathcal{F}$ fixed by $\R$-translation
is a equal to $\R\times V\subset\R\times M$, and

\item all other leaves of the foliation are
$\tilde J$-holomorphic hypersurfaces which are asymptotically cylindrical over some collection of
components of $V$,
and which project by $\pi$ to embedded submanifolds smoothly foliating $M\setminus V$.
\end{itemize}
We will refer to the stable Hamiltonian hypersurface $V$ as the
\emph{binding set} of the foliation $\F$.
For brevity, we will refer to $\R$-invariant, asymptotically cylindrical, $\tilde J$-holomorphic
foliations simply as \emph{holomorphic foliations} from now on.

Now assuming the manifold $M$ is equipped with both
a holomorphic foliation 
$\F$ with binding $V$ and an $f$-admissible fibration $(\Sigma, p, Y, f)$,
we will say that \emph{$\F$ is compatible with $(\Sigma, p, Y, f)$}
if:
\begin{itemize}
\item $p^{-1}(\crit(f))$ is equal to the binding $V$ of the foliation,

\item all other leaves of the foliation are diffeomorphic to $\R\times Y$ and admit smooth parametrizations of
the form
\[
(s, y)\in \R\times Y\mapsto (a(s, y), m(s, y))\in\R\times M
\]
where $p(m(s, y))=\gamma(s)$ for some solution $\gamma$ to the gradient flow equation 
\[
\dot\gamma(s)=\nabla f(\gamma(s))
\]
with respect to an appropriate metric $g_{\Sigma}$ on $\Sigma$, and where
$a(s, y)$ satisfies
\[
\lim_{s\to\pm\infty}a(s, y)=\sign(\lim_{s\to\pm\infty}p(m(s, y)))\pm\infty.
\]
\end{itemize}
Examples of such foliations are
constructed in Sections \ref{holhyps} and \ref{prequantcurves}.
We observe that the assumption that leaves of the foliation project to gradient flow lines implies that
with respect to a trivialization $\Phi$ of
$\xi_{V}^{\perp}$ in the preferred homotopy class determined by the fibration, the windings
$\wind^{\Phi}_{\infty}(\tilde V, \gamma)$ vanish for  any leaf $\tilde V$ of the foliation not fixed by the $\R$-action
and any periodic orbit $\gamma$ lying in the one of the asymptotic limits of $\tilde V$.

For the remainder of the section we will assume that $M$ is a manifold equipped with a
nondegenerate stable Hamiltonian structure $(\lambda, \omega)$,
a compatible almost complex structure $\tilde J$ on $\R\times M$,
an $f$-admissible fibration $(\Sigma, p, Y, f)$, and
a holomorphic foliation $\F$ with binding $V$
which is compatible with fibration $(\Sigma, p, Y, f)$.
We will denote by $\Phi$ a trivialization of the symplectic normal bundles
$\xi_{V}^{\perp}$ to $V=p^{-1}(\crit(f))$ in the preferred homotopy class of trivializations determined by
the fibration.

We are interested in understanding the behavior of
punctured pseudoholomorphic curves whose asymptotic limits are all contained in the bindings
of the foliation.
We will let
$\tilde u=(a, u):\Sigma'\setminus\Gamma^{+}\cup\Gamma^{-}\to\R\times M$
be a punctured pseudoholomorphic map and assume that at a puncture
$z\in\Gamma=\Gamma^{+}\cup\Gamma^{-}$
the map $u$ is asymptotic to a periodic orbit $\gamma_{z}^{m_{z}}$
lying in $Y_{w_{z}}=p^{-1}(w_{z})$ for some $w_{z}\in\crit(f)$.
We observe that the assumption that all asymptotic limits
of $u$ lie in $V=p^{-1}(\crit(f))$ implies that the map
$v=p\circ u:\Sigma'\setminus\Gamma\to\Sigma$
obtained by projecting $u$ to the base approaches a critical point of
$f$ as one approaches each puncture $z\in\Gamma$.  The map $v$ thus has a unique
continuous extension to a map $\bar v:\Sigma'\to\Sigma$.
The following lemma computes the degree of the map $\bar v$
in terms of the intersections
and relative normal windings of the map $\tilde u$ with any of the
stable Hamiltonian hypersurfaces $Y_{w}$.

\begin{lemma}\label{l:degree-computation}
Assume the map $u:\Sigma'\setminus\Gamma\to M$ does not have image
contained in $V=p^{-1}(\crit(f))$.  Then,
given a point $w\in\crit(f)$, the degree of the map
$\bar v$ defined above is given by the formula
\[
\deg(\bar v)=
\sign(w)
\bp{
\tilde u\cdot (\R\times Y_{w})
-\sum_{\stackrel{z\in\Gamma^{+}}{\gamma_{z}\subset Y_{w}}}\wind^{\Phi}_{rel}((\tilde u; z), Y_{w})
+\sum_{\stackrel{z\in\Gamma^{-}}{\gamma_{z}\subset Y_{w}}}\wind^{\Phi}_{rel}((\tilde u; z), Y_{w})
}
\]
and thus satisfies
\begin{equation}\label{e:degree-inequality}
\begin{aligned}
\sign(w)\deg(\bar v)
&\ge
-\sum_{\stackrel{z\in\Gamma^{+}}{\gamma_{z}\subset Y_{w}}}\fl{\mu_{N}^{\Phi}(\gamma_{z}^{m_{z}})/2}
+\sum_{\stackrel{z\in\Gamma^{-}}{\gamma_{z}\subset Y_{w}}}\ceil{\mu_{N}^{\Phi}(\gamma_{z}^{m_{z}})/2} \\
&=
\sum_{\stackrel{z\in\Gamma^{\pm}}{\gamma_{z}\subset Y_{w}}}\mp\alpha_{N}^{\Phi; \mp}(\gamma_{z}^{m_{z}})
\end{aligned}
\end{equation}
with equality occurring if and only if $u$ does not intersect
$Y_{w}$ and
$\wind_{rel}^{\Phi}((\tilde u; z), Y_{w})=\alpha_{N}^{\Phi; \mp}(\gamma_{z}^{m_{z}})$ for
every $z\in\Gamma^{\pm}$ with $\gamma_{z}\subset Y_{w}$.
\end{lemma}

\begin{proof}
Given a point $w\in\crit(f)$ it's clear from the definition of the map
$\bar v$ that points $z\in\Sigma'$ with $\bar v(z)=w$ coincide with intersection points of the map
$\tilde u$ with $\R\times Y_{w}$ (or equivalently, intersection points of
$u$ with $Y_{w}$) and with punctures $z\in\Gamma$ for which the periodic orbit
$\gamma_{z}^{m_{z}}$ is contained in $Y_{w}$.
Since we've observed in Corollary \ref{c:finite-intersections} this number is finite,
we can compute the degree of the map
$\bar v:\Sigma'\to\Sigma$ by summing the local degree of the map
at each point in $\bar v^{-1}(w)$.

Assume for the moment that $\sign(w)=1$.
Since in this case the identification of $\xi_{Y_{w}}^{\perp}$ with $T\Sigma_{w}$ via the map
$p_{*}$ is orientation preserving, it's clear that the local index of the map
$v=p\circ u$ agrees with the intersection number of $\tilde u$ with $\R\times Y_{w}$.
Summing over all such intersection points leads to the first term in the
given formula for the degree.
To compute the local degree at a positive puncture we
first choose positive holomorphic cylindrical coordinates
$(s, t)\in [0, \infty)\times S^{1}$ on a deleted neighborhood of the puncture
$z\in\Gamma^{+}$.
Because in such a coordinate system the loop $t\mapsto (s, t)$ for fixed $s$ encircles the puncture in the
clockwise direction, the local degree of the map
$\bar v$ at $z$ can be computed by identifying a neighborhood of $w$ with
$T_{w}\Sigma$ and computing
\[
-\wind\bar v(s, t)=-\wind (p\circ u)(s, t).
\]
But considering the asymptotic representation of the normal component of the map $\tilde u$ from
Theorem \ref{t:normal-asymptotics} along with the definition of the normal relative winding, 
we have that 
\[
\wind (p\circ u)(s, t)=\wind^{\Phi}_{rel}((\tilde u; z), \R\times Y_{w})
\]
which shows the local degree of the $\bar v$ at $z$ is
\[
-\wind^{\Phi}_{rel}((\tilde u; z), \R\times Y_{w}).
\]

At negative punctures, we argue similarly but instead choose negative
cylindrical coordinates $(s, t)\in (-\infty, 0]\times S^{1}$ on a deleted neighborhood of the puncture.
Since the loop $t\mapsto (s, t)$ encircles the puncture in the counterclockwise direction,
an argument analogous to that given above tells us that 
the local degree is now given by
\[
\wind\bar v(s, t)=\wind (p\circ u)(s, t)=\wind^{\Phi}_{rel}((\tilde u; z), \R\times Y_{w}).
\]
The claimed formula for the degree of $\bar v$ now follows in the case that
$\sign(w)=1$.
The case when $\sign(w)=-1$ is identical with the exception of the fact that the map
$p_{*}:\xi_{Y_{w}}^{\perp}\to T_{w}\Sigma$ is now orientation reversing
which introduces a factor of $-1$ into each
of the computations.

Finally, the inequality \eqref{e:degree-inequality} and the claim about when equality is achieved is
an immediate consequence of local positivity of intersections and the bounds
\eqref{e:wind-rel-pos} and \eqref{e:wind-rel-neg}.
\end{proof}

Finally, we close this section with our main theorem which, given some assumptions on the normal
Conley--Zehnder indices of the map $\tilde u$, bounds the degree of the map
$\bar v$ in certain cases, and shows that in the event that $\deg(\bar v)=0$, the image of the
map $\tilde u$ must be contained in a leaf of foliation.
The motivation for the assumed bounds is the fact that in the particular
examples of holomorphic foliations
constructed in Sections \ref{holhyps} and \ref{prequantcurves} above,
the normal linearized flow along the orbits contained in 
$V=p^{-1}(\crit(f))$ can be expressed in terms of the Hessian of the function $f$.
In particular, one can show that given a number $T>0$, every periodic orbit
$\gamma\in Y_{w}$ with period less than $T$ has normal Conley--Zehnder index given by the formula
\[
\mu_N^{\Phi}(\gamma)=\sign(w)\bp{\ind_{w}(f)-1}\in\br{-1, 0, 1}
\]
provided the function $f$ is sufficiently $C^{2}$ close to a constant.

\begin{thm}\label{posdegree}
Assume that 
\begin{equation}\label{e:cz-condition}
\mu_N^{\Phi}(\gamma_{z}^{m_{z}})\in\br{-1, 0, 1}
\end{equation}
for every $z\in\Gamma$.
Then:
\begin{itemize}
\item If $\sign(w)=1$ for all $w\in\crit(f)$ then the map $\bar v$ has nonnegative degree and
has degree zero
precisely when the image of the map $\tilde u$ is contained in a leaf of the foliation $\F$.

\item If $\sign(w)=-1$ for all $w\in\crit(f)$ then the map $\bar v$ has nonpositive degree and
has degree zero
precisely when the image of the map $\tilde u$ is contained in a leaf of the foliation $\F$.

\item If the map $\sign:\crit(f)\to\br{-1, 1}$ is surjective then the map $\bar v$ has degree zero and
the image of the map $\tilde u$ is contained in a leaf of the foliation $\F$.
\end{itemize}
\end{thm}

\begin{proof}
The condition \eqref{e:cz-condition} implies that
\[
\alpha_{N}^{\Phi; -}(\gamma_{z}^{m_{z}})=\fl{\mu_{N}^{\Phi}(\gamma_{z}^{m_{z}})/2}\le \fl{1/2}=0
\]
for all $z\in \Gamma^{+}$ and
\[
\alpha_{N}^{\Phi; +}(\gamma_{z}^{m_{z}})=\ceil{\mu_{N}^{\Phi}(\gamma_{z}^{m_{z}})/2}\ge \ceil{-1/2}=0
\]
for all $z\in\Gamma^{-1}$, which together are equivalent to
\begin{equation}\label{e:alpha-inequality}
\mp\alpha_{N}^{\Phi; \mp}(\gamma_{z}^{m_{z}})\ge 0
\end{equation}
for all $z\in\Gamma^{\pm}$.

Assuming that $\sign(w)=1$ for all $w\in\crit(f)$, then 
\eqref{e:degree-inequality} from Lemma \ref{l:degree-computation}
with \eqref{e:alpha-inequality}
shows that $\deg\bar v\ge 0$ and that
$\deg\bar v=0$ precisely when
\begin{equation}\label{e:degree-zero-condition}
\begin{gathered}
\text{$u(\Sigma\setminus\Gamma)$ and $Y_{w}$ are disjoint for all $w\in\crit(f)$, and} \\
\text{$\wind_{rel}^{\Phi}((C; z), Y_{w_{z}})=\mp\alpha^{\Phi;\mp}(\gamma_{z}^{m_{z}})=0$ for all $z\in\Gamma^{\pm}$}.
\end{gathered}
\end{equation}
Similarly if  $\sign(w)=-1$ for all $w\in\crit(f)$, the same argument shows that
$\deg\bar v\le 0$ and that
$\deg\bar v=0$ precisely when \eqref{e:degree-zero-condition} holds.

In the third case that there exist points $w_{+}$, $w_{-}\in\crit(f)$ with
$\sign(w_{+})=1$ and $\sign(w_{-})=-1$, applying
\eqref{e:degree-inequality} and \eqref{e:alpha-inequality}
with $w=w_{+}$
yields $\deg(\bar v)\ge 0$
while applying 
\eqref{e:degree-inequality} and \eqref{e:alpha-inequality} with $w=w_{-}$ yields
$\deg(\bar v)\le 0$.  We conclude that
$\deg(\bar v)=0$ which once again happens precisely when
\eqref{e:degree-zero-condition} holds.

To complete the proof of all three cases, it remains to show that
$\deg(\bar v)=0$ precisely when the image of the map $\tilde u$ is contained in a leaf of the foliation
$\F$.
Assume that $\deg(\bar v)=0$.  Since $\F$ is a foliation, and we assume that the image of $u$ is not
contained in the binding $V=\crit(f)$, there exists at least one leaf $\tilde Y$ not fixed by the $\R$-action
for which $\tilde u$ intersects $\tilde Y$.
Assume $\tilde Y$ projects via $p\circ\pi$ to a gradient flow line between critical points
$w_{1}$ and $w_{2}$ so that $\tilde Y$ is asymptotically cylindrical over $Y_{w_{1}}\cup Y_{w_{2}}$,
and recall that we've noted above that with respect to the global trivialization $\Phi$ for
$\xi_{V}^{\perp}$ we've chosen,
$\wind_{\infty}^{\Phi}(\tilde Y, \gamma)=0$ for all periodic orbits $\gamma\in Y_{w_{1}}\cup Y_{w_{2}}$.
We would now like to apply Corollary \ref{c:nonintersection-conditions} above to prove
the image of $\tilde u$ is contained in $\tilde Y$.
Assume to the contrary that the image of $\tilde u$ is not contained in $\tilde Y$.
By the assumptions that all asymptotic limits of the map $\tilde u$ are periodic orbits
in the binding $V=p^{-1}(\crit(f))$ we know that none of these limits intersect
$\pi(\tilde Y)$ since $\tilde Y$ is assumed to project to an embedding in $M\setminus V$.
Moreover, since
we have already observed that in each case, $\deg(\bar v)=0$ is true precisely when
\eqref{e:degree-zero-condition} holds, we can conclude that
$u(\Sigma'\setminus\Gamma)$ does not intersect the asymptotic limit set
$Y_{w_{1}}\cup Y_{w_{2}}$ and that for each $z\in\Gamma^{\pm}$ with
$\gamma_{z}^{m_{z}}\subset Y_{w_{1}}\cup Y_{w_{2}}$
we'll have that 
$\wind^{\Phi}_{rel}((u; z), Y_{w_{1}}\cup Y_{w_{2}})=0=m_{z} \wind_{\infty}^{\Phi}(\tilde Y, \gamma_{z})$.
Corollary \ref{c:nonintersection-conditions} now lets us conclude that the image of
$\tilde u$ is disjoint from $\tilde Y$ in contradiction to the assumption that they intersect.
We thus conclude that the image of $\tilde u$ is contained in $\tilde Y$ as desired.
\end{proof}

\newpage

\bibliographystyle{alpha}

\end{document}